\documentclass[a4paper]{amsart}
\usepackage{mathtools,amsfonts,amssymb,mathrsfs}
\usepackage[a4paper,hmargin=26mm,vmargin=28mm]{geometry}
\usepackage{bm,multicol,mathbbol}
\usepackage{booktabs}
\usepackage[svgnames]{xcolor}
\usepackage[pagebackref=true]{hyperref}
\usepackage[hyperpageref]{backref}
\hypersetup{pdftitle={Quiver Schur algebras for linear quivers},
  pdfauthor={Jun Hu and Andrew Mathas},
  colorlinks = true,
  linkcolor =ForestGreen,
  anchorcolor = red,
  citecolor = blue,
  urlcolor = blue
}
\renewcommand*{\backref}[1]{}
\renewcommand*{\backrefalt}[4]{%
  [{\tiny\ifcase #1 Not cited.\relax\or Page~#2.\else Pages #2.\fi}]%
}

\DeclareFontFamily{OT1}{pzc}{}
\DeclareFontShape{OT1}{pzc}{m}{it}{<-> s * [1.20] pzcmi7t}{}
\DeclareMathAlphabet{\mathpzc}{OT1}{pzc}{m}{it}

\DeclareSymbolFontAlphabet{\mathbb}{AMSb}
\usepackage{todonotes}

\usepackage{cite}
\newcommand\Comment[2][Hu]{\space\par\medskip\noindent%
   \fbox{\begin{minipage}{\textwidth}\textbf{Comment\ifx\relax#1\else---#1\fi}\newline%
        #2\end{minipage}}\medskip
}
\usepackage{tikz}
\usetikzlibrary{matrix}

\newcount\tableauRow\newcount\tableauCol
\newcommand\Tableau[2][-4]{
  \begin{tikzpicture}[scale=0.4,draw/.append style={thick,black},baseline=#1mm]
    \tableauRow=0
    \foreach \Row in {#2} {
       \tableauCol=1
       \foreach\k in \Row {
          \draw(\the\tableauCol,\the\tableauRow)+(-.5,-.5)rectangle++(.5,.5);
          \draw(\the\tableauCol,\the\tableauRow)node{\k};
          \global\advance\tableauCol by 1
       }
       \global\advance\tableauRow by -1
    }
  \end{tikzpicture}
}
\def\Bitab(#1|#2){\Bigg(%
\hspace*{1mm}\Tableau{#1}\hspace*{2mm}\Bigg|\hspace*{2mm}\Tableau{#2}\hspace*{1mm}\Bigg)
}
\def\Tritab(#1|#2|#3){\Bigg(%
  \hspace*{1mm}\Tableau{#1}\hspace*{2mm}\Bigg|\hspace*{2mm}\Tableau{#2}\hspace*{2mm}%
       \Bigg|\hspace*{2mm}\Tableau{#3}\hspace*{1mm}\Bigg)
}

\newcommand\ShadedTableau[2][\relax]{
  \begin{tikzpicture}[scale=0.4,draw/.append style={thick,black},baseline=-4mm]
    \ifx\relax#1\relax%
    \else % shade the boxes in #1
      \foreach\box in {#1} { \filldraw[blue!30]\box+(-.5,-.5)rectangle++(.5,.5); }
    \fi
    \newcount\TabRow\newcount\TabCol
    \TabRow=0
    \foreach \Row in {#2} {
       \TabCol=1
       \foreach\k in \Row {
          \draw(\the\TabCol,\the\TabRow)+(-.5,-.5)rectangle++(.5,.5);
          \draw(\the\TabCol,\the\TabRow)node{\k};
          \global\advance\TabCol by 1
       }
       \global\advance\TabRow by -1
    }
  \end{tikzpicture}
}

\newcommand\YoungDiagram[2][\relax]{
  \begin{tikzpicture}[scale=0.5,draw/.append style={thick,black},baseline=-2mm]
    \ifx\relax#1\relax%
    \else % shade the boxes in #1
    \foreach\box in {#1} {
      \filldraw[blue!30]\box rectangle ++(1,1);
    }
    \fi
    \newcount\tableauRow
    \tableauRow=0
    \foreach \diagramCol in {#2} {
       \draw(1,\the\tableauRow)grid ++(\diagramCol,1);
       \global\advance\tableauRow by -1
    }
  \end{tikzpicture}
}
\def\TriDiagram(#1|#2|#3){\Bigg(%
  \hspace*{1mm}\YoungDiagram{#1}\hspace*{2mm}\Bigg|\hspace*{2mm}\YoungDiagram{#2}\hspace*{2mm}%
    \Bigg|\hspace*{2mm}\YoungDiagram{#3}\hspace*{1mm}\Bigg)
}

\let\<=\langle
\let\>=\rangle
\def\({\big(}
\def\){\big)}

\def\Sum{\displaystyle\sum}
\def\bijection{\overset{\sim}{\longrightarrow}}
\let\surjection\twoheadrightarrow

\def\trans{\text{tr}}  % matrix transpose

\newcommand{\C}{\mathbb{C}}
\newcommand{\N}{\mathbb N}

\newcommand{\Q}{\mathbb Q}
\newcommand{\Z}{\mathbb Z}
\newcommand{\Sym}{\mathfrak S}

\def\map#1#2{\,{:}\,#1\!\longrightarrow\!#2}
\newcommand\SetBox[2][34]{\Big\{\vcenter{\hsize#1mm\centering#2}\Big\}}
{\catcode`\|=\active
  \gdef\set#1{\mathinner{\lbrace\,{\mathcode`\|"8000%
                                   \let|\midvert #1}\,\rbrace}}
}
\def\midvert{\egroup\,\mid\,\bgroup}

\def\pmod#1{\text{ }(\text{mod } #1)\,}

\def\D{\mathcal D}
\def\Ocal{\mathcal{O}}

\renewcommand\O[1][\relax]{\Ocal^\Lambda_{#1}}
\newcommand\dualO[1][\beta]{\widetilde{\Ocal}_\Lambda^{#1}}
\newcommand\gl[1][N]{\mathfrak{gl}_{#1}(\C)}
\newcommand\h{\mathfrak h}
\newcommand\p{\mathfrak p}
\newcommand\q{\mathfrak q}
\newcommand\levi[1][\psi]{\mathfrak{l}_{#1}}
\newcommand\SL[1][e]{\widehat{\mathfrak{sl}}_{#1}}

\let\eps\varepsilon

\def\eR{e_\beta^\Ocal}

\def\Fun{\mathsf{F}}
\def\Gun{\mathsf{G}}
\def\Hun{\mathsf{H}}
\def\Oun{\mathsf{O}}
\newcommand\RFun{\mathsf{Rin}^\Lambda_\beta}
\newcommand\SFun[1][n]{\Fun^\Lambda_{#1}}
\newcommand\SpFun[1][\beta]{\Fun_{\Lambda}^{#1}}

\newcommand\SOFun[1][\beta]{\Fun^{\Ocal}_{#1}}

\def\DotGun{\dot\Gun^\Lambda_n}
\def\DotFun{\dot\Fun^\Lambda_n}
\def\Fomega{\dot\Fun^\omega_n}
\def\Gomega{\dot\Gun_\omega^n}

\newcommand\OFun[1][\beta]{\Fun^{\Ocal}_{#1}}
\newcommand\EFun[1][\beta]{\mathsf{E}^\Ocal_{#1}}
\def\UnEFun{\underline{\mathsf{E}}^\Lambda_{\Ocal}}
\newcommand\EFlat{\mathsf{E}^\flat_\beta}

\def\EDJM{\underline{\mathsf{E}}^\Lambda_{\text{DJM}}}
\def\SWEquiv{\mathsf{E}^{\text{SW}}_n}
\def\SWFun{\mathsf{F}^{\text{SW}}_n}

\newcommand\HH{\mathcal{H}}   % \HH and \RR define the 'default' symbol
\newcommand\RR{\mathcal{R}}
\renewcommand\H[1][n]{\HH^\Lambda_{#1}}
\newcommand\Rmu[1][\bmu]{(\R)^{\gdom#1}}
\newcommand\R[1][n]{\RR^\Lambda_{#1}}
\newcommand\Rp[1][n]{\RR^{\Lambda'}_{#1}}
\newcommand\UnR[1][n]{\underline{\RR}^\Lambda_{#1}}
\newcommand\RO[1][\beta]{\RR^{\Ocal}_{#1}}
\newcommand\UnRO[1][\beta]{\underline{\RR}^{\Ocal}_{#1}}
\newcommand\RFlat[1][\beta]{\RR^\flat_{#1}}
\newcommand\UnRFlat[1][\beta]{\underline{\RR}^\flat_{#1}}

\def\AFlat{A_\flat}

\def\SS{\mathcal{S}}  % \SS defines the 'default' symbol for all Schur algebras
\newcommand\Sch[1][\beta]{\SS^\Lambda_{#1}}

\newcommand\Schp[1][\beta']{\SS^{\Lambda'}_{#1}}
\newcommand\UnS[1][n]{\underline{\SS}^\Lambda_{#1}}
\newcommand\dualS[1][\beta']{\SS_{\Lambda'}^{#1}}
\newcommand\dualSp[1][\beta]{\SS_{\Lambda}^{#1}}
\newcommand\ddS[1][\beta]{\SS_{\Lambda}^{#1}}
\newcommand{\ScDJM}[1][n]{\underline{\SS}^{\text{DJM}}_{#1}}
\newcommand{\ScSW}[1][n]{\SS^{\text{SW}}_{#1}}
\newcommand{\UnScSW}[1][n]{\underline{\SS}^{\text{SW}}_{#1}}

\newcommand\DotS[1][n]{\dot{\SS}^\Lambda_{#1}}
\newcommand\SO[1][\beta]{\SS^{\Ocal}_{#1}}

\newcommand\UnSO[1][\beta]{\underline{\SS}^\Ocal_{#1}}

\newcommand\UnP[1][\bmu]{\underline{P}^{#1}}
\newcommand\unL{\underline{L}}
\newcommand\UnL[1][\blam]{\underline{L}^{#1}}
\newcommand\UnM[1][\p]{\underline{M}_{#1}}
\newcommand\unM{\underline{M}}
\newcommand\DotG[1][n]{\dot{G}^\Lambda_{#1}}
\newcommand\DelO[1][\blam]{\Delta^{#1}_{\Ocal}}

\newcommand\NabO[1][\blam]{\nabla^{#1}_{\Ocal}}

\newcommand\UnDelO[1][\blam]{\underline{\Delta}^{#1}_{\Ocal}}

\newcommand\LO[1][\bmu]{L^{#1}_{\Ocal}}
\newcommand\PO[1][\bmu]{P^{#1}_{\Ocal}}

\newcommand\UnLO[1][\bmu]{\underline{L}^{#1}_{\Ocal}}

\newcommand\dualLO[1][\blam]{L_{#1}^{\Ocal}}
\newcommand\dualUnLO[1][\blam]{\underline{L}_{#1}^{\Ocal}}
\newcommand\UnPO[1][\bmu]{\underline{P\!}\,^{#1}_{\Ocal}}
\newcommand\YO[1][\bmu]{Y^{#1}_{\Ocal}}

\newcommand\SFlat[1][\beta]{\SS^\flat_{#1}}
\newcommand\UnSFlat[1][\beta]{\underline{\SS}^\flat_{#1}}
\newcommand\PFlat[1][\bmu]{P_\flat^{#1}}

\newcommand\YFlat[1][\blam]{Y_\flat^{#1}}

\newcommand\UnPFlat[1][\bmu]{\underline{P}_\flat^{#1}}
\newcommand\LFlat[1][\bmu]{L_\flat^{#1}}
\newcommand\LDJM[1][\bmu]{\underline{L}^{#1}_{\text{DJM}}}
\newcommand\WDJM[1][\blam]{\underline{\Delta}^{#1}_{\text{DJM}}}
\newcommand\Mmu[1][\bmu]{\underline{M}^{#1}}
\newcommand\Nmu[1][\bmu]{\underline{N}_{#1}}

\newcommand\LSW[1][\bmu]{L_{#1}^{\text{SW}}}
\newcommand\PSW[1][\bmu]{P_{#1}^{\text{SW}}}
\newcommand\WSW[1][\blam]{\Delta_{#1}^{\text{SW}}}
\newcommand\YSW[1][\bmu]{Y_{#1}^{\text{SW}}}
\newcommand\UnLSW[1][\bmu]{\underline{L}_{#1}^{\text{SW}}}
\newcommand\UnYSW[1][\bmu]{\underline{Y}_{#1}^{\text{SW}}}
\newcommand\UnWSW[1][\blam]{\underline{\Delta}_{#1}^{\text{SW}}}

\def\P{\mathscr{P}}
\newcommand{\Klesh}[1][n]{\mathcal{K}^{\Lambda}_{#1}}
\newcommand{\Parts}[1][n]{\P^{\Lambda}_{#1}}

\newcommand\DotParts[1][n]{\dot{\P}^\Lambda_{#1}}
\def\tomega{\t^\omega}
\newcommand\blam{{\boldsymbol\lambda}}

\newcommand\bmu{{\boldsymbol\mu}}
\newcommand\bump{{\boldsymbol\mu'}}
\newcommand\bnu{{\boldsymbol\nu}}
\newcommand\bsig{{\boldsymbol\sigma}}
\newcommand\btau{{\boldsymbol\tau}}

\newcommand\bi{\mathbf{i}}
\newcommand\bj{\mathbf{j}}
\def\a{{\mathfrak a}}
\def\b{{\mathfrak b}}

\def\s{{\mathfrak s}}

\def\t{{\mathfrak t}}
\def\unt{\underline{\t}}
\def\T{{\mathsf T}}
\let\uold=\u
\def\u{{\mathfrak u}}
\def\v{{\mathfrak v}}
\def\rest#1{{}_{\downarrow #1}}
\def\ilam{\bi^\blam}
\def\imu{\bi^\bmu}
\def\ilmu{\bi_\bmu}
\def\tlam{\t^\blam}
\def\tllam{\t_\blam}
\def\tmu{\t^\bmu}
\def\tlmu{\t_\bmu}
\def\tlmup{\t_{\bmu'}}
\def\tnu{\t^\bnu}

\def\Std{\mathop{\rm Std}\nolimits}
\newcommand\Col[1][]{\mathop{\rm Col}\nolimits^\Lambda_{#1}}
\def\SStd{\mathop{\rm Std}\nolimits^2}
\def\Tcal{\mathcal{T}}

\let\gedom=\trianglerighteq
\let\gdom=\vartriangleright
\let\Gdom=\blacktriangleright
\DeclareMathOperator\Gedom{\,{\underline{\kern-.1ex{\blacktriangleright}\kern-0.1ex}}\,}
\let\ledom=\trianglelefteq
\let\ldom=\vartriangleleft

\newcommand{\charge}{{\boldsymbol{\kappa}}}

\newcommand\Fock{\mathfrak{F}^\Lambda}

\newcommand\ch[1][q]{\mathop{\rm ch}\nolimits_{#1}}
\newcommand\Dim[1][q]{\mathop{\mathpzc{dim}}\nolimits_{#1}}
\DeclareMathOperator\End{End}
\DeclareMathOperator\Ext{Ext}
\DeclareMathOperator\Hom{Hom}

\DeclareMathOperator\Ker{Ker}
\DeclareMathOperator\Mod{{-}Mod}
\DeclareMathOperator\Proj{Proj}
\DeclareMathOperator\Rep{Rep}

\DeclareMathOperator\Shape{Shape}

\DeclareMathOperator\codeg{codeg}
\DeclareMathOperator\col{col}
\DeclareMathOperator\comp{comp}

\newcommand\defect[1][]{\mathop{\rm def}\nolimits_{#1}}

\DeclareMathOperator\mindeg{mindeg}

\DeclareMathOperator\rad{rad}
\DeclareMathOperator\res{res}

\DeclareMathOperator\soc{soc}
\DeclareMathOperator\wt{\omega}
\def\sgn{\mathtt{sgn}}
\def\mz{\mathcal{Z}}
\def\RBasis{\Theta^\Ocal_\beta}
\newcommand\thetads{\theta^{(d,s)}_{\blam\bmu}}
\newcommand\varthetads{\vartheta^{(d,s)}_{\blam\bmu}}

\def\op{{\text{op}}}

\DeclareMathOperator\Gr{\mathpzc{Gr}}
\DeclareMathOperator\ZEnd{\mathpzc{End}}

\DeclareMathOperator\ZHom{\mathpzc{Hom}}

\newenvironment{Notation}{%
  \small\noindent%
  \multicolsep=1mm%
  \columnsep=4mm%
  \begin{center}\parindent=0pt%
  \begin{multicols}{2}\noindent%
}{\end{multicols}\end{center}}

\def\notation#1#2#3{\rlap{\hyperref[#1]{#2}}\hspace*{14mm} \hbox to 47mm{#3\hfill}}

\usepackage{aliascnt}
\def\NewTheorem#1{%
  \newaliascnt{#1}{equation}
  \newtheorem{#1}[#1]{#1}
  \aliascntresetthe{#1}
  \expandafter\def\csname #1autorefname\endcsname{#1}
}
\def\equationautorefname~#1\null{(#1)\null}

\newcounter{main}
\theoremstyle{plain}
\newtheorem{THEOREM}[main]{Theorem}

\swapnumbers
\numberwithin{equation}{section}
\NewTheorem{Assumption}
\NewTheorem{Proposition}
\NewTheorem{Theorem}
\NewTheorem{Corollary}
\NewTheorem{Conjecture}
\NewTheorem{Lemma}
\NewTheorem{Definition}
\theoremstyle{definition}
\theoremstyle{remark}
\NewTheorem{Remark}

\newaliascnt{Example}{equation}
\aliascntresetthe{Example}

\newenvironment{Example}%
  {\refstepcounter{Example}\trivlist
   \item[\hskip\labelsep\theequation.~\textbf{Example}\space]
   \ignorespaces
  }{\unskip\nobreak\hfil%
    \penalty50\hskip2em\hbox{}\nobreak\hfil$\Diamond$%
    \parfillskip=0pt\finalhyphendemerits=0\penalty-100\endtrivlist
}

\newcounter{MyCase}
\numberwithin{MyCase}{equation}

\title{Quiver Schur algebras for linear quivers}
\subjclass[2000]{20G43, 20C08, 20C30, 05E10}
\keywords{Cyclotomic Hecke algebras, Schur algebras, quasi-hereditary and graded
cellular algebras, Khovanov--Lauda--Rouquier algebras}
\author{Jun Hu}\address{School of Mathematics,
Beijing Institute of Technology,
Beijing, 100081, P.R.~China}
\address{School of Mathematics and Statistics F07, University of Sydney, NSW 2006, Australia}

\def\Email#1{\email{\href{mailto:#1}{#1}}}
\Email{junhu404@bit.edu.cn}

\author{Andrew Mathas}
\address{School of Mathematics and Statistics, University of Sydney, NSW 2006, Australia}
\Email{andrew.mathas@sydney.edu.au}

\begin{document}
\bibliographystyle{andrew}

\begin{abstract}
We define a graded quasi-hereditary covering of the cyclotomic quiver
Hecke algebras $\R$ of type~$A$ when $e=0$ (the linear quiver) or $e>n$.
We prove that these algebras are quasi-hereditary graded cellular
algebras by giving explicit homogeneous bases for them. When $e=0$ we
show that the KLR grading on the quiver Hecke algebras is compatible
with the Koszul grading on the blocks of parabolic category~$\O$ given
by  Backelin, building on the work of Beilinson, Ginzburg and Soergel.
As a consequence, $e=0$ our cyclotomic quiver Schur
algebras are Koszul over fields of characteristic zero.  Finally, we
give an LLT-like algorithm for computing the graded decomposition
numbers of the cyclotomic quiver Schur algebras in characteristic zero.
\end{abstract}
\maketitle\setcounter{tocdepth}{1}\tableofcontents

\section{Introduction}
  Khovanov and Lauda~\cite{KhovLaud:diagI,KhovLaud:diagII} and
  Rouquier~\cite{Rouq:2KM} have introduced a remarkable family of $\Z$-graded
  algebras that are now known to categorify the canonical bases of Kac-Moody
  algebras~\cite{BK:GradedDecomp, BrundanStroppel:KhovanovIII,%
                 VaragnoloVasserot:CatAffineKLR}.
  Brundan and Kleshchev~\cite{BK:GradedKL} initiated the study of
  `cyclotomic' quotients of these algebras by showing that they are isomorphic
  to the degenerate and non-degenerate cyclotomic Hecke algebras of
  type~$G(\ell,1,n)$; see also \cite{Rouq:2KM}.

  This paper defines and studies certain graded quasi-hereditary
  covers~$\Sch[n]$ of the cyclotomic quiver Hecke algebras of the linear
  quiver. These algebras are graded analogues of the cyclotomic Schur
  algebras~$\ScDJM$ of
  type~$G(\ell,1,n)$~\cite{DJM:cyc,BK:HigherSchurWeyl}. This paper
  studies the cyclotomic quiver Schur algebras for the linear quiver and
  `large' cyclic quiver. Let $\Parts$ be the poset of multipartitions
  of~$n$ ordered by dominance. The first main result of this paper is
  the following.

  \begin{THEOREM}\label{THEOREMA}
    Suppose that $e=0$ or $e>n$ and let $\mz=K$ be an arbitrary
    field. The algebra $\Sch[n]$ is a quasi-hereditary graded cellular algebra
    with graded standard modules $\set{\Delta^\bmu|\bmu\in\Parts}$ and
    irreducible modules $\set{L^\bmu|\bmu\in\Parts}$. Moreover, there is an
    equivalence of (ungraded) highest weight categories
    $$\SFun\map{\Sch[n]\Mod}\ScDJM\Mod$$
    that preserves the labelling of the standard modules and simple modules.
  \end{THEOREM}

  In fact, the quiver Schur algebra $\Sch[n]$ is defined over an arbitrary integral
  domain.

  Like the cyclotomic Schur algebras of~\cite{DJM:cyc}, the quiver Schur algebra
  $\Sch[n]$ is defined to be the endomorphism algebra of a direct sum of
  ``graded permutation modules''; see \autoref{D:QuiverSchur}. Incorporating the
  grading into this picture is surprisingly difficult, not least because many of
  the structural results in this paper fail when $1<e\leq n$. After we have defined
  the quiver Schur algebras, and shown that they are quasi-hereditary
  (\autoref{T:quasi}), most of the work in the first six sections of the paper
  is geared towards proving the non-trivial result that the cyclotomic quiver
  Schur algebras are Morita equivalent, as ungraded algebras, to the cyclotomic
  $q$-Schur algebras (\autoref{T:DJMEquivalence}). Constructing a graded
  analogue of the Schur functor (\autoref{P:SchurFunctor}) is also not
  completely straightforward.

  If $e=0$ and we work over the field of complex numbers then Brundan
  and Kleshchev have shown that the degenerate cyclotomic Schur algebras
  are Morita equivalent to a sum of certain integral blocks $\O[\beta]$
  of parabolic category~$\O$ for the Lie algebra of the general linear
  group~\cite{BK:HigherSchurWeyl}. By results of
  Backelin~\cite{Backelin:Koszul}, and Beilinson, Ginzburg and
  Soergel~\cite{BGS:Koszul}, the blocks of parabolic category~$\Ocal$
  can be endowed with a Koszul grading. By \cite{BK:HigherSchurWeyl},
  the endomorphism algebra of a prinjective generator of the sum of
  these blocks is Morita equivalent to the degenerate cyclotomic Hecke
  algebra $\H$ of type~$A$.  Therefore, the Koszul grading on the
  blocks~$\O[\beta]$ induce a grading on the module category of~$\H$.
  This gives two ostensibly different gradings on the degenerate
  cyclotomic Hecke algebra~$\H$: one coming from parabolic category~$\O$
  and the KLR grading given by the Brundan-Kleshchev isomorphism
  $\H\cong\R$~\cite{BK:GradedKL} when $e=0$.

  \begin{THEOREM}\label{THEOREMB}
    Suppose that $e=0$ and $\mz=\C$ is the field of complex numbers. Then
    category~$\O$ and the quiver Hecke algebra~$\R$ induce graded Morita
    equivalent gradings on~$\H\Mod$.
  \end{THEOREM}

  To prove \autoref{THEOREMB} we first show that the cyclotomic quiver
  Schur algebra $\Sch[\beta]$ is graded Morita equivalent to the
  (non-isomorphic) quiver Schur algebras recently constructed by
  Stroppel and Webster~\cite{StroppelWebster:QuiverSchur}. This allows
  us to show that the prinjective modules of $\Sch[\beta]$ are
  \textit{rigidly graded} and using this we can construct an explicit
  isomorphism between the basic algebras of $\O$ and $\Sch[\beta]$.

  Building on \autoref{THEOREMB}, in \autoref{S:GradedDec} we prove a graded
  analogue of \cite[Theorem~C]{BK:HigherSchurWeyl}, thus lifting Brundan and
  Kleshchev's ``higher Schur-Weyl duality'' to the graded setting.

  {\samepage
  \begin{THEOREM}\label{THEOREMC}
    Suppose that $e=0$, $\beta\in Q^+$,  and $\mz=\C$. Then there are
    graded Schur functors $\OFun[\beta]\map{\O[\beta]\Mod}\R\Mod$ and
    $\SFun[\beta]\map{\Sch[\beta]\Mod}\R\Mod$ and a graded equivalence
    $\EFun[\beta]\map{\O[\beta]}\Sch[\beta]\Mod$ such that the following
    diagram commutes:
    $$\begin{tikzpicture}
     \matrix[matrix of math nodes,row sep=1cm,column sep=16mm]{
     |(O)| \O[\beta]&|(S)|\Sch[\beta]\Mod\\
                    &|(H)|\R[\beta]\Mod\\
     };
     \draw[->](O) -- node[above]{$\EFun[\beta]$} (S);
     \draw[->](S) -- node[right]{$\SFun[\beta]$} (H);
     \draw[->](O) -- node[below]{$\OFun[\beta]$} (H);
    \end{tikzpicture}$$
    In particular, $\Sch[\beta]\Mod$ is Koszul.
  \end{THEOREM}
  }

  This result can be interpreted as saying that the KLR grading on $\H$ induces
  the Koszul grading on parabolic category~$\O$. As a consequence, via our
  graded Schur algebra, we obtain a very explicit and new combinatorial
  description of the grading on parabolic category~$\O$.

  Our description of the grading on parabolic category~$\O$ gives new
  information. For example we use it in \autoref{S:LLT} to give a fast algorithm
  for computing the graded decomposition numbers of $\Sch[n]$ and parabolic
  category~$\O$, which are certain parabolic Kazhdan-Lusztig polynomials, that
  is similar in spirit to the LLT algorithm for the Hecke algebras of
  type~$A$~\cite{LLT}. What is really interesting about our ``LLT algorithm'' is
  that it computes the graded decomposition numbers of the quiver Schur algebras
  when $e=0$. In contrast, the extension of the LLT algorithm to the $q$-Schur
  algebras~\cite{LLT:Schur} is non-trivial because it requires first computing
  the action of the bar involution on the Fock space.

  If $\Lambda$ is a dominant weight of level~$2$ then we show in
  \autoref{C:BrundanStroppel} that $\Sch[n]$ is isomorphic, as a graded
  algebra, to the corresponding quasi-hereditary cover of the Khovanov's
  diagram algebra as introduced by Brundan and
  Stroppel~\cite{BrundanStroppel:KhovanovI}. This is quite surprising
  because the definition of these two algebras is very different. The
  key is to show that if $\Lambda$ is a weight of level~$2$ then
  $\Sch[n]$ is a positively graded basic
  algebra (\autoref{T:L2Positivity}).  This positivity result, and the
  uniqueness of Koszul gradings and \autoref{THEOREMC}, provides the
  bridge to Brundan and Stoppel's algebra.

\section*{Acknowledgements}

  We thank the referee for their careful reading of this manuscript.
  Both authors were supported by the Australian Research Council. The
  first author was also supported by the National Natural Science
  Foundation of China.

\section*{Index of notation}

\begin{Notation}
  \notation{S:graded}         {$A, M$}{A graded algebra or module}
  \notation{S:graded}         {$\underline{A}, \underline{M}$}%
                                  {An ungraded algebra or module}
  \notation{C:RSimples}       {$D^\bmu$}{Simple $\R$-module}
  \notation{S:tableaux}       {$d(\t),d'(\t)$}{Permutations: $\t=\tmu d(\t)=\tlmu d'(\t)$}
  \notation{D:dlammu}{$d_{\blam\bmu}(q)$}{$[\Delta^\blam:L^\bmu]_q$ for $\Sch[\beta]$}
  \notation{T:CatODecompNumbers}{$d^\Ocal_{\blam\bmu}(q)$}{$[\DelO:\LO]_q$ for $\SO$}
  \notation{S:tableaux}       {$\deg\t$}{Tableau degree}
  \notation{S:tableaux}       {$\codeg\t$}{Tableau codegree}
  \notation{E:defect}         {$\defect\beta$}{Defect of~$\beta\in Q^+$}
  \notation{E:GrDim}          {$\Dim M$}{Graded dimension of $M$}
  \notation{E:WeylBasis}      {$\Delta^\blam,\nabla^\blam$}{Weyl and costandard modules}
  \notation{T:DualSchur}      {$\Delta_\blam,\nabla_\blam$}{Sign dual (co)standard modules}
  \notation{D:psis}{$e^\bmu,e_\bmu$}{KLR idempotents $e(\bi^\bmu), e(\bi_\bmu)$}
  \notation{L:BarInvolution}  {$e_{\blam\bmu}(q)$}{Inverse decomposition number}
  \notation{D:Emu}            {$E^\bmu$}{Graded exterior powers}
  \notation{S:GradedSchurFunctors}  {$\SFun,\SFun[\beta]$}{Graded Schur functors}
  \notation{D:Gmu}            {$G^\bmu,G_\bmu$}  {Graded permutation modules}
  \notation{D:QuiverSchur}    {$G^\Lambda_n$}{$\bigoplus_{\bmu\in\Parts}G^\bmu$}
  \notation{S:graded}         {$\End_A$}{Degree preserving endomorphisms}
  \notation{S:graded}         {$\ZEnd_A$}{All $A$-module endomorphisms}
  \notation{S:Fock}           {$\Fock$}{Combinatorial Fock space}
  \notation{D:HeckeAlgebras}  {$\H,\H[\beta]$}{Cyclotomic Hecke algebras}
  \notation{S:graded}         {$\Hom_A$}{Degree preserving maps in $A\Mod$}
  \notation{S:graded}         {$\ZHom_A$}{All $A$-module homomorphisms}
  \notation{S:tableaux}       {$\imu,\ilmu$}{$\res(\tmu)$ and $\res(\tlmu)$}
  \notation{L:defect}         {$I^\beta$}{$\set{\bi\in I^n|\sum_{r=1}^\ell\alpha_{i_r}=\beta}$}
  \notation{E:multicharge}    {$\charge$}{Multicharge determining $\Lambda=\Lambda(\charge)$}
  \notation{S:ContDuality}    {$\Klesh$}{Restricted multipartitions for $\R$}
  \notation{S:KLR}            {$\Lambda$}{Dominant weight determined by $\charge$}
  \notation{T:quasi}          {$L^\bmu$}{Simple $\Sch[n]$-module}
  \notation{S:tableaux}       {$\bmu'$}{Conjugate multipartition}
  \notation{S:CatO}           {$\O[\beta]$}{Parabolic category~$\Ocal$}
  \notation{S:KLR}            {$P^+$}{Positive weight lattice}
  \notation{S:Young}          {$P^\bmu$}{Projective cover of $L^\bmu$}
  \notation{S:HigherSchurWeyl}{$P_{x,y}$}{Kazhdan-Lusztig polynomial}
  \notation{S:tableaux}       {$\Parts$}{Multipartitions of $n$}
  \notation{E:Rblocks}        {$\Parts[\beta]$}{$\set{\bmu\in\Parts|\imu\in I^\beta}$}
  \notation{T:CatODecompNumbers}{$p_{\blam\bmu}(q)$}{Inverse graded decomposition numbers}
  \notation{T:CatODecompNumbers}{$p^\Ocal_{\blam\bmu}(q)$}{Inverse graded decomposition numbers}
  \notation{E:Psi}            {$\Psi^{\bmu\blam}_{\s\t}$}{Basis elements of $\Sch[n]$}
  \notation{D:psis}           {$\psi_{\s\t},\psi'_{\s\t}$}{Basis elements of $\R$}
  \notation{S:GradedSchurFunctors}  {$\Psi^\bmu$}{Identity map on $G^\bmu$}
  %\notation{T:EBasis}         {$\theta^{\bmu\blam}_{\s\t}$}{Basis elements of $E^\bmu$}
  %\notation{L:dualEBasis}     {$\theta_{\bmu\blam}^{\u\v}$}{Basis elements of $E^\bmu$}
  \notation{S:KLR}            {$Q^+$}{Positive root lattice}
  \notation{S:blocks}         {$Q^+_n$}{$\set{\beta\in Q^+|\Parts[\beta]\ne0}$}
  \notation{S:tableaux}       {$\res$}{Residue sequence of tableau}
  \notation{D:QuiverRelations}{$\R$}{Cyclotomic quiver Hecke algebra}
  \notation{E:Rblocks}        {$\R[\beta]$}{A block of $\R$}
  \notation{E:sgn}            {$\sgn$}{Sign isomorphism}
  \notation{D:QuiverSchur}    {$\Sch[n]$}{Cyclotomic quiver Schur algebra}
  \notation{T:CatOKoszul}     {$\SO$}{Schur algebra of $\O[\beta]$}
  %\notation{S:SchurFunctors} {$\DotS$}{Extended graded Schur algebra}
  \notation{T:SBlocks}        {$\Sch$}{A block of $\Sch[n]$}
  \notation{S:SignDual}        {$\dualS$}{Sign-dual quiver Schur algebra}
  \notation{S:ContDuality}    {$S^\bmu,S_\bmu$}{Graded (dual) Specht modules}
  \notation{S:tableaux}       {$\Std(\Parts)$}{Standard tableaux}
  \notation{L:spanning}       {$\Std^\bmu(\Parts)$}{$\set{\t|\t\gedom\tmu\text{ and }\res(\t)=\imu}$}
  \notation{L:spanning}       {$\Std_\bmu(\Parts)$}{$\set{\t|\tlmu\gedom\t\text{ and }\res(\t)=\ilmu}$}
  \notation{T:SCellular}      {$\Tcal^\blam$}{$\set{(\bmu,\s)|\s\in\Std^\bmu(\blam)}$}
  \notation{T:DualSchur}      {$\Tcal_\blam$}{$\set{(\bmu,\s)|\s\in\Std_\bmu(\blam)}$}
  \notation{S:Tilting}        {$T^\bmu$}{Tilting module}
  \notation{T:trace}          {$\tau_\beta$}{Trace form on $\R[\beta]$}
  \notation{S:tableaux}       {$\tmu,\tlmu$}{Initial and final $\bmu$-tableaux}
  \notation{E:TlamWeight}     {$\wt(\T^\blam)$}{The weight of $\T^\blam$}
  \notation{D:psis}           {$y^\bmu,y_\bmu$}{$\psi_{\tmu\tmu}=e^\bmu y^\bmu$,
                                                $\psi_{\tlmu\tlmu}'=e_\bmu y_\bmu$}
  \notation{D:YoungModule}    {$Y^\bmu,Y_\bmu$}{Young modules}
  \notation{S:cellular}       {$\mz$}{Commutative ring}
  \notation{S:Young}          {$Z^\bmu$}{Graded symmetric power}
  \notation{S:tableaux}       {$\gdom,\Gedom$}{Dominance orderings}
  \notation{D:dlammu}         {$[M{:}L^\bmu]_q$}{Graded decomposition number}
  \notation{D:dlammu}         {$[N{:}D^\bmu]_q$}{Graded decomposition number}
  \notation{E:dual}           {$\circledast$}{$\ZHom_{\mz}(?,\mz)$-dual}
  \notation{S:Fock}           {$\#$}{$\ZHom_A(?,A)$-dual}
\end{Notation}

\section{Graded representation theory}
In this chapter we set our notation and give the reader some quick
reminders about graded modules and graded algebras, by which we will
always mean $\Z$-graded modules and $\Z$-graded algebras. Expert readers
may wish to skip this chapter.

\subsection{Graded modules and algebras}\label{S:graded} Throughout this paper, $\mz$
will be an integral domain. Unless otherwise stated, all modules and algebras
will be free and of finite rank over their base ring. We also
assume that if $A$ is a $K$-algebra then $A$ is split over~$K$. As we work
mainly with (graded) cellular algebras, which are always split, there is no loss
of generality in assuming this.

In this paper a \textbf{graded
$\mz$-module} is a $\Z$-graded $\mz$-module $M$. That is, as $\mz$-module,~$M$
has a direct sum decomposition
$$M=\bigoplus_{d\in\Z}M_d$$
By assumption, $M$ is $\mz$-free and of finite rank so $M_d\ne0$ for only finitely
many~$d$.  If $M$ is a graded $\mz$-module let $\underline{M}$ be the
ungraded $\mz$-module obtained by forgetting the grading on~$M$.

If $m\in M_d$, for $d\in\Z$, then $m$ is \textbf{homogeneous} of
\textbf{degree} $d$ and we set $\deg m=d$. If $M$ is a graded
$\mz$-module and $s\in\Z$ let $M\<s\>$ be the graded $\mz$-module
obtained by shifting the grading on~$M$ up by~$s$. That is,
$M\<s\>_d=M_{d-s}$, for $d\in\Z$. Let $q$ be an indeterminate. If
$\mz=K$ is a field then \textbf{graded dimension} of $M$ is the
Laurent polynomial
\begin{equation}\label{E:GrDim}
\Dim M=\sum_{d\in\Z}(\dim_K M_d)\,q^d\in\N[q,q^{-1}].
\end{equation}
In particular, $\dim_K M=(\Dim M)\!\bigm|_{q=1}$.

If $M$ is a graded module and if
$f(q)=\sum_{d\in\Z}f_dq^d\in\N[q,q^{-1}]$ is a Laurent polynomial
with non-negative coefficients $\{f_d\}_{d\in\Z}$ then define
$$f(q)M=\bigoplus_{d\in\Z}M\<d\>^{\oplus f_d}.$$
Then $f(q)M$ is again free and of finite rank and $\Dim (f(q)M)=f(q)\Dim M$.

A \textbf{graded $\mz$-algebra} is a unital associative
$\mz$-algebra $A=\bigoplus_{d\in\Z}A_d$ that is a graded
$\mz$-module such that $A_dA_e\subseteq A_{d+e}$, for all
$d,e\in\Z$. It follows that $1\in A_0$ and that $A_0$ is a graded
subalgebra of $A$.  A graded (right) $A$-module is a graded
$\mz$-module $M$ such that $\underline{M}$ is an
$\underline{A}$-module and $M_dA_e\subseteq M_{d+e}$, for all
$d,e\in\Z$. Graded submodules, graded left $A$-modules and so on are
all defined in the obvious way.

Let $A\Mod$ be the category of
finitely generated \textit{graded} $A$-modules with degree preserving maps. Then
$$\Hom_A(M,N)=\set{f\in\Hom_{\underline A}(\underline M,\underline N)|
               f(M_d)\subseteq N_d\text{ for all }d\in\Z},$$
for all $M,N\in A\Mod$. The elements of $\Hom_A(M,N)$ are homogeneous
maps of degree~$0$. More generally, for each $d\in\Z$ set
$$\Hom_A(M,N)_d=\Hom_A(M,N\<-d\>)\cong\Hom_A(M\<d\>,N).$$
Thus, $\Hom_A(M,N)=\Hom_A(M,N)_0$. If $f\in\Hom_A(M,N)_d$ then $f$
is \textbf{homogeneous} of degree $d$ and we set $\deg f=d$. Define
$$\ZHom_A(M,N)=\bigoplus_{d\in\Z}\Hom_A(M,N)_d=\bigoplus_{d\in\Z}\Hom_A(M,N\<-d\>).$$
Any map is a sum of its homogeneous components so
$\Hom_{\underline A}(\underline M,\underline
N)\cong\ZHom_A(M,N)$ as a $\mz$-module. Set
$\End_A(M)=\Hom_A(M,M)$ and $\ZEnd_A(M)=\ZHom_A(M,M)$.

If $r\ge0$ and $M$ and~$N$ are graded $A$-modules let $\Ext_A^r(M,N)$ be the
space of $r$-fold extensions of $M$ by~$N$ in the category $A\Mod$ of (graded)
$A$-modules.
% and set
% $$\ZExt_A^r(M,N)=\bigoplus_{d\in\Z}\Ext_A^r(M\<d\>,N).$$
% Once again, $\Ext^r_{\underline A}(\underline M,\underline N)\cong\ZExt_A^r(M,N)$
% as a $\mz$-module, for all $r\ge0$.
%
We emphasize that $\Hom_A$ and $\Ext_A$ are the spaces of homomorphisms and
extensions in the category $A\Mod$ of finitely generated (graded) $A$-modules.
These should not be confused with $\Hom_{\underline{A}}$ and
$\Ext_{\underline{A}}$ in the (ungraded) category $\underline{A}\Mod$.

Now suppose that $A$ comes equipped with a homogeneous
anti-isomorphism $\star$. Then the \textbf{graded dual} of
the graded $A$-module $M$ is the graded $A$-module
\begin{equation}\label{E:dual}
M^\circledast =
\ZHom_\mz(M,\mz)=\bigoplus_{d\in\Z}\Hom_\mz(M\<d\>,\mz)
\end{equation}
where $\mz$ is concentrated in degree zero and where the action of
$A$ on~$M^\circledast$ is given by $(fa)(m)=f(ma^\star)$ for all
$f\in M^\circledast$, $a\in A$ and~$m\in M$. The module $M$ is
\textbf{self dual} if $M\cong M^\circledast$ as graded $A$-modules.
If $\mz=K$ is a field then, as a vector space,
$M^\circledast_d=\Hom_\mz(M_{-d},K)$, so that
$\Dim M^\circledast=\overline{\Dim M}$, where the bar involution
$\rule[1.5ex]{.6em}{.1ex}\map{\Z[q,q^{-1}]}{\Z[q,q^{-1}]}$ is the $\Z$-linear
map determined by $q^k\mapsto q^{-k}$, for all $k\in\Z$.

If $m$ is an $\underline{A}$-module then a \textbf{graded lift} of
$m$ is an $A$-module $M$ such that $\underline{M}\cong m$ as
$A$-modules. In general, there is no guarantee that an $\underline
A$-module will have a graded lift but if $A$ is a graded Artin algebra (e.g., a finite dimensional algebra over a field), it is easy to see that if a finitely generated indecomposable $\underline{A}$-module has a graded lift then this
lift is unique up to isomorphism and grading shift; see for example
\cite[Lemma~2.5.3]{BGS:Koszul}. In this case, the irreducible and projective
indecomposable $\underline{A}$-modules always have graded lifts; see
\cite{GordonGreen:GradedArtin}.

Suppose that $M$ is a graded $A$-module and that
$X=\set{X^\mu|\mu\in\P}$ is a collection of $A$-modules such that
$\set{\underline{X}^\mu|\mu\in\P}$ are pairwise non-isomorphic
$\underline{A}$-modules. Then $M$ has a $X$-module
\textbf{filtration} if there exists a filtration
$$M=M_0\supset M_1\supset \dots\supset M_s=0$$
such that there exist $\mu_r\in\P$ and $d_r\in\Z$ with $M_r/M_{r+1}\cong
X^{\mu_r}\<d_r\>$, for $0\le r<s$. The
\textbf{graded multiplicity} of $X^\mu$ in $M$ is the Laurent polynomial
\begin{equation}\label{E:filtMult}
  (M:X^\mu)_q=\sum_{\substack{0\leq r\leq s-1\\\mu_r=\mu}}q^{d_r}\,\in\N[q,q^{-1}].
\end{equation}
In general, this multiplicity will depend upon the choice of filtration but for
many modules, such as irreducible modules and Weyl modules, the
Laurent polynomial $(M:X^\mu)_q$ will be independent of this choice.
% We set $[M:X^\mu]_q=(M:X^\mu)_q$ when this multiplicity is independent of the
% choice of filtration.

\subsection{Cellular algebras}\label{S:cellular}
Many of the algebras considered in this paper are (graded) cellular
algebras so we quickly recall the definition and some of the
important properties of these algebras. Cellular algebras were
defined by Graham and Lehrer~\cite{GL} with their natural extension
to the graded setting given in~\cite{HuMathas:GradedCellular}.

\begin{Definition}[Graded cellular algebra~\cite{GL,HuMathas:GradedCellular}]\label{graded cellular def}\label{D:cellular}
  Suppose that $A$ is a $\Z$-graded $\mz$-algebra that is free of
  finite rank over $\mz$. A \textbf{graded cell datum} for $A$ is an
  ordered quadruple $(\P,T,B,\deg)$, where $(\P,\gdom)$ is the
  \textbf{weight poset}, $T(\lambda)$ is a finite set for
  $\lambda\in\P$, and
  $$B\map{\coprod_{\lambda\in\P}T(\lambda)\times T(\lambda)}A;
     (\s,\t)\mapsto b_{\s\t},
     \quad\text{and}\quad
     \deg\map{\coprod_{\lambda\in\P}T(\lambda)}\Z$$
  are two functions such that  $B$ is injective and
  \begin{enumerate}
    \item[(GC$_d$)] If $\lambda\in\P$ and $\s,\t\in\T(\lambda)$ then $b_{\s\t}$ is homogeneous
    of degree $\deg b_{\s\t}=\deg\s+\deg\t$.
    \item[(GC$_1$)] $\set{b_{\s\t}|\s,\t\in T(\lambda) \text{ for } \lambda\in\P}$ is a
      $\mz$-basis of $A$.
    \item[(GC$_2$)] If $\s,\t\in T(\lambda)$, for some $\lambda\in\P$, and $a\in A$ then
    there exist scalars $r_{\t\v}(a)$, which do not depend on $\s$, such that
      $$b_{\s\t} a=\sum_{\v\in T(\lambda)}r_{\t\v}(a)b_{\s\v}\pmod
      {A^{\gdom\lambda}},$$
      where $A^{\gdom\lambda}$ is the $\mz$-submodule of $A$ spanned by
      $\set{b^\mu_{\a\b}|\mu\gdom\lambda\text{ and }\a,\b\in T(\mu)}$.
    \item[(GC$_3$)] The $\mz$-linear map $\star\map AA$ determined by
      $(b_{\s\t})^\star=b_{\t\s}$, for all $\lambda\in\P$ and
      all $\s,\t\in T(\lambda)$, is a homogeneous anti-isomorphism of $A$.
  \end{enumerate}
  A \textbf{graded cellular algebra} is a graded algebra that has a graded
  cell datum. The basis $\{b_{\s\t}\}$ %|\lambda\in\P\text{ and } \s,\t\in T(\lambda}$
  is a \textbf{graded cellular basis} of~$A$.
\end{Definition}

If we omit the degree assumption (GC$_d$) then we recover Graham and
Lehrer's~\cite{GL} definition of an (ungraded) cellular algebra.

Fix a graded cellular algebra $A$ with graded cellular basis
$\{b_{\s\t}\}$. If $\lambda\in\P$ then the graded \textbf{cell
module} is the $\mz$-module $\Delta^\lambda$ with basis
$\set{b_\t|\t\in T(\lambda)}$ and with $A$-action
$$b_{\t} a=\sum_{\v\in T(\lambda)}r_{\t\v}(a)b_{\v},$$ where the scalars $r_{\t\v}(a)\in \mz$
are the same scalars appearing in (GC$_2)$. One of the key
properties of the graded cell modules is that by
\cite[Lemma~2.7]{HuMathas:GradedCellular} they come equipped with a
homogeneous bilinear form $\<\ ,\ \>$ of degree zero that is
determined by the equation
$$\<b_\t, b_\u\>b_{\v}\equiv b_{\t}b_{\u\v},$$
for $\s, \t, \u, \v\in T(\lambda)$. The \textbf{radical} of this form,
$$\rad \Delta^\lambda=\set{x\in \Delta^\lambda|\<x,y\>=0
            \text{ for all }y\in \Delta^\lambda},$$
is a graded $A$-submodule of $\Delta^\lambda$ so that $L^\lambda=\Delta^\lambda/\rad
\Delta^\lambda$ is a graded $A$-module.

As in \autoref{S:graded}, if $M$ is an $A$-module let $M^\circledast$ be
the (graded) dual of~$M$.

\begin{Theorem}[\protect{\cite[Theorem~2.10]{HuMathas:GradedCellular}}]\label{T:CellularSimples}
  Suppose that $\mz$ is a field and that $A$ is a graded cellular algebra. Then:
  \begin{enumerate}
    \item If $L^\lambda\ne0$, for $\lambda\in\P$, then $L^\lambda$ is an absolutely
      irreducible graded $A$-module and $(L^\lambda)^\circledast\cong L^\lambda$.
    \item $\set{L^\lambda\<k\>|\lambda\in\P, L^\lambda\ne0 \text{ and }k\in\Z}$ is a
      complete set of pairwise non-isomorphic irreducible (graded) $A$-modules.
  \end{enumerate}
\end{Theorem}

Suppose that $\mz=K$ is a field and let $M$ be a (graded) $A$-module
and $L^\mu$ be a graded simple $A$-module, for $\mu\in\P$. Define
\begin{equation}\label{E:multiplicity}
[M:L^\mu]_q = \sum_{d\in\Z}[M:L^\mu\<d\>]q^d
\end{equation}
to be the graded multiplicity of $L^\mu$ in $M$. By the
Jordon-H\"older theorem, $[M:L^\mu]_q$ depends only on~$M$ and
$L^\mu$ and not on the choice of composition series for~$M$.
Moreover, $[M:L^\mu]_q\in\N[q,q^{-1}]$ and
$[M:L^\mu]_{q=1}=[\underline M:\underline L^\mu]$ is the usual
decomposition multiplicity of $\underline{L}^\mu$ in
$\underline{M}$.

\begin{Corollary}[\protect{\cite[Lemma~2.13]{HuMathas:GradedCellular}}]\label{C:triangular}
  Suppose that $\mz$ is a field and that $\lambda,\mu\in\Lambda$ with $L^\mu\ne0$. Then
  $[\Delta^\mu:L^\mu]_q=1$ and $[\Delta^\lambda:L^\mu]_q\ne0$ only if $\lambda\gedom\mu$.
\end{Corollary}

Let $\P_0=\set{\mu\in\P|L^\mu\ne0}$. Then
$\mathbf{D}_A(q)=([\Delta^\lambda:L^\mu]_q)_{\lambda\in\P,\mu\in\P_0}$ is the
\textbf{decomposition matrix} of $A$. For each $\mu\in\P_0$ let $P^\mu$ be the
projective cover of~$L^\mu$ in~$A\Mod$. Then
$\mathbf{C}_A(q)=([P^\lambda:L^\mu]_q)_{\lambda,\mu\in\P_0}$ is the
\textbf{Cartan matrix} of~$A$.

If $M=(m_{ij})$ is a matrix let $M^\trans=(m_{ji})$ be its transpose. We will need
the following fact.

\begin{Corollary}[\protect{Brauer-Humphreys reciprocity~\cite[Theorem 2.17]{HuMathas:GradedCellular}}]%
  \leavevmode\newline\label{C:CartanSymmetric}%
  Suppose that $\mz=K$ is a field, $\lambda\in\P$ and $\mu\in\P_0$. Then
  $P^\mu$ has a
  cell filtration in which $\Delta^{\lambda}$ appears with graded multiplicity
  $(P^\mu:\Delta^\lambda)_q=[\Delta^\lambda:L^\mu]_q$. Consequently,
  $\mathbf{C}_A(q)=\mathbf{D}_A(q)^{\trans}\mathbf{D}_A(q)$ is a symmetric matrix.
\end{Corollary}

Finally, we note the following criterion  for a cellular algebra to be
quasi-hereditary. In particular, this implies that $A\Mod$ is a highest weight
category. The definitions of these objects can be found, for example, in
\cite[Appendix]{Donkin:book}. Alternatively, the reader can take the following
result to be the definition of a (graded split) quasi-hereditary algebra (with a
graded duality).

\begin{Corollary}[\protect{\cite[Remark 3.10]{GL}}]\label{C:quasiH}
  Suppose that $A$ is a graded cellular algebra. Then $A$ is a split
  quasi-hereditary algebra, with standard modules $\set{\Delta^\mu|\mu\in\P}$,
  if and only if $L^\mu\ne0$ for all $\mu\in\P$.
\end{Corollary}

\subsection{Basic algebras and graded Morita
equivalences}\label{S:Basic} Let $\mz=K$ be a field. Recall that a
finite dimensional \textit{ungraded} split $K$-algebra $\underline{B}_0$
is a \textbf{basic algebra} if every irreducible
$\underline{B}_0$-module is one dimensional. It is well-known that, up
to isomorphism, every finite dimensional (ungraded) split $K$-algebra
$\underline{B}$ is Morita equivalent to a unique
basic algebra $\underline{B}_0$. In fact, if
$\{\underline{P}_1,\dots,\underline{P}_z\}$ is a complete set of
pairwise non-isomorphic projective indecomposable
$\underline{B}$-modules then the basic algebra of~$\underline{B}$ is
isomorphic to $\End_{\underline{B}}(\underline{P}_1\oplus\dots\oplus
\underline{P}_r)^\op$.  These facts can be found, for example, in
\cite[\S2.2]{Benson:I}.

We need analogues of these results for graded categories.

A \textbf{graded category} is any category whose objects are finite
dimensional $\Z$-graded $K$-modules and whose morphisms are homogeneous
maps of degree zero, where $K$ is a field. If $\mathcal{C}$ is a graded
category let $\sigma_{\mathcal{C}}$ be the shift functor that sends a
module $M\in\mathcal{C}$ to $M\<1\>$. If $\mathcal{D}$ is another graded
category then a graded functor $F\map{\mathcal{C}}\mathcal{D}$ is any
functor such that $F\circ\sigma_{\mathcal{C}}=\sigma_{\mathcal{D}}\circ
F$.  Similarly, a graded equivalence is an equivalence given by a graded
functor.

Let $A$ and $B$ be two finite dimensional $\Z$-graded $K$-algebras.
Following~\cite[\S5]{GordonGreen:GradedArtin}, the $K$-algebras $A$ and $B$
are \textbf{graded Morita equivalent} if
there is a graded equivalence of graded module categories
$A\Mod\cong B\Mod$. Equivalently, by the results
of~\cite[\S5]{GordonGreen:GradedArtin}, $A$ and $B$
are graded Morita equivalent if and only if there is an (ungraded)
Morita equivalence $\underline{E}:\underline{A}\Mod\cong
\underline{B}\Mod$ and a graded functor $G: A\Mod\rightarrow
B\Mod$ such that the following diagram commutes:
$$\begin{tikzpicture}
     \matrix[matrix of math nodes,row sep=1cm,column sep=16mm]{
       |(O)| A\Mod & |(S)| B\Mod\\
       |(P)| \underline{A}\Mod  &|(H)|\underline{B}\Mod\\
     };
     \draw[->](O) -- node[above]{$G$} (S);
     \draw[->](S) -- node[right]{Forget} (H);
     \draw[->](O) -- node[left]{Forget} (P);
     \draw[->](P) -- node[below]{$\underline{E}$} (H);
    \end{tikzpicture}$$
where the vertical functors are the natural forgetful functors.
Let $\{P_1,\dots,P_z\}$ be a complete set of pairwise non-isomorphic
graded projective indecomposable $A$-modules such that $P_i\not\cong
P_j\<k\>$ for any $i\neq j$ and $k\in\Z$. The \textbf{graded basic algebra} of~$A\Mod$ is the
endomorphism algebra
$$ \AFlat=\ZEnd_{A}\bigl(P_1\oplus\dots\oplus P_r\bigr)^\op.$$
By construction,
%every irreducible $\AFlat$-module is one dimensional so $\AFlat$ is a basic algebra. Moreover,
$\AFlat$ is naturally $\Z$-graded and, on
forgetting the grading, $\underline{A}_\flat$ is the basic algebra of
$\underline{A}$. Unlike the ungraded case, two graded
Morita equivalent graded basic algebras need not be isomorphic as graded
algebras because, for example, we can change the degree shifts on
$P_1,\dots,P_r$. This is discussed in more detail after
\cite[Corollary~5.10]{GordonGreen:GradedArtin}.

\subsection{Schur functors}\label{S:SchurFunctors}
Several places in this paper rely on Auslander's theory of ``Schur
functors'', or quotient functors. We briefly recall how this works in the
graded setting following \cite[\S3.1]{BDK}.

Let $A$ be a finite dimensional graded algebra (with $1$) that is split over a
field $K$ and let $A\Mod$ be the category of finite dimensional graded right
$A$-modules. Suppose that $e\in A$ is a non-zero homogeneous idempotent and consider the subalgebra $eAe$ of~$A$.  Then $eAe$ is a graded
algebra with identity element~$e$. (In all of our applications,~$A$ will be a
quasi-hereditary graded cellular algebra.)

Define functors $\Fun\map{A\Mod}{eAe}\Mod$ and $\Gun\map{eAe\Mod}A\Mod$ by
\begin{equation}\label{E:SchurFunctors}
  \Fun(M)=Me\cong\Hom_A(eA,M)\qquad\text{and}\qquad\Gun(N)=N\otimes_{eAe}eA,
\end{equation}
for $M\in A\Mod$ and $N\in eAe\Mod$, together with the obvious action on morphisms.
Both of these functors are graded because they respect the gradings on
both categories. In general, however, these functors will not be
equivalences between the (graded) module categories of~$A$ and~$eAe$.

To define a graded equivalence between $eAe\Mod$ and a subcategory of $A\Mod$ we need
to work a little harder. Suppose that $M$ is a graded $A$-module and define
$\Oun_e(M)$ to be the largest graded submodule $M'$ of~$M$ such that $\Fun(M')=0$
and define $\Oun^e(M)$ to be the smallest graded submodule $M''$ of~$M$ such
that $\Fun(M/M'')=0$. Any $A$-module (degree preserving) homomorphism
$M\longrightarrow N$ sends $\Oun_e(M)$ to $\Oun_e(N)$ and $\Oun^e(M)$ to
$\Oun^e(N)$, so~$\Oun_e$ and~$\Oun^e$ define graded functors on the category of
$A$-modules.

%\begin{Lemma}[\protect{\cite[Corollary~3.1c]{BDK}}]\label{L:FullyFaithful}
%  Suppose that $M$ and $N$ are graded $A$-modules such that $\Oun^e(M)\cong M$
%  and $\Oun_e(N)=0$.  Then $\ZHom_A(M,N)\cong\ZHom_{eAe}(\Fun M,\Fun N)$.
%\end{Lemma}

Let $\mathcal{A}_e$ be the full subcategory of $A\Mod$ with objects all
graded $A$-modules $M$ such that $\Oun_e(M)=0$ and $\Oun^e(M)=M$. It is
easy to check that any $A$-module homomorphism $M\longrightarrow N$
induces a well-defined map $M/\Oun_e(M)\longrightarrow N/\Oun_e(N)$
so that there is an exact graded functor $$\Hun\map{A\Mod}A\Mod; M\mapsto
M/\Oun_e(M).$$ By \cite[Lemma~3.1a]{BDK}, the functors
$\Hun\circ\Gun\circ\Fun$ and $\Fun\circ\Hun\circ\Gun$ are isomorphic
to the identity functors on $\mathcal{A}_e$ and on $eAe\Mod$
respectively. This implies the following.

\begin{Theorem}[\protect{\cite[Theorem~3.1d]{BDK}}]\label{T:SchurEquiv}
  The restrictions of the functors $\Fun$ and $\Hun\circ\Gun$ induce mutually inverse
  graded equivalences of categories between $\mathcal{A}_e$ and $eAe\Mod$.
\end{Theorem}

In \cite{BDK} this result is proved only for ungraded algebras, however,
the proof there generalizes without change to graded module categories.

\subsection{Koszul algebras}\label{S:Koszul}
In this section we recall the definition of Koszul algebras and the
properties of these algebras that we will need. Throughout this
section we work over a field $K$.

Let $A=\bigoplus_{d\in\Z}A_d$ be a finite dimensional graded $K$-algebra. Then $A$
is \textbf{positively graded} if $A_d=0$ whenever $d<0$. That is, all of the
homogeneous elements of~$A$ have \textit{non-negative} degree.

% \begin{Lemma}[\protect{\cite[Proposition~2.3.1]{BGS:Koszul}}]\label{L:Extensions}
%   Suppose that $A$ is positively graded and that~$A_0$ is semisimple. Then the
%   following are equivalent.
%   \begin{enumerate}
%     \item $\Ext_A^{1}(A_0, A_0\<k\>)\neq 0$ only if $k=1$;
%     \item $A$ is generated by $A_0$ and $A_1$.
%   \end{enumerate}
% \end{Lemma}

Suppose that $M$ is a (graded) $A$-module. A \textbf{linear projective
resolution} of~$M$ is a (graded) projective resolution of~$M$ of the form
$$\cdots\to P^2\to P^1\to P^0 \to M \to 0$$
such that $P^d=P^d_d A$ is generated by its elements of degree~$d$, for $d\ge0$.

\begin{Definition}[\protect{\cite[Definition~1.2.1]{BGS:Koszul}}]%
  \label{D:Koszul}%
  A \textbf{Koszul algebra} is a positively graded algebra
  $A=\bigoplus_{d\ge0}A_d$ such that $A_0$ is semisimple and $A_0$ has a
  linear projective resolution.
\end{Definition}

More generally, if $A$ is a graded algebra then $A\Mod$ is
\textbf{Koszul} if it is graded Morita equivalent to the module category
of a Koszul algebra.  By definition, if the category $A\Mod$ is Koszul
then $A$ is not necessarily a Koszul algebra. For example, by
\autoref{THEOREMC}, in characteristic zero the category $\Sch[n]\Mod$ is
Koszul when $e=0$, however, the graded Schur algebra $\Sch[n]$ is not
usually positively graded so it is not a Koszul algebra in general.

Koszul algebras play a key role in \autoref{S:CatO} where we prove
\autoref{THEOREMC}. For us, one of the most important properties of a
Koszul algebra is that their grading determines the radical and socle
filtrations of certain modules. To make this statement precise, let $A$
be a Koszul algebra and suppose that $M=\bigoplus_{d=a}^zM_d$ is a
finite dimensional graded $A$-module. For each $d\in\Z$ let
$\Gr_dM=\bigoplus_{k\ge d}M_k$. As $A$ is positively graded, $\Gr_{d}M$
is an $A$-submodule of $M$. The \textbf{grading filtration} of~$M$ is the filtration
$$M=\Gr_{a}M\supseteq \Gr_{a+1}M\supseteq\dots\supseteq \Gr_{z}M\supset0.$$
As $A$ is Koszul, $A_0$ is semisimple so the quotients
$\Gr_{d}M/\Gr_{d+1}M$ are semisimple for all $d\in\Z$.

Let $M=\rad^0M\supset \rad^1M\supset \rad^2 M\supset\dots\supset \rad^r M=0$
be the \textbf{radical filtration} of~$M$ so that $\rad^1M=\rad M$ and
$\rad^{i+1}M=\rad(\rad^i M)$ for each $i\geq 1$. Similarly, let
$M=\soc^sM\supset \soc^{s-1} M\supset\dots\supset \soc^1
M\supset0$ be the \textbf{socle filtration} of $M$ where $\soc^1 M=\soc
M$ and $\soc^{i+1}M$ is the inverse image of $\soc(M/\soc^iM)$ under the
natural projection $M\surjection M/\soc^iM$. By \cite[\S4]{NV:graded},
the radical and socle filtrations of finite dimensional modules have the
same length~$r=\ell\ell(M)=s$, which is the \textbf{Loewy length}
of~$M$.

An $A$-module $M$ is \textbf{rigid} if its socle and radical filtrations
coincide. That is, $\rad^rM=\soc^{\ell\ell(M)-r}$, for
$0\le r\le\ell\ell(M)$.  The module $M$ is \textbf{rigidly graded} if $M$ is
rigid and $\Gr_r M=\rad^rM$, for $0\le r\le M$.

For the following result is essentially~\cite[Proposition~2.4.1]{BGS:Koszul},
the proof of which is elementary.

\begin{Proposition}[\protect{\cite[Proposition~2.4.1]{BGS:Koszul}}]\label{P:Rigidity}
  Suppose that $A$ is a Koszul algebra and that $M$ is a finite dimensional
  $A$-module.
  \begin{enumerate}
    \item If $M/\rad M$ is irreducible then the radical filtration of $M$
    coincides with the grading filtration of~$M$, up to constant shift.
    That is, there exists $d\in\Z$ such that $\rad^k M=\Gr_{k+d}M$,
    for all $k$.
    \item If $\soc M$ and $M/\rad M$ are irreducible and $M/\rad M$
    concentrated in degree zero then $M$ is rigidly graded.
  \end{enumerate}
\end{Proposition}

\begin{proof}
  By \cite[Corollary~2.3.3]{BGS:Koszul} any Koszul algebra~$A$ is a
  \textit{quadratic algebra}. That is, $A$ is generated by~$A_0$
  and~$A_1$ subject only to relations in degree~$2$. Therefore, part~(a)
  is a special case of \cite[Proposition~2.4.1]{BGS:Koszul}. In
  turn, part~(a) implies that if both $M/\rad M$ and $\soc M$ are irreducible
  then the radical, socle and grading filtrations of~$M$ coincide up to
  a constant shift. In particular, $M$ is rigid. Finally, if~$M/\rad M$ is
  concentrated in degree zero then the radical and grading filtrations
  of~$M$ coincide exactly, so~$M$ is rigidly graded.
\end{proof}

If $A$ is a Koszul algebra then the
\textbf{Koszul dual} of~$A$ is the algebra
\begin{equation}\label{E:KoszulDual}
  E(A):=\Ext^\bullet_{\underline{A}}(A_0,A_0),
\end{equation}
which we consider as a positively graded algebra under Yoneda extension.
By \cite[Theorem~1.2.5]{BGS:Koszul}, if~$A$ is Koszul then $E(A)$ is a Koszul
and, moreover, $E(E(A))\cong A$, whenever $A_d$ is a finitely generated
as a left $A_0$-module, for $d\ge0$. We use Koszul duality implicitly in
\autoref{T:CatOKoszul} to define the Koszul algebras that we use to prove
\autoref{THEOREMC}.

All of the Koszul algebras that we consider in this paper will be
quasi-hereditary, where there is a strengthening of these ideas.
Following~\cite{AgostonDlabLukacs:QuasiKoszul}, an algebra $A$ is
\textbf{standard Koszul} if it is a positively graded (split)
quasi-hereditary algebra such that all of its standard and costandard
modules have linear projective resolutions.
By~\cite[Theorems 1 and 3]{AgostonDlabLukacs:QuasiKoszul}, any standard Koszul algebra is
Koszul and if $A$ is standard Koszul then so is~$E(A)$. Further, if the simple
$A$-modules are indexed by the poset $(X,\ge)$ then the simple $E(A)$-modules
are indexed by the opposite poset $(X,\le)$.

Suppose that $A$ is a (finite dimensional) standard Koszul basic
$K$-algebra with graded simple modules $\set{L^\mu|\mu\in X}$,
concentrated in degree zero, and graded standard modules
$\set{\Delta^\lambda|\lambda\in X}$ such that the natural map
$\Delta^\lambda\rightarrow L^\lambda$ is homogeneous of degree zero, for
$\lambda\in X$. Let $P^\mu$ be the projective cover of~$L^\mu$ and let
$\iota_\mu$ be the degree $0$ homogeneous idempotent such that $P^\mu=\iota_\mu A$, for $\mu\in
X$. Define two matrices $D_A(q)=(
[\Delta^\lambda:L^\mu]_q)_{\lambda,\mu\in X}$ and $C_A(q)=( \Dim
\iota_\lambda A\iota_\mu )_{\lambda,\mu\in X}$, where the rows and
columns of these matrices are ordered in a way that is compatible with
the partial order on~$X$. Then $D_A(q)$ is the decomposition matrix
of~$A$ and $C_A(q)$ its Cartan matrix.

Let $\set{L_E^\mu|\mu\in X}$ and $\set{\Delta_E^\lambda|\lambda\in X}$
be the graded simple and standard modules of $E(A)$. Then we have
matrices $D_{E(A)}(q)$ and $C_{E(A)}(q)$ as above.

We have not found the next result in the literature, even though we
think it is well-known.

\begin{Lemma}\label{L:KoszuldualDecomp}
  Suppose that $A$ is a standard Koszul $K$-algebra. Then
  $D_A(q)^{-1} = D_{E(A)}(-q)^{\trans}$.
\end{Lemma}

\begin{proof} Without loss of generality we may assume that $A$ and
  $E(A)$ are both basic algebras.
  By standard arguments, $C_A(q)=D_A(q)^\trans D_A(q)$ and
  $C_{E(A)}(q)=D_{E(A)}(q)^\trans D_{E(A)}(q)$; compare with
  \autoref{C:CartanSymmetric}. On the other hand, since $A$ and $E(A)$
  are basic algebras, the matrices $C_A(q)$ and $C_{E(A)}(q)$
  coincide with the \textit{Hilbert polynomials} of $A$ and $E(A)$, respectively,
  as defined in \cite[\S2.11]{BGS:Koszul}. By the \textit{numerical
  condition for Koszulity} given in \cite[Lemma~2.11.1]{BGS:Koszul},
  $C_A(q)C_{E(A)}(-q)^\trans=I_X$ where $I_X$ is
  the $|X|\times|X|$ identity matrix. Expanding,
  $$
  D_A(q)^\trans D_A(q)D_{E(A)}(-q)^\trans D_{E(A)}(-q)=I_X.
  $$
  Hence,
  $D_A(q)D_{E(A)}(-q)^\trans=\(D_{E(A)}(-q)D_A(q)^\trans\)^{-1}$.
  As noted above, the quasi-hereditary structures on~$A$ and~$E(A)$ are
  governed by opposite posets, so the matrix $D_A(q)D_{E(A)}(-q)^\trans$
  is upper unitriangular whereas the matrix
  $\(D_{E(A)}(-q)D_A(q)^\trans\)^{-1}$ is lower unitriangular.
  Therefore,
  $$D_A(q)D_{E(A)}(-q)^\trans=I_X=\(D_{E(A)}(-q)D_A(q)^\trans\)^{-1}.$$
  Hence, $D_{E(A)}(-q)^\trans=D_A(q)^{-1}$ as claimed.
\end{proof}

\section{Cyclotomic Quiver Hecke algebras and combinatorics}
In this chapter we recall the facts about the cyclotomic quiver Hecke algebras
of type~$A$ and the cyclotomic Hecke algebras of type $G(\ell,1,n)$ that
are needed in this paper.

\subsection{Cyclotomic quiver Hecke algebras}\label{S:KLR}
Khovanov and Lauda~\cite{KhovLaud:diagI,KhovLaud:diagII} and
Rouquier~\cite{Rouq:2KM} introduced (cyclotomic) quiver Hecke algebras for
arbitrary oriented quivers.  In this paper we consider mainly the linear quiver
of type~$A_\infty$.

Fix a non-negative integer $n$ and an
integer $e\in\{0,2,3,4\dots\}$. Let $\Gamma_e$ be the quiver with vertex
set $I=\Z/e\Z$ and edges $i\longrightarrow i+1$, for all $i\in I$.
Following \cite[Chapter~1]{Kac},  attach to $\Gamma_e$ the standard Lie
theoretic data of a Cartan matrix $(a_{ij})_{i,j\in I}$, simple roots
$\set{\alpha_i|i\in I}$, fundamental weights $\set{\Lambda_i|i\in I}$,
the positive weight lattice $P^+=\bigoplus_{i\in I}\N\Lambda_i$ and the
positive root lattice $Q^+=\bigoplus_{i\in I}\N\alpha_i$. Let
$(\cdot,\cdot)$ be the usual invariant form, normalised so that
$$(\alpha_i,\alpha_j)=a_{ij}\qquad\text{and}\qquad
          (\Lambda_i,\alpha_j)=\delta_{ij},\qquad\text{for }i,j\in I.$$

Fix, once and for all, a \textbf{multicharge}
$\charge=(\kappa_1,\dots,\kappa_\ell)\in\Z^\ell$ and define
$\Lambda=\Lambda(\charge)=\Lambda_{\bar\kappa_1}+\dots+\Lambda_{\bar\kappa_\ell}$,
where $\bar\kappa=\kappa\pmod e$. Equivalently, $\Lambda$ is the
unique element of $P^+$ such that
\begin{equation}\label{E:multicharge}
(\Lambda,\alpha_i) = \#\set{1\le l\le\ell | \kappa_l\equiv i\pmod e},
     \qquad\text{ for all }i\in I.
\end{equation}
All of the bases for the modules and algebras in this paper depend implicitly
on the choice of~$\charge$ even though the algebras themselves depend only on~$\Lambda$.

Let $\Sym_n$ be the symmetric group on~$n$ letters and let
$s_r=(r,r+1)$, for $1\le r<n$. Then $\{s_1,s_2,\dots,s_{n-1}\}$ is the standard set
of Coxeter generators for~$\Sym_n$. The group $\Sym_n$ acts from
the left on $I^n$ by place permutations. More explicitly, if $1\le
r<n$ and $\bi=(i_1,\dots,i_n)\in I^n$ then
$s_r\cdot \bi=(i_1,\dots,i_{r-1},i_{r+1},i_r,i_{r+2},\dots,i_n)\in I^n$.

Fix $\beta\in Q^+$ with $\sum_{i\in I}(\Lambda_i,\beta)=n$ and let
$$I^\beta=\set{\bi\in I^n|\alpha_{i_1}+\dots+\alpha_{i_n}=\beta}.$$
Then $I^\beta$ is a $\Sym_n$-orbit of $I^n$ and every $\Sym_n$-orbit can be
written uniquely in this way for some $\beta\in Q^+$.

\begin{Definition}\label{D:QuiverRelations}
  Suppose that $n\ge0$, $e\in\{0,2,3,4,\dots\}$ and $\beta\in Q_n^{+}$. Define
 $\R[\beta]$ to be the unital associative $\mz$-algebra with generators
  $$\{\psi_1,\dots,\psi_{n-1}\} \cup
  \{ y_1,\dots,y_n \} \cup \set{e(\bi)|\bi\in I^\beta} $$
  and relations
  \bgroup
      \setlength{\abovedisplayskip}{1pt}
      \setlength{\belowdisplayskip}{1pt}
  \begin{align*}
    y_1^{(\Lambda,\alpha_{i_1})}e(\bi)&=0,
    & e(\bi) e(\bj) &= \delta_{\bi\bj} e(\bi),
    &{\textstyle\sum_{\bi \in I^{\beta}}} e(\bi)&= 1,\\
    y_r e(\bi) &= e(\bi) y_r,
    &\psi_r e(\bi)&= e(s_r{\cdot}\bi) \psi_r,
    &y_r y_s &= y_s y_r,
  \end{align*}
  \begin{align}
    \psi_r y_{r+1} e(\bi)&=(y_r\psi_r+\delta_{i_ri_{r+1}})e(\bi),&
    y_{r+1}\psi_re(\bi)&=(\psi_r y_r+\delta_{i_ri_{r+1}})e(\bi),\label{E:ypsi}\\
    \psi_r y_s  &= y_s \psi_r,&&\text{if }s \neq r,r+1,\label{E:ypsiCommute}\\
    \psi_r \psi_s &= \psi_s \psi_r,&&\text{if }|r-s|>1,\notag
\end{align}
\begin{align*}
  \psi_r^2e(\bi) &= \begin{cases}
       0,&\text{if }i_r = i_{r+1},\\
       (y_{r+1}-y_r)e(\bi),&\text{if  }i_r\rightarrow i_{r+1},\\
       (y_r - y_{r+1})e(\bi),&\text{if }i_r\leftarrow i_{r+1},\\
       (y_{r+1} - y_{r})(y_{r}-y_{r+1}) e(\bi),&\text{if }i_r\rightleftarrows i_{r+1}\\
      e(\bi),&\text{otherwise},
\end{cases}\\
\psi_{r}\psi_{r+1} \psi_{r} e(\bi) &= \begin{cases}
    (\psi_{r+1} \psi_{r} \psi_{r+1} +1)e(\bi),\hspace*{7mm} &\text{if }i_r=i_{r+2}\rightarrow i_{r+1} ,\\
  (\psi_{r+1} \psi_{r} \psi_{r+1} -1)e(\bi), &\text{if }i_r=i_{r+2}\leftarrow i_{r+1},\\
  \rlap{$\big(\psi_{r+1} \psi_{r} \psi_{r+1} +y_r -2y_{r+1}+y_{r+2}\big)e(\bi)$,}\\
           &\text{if }i_r=i_{r+2} \rightleftarrows i_{r+1},\\
  \psi_{r+1} \psi_{r} \psi_{r+1} e(\bi),&\text{otherwise.}
\end{cases}
\end{align*}
\egroup
for $\bi,\bj\in I^\beta$ and all admissible $r$ and $s$.

The \textbf{cyclotomic quiver Hecke algebra}, or \textbf{cyclotomic
Khovanov-Lauda--Rouquier algebra}, of weight $\Lambda$ and type
$\Gamma_e$ is the algebra $\R=\bigoplus_{\beta\in Q_n^{+}}\R[\beta]$.
The algebras $\R=\R(\mz)$, and $\R[\beta]=\R[\beta](\mz)$ for $\beta\in Q^+$
are $\Z$-graded with degree function determined by
$$\deg e(\bi)=0,\qquad \deg y_r=2\qquad\text{and}\qquad
         \deg \psi_s e(\bi)=-a_{i_s,i_{s+1}},$$
for $1\le r\le n$, $1\le s<n$ and $\bi\in I^n$.
\end{Definition}

Inspecting the relations in \autoref{D:QuiverRelations},
there is a unique anti-isomorphism~$\star$ of~$\R$ that fixes
each of the generators of~$\R$. Thus $\star$ is homogeneous of
order~$2$. Hence, by twisting with $\star$ we can define the
graded dual $M^\circledast$ of an $\R$-module
$M^\circledast=\ZHom_\mz(M,\mz)$ as in~\autoref{E:dual}.

In this paper we will mainly be concerned with the special cases when either
$e=0$ or $e>n$. The presentation of $\R$ depends on the orientation
of~$\Gamma_e$, however, different orientations of
$\Gamma_e$ yield isomorphic algebras; see, for example,
\cite[Proposition~3.12]{Rouq:2KM}.

\subsection{Cyclotomic Hecke algebras}
Recall that $\Lambda\in P^+$ and that we have fixed an integer
$e\in\{0,2,3,4,\dots\}$. We now define the `integral' cyclotomic Hecke algebras
$\H$ of type $G(\ell,1,n)$, where $\ell=\sum_{i\in I}(\Lambda,\alpha_i)$ is the
\textbf{level} of $\Lambda$.

Fix an integral domain $\mz$ that contains an element $\xi=\xi(e)$ such that
one of the following holds:
\begin{itemize}
  \item $e>0$ and $\xi$ is a primitive $e$th root of unity in $\mz$.
  \item $e=0$ and $\xi$ is not a root of unity.
  \item $\xi=1$ and $e$ is the characteristic of $\mz$.
\end{itemize}
Define
$\delta_{\xi1}=1$ if $\xi=1$ and $\delta_{\xi1}=0$ otherwise. For $k\in\Z$ set
\begin{equation}\label{E:xi}
\xi^{(k)}=\begin{cases} \xi^k,&\text{if }\xi\ne1,\\
                          k,&\text{if }\xi=1.
                        \end{cases}
\end{equation}
The definition of $\xi=\xi(e)$ above ensures that $\xi^{(i)}=\xi^{(i+e)}$.
Hence,~$\xi^{(i)}$ is well-defined for all $i\in I=\Z/e\Z$.

\begin{Definition}\label{D:HeckeAlgebras}
  The  (integral) \textbf{cyclotomic Hecke algebra} $\H=\H(\mz,\xi)$ of type $G(\ell,1,n)$
  is the unital associative $\mz$-algebra
    with generators $L_1,\dots,L_n$, $T_1,\dots,T_{n-1}$ and relations
    \begin{align*}
      \prod_{i\in I}(L_1-\xi^{(i)})^{(\Lambda,\alpha_i)}&=0,
    & L_rL_t&=L_tL_r,\\
     (T_r+1)(T_r-\xi )&=0,
    &T_rL_r+\delta_{\xi 1}&=L_{r+1}(T_r-\xi +1),\\
    T_sT_{s+1}T_s&=T_{s+1}T_sT_{s+1},\\
    T_rL_t&=L_tT_r,&\text{if }t\ne r,r+1,\\
    T_rT_s&=T_sT_r,&\text{if }|r-s|>1,
  \end{align*}
  where $1\le r<n$, $1\le s<n-1$ and $1\le t\le n$.
\end{Definition}

It is well-known that $\H$ decomposes into a direct sum of
simultaneous generalized eigenspaces for the elements
$L_1,\dots,L_n$~(cf. \cite{Groj:control}). Moreover, the possible
eigenvalues for $L_1,\dots,L_n$ belong to the set
$\set{\xi^{(i)}|i\in I}$. Hence, the generalized eigenspaces for
these elements are indexed by $I^n$. For each $\bi\in I^n$ let
$e(\bi)$ be the corresponding idempotent in~$\H$ (or zero if the
corresponding eigenspace is zero).

\begin{Theorem}[\protect{Brundan-Kleshchev~\cite[Theorem~1.1]{BK:GradedKL}}]\label{T:BKiso}
Suppose that $\mz=K$ is a field, $\xi\in K$ as above, and that
$\Lambda=\Lambda(\charge)$. Then there is an isomorphism of algebras
$\UnR\cong\H$ that sends $e(\bi)\mapsto e(\bi)$, for all $\bi\in
I^n$ and
\begin{align*}
y_r&\mapsto\begin{cases}
  \Sum_{\bi\in I^n}(1 - \xi^{-i_r}L_r)e(\bi),&\text{if }\xi\ne1,\\[5mm]
  \Sum_{\bi\in I^n}(L_r - i_r)e(\bi),&\text{if }\xi=1.
\end{cases}\\
\psi_s&\mapsto \sum_{\bi\in I^n}(T_s+P_s(\bi))Q_s(\bi)^{-1}e(\bi),
\end{align*}
where $P_s(\bi),Q_s(\bi)\in\mz[y_s,y_{s+1}]$, for $1\le r\le n$ and $1\le s<n$.
\end{Theorem}

By \cite[Theorem~1.1]{BK:GradedKL}, the inverse
isomorphism $\H\bijection\UnR$ is determined by
\begin{align}\label{E:InverseBKIsoL}
  L_r&\mapsto\begin{cases}
      \Sum_{\bi\in I^n}\xi^{i_r}(1-y_r)e(\bi),&\text{if }\xi\ne1,\\
      \Sum_{\bi\in I^n}(i_r+y_r)e(\bi),&\text{if }\xi=1.\\
    \end{cases}\\[5pt]\label{E:InverseBKIsoT}
    T_s&\mapsto \sum_{\bi\in I^n}(\psi_sQ_s(\bi)-P_s(\bi))e(\bi),
\end{align}
for $1\le r\le n$ and $1\le s<n$.

Henceforth, we identify the algebras $\UnR$ and $\H$ under this
isomorphism. In particular, we will not distinguish between the
homogeneous generators of $\UnR$ and their images in $\H$ under the
isomorphism of \autoref{T:BKiso}.

Even though we will not distinguish between $\UnR$ and $\H$ we will
usually write~$\R$ when we are working with graded representations
and~$\H$ for ungraded representations.

\subsection{Tableau combinatorics}\label{S:tableaux}
This section sets up the tableau combinatorics that will be used throughout this
paper. Recall that a \textbf{partition} of~$m$ is a weakly decreasing sequence
$\mu=(\mu_1\ge\mu_2\ge\dots)$ of non-negative integers that sum to~$m$. Set
$|\mu|=m$.

% When writing partitions we will
% omit trailing zeros and we use exponential notation for repeated
% parts. So, for example, $(5,4^3,2,1^2)$ represents the partition
% $(5,4,4,4,2,1,1,0,0,0,\dots)$.

A \textbf{multipartition} of~$n$ is an $\ell$-tuple
$\bmu=(\mu^{(1)}|\dots|\mu^{(\ell)})$ of partitions such that
$|\mu^{(1)}|+\dots+|\mu^{(\ell)}|=n$. We \textit{identify} a multipartition with
its diagram
$$\bmu=\set{(r,c,l)|1\le c\le\mu_r^{(l)}\text{ for }r\ge1\text{ and }1\le l\le\ell},$$
which we think of as an $\ell$-tuple of boxes in the plane. For example,
$$(3,2|2,1^2|3,1)=\TriDiagram(3,2|2,1,1|3,1).$$
The partitions $\mu^{(1)},\dots,\mu^{(\ell)}$ are the
\textbf{components} of $\bmu$ and we identify $\mu^{(l)}$ with the
subdiagram
$$\set{(r,c,l)|1\le c\le \mu_r^{(l)} \text{ for } r\ge1}$$
of~$\bmu$. A \textbf{node} is any triple
$A=(r,c,l)\in\N^2\times\{1,2,\dots,\ell\}$. In particular, the
elements of (the diagram of)~$\bmu$ are nodes.

Let $\Parts$ be the set of multipartitions of~$n$. Then $\Parts$ is a poset
under the \textbf{dominance order} $\gedom$ where $\blam\gedom\bmu$, for
multipartitions~$\blam$ and~$\bmu$ of~$n$,  if
$$\sum_{k=1}^{l-1}|\lambda^{(k)}|+\sum_{j=1}^i\lambda^{(l)}_j
     \ge\sum_{k=1}^{l-1}|\mu^{(k)}|+\sum_{j=1}^i\mu^{(l)}_j,$$
for $1\le l\le\ell$ and $i\ge1$. If $\blam\gedom\bmu$ and $\blam\ne\bmu$
then we write $\blam\gdom\bmu$.

Suppose that  $\bmu\in\Parts$ is a multipartition of~$n$. A
\textbf{$\bmu$-tableau} is bijection $\t\map{\bmu}\{1,2,\dots,n\}$. We
think of a $\bmu$-tableau $\t=(\t^{(1)},\dots,\t^{(\ell)})$ as a labelling of
(the diagram of) $\bmu$, where $\t^{(r)}$ is the restriction of $\t$ to
$\mu^{(r)}$. In this way, we talk of the rows, columns and components of a
tableau~$\t$. For example,
$$\Tritab({1,2,3},{4,5}|{6,7},{8},{9}|{10,11,12},{13})\text{ and }
  \Tritab({9,11,13},{10,12}|{5,8},{6},{7}|{1,3,4},{2})$$
are two $(3,2|2,1^2|3,1)$-tableaux. If $\t=(\t^{(1)},\dots,\t^{(\ell)})$ is a
$\bmu$-tableau then define $\Shape(\t)=\bmu$, so that
$\Shape(\t^{(r)})=\mu^{(r)}$, for $1\le r\le\ell$. If $\t^{-1}(k)=(r,c,l)$, then
we set $\comp_{\t}(k)=l$.

A $\bmu$-tableau $\t$ is \textbf{standard} if its entries increase
along the rows and down the columns of each component. For example,
the two tableaux above are standard. If $\t$ is a standard tableau
let $\t\rest k$ be the subtableau of~$\t$ that contains
$1,2,\dots,k$. Then a tableau~$\t$ is standard if and only if
$\Shape(\t\rest k)$ is a multipartition for $1\le k\le n$. The
dominance order induces a partial order on the set of tableaux where
$\s\gedom\t$ if
$$\Shape(\s\rest k)\gedom\Shape(\t\rest k),\qquad\text{ for }1\le k\le n,$$
for $\s\in\Std(\blam)$ and $\t\in\Std(\bmu)$, where $\blam,\bmu\in\Parts$.
Again we write $\s\gdom\t$ if $\s\gedom\t$ and $\s\ne\t$.
Let $\Std(\bmu)$ be the poset of standard
$\bmu$-tableau and set $\SStd(\bmu)=\Std(\bmu)\times\Std(\bmu)$,
$\Std(\Parts)=\bigcup_{\bmu\in\Parts}\Std(\bmu)$ and
$\SStd(\Parts)=\bigcup_{\bmu\in\Parts}\SStd(\bmu)$.

We extend the dominance order to $\SStd(\Parts)$ by declaring that
$(\s,\t)\Gedom(\u,\v)$ if~$\s\gedom\u$ and $\t\gedom\v$. We write $(\s,\t)\Gdom(\u,\v)$
if $(\s,\t)\Gedom(\u,\v)$ and $(\s,\t)\ne(\u,\v)$.

If $\bmu\in\Parts$ let $\bump=({\mu^{(\ell)}}',\dots,{\mu^{(1)}}')\in\Parts$
be the \textbf{conjugate} multipartition, which is obtained from~$\bmu$
by reversing the order of its components and then swapping the rows and
columns in each component. Similarly, the \textbf{conjugate} of
the~$\bmu$-tableau~$\t$ is the $\bump$-tableau~$\t'$ that is obtained
from~$\t$ by reversing the order of its components and then swapping its rows and
columns in each component.  The reader is invited to check that
$\blam\gedom\bmu$ if and only if~$\bmu'\gedom\blam'$ and that
$\s\gedom\t$ if and only if $\t'\gedom\s'$, for $\blam,\bmu\in\Parts$
and for $\s,\t\in\Std(\Parts)$.

Fix a multipartition $\bmu\in\Parts$. Define $\tmu$ to be the unique
standard $\bmu$-tableau such that $\tmu\gedom\t$, for all
$\t\in\Std(\bmu)$. More explicitly, $\tmu$ is the $\bmu$-tableau
that has the numbers $1,2,\dots,n$ entered in order, from left to
right, and then top to bottom, along the rows of the components
$\mu^{(1)},\dots,\mu^{(\ell)}$ of~$\bmu$. Define
$\tlmu=(\t^\bump)'$. By the last paragraph $\tlmu$ is the unique
$\bmu$-tableau such that $\t\gedom\tlmu$, for all
$\tlmu\in\Std(\bmu)$. The numbers $1,2,\dots,n$ are entered in order
down the columns of the components $\mu^{(\ell)},\dots,\mu^{(1)}$
of~$\bmu$. The two tableaux displayed above are $\tmu$ an
$\tlmu$, respectively, for $\bmu=(3,2|2,1^2|3,1)$.

Recall from \autoref{S:KLR} that we have fixed a multicharge
$\charge\in\Z^\ell$. The \textbf{residue} of the node $A=(r,c,l)$ is
$\res(A) = \kappa_l+c-r \pmod e$. (If $e=0$ then we adopt the convention that
$i\equiv i\pmod 0$, for $i\in\Z$.) Thus, $\res(A)\in I$. A node $A$ is an
\textbf{$i$-node} if $\res(A)=i$. If $\t$ is a $\bmu$-tableau
and $1\le k\le n$ then the \textbf{residue} of $k$ in~$\t$ is $\res_\t(k)=\res(A)$,
where $A\in\bmu$ is the unique node such that $\t(A)=k$. The
\textbf{residue sequence} of~$\t$ is
$$\res(\t) = (\res_\t(1),\res_\t(2),\dots,\res_\t(n))\in I^n.$$
As two important special cases we set $\bi^\mu=\res(\tmu)$ and $\bi_\bmu=\res(\tlmu)$, for
$\bmu\in\Parts$.

% All of the main results in this paper rely on the assumption that $e=0$ or
% $e>n$. This assumption enters our arguments through the following elementary
% lemma which we leave to the reader.
%
% \begin{Lemma}\label{L:distinct} Suppose that $e=0$ or $e>n$ and that
% $\bmu=(\mu^{(1)},\dots,\mu^{(\ell)})\in\Parts$ and $1\le l\le\ell$. Let $A\ne
% B$ be nodes which either belong to~$\mu^{(l)}$ or are addable nodes
% of~$\mu^{(l)}$. Then $\res(A)\ne\res(B)$. \end{Lemma}

Following Brundan, Kleshchev and Wang~\cite[Definition.~3.5]{BKW:GradedSpecht}
we now define the degree and codegree of a standard tableau. Suppose that
$\bmu\in\Parts$. A node $A$ is an \textbf{addable node} of~$\bmu$ if
$A\notin\bmu$ and $\bmu\cup\{A\}$ is the (diagram of) a multipartition of~$n+1$.
Similarly, a node $B$ is a \textbf{removable node} of~$\bmu$ if $B\in\bmu$ and
$\bmu\setminus\{B\}$ is a multipartition of~$n-1$. Given any two nodes
$A=(r,c,l), B=(r',c',l')$, say that $B$ is strictly below $A$, or $A$ is
strictly above $B$, if either $l'>l$ or $l'=l$ and $r'>r$. Suppose that $A$ is
an $i$-node and define integers
  \begin{align*}
  d_A(\bmu)&=\#\SetBox{addable $i$-nodes of $\bmu$\\[-2mm] strictly below $A$}
                -\#\SetBox[36]{removable $i$-nodes of $\bmu$\\[-2mm] strictly below $A$},\\
  \intertext{and}
  d^A(\bmu)&=\#\SetBox{addable $i$-nodes of $\bmu$\\[-2mm] strictly above $A$}
                -\#\SetBox[36]{removable $i$-nodes of $\bmu$\\[-2mm] strictly above $A$}.
  \end{align*}
  If $\t$ is a standard $\bmu$-tableau then its \textbf{degree} and
  \textbf{codegree} are defined inductively by setting $\deg\t=0=\codeg\t$ when
  $n=0$. If $n>0$ then define
  $$ \deg\t=\deg\t\rest{(n-1)}+d_A(\bmu)\quad\text{and}\quad
     \codeg\t= \codeg\t\rest{(n-1)}+d^A(\bmu),
  $$
where $A=\t^{-1}(n)$ is the node containing~$n$.

The definitions of the residue, degree and codegree of a tableau all depend on
the choice of multicharge~$\charge$. We write $\res_\t^\charge$,
$\deg^\charge\t$ and $\codeg^\charge\t$ when we want to emphasize this choice.

Fix $\beta\in Q^+$ and set $\Parts[\beta]=\set{\blam\in\Parts|\ilam\in I^\beta}$. The
\textbf{defect} of $\beta$ is the integer
\begin{equation*}\label{E:defect}
  \defect[\Lambda]\beta =(\Lambda,\beta) - \frac12(\beta,\beta).
\end{equation*}
When $\Lambda$ is clear we write $\defect\beta=\defect[\Lambda]\beta$.
The defect of $\beta\in Q^+$ is closely related to the degree and codegree of the
corresponding tableaux.

\begin{Lemma}[\protect{\!\cite[Lemma~3.12]{BKW:GradedSpecht}}]\label{L:defect}
  Suppose that $\beta\in Q^+$ and $\s\in\Std(\bmu)$, for $\bmu\in\Parts[\beta]$. Then
  $\deg\t+\codeg\t=\defect\beta$.
\end{Lemma}

\subsection{Standard homogeneous bases}

We are now ready to define some bases for the cyclotomic quiver
Hecke algebra $\R$. Recall from the last section that $\Sym_n$ is
the symmetric group on~$n$ letters and that
$\{s_1,s_2,\dots,s_{n-1}\}$ is the standard set of Coxeter
generators for~$\Sym_n$. If $w\in\Sym_n$ then the \textbf{length}
of~$w$ is the integer
$$\ell(w)=\min\set{k|w=s_{r_1}\dots s_{r_k}\text{ for some }1\le r_1,\dots,r_k<n}.$$
A \textbf{reduced expression} for $w$ is a word $w=s_{r_1}\dots
s_{r_k}$ such that $k=\ell(w)$. It is a general fact from the theory
of Coxeter groups that any reduced expression for~$w$ can be
transformed into any other reduced expression using just the braid
relations $s_rs_t=s_ts_t$, if $|r-t|>1$, and
$s_rs_{r+1}s_r=s_{r+1}s_rs_{r+1}$, for $1\le r<n-1$.

Hereafter, unless otherwise stated, we fix a reduced expression
$w=s_{r_1}\dots s_{r_k}$ for each element $w\in\Sym_n$, with $1\le
r_1,\dots,r_k<n$. We  define $\psi_w=\psi_{r_1}\dots\psi_{r_k}$.  By
\autoref{D:QuiverRelations}, the generators $\psi_r$, for
$1\le r<n$, do not satisfy the braid relations. Therefore, the
element $\psi_w\in\R$ depends upon our choice of reduced expression
for~$w$.

The symmetric group $\Sym_n$ acts from the right on the set of tableaux by
composition of maps. If $\t\in\Std(\bmu)$ define two permutations~$d(\t)$
and~$d'(\t)$ in~$\Sym_n$ by $\t=\tmu d(\t)$ and $\t=\tlmu d'(\t)$. Conjugating
either of the last two equations shows that $d'(\t)=d(\t')$. Let
$w_\bmu=d(\t_\bmu)$. Then it is easy to check that $w_\bmu=d(\t)d'(\t)^{-1}$ and
$\ell(w_\bmu)=\ell(d(\t))+\ell(d'(\t))$, for all $\t\in\Std(\bmu)$.

Recall from \autoref{S:tableaux} that
$\bi^\bmu=\res(\tmu)$ and that $\bi_\bmu=\res(\tlmu)$.

\begin{Definition}[\protect{\cite[Definitions~4.9, 5.1 and 6.9]{HuMathas:GradedCellular}}]\label{D:psis}
Suppose that $\bmu\in\Parts$. Define non-negative integers $d^\bmu_1,\dots,d^\bmu_n$ and
$d_\bmu^1,\dots,d_\bmu^n$ recursively by requiring that
$$d^\bmu_1+\dots+d^\bmu_k=\deg(\tmu_{\downarrow k})\quad\text{and}\quad
  d_\bmu^1+\dots+d_\bmu^k=\codeg(\tlmu\rest k),$$
  for $1\le k\le n$. Now set $e^\bmu=e(\bi^\bmu)$, $e_\bmu=e(\bi_\bmu)$,
$$ y^\bmu=y_1^{d^\bmu_1}\dots y_n^{d^\bmu_n}\quad\text{and}\quad y_\bmu=y_1^{d_\bmu^1}\dots y_n^{d_\bmu^n}.$$
For a pair of tableaux $(\s,\t)\in\SStd(\bmu)$ define
$$\psi_{\s\t}=\psi_{d(\s)}^\star e^\bmu y^\bmu\psi_{d(\t)}\quad\text{and}\quad
  \psi'_{\s\t}=\psi_{d'(\s)}^\star e_\bmu y_\bmu\psi_{d'(\t)}.$$
\end{Definition}

\begin{Remark}\label{R:warning}
We warn the reader that the element $\psi'_{\s\t}$ is equal to the element
$\psi'_{\s'\t'}$ in the notation of
\cite{HuMathas:GradedCellular,HuMathas:GradedInduction} so care should be taken
when comparing the results in this paper with those
in~\cite{HuMathas:GradedCellular,HuMathas:GradedInduction}. We have changed
notation because \autoref{D:psis} makes several subsequent definitions and
results more intuitive. For example, see \autoref{C:BlockBasis} and
\autoref{P:signedPsi} below.
\end{Remark}

In general, the elements $\psi_{\s\t}$
and $\psi'_{\s\t}$ depend upon the choice of reduced expression that we fixed in
\autoref{D:psis} because $\psi_1,\dots,\psi_{n-1}$ do not satisfy the braid
relations. Similarly, $\psi_{d(\s)^{-1}}$ and $\psi_{d(\s)}^\star$
will generally be different elements of~$\R$.

It follows from \autoref{D:psis} and the relations that if $(\s,\t)\in\SStd(\Parts)$ then
\begin{equation}\label{E:weights}
e(\bi)\psi_{\s\t}e(\bj)=\delta_{\bi,\res(\s)}\delta_{\bj,\res(\t)}\psi_{\s\t}
\quad\text{and}\quad
  e(\bi)\psi'_{\s\t}e(\bj)=\delta_{\bi,\res(\s)}\delta_{\bj,\res(\t)}\psi'_{\s\t},
\end{equation}
for all $\bi,\bj\in I^n$. More importantly we have the following.

\begin{Theorem}[\protect{%
  Hu-Mathas~\cite[Theorems~5.8 and 6.11]{HuMathas:GradedCellular},
  Li~\cite{Li:PhD}}] \label{T:PsiBases}
  Let $\mz$ be an integral domain. Then:
  %such that $e$ is invertible in $\mz$ whenever $e\ne0$ and $e$ is not prime. Then:
  \begin{enumerate}
  \item $\set{\psi_{\s\t}|(\s,\t)\in\SStd(\Parts)}$ is a graded cellular basis
  of $\H$ with weight poset $(\Parts,\gedom)$ and degree function
  $\deg\psi_{\s\t}=\deg\s+\deg\t$.
  \item $\set{\psi'_{\s\t}|(\s,\t)\in\SStd(\Parts)}$ is a graded cellular basis
  of $\H$ with weight poset $(\Parts,\ledom)$ and degree function
  $\deg\psi'_{\s\t}=\codeg\s+\codeg\t$.
  \end{enumerate}
\end{Theorem}

This result was proved under some assumptions on $e$ and $\mz$ in
\cite{HuMathas:GradedCellular}. Li~\cite{Li:PhD} has proved that $\R(\Z)$ is
free as a $\Z$-module with basis~$\{\psi_{\s\t}\}$. This implies
that the results of \cite{HuMathas:GradedCellular} extend to an
arbitrary integral domain because all of the arguments in
\cite{HuMathas:GradedCellular} can now be carried out over $\Z$ using
the embedding $\R[\Z]\hookrightarrow\R[\Q]$. The results involving the
$\psi'$-basis also hold over~$\Z$ in view of \autoref{P:signedPsi}
below.

The $\psi$-basis and the $\psi'$-basis are dual to each other in the
following sense.

\begin{Lemma}[\protect{\cite[Corollary 3.10]{HuMathas:GradedInduction}}]
  \label{L:PsiDual}
  Suppose that $(\s,\t),(\u,\v)\in\SStd(\Parts)$. Then:
  \begin{enumerate}
    \item $\psi_{\s\t}\psi'_{\u\v}\ne0$ only if $\res(\t)=\res(\u)$ and $\u\gedom\t$.
    \item $\psi_{\u\v}'\psi_{\s\t}\ne0$ only if $\res(\s)=\res(\v)$ and $\v\gedom\s$
  \end{enumerate}
\end{Lemma}

We need the following dominance results. Recall from \autoref{S:tableaux}
that $(\s,\t)\Gedom(\u,\v)$ if $\s\gedom\u$ and $\t\gedom\v$.

\begin{Lemma}[\protect{\cite[Corollary~3.11]{HuMathas:GradedInduction}}]\label{L:yAction}
  Suppose that $(\s,\t)\in\SStd(\Parts)$ and $1\le r\le n$. Then
  $$\psi_{\s\t}y_r=\sum_{\substack{(\u,\v)\in\Parts\\(\u,\v)\Gdom(\s,\t)}}a_{\u\v}\psi_{\u\v}
    \quad\text{and}\quad
    \psi'_{\s\t}y_r=\sum_{\substack{(\u,\v)\in\Parts\\(\s,\t)\Gdom(\u,\v)}}b_{\u\v}\psi'_{\u\v},$$
  for some scalars $a_{\u\v},b_{\u\v}\in \mz$.
\end{Lemma}

The next result strengthens \cite[Lemma~5.7]{HuMathas:GradedCellular}.

\begin{Lemma}\label{L:DiffRedExps}
  Suppose that $\psi_{\s\t}$, $\hat\psi_{\s\t}$ and
  $\psi'_{\s\t}$, $\hat\psi'_{\s\t}$ are defined using possibly
  different reduced expressions for $d(\s), d(\t)$ and $d'(\s), d'(\t)$, where
  $\s,\t\in\Std(\blam)$ for some $\blam\in\Parts$. Then
  $$\psi_{\s\t}-\hat\psi_{\s\t} =\sum_{(\u,\v)\Gdom(\s,\t)} s_{\u\v}\psi_{\u\v}
  \quad\text{ and }\quad
    \psi'_{\s\t}-\hat\psi'_{\s\t} =\sum_{(\s,\t)\Gdom(\u,\v)} t_{\u\v}\psi'_{\u\v}
  $$
  where $s_{\u\v}\ne0$ only if $\res(\u)=\res(\s)$,
  $\res(\v)=\res(\t)$ and $\deg\u+\deg\v=\deg\s+\deg\t$ and
  $t_{\u\v}\ne0$ only if $\res(\u)=\res(\s)$,
  $\res(\v)=\res(\t)$ and $\codeg\u+\codeg\v=\codeg\s+\codeg\t$.
\end{Lemma}

\begin{proof}
  By \cite[Theorem~3.9]{HuMathas:GradedInduction} the transition matrices between the
  $\psi$-basis and the (non-homogeneous) standard basis of~$\H$ from
  \cite{DJM:cyc} is triangular with respect to strong dominance partial
  order~$\Gedom$. The same remark applies to the $\hat\psi$-basis, which is
  defined using possibly different choices of reduced expressions.
  Moreover, the transition matrices for the $\psi$-basis and the
  $\hat\psi$-basis have the same elements on the diagonal. Applying this
  result twice to rewrite $\hat\psi_{\s\t}$ in terms of the
  $\psi$-basis, via the standard basis, proves the first statement. The
  second statement can be proved similarly.
\end{proof}

\subsection{The blocks of $\R$}\label{S:blocks}
We now show how \autoref{T:PsiBases} restricts to give a basis for the blocks, or the
indecomposable two-sided ideals, of $\R$. By
\autoref{D:QuiverRelations}, if $\beta\in Q^+$ then $\R[\beta]=e_\beta\R
e_\beta$, where $e_\beta=\sum_{\bi\in I^\beta}e(\bi)$.  Set
$Q^+_n=\set{\beta\in Q^+|e_\beta\ne0\text{ in }\R}$. By
\cite[Theorem~2.11]{LM:AKblocks} and
\cite[Theorem~1]{Brundan:degenCentre}, if $\mz=K$ is a field
and~$\beta\in Q_n^{+}$ then $\R[\beta]$ is a (non-zero) block of~$\R$
and $e_{\beta}$ is a central primitive idempotent.  That is,
\begin{equation}\label{E:Rblocks}
  \R=\bigoplus_{\beta\in Q^+_n}\R[\beta]
\end{equation}
is the decomposition of $\R$ into blocks. \autoref{T:BKiso} implies that
$\UnR[\beta]\cong\H[\beta]$, where $\H[\beta]=e_\beta\H e_\beta$.

Recall that
$\Parts[\beta]=\set{\blam\in\Parts|\ilam\in I^\beta}$. Combining \autoref{T:PsiBases},
\autoref{E:weights} and \autoref{E:Rblocks} we obtain the following.

\begin{Corollary}[\cite{HuMathas:GradedCellular}]\label{C:BlockBasis}
  Suppose that $\mz=K$ is a field and that $\beta\in Q^+_n$. Then
  $$\set{\psi_{\s\t}|\s,\t\in\Std(\blam) \text{ for }\blam\in\Parts[\beta]}
                           \text{ and }
    \set{\psi'_{\s\t}|\s,\t\in\Std(\blam) \text{ for }\blam\in\Parts[\beta]}
  $$
  are graded cellular bases of $\R[\beta]$. In particular, $\R[\beta]$ is a
  graded cellular algebra.
\end{Corollary}

\subsection{Trace forms and graded duality}\label{S:ContDuality}
Recall that a \textbf{trace form} on $A$ is a map $\tau\map A\mz$
such that $\tau(ab)=\tau(ba)$, for all $a,b\in A$. The trace form
$\tau$ is \textbf{non-degenerate} if whenever $a\in A$ is non-zero
then $\tau(ab)\ne0$ for some $b\in A$. An algebra $A$ is a
\textbf{symmetric} algebra if it has a non-degenerate trace form.

\begin{Theorem}[\protect{\cite[Theorem~6.17]{HuMathas:GradedCellular}}]\label{T:trace}
  Suppose that $\beta\in Q^+_n$ and that $\mz=K$ is a field. Then there is a
  non-degenerate homogeneous trace form $\tau_\beta\map{\R[\beta]} K$ of degree
  $-2\defect\beta$ such that $\tau_\beta(\psi_{\s\t}\psi'_{\v\u})\ne0$ only if
  $(\u,\v)\Gedom(\s,\t)$, for $(\s,\t),(\u,\v)\in\SStd(\Parts)$.  Moreover,
  $\tau_\beta(\psi_{\s\t}\psi'_{\t\s})\ne0$, for all $(\s,\t)\in\SStd(\Parts)$.
  Consequently, $\R[\beta]$ is a graded symmetric algebra.
\end{Theorem}

By the results in \autoref{S:cellular} the two cellular bases
$\{\psi_{\s\t}\}$ and $\{\psi'_{\u\v}\}$ both determine cell modules for
$\R$. Suppose that $\bmu\in\Parts$. The \textbf{Specht module}~$S^\bmu$
is the cell module of $\R$ indexed by~$\bmu$ determined by the
$\psi$-basis and the \textbf{dual Specht module} $S_\bmu$ is the cell
module indexed by $\bmu$ determined by the $\psi'$-basis. In more
detail, where we use the notation of~\autoref{S:cellular}, as a
$\mz$-module $S^\bmu$ has homogeneous basis
$\set{\psi_\t|\t\in\Std(\bmu)}$, with $\deg\psi_\t=\deg\t$, and the
$\R$-module structure on~$S^\bmu$ is determined by requiring that for
any $\s\in\Std(\bmu)$ the map
$$S^\bmu\<\deg\s\>\rightarrow\R/\Rmu; \psi_\t\mapsto\psi_{\s\t}+\Rmu,$$
is an $\R$-module isomorphism. Similarly, $S_\bmu$ has homogeneous basis
$\set{\psi'_\t|\t\in\Std(\bmu)}$, with $\deg\psi'_\t=\codeg\t$, and
where the $\R$-action is determined in the exactly same way except that
we use the $\psi'$-basis of~$\R$.

The modules $S^\bmu$ and
$S_\bmu$ are dual to each other in the following sense.

\begin{Proposition}[\protect{\cite[Proposition~6.19]{HuMathas:GradedCellular}}]\label{P:SpechtDuality}
  Suppose that $\bmu\in\Parts[\beta]$, where $\beta\in Q^+_n$. Then
  $S^\bmu\cong S_\bmu^\circledast\<\defect\beta\>$ as graded $\R$-modules.
\end{Proposition}

We warn the reader that the module $S_\bmu$ is denoted $S_{\bmu'}$ in
\cite[\S6]{HuMathas:GradedCellular}. This change in notation is a consequence
of \autoref{R:warning}. The notation for Specht modules and dual Specht
modules in this paper is compatible with~\cite{KMR:UniversalSpecht}.

As in \autoref{S:cellular}, define $D^\bmu=S^\bmu/\rad S^\bmu$. Let
$\Klesh=\set{\bmu\in\Parts|\underline{D}^\bmu\ne0}$ be the set of
\textbf{Kleshchev multipartitions}. By \cite{Ariki:class,BK:GradedDecomp,ArikiJaconLecouvey:ModularBranching}, there
is a recursive description of the Kleshchev multipartitions. Observe
that $\Klesh$ depends on the choice of the $\ell$-tuple $(\overline{\kappa}_1,\cdots,\overline{\kappa}_{\ell})$ and not
just on~$\Lambda$.

In the graded setting the irreducible $\R$-modules were first
constructed by Brundan and Kleshchev~\cite[Theorem~4.11,
Theorem~5.10]{BK:GradedDecomp}. By
\cite[Corollary~5.1]{HuMathas:GradedCellular}, Brundan and Kleshchev's
irreducible modules coincide exactly with the irreducible $\R$-modules
constructed using the $\psi$-basis and the cellular algebra framework of
\autoref{T:CellularSimples}.

\begin{Proposition}[\protect{\cite{BK:GradedDecomp,HuMathas:GradedCellular}}]
  \label{C:RSimples}
  Suppose that $\mz=K$ is a field. Then
  $$\set{D^\bmu\<d\>|\bmu\in\Klesh\text{ and }d\in\Z}$$
  is a complete set of pairwise non-isomorphic irreducible graded $\R$-modules.
\end{Proposition}

When $e=0$ it is well-known that $\Klesh$ is in bijection with the set
of \textbf{FLOTW} multipartitions from~\cite{FLOTW}, which have a
particularly simple description. In fact, these two sets of
multipartitions coincide. As with \autoref{C:RSimples}, the main point is
to match the labelling of the simple modules with the cellular algebra structure
of \autoref{T:CellularSimples}.

\begin{Corollary}\label{C:Restricted} Suppose that $e=0$ or $e>n$,
  $\kappa_1\geq\kappa_2\geq\cdots\geq\kappa_{\ell}$ and $\bmu\in\Parts$. Then
  $\bmu=(\mu^{(1)},\dots,\mu^{(\ell)})\in\Klesh$ if and only if
  $\mu^{(l)}_{r+\kappa_l-\kappa_{l+1}}\le\mu^{(l+1)}_r$,
  for $1\le l<\ell$ and $r\ge1$.
\end{Corollary}

\begin{proof}
  This result is part of the folklore for these algebras when $e=0$.
  As far as we are aware, however, it is not stated explicitly in
  the literature. This said, when $e=0$ the result does follow from
  \cite[Remark~3.4(2) and Theorem~5.3]{BK:GradedDecomp} and it is
  implicit in the papers
  \cite{Vazirani:Parametrisation,Ariki:branching,BK:DegenAriki,BK:HigherSchurWeyl}.

  To prove the result directly, it is enough to observe that if
  $\bmu\in\Parts$ then $\bmu$ has at most one addable or removable
  $i$-node in each component, for any $i\in I$, because $e=0$ or $e>n$.
  Therefore, in view of \cite[(11)]{FLOTW}, up to height~$n$ the crystal
  graphs determining the sets of FLOTW and Kleshchev multipartitions
  coincide exactly, which gives the result.
\end{proof}

\subsection{The sign isomorphism}\label{S:SignDual0}
Following~\cite[\S3.2]{KMR:UniversalSpecht} we now introduce an analogue of
the sign involution of the symmetric groups for the quiver Hecke algebras.
Unlike the case of the symmetric groups, this map is generally not an
automorphism of~$\R$.

In \autoref{S:KLR} we fixed the multicharge
$\charge=(\kappa_1,\dots,\kappa_\ell)\in\Z^\ell$ that
determines~$\Lambda=\Lambda(\charge)$. Define
$\charge'=(-\kappa_\ell,\dots,-\kappa_1)\in\Z^\ell$ and let
$\Lambda'=\Lambda(\charge')$.  Then $\Lambda'\in P^+$. More precisely, if
$\Lambda=\sum_{i\in I}l_i\Lambda_i$, for $l_i\in\N$, then $\Lambda'=\sum_{i\in
I}l_i\Lambda_{-i}$. Similarly, if $\beta=\sum_{i\in I}b_i\alpha_i$, for some
$b_i\in\N$, define $\beta'=\sum_{i\in I}b_i\alpha_{-i}$. Then $\beta'\in Q^+$
and $\defect[\Lambda']\beta'=\defect[\Lambda]\beta$.

As noted in \cite[\S3.2]{KMR:UniversalSpecht}, the relations
of~$\R[\beta]$ given in \autoref{D:QuiverRelations} imply
that there is a unique degree preserving isomorphism of graded algebras
$\sgn\map{\R[\beta]}\Rp[\beta']$ such that
\begin{equation}\label{E:sgn}
   e(\bi)\mapsto e(-\bi), \quad y_r\mapsto -y_r\quad\text{and}\quad
   \psi_s\mapsto-\psi_s,
\end{equation}
for $\bi\in I^{\beta},\ 1\leq r\leq n$, and $1\leq s<n$. The map $\sgn$ induces
a graded equivalence
\begin{equation}\label{E:sgnEquiv}
  \Rp[\beta']\Mod\bijection\R[\beta]\Mod
\end{equation}
in the sense of \autoref{S:Basic}. This equivalence sends an
$\Rp[\beta']$-module $M$ to the $\R[\beta]$-module $M^\sgn$ where $M^\sgn$ is
equal to $M$ as a graded vector space and where the $\R[\beta]$-action on $M^\sgn$ is
given by $m\cdot a=m\,\sgn(a)$, for $a\in\R[\beta]$ and $m\in M^\sgn$.

In \autoref{S:tableaux} we defined residues, degrees and codegrees for
tableaux, all as functions of the multicharge~$\charge$. The same
definitions, but with respect to the multicharge~$\charge'$, give analogous statistics
for $\Rp[\beta']$. To distinguish these definitions from the
previous ones set $\res'=\res^{\charge'}$, $\deg'=\deg^{\charge'}$ and
$\cramped{\codeg'=\codeg^{\charge'}}$. In particular, if $A=(r,c,l)$ is a node
then $\res'(A)=\kappa'_l+c-r\pmod e$.  Residue sequences, degrees
and codegrees are now defined exactly as before using~$\res'$. In this way
the set of multipartitions
$\mathscr{P}^{\Lambda'}_{\beta'}=\set{\bmu\in\Parts|\res'(\tmu)\in I^{\beta'}}$
is attached to $\Rp[\beta']$ in the same way that $\Parts[\beta]$ is attached
to~$\R[\beta]$.

Recall that $d'(\t)\in\Sym_n$ is the permutation determined by $\t=\tlmu d'(\t)$
and that $\t'$ is the tableau that is conjugate to~$\t$.

\begin{Lemma}\label{L:SgnResDegCodeg}
  Suppose that $\beta\in Q^+_n$. Then
  $\mathscr{P}^{\Lambda'}_{\beta'}=\set{\bmu'|\bmu\in\Parts[\beta]}$. Moreover,
  if $\t\in\Std(\Parts[\beta])$ then $\t'\in\Std(\mathscr{P}^{\Lambda'}_{\beta'})$
  and $\res'(\t')=-\res(\t)$, $\deg'\t'=\codeg\t$, $\codeg'\t'=\deg\t$,
  and $d'(\t)=d(\t')$.
\end{Lemma}

\begin{proof}By definition, if $\t\in\Std(\bmu)$ then $\t=\tmu d(\t)$.
  Conjugating shows that $\t'=\t_{\bmu'} d(\t)$ and hence that $d(\t')=d'(\t)$.
  By definition, $\kappa'_l=-\kappa_{\ell-l+1}$ for $1\le l\le\ell$.
  Therefore, if the node $A=(r,c,l)\in\bmu$ has residue $i=\res(A)$ then the
  ``conjugate node'' $A'=(c,r,\ell-l+1)\in\bmu'$ has residue~$\res'(A')=-i$.
  Consequently, $\bmu\in\Parts[\beta]$ if and only if
  $\bmu'\in\mathscr{P}^{\Lambda'}_{\beta'}$ and if $\t\in\Std(\Parts[\beta])$
  then $\res'(\t')=-\res(\t)$. As $\res'(\t')=-\res(\t)$, it is now easy
  to see that conjugation interchanges the definitions of degrees and codegrees
  and hence that $\deg'\t'=\codeg\t$ and $\codeg'\t'=\deg\t$.
\end{proof}

Using the residue and degree functions defined using the
multicharge~$\charge'$, \autoref{D:psis} gives $\psi$ and $\psi'$-bases
for~$\Rp[\beta']$. We emphasize that here and below the bases for
$\Rp[\beta']$ will always be defined using~$\charge'$. The bases
of~$\R[\beta]$ and~$\Rp[\beta']$ are interchanged by the $\sgn$
automorphism. More precisely, we have the following.

\begin{Proposition}\label{P:signedPsi}
  Suppose that $\beta\in Q^+_n$, $\bmu\in\Parts[\beta]$, that $\sgn\map{\R[\beta]}{\Rp[\beta']}$ is the
  sign isomorphism and $(\s,\t)\in\SStd(\bmu)$. Then
  $$\sgn(\psi_{\s\t})=\varepsilon_{\s\t}\psi'_{\s'\t'}\quad\text{and}\quad
    \sgn(\psi'_{\s\t})=\varepsilon'_{\s\t}\psi_{\s'\t'},$$
  where $\varepsilon_{\s\t}=(-1)^{\deg\tmu+\ell(d(\s))+\ell(d(\t))}$
  and $\varepsilon'_{\s\t}=(-1)^{\codeg\tlmu+\ell(d'(\s))+\ell(d'(\t))}$.
\end{Proposition}

\begin{proof}By \autoref{L:SgnResDegCodeg},
  $\deg\tnu=\codeg'\t_{\bnu'}$ for all $\bnu\in\Parts[\beta]$. This
  implies that $\deg(\tmu_{\downarrow k})=\codeg(\t_{\bmu'\downarrow k})$, for
  $1\le k\le n$. Therefore, $\sgn(y^\bmu)=(-1)^{\deg\tmu}y_{\bmu'}$
  by \autoref{D:psis}.  Hence,
  $\psi_{\tmu\tmu}=(-1)^{\deg\tmu}\psi'_{\t_{\bmu'}\t_{\bmu'}}$ since
  $\res(\tmu)=-\res'(\t_{\bmu'})$ by \autoref{L:SgnResDegCodeg}. It
  follows that $\sgn(\psi_{\s\t})=\varepsilon_{\s\t}\psi'_{\s'\t'}$
  since $d'(\s')=d(\s)$ and $d'(\t')=d(\t)$ by
  \autoref{L:SgnResDegCodeg}. That
  $\sgn(\psi'_{\s\t})=\varepsilon'_{\s\t}\psi_{\s'\t'}$ is proved
  similarly.
\end{proof}

The involution $\sgn$ sends the $\psi$-basis of $\R[\beta]$ to the
$\psi'$-basis of~$\Rp$, up to sign. Therefore, $\sgn$ induces an
isomorphism between the Specht modules of $\R[\beta]$, which are
constructed using the $\psi$-basis, and the dual Specht modules of
$\Rp$, which are constructed using the $\psi'$-basis. That is, we have
the following.

\begin{Corollary}[\protect{\cite[Theorem~8.5]{KMR:UniversalSpecht}}]
  \label{C:SignedSpechts}
  Suppose that $\mz$ is an integral domain. Then
  $S^\bmu\cong(S_{\bmu'})^\sgn$ and $S_\bmu\cong(S^{\bmu'})^\sgn.$
\end{Corollary}

Note that no degree shifts are required in \autoref{C:SignedSpechts}
because $\sgn$ is homogeneous. Alternatively, \autoref{L:SgnResDegCodeg}
shows that the combinatorics giving the degrees on both sides of these
isomorphisms agree.

\section{Graded Schur algebras}
In this chapter we introduce the cyclotomic quiver Schur algebras.  We
will show that they are quasi-hereditary graded cellular algebras.
Unless otherwise stated, the following assumption will be in force for
the rest of this paper.

\begin{Assumption}\label{A:StandingAssumption}
  We assume that $e=0$ or $e>n$.
  % and that $\mz$ is an integral domain in which
  %$e$ is invertible if $e\ne0$ and $e$ is not prime.
\end{Assumption}

Unless otherwise stated we work over an arbitrary integral domain~$\mz$.

\subsection{Permutation modules}
Following the recipe in~\cite{DJM:cyc} we will define the graded
cyclotomic Schur algebra to be the algebra of graded
$\R$-endomorphisms of a particular $\R$-module $G_n^\Lambda$. In
this section we introduce and investigate the summands
of ``graded tensor space''~$G_n^\Lambda$, as defined below.

\begin{Definition}\label{D:Gmu}
  Suppose that $\bmu\in\Parts$. Define $G^\bmu$ and $G_\bmu$ to be the $\R$-modules:
  $$G^\bmu=\psi_{\tmu\tmu}\R\<-\deg\tmu\>\quad\text{and}\quad
    G_\bmu=\psi'_{\tlmu\tlmu}\R\<-\codeg\tlmu\>.$$
  Set $G^\Lambda_n=\bigoplus_\bmu G^\bmu$ and $G_\Lambda^n=\bigoplus_\bmu G_\bmu$.
\end{Definition}

The degree shifts appear in \autoref{D:Gmu} because we want $G^\bmu$ to
have a graded Specht filtration in which~$S^\bmu$ has graded multiplicity one
and we want $G_\bmu$ to have a graded dual Specht filtration in which~$S_\bmu$
appears with graded multiplicity one.

The modules $G^\bmu$ and $G_\bmu$ are closely related. To explain this recall
the isomorphism $\sgn\map\R\Rp$ from \autoref{E:sgn}.

\begin{Lemma}\label{L:Gsgn}
  Suppose that $\bmu\in\Parts[\beta]$, for $\beta\in Q^+_n$. Then
  $$G^\bmu\cong(G_{\bmu'})^\sgn \quad\text{and}\quad G_\bmu\cong(G^{\bmu'})^\sgn$$
  as graded $\R$-modules.
\end{Lemma}

\begin{proof}
  This result is a consequence of \autoref{D:Gmu} and
  \autoref{P:signedPsi}. We give the details because applying
  the~$\sgn$ involution is an important theme throughout this paper.
  For the proof, let~$\psi_{\tmu\tmu}\in\R$ and
  $\psi'_{\tlmup\tlmup}\in\Rp$ be the corresponding basis elements of $\R$
  and $\Rp$, respectively. We remind the reader that we have to use the
  definitions of \autoref{S:SignDual0} when defining the $\psi'$-basis of $\Rp$
  and when defining $\Rp$-modules. By definition,
  $G^\bmu=\psi_{\tmu\tmu}\R\<-\deg\tmu\>$
  and $G_{\bmu'}=\psi'_{\tlmup\tlmup}\Rp\<-\codeg'\tlmup\>$.
  By \autoref{L:SgnResDegCodeg}, $\deg\tmu=\codeg'\tlmup$
  so it is enough to show that $\psi_{\tmu\tmu}\R\cong(\psi'_{\tlmup\tlmup}\Rp)^\sgn$
  as $\R$-modules. By \autoref{P:signedPsi},
  $\sgn(\psi_{\tmu\tmu})=\pm\psi'_{\tlmup\tlmup}$, so the algebra isomorphism
  $\sgn\map\R\Rp$ restricts to give an isomorphism of $\mz$-modules
  $\psi_{\tmu\tmu}\R\longrightarrow\psi'_{\tlmup\tlmup}\Rp$. Hence,
  $\psi_{\tmu\tmu}\R\cong(\psi'_{\tlmup\tlmup}\Rp)^\sgn$ as
  $\R$-modules as required.
\end{proof}

As a consequence, any result that we prove for $G^\bmu$ immediately
translates into an equivalent ``sign dual'' result for $G_{\bmu'}$. Our
first aim is to give a basis for these modules. If $\bmu,\blam\in\Parts$
define
\begin{align*}
\Std^\bmu(\blam)&=\set{\s\in\Std(\blam)|\s\gedom\tmu\text{ and }\res(\s)=\imu},\\
\Std_\bmu(\blam)&=\set{\s\in\Std(\blam)|\tlmu\gedom\s\text{ and }\res(\s)=\ilmu},
\end{align*}
and set $\Std^\bmu(\Parts)=\bigcup_{\blam\in\Parts}\Std^\bmu(\blam)$
and $\Std_\bmu(\Parts)=\bigcup_{\blam\in\Parts}\Std_\bmu(\blam)$. We next show that these
sets parameterise bases (or, more accurately, Specht filtrations),
for~$G^\bmu$ and~$G_\bmu$.  Observe that $\t\in\Std^\bmu(\blam)$ if and
only if $\t'\in\Std_{\bmu'}(\blam')$.  This is a combinatorial
reflection of \autoref{L:Gsgn} in these definitions.

\begin{Lemma}\label{L:spanning}
  Suppose that $\bmu\in\Parts$. Then
  \begin{enumerate}
    \item $G^\bmu$ is spanned by
      $\set{\psi_{\tmu\tmu}\psi'_{\u\v}|\u\in\Std^\bmu(\blam), \v\in\Std(\blam)
          \text{ for }\blam\in\Parts}$.
    \item $G_\bmu$ is spanned by
      $\set{\psi'_{\tlmu\tlmu}\psi_{\s\t}\,|\s\in\Std_\bmu(\blam), \t\in\Std(\blam)
          \text{ for }\blam\in\Parts}$.
  \end{enumerate}
\end{Lemma}

\begin{proof}By definition, $e^\bmu y^\bmu=\psi_{\tmu\tmu}$ so, by
  \autoref{T:PsiBases}, $G^\bmu$ is spanned by the elements of the form
  $\psi_{\tmu\tmu}\psi'_{\u\v}$, for $(\u,\v)\in\SStd(\Parts)$.
  By \autoref{L:PsiDual}(a), $\psi_{\tmu\tmu}\psi'_{\u\v}\ne0$ only if
  $\u\gedom\tmu$ so this proves~(a). Part~(b) can be proved by repeating
  this argument with the roles of the $\psi$-basis and $\psi'$-basis
  interchanged, since $G_\bmu$ is spanned by the elements
  $\{\psi'_{\tlmu\tlmu}\psi_{\u\v}\}$. Alternatively, apply \autoref{L:Gsgn}.
\end{proof}

For any positive integer $m\le n$ set $s_{m,m}=1$ and
$\psi_{m,m}=1$. If $1\le r<m$ let
$$s_{r,m}=s_r\dots s_{m-1}\quad \text{ and }\quad\psi_{r,m}=\psi_r\dots\psi_{m-1},$$
and set $s_{m,r}=s_{r,m}^{-1}$ and $\psi_{m,r}=\psi_{r,m}^\star$. To
show that the elements in \autoref{L:spanning} give bases of
$G^\bmu$ and $G_\bmu$ we need the following technical lemma. This
result does not require the assumption that $e=0$ or $e>n$.

\begin{Lemma} \label{L:bump}
  Suppose that $e\in\{0,2,3,4\dots\}$, $\bi\in I^n$ and that there exists an integer $r$, with
  $1\le r\le n$, and non-negative integers $d_r,\dots,d_n,d_{n+1}$ such that
  $$d_r\ge d_s\ge d_t\ge d_{n+1}\text{ whenever }r\le s\le t<n\text{  and }i_r=i_s=i_t.$$
  Then $\psi_{n,r}y_r^{d_r}\dots  y_n^{d_n}e(\bi)\in y_r^{d_{r+1}}\dots y_n^{d_{n+1}}e(s_{n,r}\bi)\R$.
\end{Lemma}

\begin{proof}We argue by downwards induction on $r$. If $r=n$ then $\psi_{n,r}=1$ and there is
  nothing to prove since, by assumption, $d_r\ge d_{n+1}$. Suppose then that $r<n$. We divide the
  proof into two cases.

First suppose that $i_r\neq i_{r+1}$. Then, using \autoref{E:ypsi}, we have that
\begin{align*}
  \psi_{n,r}y_r^{d_r}y_{r+1}^{d_{r+1}}\dots y_n^{d_n}e(\bi)
    &=\psi_{n,r+1}\psi_r y_r^{d_r}y_{r+1}^{d_{r+1}}\dots y_n^{d_n}e(\bi)\\
    &=\psi_{n,r+1} y_{r+1}^{d_r}y_r^{d_{r+1}}\psi_r y_{r+2}^{d_{r+2}}\dots y_n^{d_n}e(\bi)\\
    &=y_r^{d_{r+1}}\psi_{n,r+1} y_{r+1}^{d_r} y_{r+2}^{d_{r+2}}\dots y_n^{d_n}e(s_r\bi)\psi_r.
\end{align*}
Therefore, $\psi_{n,r}y_r^{d_r}y_{r+1}^{d_{r+1}}\dots y_n^{d_n}e(\bi)
               \in y_r^{d_{r+1}}y_{r+1}^{d_{r+2}}\dots y_n^{d_{n+1}}e(s_{n,r}\bi)\R$
by induction because the sequence $s_r\bi$ and the non-negative integers
$d_r,d_{r+2},\dots,d_{n+1}$ satisfy the assumptions of the Lemma.

Now consider the remaining case when $i_r=i_{r+1}$. A quick calculation using
\autoref{E:ypsi} shows that $\psi_r$ commutes with any symmetric polynomial in
$y_r$ and  $y_{r+1}$, so $\psi_r(y_ry_{r+1})=(y_ry_{r+1})\psi_r$. By assumption,
$d_r\geq d_{r+1}\geq d_{n+1}$, so
\begin{align*}
  \psi_{n,r}y_r^{d_r}\dots  y_n^{d_n}e(\bi)
  &=\psi_{n,r+1}\psi_r(y_ry_{r+1})^{d_{r+1}}y_{r+2}^{d_{r+2}}\dots y_n^{d_n} e(\bi)y_r^{d_r-d_{r+1}}\\
  &=\psi_{n,r+1}(y_ry_{r+1})^{d_{r+1}}y_{r+2}^{d_{r+2}}\dots y_n^{d_n} e(s_r\bi)\psi_ry_r^{d_r-d_{r+1}}\\
  &=y_r^{d_{r+1}}\psi_{n,r+1}y_{r+1}^{d_{r+1}}y_{r+2}^{d_{r+2}}\dots y_n^{d_n} e(s_r\bi)\psi_ry_r^{d_r-d_{r+1}}.
\end{align*}
Since $i_r=i_{r+1}$ the sequence $s_r\bi$ and the integers $d_{r+1},\dots,d_n$ again
satisfy the assumptions of the Lemma. Hence, the result follows by induction.
\end{proof}

Before we can give bases for $G^\bmu$ and $G_\bmu$ we need to introduce
a special choice of reduced expression. Recall the definition of
$\psi_{r,n}$ and $\psi_{n,r}$ from above
\autoref{L:bump}, for $1\leq r\leq n$. We remind the reader that, by convention, the
symmetric group $\Sym_n$ acts from the right on $\{1,2,\cdots,n\}$.  It
is well-known and easy to prove that
$$\Sym_n = \bigsqcup_{r=1}^ns_{r,n}\Sym_{n-1}  \qquad\text{(disjoint union)},$$
and that $\ell(s_{r,n}w)=\ell(s_{r,n})+\ell(w)=\ell(w)+n-r$, for all $w\in\Sym_{n-1}$.
Hence, we have the following:

\begin{Lemma}\label{L:reduced}
  Suppose that $w\in\Sym_n$. Then there exist unique integers $r_2,\dots,r_n$, with
  $1\le r_k\le k$, such that $w=s_{r_n,n}\dots s_{r_2,2}$ and
  $\ell(w)=\ell(s_{r_n,n})+\dots+\ell(s_{r_2,2})$.
\end{Lemma}

The factorization $w=s_{r_n,n}\dots s_{r_2,2}$ in \autoref{L:reduced} gives a reduced
expression for~$w$. As a temporary notation, define
$\hat\psi_w=\psi_{r_n,n}\dots\psi_{r_2,2}\in\R$ and if $(\s,\t)\in\SStd(\blam)$ let
$\hat\psi_{\s\t}=\hat\psi_{d(\s)}^\star e^\blam y^\blam\hat\psi_{d(\t)}$ and let
$\hat\psi'_{\s\t}=\hat\psi_{d'(\s)}^\star e_{\blam}y_{\blam}\hat\psi_{d'(\t)}$.

Observe that the choice of reduced expression used to define
$\hat\psi_{\s\t}$ is compatible with the natural embeddings
$\Sym_m\hookrightarrow\Sym_n$, for $1\le m\le n$. More precisely,
if~$n$ appears in~$\t$ in the same position as~$r$ appears in $\tmu$
then $d(\t)=s_{r,n}d(\t_{\downarrow_{n-1}})$ and
$\ell(d(\t))=n-r+\ell(d(\t_{\downarrow_{n-1}}))$. Consequently,
$\hat\psi_{d(\t)}=\psi_{r,n}\hat\psi_{\t_{\downarrow_{n-1}}}$.
Similarly, if~$n$ appears in~$\t$ in the same position as~$r$
appears in $\tlmu$ then $d'(\t)=s_{r,n}d'(\t_{\downarrow_{n-1}})$
and $\ell(d'(\t))=n-r+\ell(d'(\t_{\downarrow_{n-1}}))$ so that
$\hat\psi'_{d'(\t)}=\psi_{r,n}\hat\psi'_{d'(\t_{\downarrow_{n-1}})}$.

The next Lemma makes heavy use of \autoref{A:StandingAssumption}.

\begin{Lemma}\label{L:inclusion}
  Suppose that $e=0$ or $e>n$ and $\s\in\Std^\bmu(\blam)$ and
  $\u\in\Std_\bmu(\blam)$, for some $\blam\in\Parts$. Then
  $\hat\psi_{\s\tlam}\in G^\bmu$ and $\hat\psi'_{\u\tllam}\in G_\bmu$.
\end{Lemma}

\begin{proof} By \autoref{P:signedPsi} and \autoref{L:Gsgn} both
  statements are equivalent so we prove only that $\hat\psi_{\s\tlam}\in
  G^\bmu$.

  We argue by induction on~$n$. If $n=1$ then $\s=\tlam$ so that
  $\hat\psi_{\s\tlam}=e^{\blam}y^{\blam}=e^\bmu y_1^{\deg\s}$.
  Similarly, $\hat\psi_{\tmu\tmu}=e^\bmu y_1^{\deg\tmu}$.  Now
  $\deg\s\ge\deg\tmu$, since by assumption $\s\gedom\tmu$ and $n=1$, so
  $\hat\psi_{\s\tlam}=y_1^{\deg\s-\deg\tmu}\hat\psi_{\tmu\tmu}\in
  G^\bmu$ as claimed.

  Now assume that $n>1$.  Let $\s_\downarrow=\s\rest{(n-1)}$,
  $\blam^{\s_\downarrow}=\Shape(\s_\downarrow)$, and
  $\bmu_\downarrow=\Shape(\tmu_{\downarrow(n-1)})$. Then
  $\s_\downarrow\in\Std^{\bmu_\downarrow}(\blam^{\s_\downarrow})$.
  Let~$r$ be the integer such that~$r$ appears in the same position
  in~$\tlam$ as $n$ does in~$\s$. Let
  $\bnu=\Shape(\s_{\downarrow(r-1)})$. By definition,
  $\hat\psi_{d(\s)}=\psi_{r,n}\hat\psi_{\s_\downarrow}$ so, recalling
  the definition of the integers $d^\blam_1,\dots,d^\blam_n$ from
  \autoref{D:psis}, we have
  $$\hat\psi_{\s\tlam}=\hat\psi_{d(\s)}^\star y^\blam e^\blam
        =\hat\psi_{\s_{\downarrow}}^\star\psi_{n,r} y^\bnu y_r^{d^\blam_r}\dots y_n^{d^\blam_n}e^\blam
        =\hat\psi_{\s_{\downarrow}}^\star y^\bnu \psi_{n,r} y_r^{d^\blam_r}\dots y_n^{d^\blam_n}e^\blam,$$
  where the last equality follows because if $r\le j<n$ then $\psi_j$ commutes
  with~$y^\bnu$ by~\autoref{E:ypsi}. We want to apply \autoref{L:bump} to the sequence
  $d_r=d^\blam_r,\dots,d_n=d^\blam_n, d_{n+1}=d^\bmu_n$, so we have to check that
  $d^\blam_r\ge d^\blam_s\ge d^\blam_t\ge d^\bmu_n$ whenever there exist $s$
  and $t$ such that $r\le s\le t<n$ and $i^\blam_r=i^\blam_s=i^\blam_t$.
  Suppose then, if possible, that $r\le s\le t<n$ and
  $i^\blam_r=i^\blam_s=i^\blam_t$.  If $\bsig\in\Parts$ then, because $e=0$ or
  $e>n$, each component of~$\bsig$ contains at most one addable or
  removable $i$-node, for all $i\in I$. Therefore,
  $$d^\bsig_m = \#\set{1\le l\le\ell|l>\comp_{\t^\bsig}(m) \text{ and }
                        \kappa_l\equiv i^\bsig_m\pmod e},$$
  for $1\le m\le n$.  (Here, as usual, $\bi^\bsig=\res(\t^\bsig)\in I^n$.)
  By assumption, $r\le s\le t<n$ so
  $\comp_{\tlam}(r)\le\comp_{\tlam}(s)\le\comp_{\tlam}(t)$ and consequently
  $d^\blam_r\ge d^\blam_s\ge d^\blam_t$ since $i^{\blam}_r=i^\blam_s=i^\blam_t$.
  Moreover, $\comp_{\tlam}(t)\le\comp_{\tmu}(n)$ since $\s\gedom\tmu$
  and $i^\blam_t=i^\s_n=i^\bmu_n$, so that $d^\blam_t\ge d^\bmu_n$.
  Therefore, $d^\blam_r\ge d^\blam_s\ge d^\blam_t\ge d^\bmu_n$ whenever
  $r\le s\le t<n$ and $i^\blam_r=i^\blam_s=i^\blam_t$. Consequently,
  \autoref{L:bump} applies and we deduce that
  $$\hat\psi_{\s\tlam}=\hat\psi_{d(\s)}^\star y^\blam e^\blam
           =\hat\psi_{\s_{\downarrow}}^\star y^\bnu y_r^{d^\blam_{r+1}}\dots
               y_{n-1}^{d^\blam_n}y_n^{d_n^\bmu}e(s_{n,r}\ilam)h,$$
  for some $h\in\R$.
  Now $y^{\blam^{\s_\downarrow}}
           =y^\bnu y_r^{d^\blam_{r+1}}\dots y_{n-1}^{d^\blam_{n}}$, by
  definition, and $y_n$ commutes with $\hat\psi_{\s{\downarrow}}$ by
  \autoref{E:ypsiCommute}.  Therefore, by induction, there
  exists $h'\in\R$ such that
   $$\hat\psi_{\s\tlam}=y_n^{d^\bmu_n}\hat\psi_{\s_{\downarrow}}^\star
         y^{\blam^{\s_{\downarrow}}} e(\bi^{\blam^{\s_\downarrow}}\vee i^\bmu_n)h
      =y_n^{d^\bmu_n}e^\bmu y_1^{d^\bmu_1}\dots y_{n-1}^{d^\bmu_{n-1}}h'h\\
     \in e^\bmu y^\bmu\R.
   $$
   This completes the proof of the Lemma.
\end{proof}

%\begin{Remark}
  If we drop \autoref{A:StandingAssumption} then it is easy to construct examples
  where \autoref{C:inclusion} fails when $0<e\le n$.
%\end{Remark}

\begin{Corollary}\label{C:inclusion}
  Suppose that $e=0$ or $e>n$ and $\s\in\Std^\bmu(\blam)$ and
  $\u\in\Std_\bmu(\blam)$, for $\blam\in\Parts$. Then for any $\t\in\Std(\blam)$, $\psi_{\s\t}\in G^\bmu$ and
  $\psi'_{\u{\t}}\in G_\bmu$.
\end{Corollary}

\begin{proof}
  We show only that $\psi_{\s\t}\in G^\bmu$. If $\psi_{\s\t}=\hat\psi_{\s\t}$
  then the result follows by \autoref{L:inclusion}. Otherwise, by \autoref{L:DiffRedExps},
  there exist $s_{\u\v}\in \mz$ such that
  $$\psi_{\s\t}=\hat\psi_{\s\t}
  +\sum_{\substack{(\u,\v)\in\SStd(\Parts)\\(\u,\v)\Gdom(\s,\t)}}s_{\u\v}\psi_{\u\v},$$
  where $s_{\u\v}\ne0$ only if $\res(\u)=\res(\s)$. Consequently, if $s_{\u\v}\ne0$
  then $\u\gedom\s\gedom\tmu$ and $\v\gedom\t$ so that $\u\in\Std^\bmu(\bnu)$,
  for some $\bnu\gedom\blam$. By induction on dominance,
  $\psi_{\u\v}=\psi_{\u\tnu}\psi_{d(\v)}$ belongs to~$G^\bmu$ whenever
  $(\u,\u)\Gdom(\s,\t)$. Moreover, $\hat\psi_{\s\t}\in G^\bmu$ by
  \autoref{L:inclusion}. Hence, $\psi_{\s\t}\in G^\bmu$ as we wanted to
  show.
\end{proof}

We can now give bases for $G^\bmu$ and $G_\bmu$. Almost everything in this paper
relies on the next result.

\begin{Theorem}\label{T:GBasis}
  Suppose that $\bmu\in\Parts$. Then
  \begin{enumerate}
    \item $\set{\psi_{\s\t}\,|\s\in\Std^\bmu(\bnu)\text{ and }\t\in\Std(\bnu),
        \text{for }\bnu\in\Parts}$ is a basis of $G^\bmu$.
    \item $\set{\psi'_{\u\v}|\u\in\Std_\bmu(\bnu)\text{ and }\v\in\Std(\bnu),
        \text{for }\bnu\in\Parts}$ is a basis of $G_\bmu$.
  \end{enumerate}
\end{Theorem}

\begin{proof}Parts (a) and (b) are equivalent by \autoref{L:Gsgn} and
  \autoref{P:signedPsi}, so it is enough to prove~(a). Suppose first that
  $\mz=K$ is a field.  By
  \autoref{C:inclusion}, $\psi_{\s\t}\in G^\bmu$ whenever $\s\in\Std^\bmu(\bnu)$
  and $\t\in\Std(\bnu)$, for some multipartition $\bnu\in\Parts$. Therefore, by
  \autoref{T:PsiBases},
  $$\dim_K G^\bmu \ge \sum_{\u\in\Std^\bmu(\bnu)}\#\Std(\bnu).$$
  On the other hand, by \autoref{L:spanning} the dimension of $G^\bmu$ is at
  most the number on the right hand side. Hence, the set in the statement of the
  theorem is a basis of~$G^\bmu$, so that the Lemma holds over any field~$K$.

  To prove the proposition when $\mz$ is not a field
  % by \autoref{A:StandingAssumption}
  it suffices to consider the case where
  $\mz=\Z$.
  % $\mz=\Z$ if $e=0$ or $e$ is a prime; or $\mz=\Z[e^{-1}]$ if $e>0$ and~$e$ is
  % not a prime. In these cases, $\mz$ is always a principal ideal domain.
  Let $G$ be the free $\Z$-submodule of $G^\bmu$ spanned by the basis elements
  $\set{\psi_{\s\t}|\s\in\Std^\bmu(\bnu), \t\in\Std(\bnu), \text{ for }\bnu\in\Parts}$.
  Then $G$ is a pure $\Z$-submodule of $\R$, by \autoref{T:PsiBases}, and hence
  a pure $\Z$-submodule of $G^{\bmu}$. As a result, $G$ is a direct summand
  of $G^{\bmu}$ as a $\Z$-module. Therefore, there is a short exact sequence
  of $\Z$-modules
  $$
  0\longrightarrow G\longrightarrow G^\bmu\longrightarrow G^\bmu/G\longrightarrow 0,
  $$
  that splits as sequence of $\Z$-modules. Therefore, for every field $K$
  there is an exact sequence
  $$
  %G\otimes_{\mz}K\longrightarrow G^\bmu\otimes_{\mz}K\rightarrow G^\bmu/G\otimes_{\mz}K\rightarrow 0.
  0\longrightarrow G\otimes_{\Z}K\longrightarrow G^\bmu\otimes_{\Z}K\longrightarrow G^\bmu/G\otimes_{\Z}K\longrightarrow 0.
  $$
  By \autoref{L:spanning}, $\dim G^\bmu\otimes_{\Z}K\leq \dim G\otimes_{\Z}K$.
  Hence the first homomorphism in the last exact sequence must be an
  isomorphism. It follows that $G^\bmu/G\otimes_{\Z}K=0$ for any field $K$ that
  is an $\Z$-algebra. Applying Nakayama's Lemma (see, for
  example,~\cite[Proposition 3.8]{AtiyahMacdonald}), $G^\bmu/G=0$. That is,
  $G^\bmu=G$. Hence, elements in the statement of the theorem are a
  basis for~$G^\bmu$ as required.
\end{proof}

\autoref{T:GBasis} has several useful corollaries. We first note that it gives
explicit formulae for the graded dimensions of these two modules:
\begin{align*}
  \Dim G^\bmu&=q^{-\deg\tmu}\Sum_{\s\in\Std^\bmu(\bnu)}
                            \sum_{\t\in\Std(\bnu)}q^{\deg\s+\deg\t},\\
  \Dim G_\bmu&=q^{-\codeg\tlmu}\Sum_{\u\in\Std_\bmu(\bnu)}
                          \sum_{\v\in\Std(\bnu)}q^{\codeg\u+\codeg\v}.
\end{align*}

\begin{Corollary}\label{C:intersection}
  Suppose that $\bmu,\blam\in\Parts$. Then
  $$\set{\psi_{\s\t}|\s\in\Std^\bmu(\bnu)\text{ and }\t\in\Std^\blam(\bnu),
             \text{for }\bnu\in\Parts}$$
  is a basis of $G^\bmu\cap(G^\blam)^\star$ and
  $$\set{\psi'_{\s\t}|\s\in\Std_\bmu(\bnu)\text{ and }\t\in\Std_\blam(\bnu),
             \text{for }\bnu\in\Parts}$$
  is a basis of $G_\bmu\cap(G_\blam)^\star$.
\end{Corollary}

\begin{proof}Suppose that $a\in G^\bmu\cap(G^\blam)^\star$ and write
  $a=\sum_{(\s,\t)\in\Parts}r_{\s\t}\psi_{\s\t}$, for $r_{\s\t}\in \mz$ and
  $(\s,\t)\in\Std(\Parts)$. Then $r_{\s\t}\ne0$ only if $\s\in\Std^\bmu(\Parts)$ by
  \autoref{T:GBasis}. Similarly, since $a^\star\in G^\blam$ we see that
  $r_{\s\t}\ne0$ only if $\t\in\Std^\blam(\Parts)$. Moreover,
  if~$\s\in\Std^\bmu(\bnu)$ and $\t\in\Std^\blam(\bnu)$ then $\psi_{\s\t}\in
  G^\bmu\cap(G^\blam)^\star$ by two more applications of \autoref{T:GBasis}.
  This proves the first claim. The second statement follows similarly.
\end{proof}

\begin{Corollary}\label{C:GSpecht}
  Suppose that $\bmu\in\Parts$.
  \begin{enumerate}
    \item Write $\Std^\bmu(\Parts)=\{\s_1,\dots,\s_m\}$, ordered so that $i\le j$
      whenever $\s_i\gedom\s_j$ and set $\bnu^i=\Shape(\s_i)$, for $1\le i,j\le m$.
      Then~$G^\bmu$ has a (graded) Specht filtration
    $$G^\bmu=G^{m} \ge G^{m-1}\ge\dots\ge G^1\ge G^{0}=0$$
    such that $G^i/G^{i-1}\cong S^{\bnu^i}\<\deg\s_i-\deg\tmu\>$, for $1\le i\le m$.
  \item Write $\Std_\bmu(\Parts)=\{\u_1,\dots,\u_l\}$, ordered so that $i\ge j$ whenever
    $\u_i\gedom\u_j$ and set $\bnu_i=\Shape(\u_i)$, for $1\le i,j\le l$. Then~$G_\bmu$
    has a (graded) dual Specht filtration
    $$G_\bmu=G_l \ge G_{l-1}\ge\dots\ge G_1\ge G_0=0$$
    such that $G_i/G_{i-1}\cong S_{\bnu_i}\<\codeg\u_i-\codeg\t_{\bmu}\>$, for $1\le i\le l$.
  \end{enumerate}
\end{Corollary}

\begin{proof} Suppose that $1\le i\le m$. Define $G^i$ to be the $\mz$-submodule of $G^\bmu$
  spanned by
  $$\set{\psi_{\s_j\t}|1\le j\le i\text{ and }\t\in\Std(\bnu^j)}.$$ Then $G^i$
  is a submodule of~$G^\bmu$ by \autoref{T:GBasis} and
  (GC$_2$) of \autoref{D:cellular}. Finally, $G^i/G^{i-1}\cong
  S^{\bnu^i}\<\deg\s_i-\deg\tmu\>$ by the construction of the cell modules given in
  \autoref{S:cellular}. More precisely, recalling that the Specht module $S^{\bnu_i}$
  has basis $\set{\psi_\t|\t\in\Std(\bnu^i)}$, the isomorphism is given by
  $\psi_\t\mapsto \psi_{\s_i\t} + G^{i-1}$, for all $\t\in\Std(\bnu^i)$. This map
  has degree $\deg\s_i-\deg\tmu$ because $\psi_{\s_i\t}$ has degree
  $\deg\s_i+\deg\t-\deg\tmu$ when considered as an element of
  $G^\bmu=\psi_{\tmu\tmu}\R\<-\deg\tmu\>$. The proof of~(b) is almost identical.
\end{proof}

In particular, note that $S^\bmu$ is a quotient of $G^\bmu$ and that $S_\bmu$ is a
quotient of~$G_\bmu$.

\begin{Corollary}\label{C:GDualBasis}
  Suppose that $\bmu\in\Parts$. Then:
  \begin{enumerate}
    \item $\set{\psi_{\tmu\tmu}\psi'_{\u\v}|\u\in\Std^\bmu(\bnu)\text{ and }\t\in\Std(\bnu),
                   \text{for }\bnu\in\Parts}$ is a basis of $G^\bmu$.
    \item $\set{\psi'_{\tlmu\tlmu}\psi_{\s\t}\,|\s\in\Std_\bmu(\bnu)\text{ and }\t\in\Std(\bnu),
                   \text{for }\bnu\in\Parts}$ is a basis of $G_\bmu$.
  \end{enumerate}
\end{Corollary}

\begin{proof}
  By \autoref{L:spanning} and \autoref{T:PsiBases}(b) the elements in (a) span $G^\bmu$, so it remains to
  show that they are linearly independent. This is a direct consequence of \autoref{T:GBasis}. The proof
  of~(b) is similar.
\end{proof}

\begin{Corollary}\label{C:GDualSpecht}
  Suppose that $\bmu\in\Parts[\beta]$. Using the notation of \autoref{C:GSpecht}:
  \begin{enumerate}
    \item $G^\bmu$ has a dual Specht filtration
      $G^\bmu=H_0 \ge H_{1}\ge\dots\ge H_{m-1}\ge H_m=0$
    such that $H_i/H_{i+1}\cong S_{\bnu^{i+1}}\<\deg\tmu+\codeg\s_{i+1}\>$, for $0\le i<m$.
  \item $G_\bmu$ has a Specht filtration
    $G_\bmu=H^0 \ge H^{1}\ge\dots\ge H^{l-1}\ge H^{l}=0$
    such that $H^i/H^{i+1}\cong S^{\bnu_{i+1}}\<\codeg\tlmu+\deg\u_{i+1}\>$, for $0\le i<l$.
  \end{enumerate}
\end{Corollary}

\begin{proof}
  We prove only~(b). Part~(a) can be proved in a similar way. Mirroring the proof of
  \autoref{C:GSpecht}, define $H^i$ to be the $\mz$-submodule of $G_\bmu$
  spanned by the elements
  $$\set{\psi'_{\tlmu\tlmu}\psi_{\u_j\t}|\t\in\Std(\bnu_j)
              \text{ and } 1+i\leq j\leq l}.$$
  This is an $\R$-submodule of $G_\bmu$ by \autoref{T:PsiBases}
  and (GC$_2$) of \autoref{D:cellular}. As in the proof of
  \autoref{C:GSpecht} it is easy to verify that
  $H^i/H^{i+1}\cong S^{\bnu^{i+1}}\<\codeg\tlmu+\deg\u_{i+1}\>$; compare with
  \cite[Corollaries 3.11, 3.12]{HuMathas:GradedInduction}. The degree shift is just
  the difference of the degrees of the basis elements of $S^{\bnu^{i+1}}$ and the
  degrees of the elements $\psi'_{\tlmu\tlmu}\psi_{\u_{i+1}\t}$.
\end{proof}

Recall from~\autoref{E:Rblocks} that $\R=\bigoplus_\beta\R[\beta]$ and that
$\R[\beta]$ carries a non-degenerate homogeneous trace form~$\tau_\beta$ of
degree $-2\defect\beta$ by \autoref{T:trace}. The following argument is
lifted from \cite[Proposition~5.13]{M:tilting}.

\begin{Theorem}\label{T:GSelfDual} Suppose that $\mz=K$ is a field and that
   $\bmu\in\Parts[\beta]$, for $\beta\in Q^+_n$. Then, as $\R$-modules,
   $$G^\bmu\cong(G^\bmu)^\circledast\<2\defect\beta\>\quad\text{and}\quad
     G_\bmu\cong(G_\bmu)^\circledast\<2\defect\beta\>.$$
\end{Theorem}

\begin{proof}Both isomorphisms can be proved similarly, so we consider
  only the first one. Using \autoref{T:GBasis} and
  \autoref{C:GDualBasis}, define a pairing
  $G^\bmu\times G^\bmu\longrightarrow \mz$ by
  $$\<\psi_{\s\t},\psi_{\tmu\tmu}\psi'_{\u\v}\>_\bmu
                         =\tau_\beta(\psi_{\s\t}\psi'_{\v\u}),$$
  for all $\s\in\Std^\bmu(\blam)$, $\t\in\Std(\blam)$, $\u\in\Std^\bmu(\bnu)$,
  $\v\in\Std(\bnu)$, for some $\blam,\bnu\in\Parts[\beta]$. By
  \autoref{T:trace}, $\tau_\beta(\psi_{\s\t}\psi'_{\t\s})\ne0$ and
  $\tau_\beta(\psi_{\s\t}\psi'_{\v\u})\ne0$ only if $(\u,\v)\Gedom(\s,\t)$.
  Therefore, the Gram matrix of $\<\ ,\ \>_\bmu$ is upper triangular with non-zero
  elements on the diagonal, so $\<\ ,\ \>_\bmu$ is non-degenerate. Recalling
  the degree shift in the definition of~$G^\bmu$ from \autoref{D:Gmu}, it
  is easy to check that $\<\ ,\ \>_\bmu$ is a homogeneous bilinear map of degree
  $-2\defect\beta$.

  We claim that $\<\ ,\ \>_\bmu$ is associative in the sense that
  $$\<\psi_{\s\t}h,\psi_{\tmu\tmu}\psi'_{\u\v}\>_\bmu
      =\<\psi_{\s\t},\psi_{\tmu\tmu}\psi'_{\u\v}h^\star\>_\bmu,$$
  for all $h\in\R$ and all $(\s,\t)$ and $(\u,\v)$ as above.  Write
  $\psi'_{\u\v}h^\star=\sum r_{\a\b}\psi'_{\a\b}$, where in the sum
  $(\a,\b)\in\SStd(\Parts[\beta])$ and $r_{\a\b}\in \mz$. Then the left hand side is equal to
  \begin{align*}
     \<\psi_{\s\t}h,\psi_{\tmu\tmu}\psi'_{\u\v}\>_\bmu
        &=\tau_\beta(\psi_{\s\t}h\psi'_{\v\u})
         =\sum_{(\a,\b)\in\SStd(\Parts[\beta])}
            r_{\a\b}\tau_\beta(\psi_{\s\t}\psi'_{\b\a}).\\
  \intertext{Now $\tau_\beta$ is a trace form, so
           $\tau_\beta(\psi_{\s\t}\psi'_{\b\a})=\tau_\beta(\psi'_{\b\a}\psi_{\s\t})$
           is non-zero only if $\a\gedom\s$ and $\res(\a)=\res(\s)$ by
           \autoref{L:PsiDual}, so that $\a\in\Std^\bmu(\Parts[\beta])$. Consequently,}
     \<\psi_{\s\t}h,\psi_{\tmu\tmu}\psi'_{\u\v}\>_\bmu
        &=\sum_{\substack{\a\in\Std^\bmu(\bnu),\b\in\Std(\bnu)\\\bnu\in\Parts[\beta]}}
                  r_{\a\b}\tau_\beta(\psi_{\s\t}\psi'_{\b\a})\\
        &=\sum_{\substack{\a\in\Std^\bmu(\bnu),\b\in\Std(\bnu)\\\bnu\in\Parts[\beta]}}
                   r_{\a\b}\<\psi_{\s\t},\psi_{\tmu\tmu}\psi'_{\a\b}\>_\bmu\\
        &=\<\psi_{\s\t},\psi_{\tmu\tmu}\psi'_{\u\v}h^{\star}\>_\bmu,
  \end{align*}
  where the last equality follows using \autoref{L:PsiDual} and
  \autoref{C:GDualBasis}.  Hence, the form $\<\ ,\ \>_\bmu$ is
  associative. Taking duals reverses the grading. Therefore, the map
  $x\mapsto\<x,?\>_\bmu$, for $x\in G^\bmu$, gives the required
  isomorphism.
\end{proof}

\subsection{Quiver Schur algebras}\label{S:SchurAlgebra}
We are now ready to define the quiver Schur algebras of type
$\Gamma_e$, which are the main objects of study in this paper. Recall that $\mz$
is an arbitrary integral domain.

\begin{Definition}\label{D:QuiverSchur} Suppose that $\Lambda\in P^+$ and let
  $G^\Lambda_n=\bigoplus_{\bmu\in\Parts}G^\bmu$. The \textbf{quiver Schur algebra} of type
  $(\Gamma_e,\Lambda)$ is the endomorphism algebra
  $$\Sch[n]=\Sch[n](\Gamma_e)=\ZEnd_{\R}(G^\Lambda_n).$$
\end{Definition}

By definition $\Sch[n]$ is a graded $\mz$-algebra. As a $\mz$-module,
$\Sch[n]$ admits a decomposition
$$\Sch[n]=\bigoplus_{\bnu,\bmu\in\Parts}\ZHom_{\R}(G^\bnu,G^\bmu).$$
By \autoref{T:trace}, $\R$ is a graded symmetric algebra, so by
\cite[61.2]{C&R}
\begin{equation}\label{E:Homs}
  \ZHom_{\R}(G^\bnu,G^\bmu)\cong G^\bmu\cap(G^\bnu)^\star
\end{equation}
as graded $\mz$-modules, where an isomorphism is given by
$\Psi\mapsto\Psi(e^\bnu y^\bnu)$. By \autoref{C:intersection},
if $\s\in\Std^\bmu(\blam)$ and $\t\in\Std^\bnu(\blam)$, for
$\blam\in\Parts[\beta]$, then $\psi_{\s\t}\in
G^\bmu\cap(G^\bnu)^\star$ so we can define a homomorphism
$\Psi^{\bmu\bnu}_{\s\t}\in\ZHom_{\R}(G^\bnu,G^\bmu)$ by
\begin{equation}\label{E:Psi}
  \Psi^{\bmu\bnu}_{\s\t}(e^\bnu y^\bnu h) = \psi_{\s\t}h, \qquad\text{for all }h\in\R.
\end{equation}
We think of $\Psi^{\bmu\bnu}_{\s\t}$ as an element of $\Sch[n]$ in the obvious way.
By definition, $\Phi^{\bmu\bnu}_{\s\t}$ is homogeneous of degree
$(\deg\s-\deg\tmu)+(\deg\t-\deg\tnu)$ since $\psi_{\s\t}$ has degree
$\deg\s+\deg\t-\deg\tnu$ when considered as an element of~$G^\bnu$.

\begin{Example}
  It is necessary to include $\bmu$ and $\bnu$ in the notation
  $\Psi^{\bmu\bnu}_{\s\t}$ because a given tableau can belong to $\Std^\bmu(\bnu)$ for many
  different~$\bmu$. The simplest example of this
  phenomenon occurs when $\t=\(\,\Tableau[-1]{{1}}\,|\,\emptyset\)$ and $\charge=(0,0)$, so that
  $\Lambda=2\Lambda_0$. Let $\bmu=(1|-)$ and $\bnu=(-|1)$. Then
  $\t\in\Std^\bmu(\bmu)\cap\Std^\bnu(\bmu)$ and
  $\psi_{\t\t}=e^\bmu y^\bmu\in G^\bmu\cap G^\bnu\cap(G^\bmu)^\star\cap(G^\bnu)^\star$
  by \autoref{C:intersection}. Therefore, the tableau $\t$ determines four different maps
  in~$\Sch[n]$:
  \begin{alignat*}{4}
    \Psi^{\bmu\bmu}_{\t\t}&\map{G^\bmu}G^\bmu; e^\bmu y^\bmu h\mapsto \psi_{\t\t}h,&\qquad&
    \Psi^{\bnu\bmu}_{\t\t}&\map{G^\bmu}G^\bnu; e^\bmu y^\bmu h\mapsto \psi_{\t\t}h,\\
    \Psi^{\bmu\bnu}_{\t\t}&\map{G^\bnu}G^\bmu; e^\bnu y^\bnu h\mapsto \psi_{\t\t}h, &&
    \Psi^{\bnu\bnu}_{\t\t}&\map{G^\bnu}G^\bnu; e^\bnu y^\bnu h\mapsto \psi_{\t\t}h.
  \end{alignat*}
  We have
  $\deg\Psi^{\bmu\bmu}_{\t\t}=0$,
  $\deg\Psi^{\bmu\bnu}_{\t\t}=1=\deg\Psi^{\bnu\bmu}_{\t\t}$ and
  $\deg\Psi^{\bnu\bnu}_{\t\t}=2$.
\end{Example}

For $\blam\in\Parts$ let
$\Tcal^\blam=\set{(\bmu,\s)|\s\in\Std^\bmu(\blam) \text{ for }
\bmu\in\Parts}$.

\begin{Theorem}\label{T:SCellular}
  Suppose that $e=0$ or $e>n$ and that $\mz$ is an integral domain.
  Then $\Sch[n]$ is a graded cellular algebra with cellular basis
  $\set{\Psi^{\bmu\bnu}_{\s\t}|(\bmu,\s),(\bnu,\t)\in\Tcal^\blam
                         \text{ and }\blam\in\Parts}$,
  weight poset $(\Parts,\gedom)$ and degree function
  $\deg\Psi^{\bmu\bnu}_{\s\t}=\deg\s-\deg\tmu+\deg\t-\deg\tnu$.
\end{Theorem}

\begin{proof}
  By \autoref{C:intersection} and \autoref{E:Homs} the maps in the statement of the
  Theorem are a basis of~$\Sch[n]$. As in \cite[\S6]{DJM:cyc}, it is now a purely formal
  argument to show that this basis is a cellular basis of~$\Sch[n]$. We have already verified
  axioms (GC$_d$) and (GC$_1$) from \autoref{S:cellular}. Axiom (GC$_3$) is a
  straightforward calculation using the fact that $\psi_{\s\t}^\star=\psi_{\t\s}$ by
  \autoref{T:PsiBases}; see \cite[Proposition~6.9]{DJM:cyc}. It remains to check
  (GC$_2$) but this follows by repeating the argument from
  \cite[Theorem~6.6(ii)]{DJM:cyc}, essentially without change, using
  \autoref{C:intersection} and \autoref{T:PsiBases}.
\end{proof}

\begin{Remark}\label{R:Semistandard}
In \cite[Theorem~6.6]{DJM:cyc} the cellular basis of the cyclotomic $q$-Schur
algebras is labelled by \textit{semistandard tableaux} of type~$\bnu$. The
tableaux in $\Tcal^\blam$ are, in fact, closely related to semistandard
tableaux. Using the notation of \cite[Definition~4.2]{DJM:cyc}, if
$(\bnu,\t)\in\Tcal^\blam$ then $\bnu(\t)$ is a semistandard $\blam$-tableau of
type~$\bnu$. In fact, \autoref{A:StandingAssumption} implies that the
entries in each row of~$\bmu$ have distinct residues. Consequently, if
$\s,\t\in\Std^\bnu(\blam)$ then $\s=\t$ if and only if $\bnu(\s)=\bnu(\t)$.
\end{Remark}

\begin{Example}
  If $\ell=2$ then $\Sch$ is positively graded by \autoref{T:L2Positivity},
  proved in the appendix. If $\ell>2$ then $\Sch$ is, in general, not
  positively graded. For example, suppose that $\Lambda=3\Lambda_0$,
  $\bmu=(1|2,1|2^2)$ and
  $$\t=\Tritab({1,6},{7}|{2,3},{4,8}|{5}).$$
  Then it is easy to check that $\t\in\Std^\bmu(2,1|2^2|1)$ and that
  $\deg\t=2<\deg\tmu=3$. So, $\deg\Psi^{\bmu\bmu}_{\t\t}=-2$.
\end{Example}

Now that we know that $\Sch[n]$ is a graded cellular algebra we can use the
general theory from \autoref{S:cellular} to construct cell modules and
irreducible $\Sch[n]$-modules.

Suppose that $\blam\in\Parts$. The \textbf{graded Weyl module} $\Delta^\blam$ is
the cell module for~$\Sch[n]$ corresponding to~$\blam$. More explicitly,
$\Delta^\blam$ is the $\Sch[n]$-module with basis
\begin{equation}\label{E:WeylBasis}
  \set{\Psi^\bnu_\t|(\bnu,\t)\in\Tcal^\blam}
\end{equation}
such that
$\(\Phi^{\blam\blam}_{\tlam\tlam}\Sch[n]+(\Sch[n])^{\gdom\blam}\)/(\Sch[n])^{\gedom\blam}
    \cong\Delta^\blam$
under the map that sends
$\Psi^{\blam\bnu}_{\tlam\t}+(\Sch[n])^{\gdom\blam}$ to $\Psi^\bnu_\t$,
for $(\bnu,\t)\in\Tcal^\blam$. (Note that
$\deg\Psi^{\blam\blam}_{\tlam\tlam}=0$, for all $\blam\in\Parts$.)

As in \autoref{S:cellular}, the graded Weyl module $\Delta^\blam$ comes equipped
with a homogeneous bilinear form $\<\  ,\ \>$ of degree zero such that
\begin{equation}\label{E:innerProd}
\<\Psi^\bmu_\s,\Psi^\bnu_\t\>\Psi^{\blam}_{\tlam}
    \equiv\Psi^{\bmu}_{\s}\Psi^{\bnu\blam}_{\t\tlam},
    \qquad\text{for $(\bmu,\s),(\bnu,\t)\in\Tcal^\blam$}.
\end{equation}
Define $L^\blam=\Delta^\blam/\rad\Delta^\blam$, where $\rad\Delta^\blam$
is the radical of this form. Set $\nabla^\blam=(\Delta^\blam)^\circledast$.

\begin{Theorem}\label{T:quasi}
  Suppose that $e=0$ or $e>n$ and that $\mz=K$ is a field. Then~$\Sch[n]$ is
  a quasi-hereditary graded cellular algebra with:
  \begin{itemize}
    \item weight poset $(\Parts,\gedom)$,
    \item graded standard modules $\set{\Delta^\blam|\blam\in\Parts}$,
    \item graded costandard modules $\set{\nabla^\blam|\blam\in\Parts}$, and,
    \item graded simple modules $\set{L^\blam\<k\>|\blam\in\Parts\text{ and }k\in\Z}$.
  \end{itemize}
  Moreover, $L^\blam\cong(L^\blam)^\circledast$ for all $\blam\in\Parts$.
\end{Theorem}

\begin{proof}
By definition, $\Psi^{\blam\blam}_{\tlam\tlam}$ is the identity map
on~$G^\blam$, so $\<\Psi^\blam_{\tlam},\Psi^\blam_{\tlam}\>=1$ by
\autoref{E:innerProd}. Consequently, $L^\blam\ne0$ for all $\blam\in\Parts$.
Therefore, $L^\blam\cong(L^\blam)^\circledast$, for $\blam\in\Parts$, and
$$\set{L^\blam\<k\>|\blam\in\Parts\text{ and }k\in\Z}$$ is a complete
set of pairwise non-isomorphic irreducible $\Sch[n]$-modules by
\autoref{T:CellularSimples}. In turn, this implies that~$\Sch[n]$ is a
quasi-hereditary algebra by \autoref{C:quasiH}, with standard and costandard
modules as stated.
\end{proof}

For each $\blam\in\Parts$ set
$\Psi^\blam=\Psi^{\blam\blam}_{\tlam\tlam}$. Then $\Psi^\blam$
(restricts to) the identity map on~$G^\blam$ and
$\sum_\blam\Psi^\blam$ is the identity element of~$\Sch[n]$. As a
$\mz$-module, every $\Sch[n]$-module~$M$ has a \textbf{weight space}
decomposition
\begin{equation}\label{E:WeightSpace} M=\bigoplus_{\blam\in\Parts} M_\blam,\quad\text{ where
  }M_\blam=M\Psi^\blam. \end{equation} In particular, if $\blam,\bnu\in\Parts$ then
$\set{\Psi^\bnu_\t|(\bnu,\t)\in\Tcal^\blam}$ is a basis of $\Delta^\blam_\bnu$ by
\autoref{E:WeylBasis}.

\begin{Remark}\label{R:Weyl}
  Although we will not need this, the
  reader can check that if $(\bnu,\t)\in\Tcal^\blam$ then we can
  identify $\Psi^\bnu_\t$ with the homomorphism $G^\bnu\to S^\blam$
  that sends $\psi_{\tnu\tnu}h$ to~$\psi_{\tlam\t}h$, for~$h\in\R$.
  In this way, $\Delta^\blam$ can be identified with a
  $\Sch[n]$-submodule of $\Hom_{\R}(G_n^\Lambda,S^\blam)$.
  By \autoref{C:GSpecht} there is a projection map
  $\pi^\blam:G^\blam\surjection S^\blam$ such that
  $\pi^\blam(\psi_{\tlam\tlam}h)=\psi_{\tlam}h$, for all $h\in\R$. By
  \autoref{T:SCellular} and the remarks after~\autoref{E:WeylBasis}, the
  weight space $\Delta^\blam_\bnu$ of the Weyl module $\Delta^\blam$ can be
  identified with the set of maps in $\ZHom_{\R}(G^\bnu,S^\blam)$ that factor
  through~$\pi^\blam$ so that the following diagram commutes:
  $$\begin{tikzpicture}
     \matrix[matrix of math nodes,row sep=1cm,column sep=16mm]{
         |(GLam)| G^\bnu & |(Gmu)| G^\blam \\
                         & |(Smu)| S^\blam\\
     };
     \draw[->,dashed] (GLam) -- node[above]{$\exists\Psi'$} (Gmu);
     \draw[->] (GLam) -- node[below] {$\Psi$} (Smu);
     \draw[->>] (Gmu) -- node[right]{$\pi^\blam$}(Smu);
  \end{tikzpicture}$$
\end{Remark}

\subsection{Graded Schur functors}\label{S:GradedSchurFunctors}
We now define an exact functor from the category of graded
$\Sch[n]$-modules to the category of graded $\R$-modules and use this to
relate the graded decomposition numbers of the two algebras. To do this
we introduce a slightly bigger version of the quiver Schur
algebra~$\Sch[n]$. The idea is to enlarge $\Sch[n]$ so that it contains
a copy of $\End_{\R}(\R)$ and then use this to construct a graded Schur
functor via \autoref{E:SchurFunctors}.

Let $\DotParts=\Parts\cup\{\omega\}$, where $\omega$ is a dummy symbol, and set
$G^\omega=\R$ and $\DotG=G^\Lambda_n\oplus G^\omega$. The
\textbf{extended quiver Schur algebra} is the algebra
$$\DotS=\ZEnd_{\R}(\DotG).$$
Suppose that $\bnu\in\Parts$. For convenience of notation, set
$\Std^\omega(\bnu)=\Std(\bnu)$ and define $e^\omega=1=y^\omega\in\R$ so
that $G^\omega=e^\omega y^\omega\R$. Let $\tomega=1$ and set
$\psi_{\tomega\tomega}=e^\omega y^\omega=1$ and define
$\deg\t^\omega=0$. We consider $\Sch[n]$ is a graded subalgebra of $\DotS$
in the obvious way.

Extending \autoref{E:Psi}, if $\bnu,\bmu\in\DotParts$ and
$\s\in\Std^\bmu(\bnu)$ and $\t\in\Std^\bnu(\bnu)$ then define
$$\Psi^{\bmu\bnu}_{\s\t}(e^\bnu y^\bnu h) = \psi_{\s\t}h, \qquad\text{for all }h\in\R.$$
Then $\Psi^{\bmu\bnu}_{\s\t}\in\DotS$ and $\deg\Psi^{\bmu\bnu}_{\s\t}=\deg\s-\deg\tmu+\deg\t-\deg\tnu$.
For each multipartition $\blam\in\Parts$ set
$\dot\Tcal^\blam=\set{(\bnu,\t)|\t\in\Std^\bnu(\blam)\text{ for }\bnu\in\DotParts}
              =\Tcal^\blam\cup\,\{\omega\}\times\Std(\blam)$.

\begin{Proposition}\label{P:SSCellular}
  The algebra $\DotS$ is a graded cellular algebra with cellular basis
  $$\set{\Psi^{\bmu\bnu}_{\s\t}|(\bmu,\s),(\bnu,\t)\in\dot\Tcal^\blam
               \text{ for }\blam\in\Parts},$$
  weight poset $(\Parts,\gedom)$ and degree function
  $$\deg\Psi^{\bmu\bnu}_{\s\t}=\deg\s-\deg\tmu+\deg\t-\deg\tnu.$$
  Moreover, if $\mz=K$ is a field then $\DotS$ is a quasi-hereditary algebra
  with standard modules $\set{\dot\Delta^\blam|\blam\in\Parts}$ and simple
  modules $\set{\dot L^\blam\<k\>|\blam\in\Parts\text{ and } k\in\Z}$.
\end{Proposition}

\begin{proof}
By definition, $\Sch[n]$ is a subalgebra of $\DotS$ and, as a $\mz$-module,
$$\DotS=\Sch[n]\oplus\ZHom_{\R}(G^\omega,G^\Lambda_n)\oplus \ZHom_{\R}(G^\Lambda_n,G^\omega)
                  \oplus \ZEnd_{\R}(G^\omega).$$
For $\bmu\in\DotParts$ there are isomorphisms of graded
$\mz$-modules $G^\bmu\cong\ZHom_{\R}(G^\omega,G^\bmu)$ given by
$\psi_{\s\t}\mapsto\Psi^{\bmu\omega}_{\s\t}$, for
$\s\in\Std^\bmu(\bnu)$ and $\t\in\Std^\omega(\bnu)$ and
$\bnu\in\Parts$. Therefore, the elements in the statement of the
proposition give a basis of $\DotS$ by \autoref{T:SCellular} and
\autoref{T:GBasis}.

Now suppose that $\mz=K$ is a field. Repeating the arguments from
\autoref{T:SCellular} and \autoref{T:quasi} shows that
$\DotS$ is a quasi-hereditary graded cellular algebra.
\end{proof}

By \autoref{P:SSCellular}, there exist Weyl modules $\dot\Delta^\blam$
and simple modules $\dot L^\blam=\dot\Delta^\blam/\rad\dot\Delta^\blam$ for
$\DotS$, for each $\blam\in\Parts$. As in \autoref{E:WeylBasis}, let
$\set{\Psi^\bnu_\t|(\bnu,\t)\in\dot\Tcal^\blam}$ be the basis
of~$\dot\Delta^\blam$.

Set $\Psi^\Lambda_n=\sum_{\bmu\in\Parts}\Psi^\bmu$ and let $\Psi^\omega$ be the
identity map on $G^\omega=\R$. Then~$\Psi^\Lambda_n$ is the identity element of
$\Sch[n]$ and $\Psi^\Lambda_n+\Psi^\omega$ is the identity element of $\DotS$. By
definition,~$\Psi^\Lambda_n$ and $\Psi^\omega$ are both idempotents in~$\DotS$
and $\Psi^\Lambda_n\DotS \Psi^\Lambda_n\cong\Sch[n]$. Therefore, by
\autoref{E:SchurFunctors}, there are exact graded functors
$$\Fomega\map{\DotS\Mod}\Sch[n]\Mod\quad\text{and}\quad
  \Gomega\map{\Sch[n]\Mod}{\DotS\Mod}$$
given by $\Fomega(M)=M\Psi^\Lambda_n$ and
$\Gomega(N)=N\otimes_{\Sch[n]}\Psi^\Lambda_n\DotS$. By
\autoref{S:SchurFunctors} there are graded functors $$
\Hun_{n,\omega}:=\Hun_{\Psi^\Lambda_n},\quad  \Oun_{n,\omega}:=\Oun_{\Psi^\Lambda_n},\quad \Oun^{n,\omega}:=\Oun^{\Psi^\Lambda_n} $$ from
$\DotS\Mod$ to~$\DotS\Mod$ such that $\Hun_{n,\omega}(M)=M/\Oun_{n,\omega}(M)$.

\begin{Lemma}\label{L:Equiv}
  Suppose that $\mz=K$ is a field. Then the functors $\Fomega$ and $\Gomega$
  induce mutually inverse graded equivalences of categories between $\DotS\Mod$ and
  $\Sch[n]\Mod$. Moreover,
  $$\Fomega(\dot\Delta^\blam)\cong\Delta^\blam\qquad\text{and}\qquad
    \Fomega(\dot L^\blam)\cong L^\blam,$$
  for all $\blam\in\Parts$.
\end{Lemma}

\begin{proof}
  Let $M$ be an $\DotS$-module. Then, extending \autoref{E:WeightSpace}, $M$ has a
  weight space decomposition
  $$M=\bigoplus_{\bmu\in\DotParts}M_\bmu,\qquad\text{ where }M_\bmu=M\Psi^\bmu.$$
  Then, essentially by definition,
  $\Fomega(M)=\bigoplus_{\blam\in\Parts}M_\blam$. That is, $\Fomega$ removes the
  $\omega$-weight space of $M$. In particular,
  $\Fomega(\dot\Delta^\bmu)=\Delta^\bmu$ and $\Fomega(\dot L^\bmu)=L^\bmu$, for
  all $\bmu\in\Parts$. The fact that $\Fomega(\dot L^\bmu)=L^\bmu$ for all
  $\bmu\in\Parts$ implies that $\Oun^{n,\omega}(M)=M$, $\Oun_{n,\omega}(M)=0$,
  for all $M\in\DotS\Mod$. Therefore, $\Hun_{n,\omega}$ is the identity
  functor and $\Gomega\cong\Hun_{n,\omega}\circ\Gomega$. Hence, the
  lemma is an application of the theory of quotient functors
  given in \autoref{T:SchurEquiv}.
\end{proof}

The identity map $\Psi^\omega$ on $\R=G^\omega$ is idempotent
in~$\DotS$ and there is a graded isomorphism of $\mz$-algebras
$\Psi^\omega\DotS\Psi^\omega\cong\R$. Therefore, by
\autoref{E:SchurFunctors}, there are functors
\begin{equation}\label{E:SequivSS}
\DotFun\map{\DotS\Mod}\R\Mod\quad\text{and}\quad
  \DotGun\map{\R\Mod}{\DotS\Mod}
\end{equation}
given by $\DotFun(M)=M\Psi^\omega=M_\omega$ and $\DotGun(N)=N\otimes_{\R}\Psi^\omega\DotS$.

\begin{Proposition}\label{P:SchurFunctor}
  Suppose that $\mz=K$ is a field and $\beta\in Q^+_n$. Then there is an exact graded functor
  $\SFun\map{\Sch[n]\Mod}\R\Mod$ given by
  $\SFun(M)=(M\otimes_{\Sch[n]} \Psi^\Lambda_n\DotS)\Psi^\omega$,
  for $M\in\Sch[n]\Mod$, such that if $\blam,\bmu\in\Parts[\beta]$ then
  $\SFun(\Delta^\blam)\cong S^\blam$, $\SFun(\nabla^\blam)\cong S_\blam\<-\defect\beta\>$ and
  $$ \SFun(L^\bmu)\cong\begin{cases}
               D^\bmu,&\text{if }\bmu\in\Klesh,\\
               0,     &\text{if }\bmu\notin\Klesh.
             \end{cases}
  $$
\end{Proposition}

\begin{proof} By definition, $\SFun=\DotFun\circ\Gomega$, so
  $\SFun$ is an exact graded functor from $\Sch[n]\Mod$ to $\R\Mod$. The
  functor $\DotFun$ is
  nothing more than projection onto the $\omega$-weight space. Hence, if
  $\blam\in\Parts$ then $\DotFun(\dot\Delta^\blam)$ is spanned by the maps
  $\set{\Psi^\omega_\t|\t\in\Std(\blam)}$, since $\Std^\omega(\blam)=\Std(\blam)$.
  The map $\Phi^\omega_\t\mapsto\psi_\t$, for $\t\in\Std(\blam)$, defines an
  isomorphism $\DotFun(\dot\Delta^\blam)\cong S^\blam$ of $\R$-modules.
  Therefore, $\SFun(\Delta^\blam)\cong S^\blam$ by \autoref{L:Equiv}.
  The functor $\SFun$ is easily seen to commute with duality (in
  $\Sch[n]$ and $\R$), so
  $\SFun(\nabla^\blam)\cong\SFun(\Delta^\blam)^\circledast\cong S_\blam\<-\defect\beta\>$
  by \autoref{P:SpechtDuality}.

  By \autoref{T:SchurEquiv}, $\SFun(L^\bmu)$ is an irreducible $\R$-module whenever
  it is non-zero. A straightforward argument by induction on the dominance
  ordering using $\SFun(\Delta^\blam)\cong S^\blam$, \autoref{C:triangular} and
  \autoref{C:RSimples} now shows that $\SFun(L^\bmu)\cong D^\bmu$ if
  $\bmu\in\Klesh$ and that $\SFun(L^\bmu)=0$ otherwise.
\end{proof}

Since $\SFun$ is exact and graded, we obtain the promised relationship
between the graded decomposition numbers of $\Sch[n]$ and $\R$.

\begin{Corollary}\label{C:submatrix} Suppose that $\mz=K$ is a field and
  that $\blam\in\Parts$ and $\bmu\in\Klesh$.  Then $[S^\blam:D^\bmu]_q
  =[\Delta^\blam:L^\bmu]_q$.  \end{Corollary}

The graded decomposition multiplicities $[\Delta^\blam:L^\bmu]_q$ are one of the
main objects of interest in this paper so we give them a special name.

\begin{Definition}\label{D:dlammu}
  Suppose that $\blam, \bmu\in\Parts$. Set
  $$d_{\blam\bmu}(q)=[\Delta^\blam:L^\bmu]_q
                    =\sum_{d\in\Z}\, [\Delta^\blam:L^\bmu\<d\>]\,q^d.$$
  Let $\mathbf{D}_{\Sch[n]}(q)=(d_{\blam\bmu}(q))_{\blam,\bmu\in\Parts}$ and
  $\mathbf{D}_{\R}(q)=(d_{\blam\bmu}(q))_{\blam\in\Parts,\bmu\in\Klesh}$
  be the \textbf{graded decomposition matrices} of~$\Sch[n]$ and $\R$, respectively.
\end{Definition}

By \autoref{C:submatrix}, $\mathbf{D}_{\R}(q)$ can be considered as a
submatrix of $\mathbf{D}_{\Sch[n]}(q)$. For future use we note the
following important property of these Laurent polynomials. This is the
general property of (graded) cellular algebras given in
\autoref{C:triangular}.

\begin{Corollary}\label{C:TriangularDecomp} Suppose that $\blam,\bmu\in\Parts$. Then
  $d_{\bmu\bmu}(q)=1$ and $d_{\blam\bmu}(q)\ne0$ only if $\blam\gedom\bmu$ and
  $\blam,\bmu\in\Parts[\beta]$ for some $\beta\in Q^+_n$. \end{Corollary}

\subsection{Blocks of quiver Schur algebras}
We now give the block decomposition of the graded Schur algebra
$\Sch[n]$. The key observation is the following double centralizer result.

Recall from \autoref{S:GradedSchurFunctors} that
$\DotG=G^\Lambda_n\oplus\R$ and $\DotS=\End_{\R}(\DotG)$.

\begin{Lemma}[A double centralizer property]\label{L:DoubleCentralizer}
  There are canonical isomorphisms of graded algebras such that
  $\DotS\cong\ZEnd_{\R}(\DotG)$ and $\R\cong\ZEnd_{\DotS}(\DotG)$. In particular, the
  functor $\dot{\Fun}_n^{\Lambda}$ is fully faithful on projectives.
\end{Lemma}

\begin{proof}
  The first isomorphism is the definition of $\DotS$ whereas the second follows
  directly from the definition of~$\DotS$ because
  $$\R\cong\ZHom_{\R}(\R,\R)
       \cong\Psi^\omega\DotS\Psi^\omega\cong\ZEnd_{\DotS}(\Psi^\omega\DotS),
  $$
  and $\Psi^\omega\DotS\cong\DotG$ as a right $\DotS$-module.
\end{proof}

In order to describe the block decomposition of $\Sch[n]$ we set
$G^\Lambda_\beta=\bigoplus_{\bmu\in\Parts[\beta]}G^\bmu$ and define
$\Sch=\ZEnd_{\R}(G^\Lambda_\beta)$ if $\beta\in Q^+_n$.  Equivalently,
$\Sch=\Psi^\beta\Sch[n]\Psi^\beta$, where
$\Psi^\beta=\sum_{\bmu\in\Parts[\beta]}\Psi^\bmu$.

The subalgebras $\Sch$ of $\Sch[n]$ are the blocks of~$\Sch[n]$. More precisely, we have
the following.

\begin{Theorem}\label{T:SBlocks}
  Suppose that $\mz=K$ is a field and $\Lambda\in P^+$. Then
  $$\Sch[n]=\bigoplus_{\beta\in Q^+_n}\Sch,$$
  is the block decomposition of $\Sch[n]$ into a direct sum of indecomposable
  two-sided ideals. Moreover, if~$\beta\in Q^+_n$ then the cellular basis
  of~$\Sch[n]$ in \autoref{T:SCellular} restricts to give a graded cellular
  basis of~$\Sch$. Consequently, $\Sch[\beta]$ is
  a quasi-hereditary graded cellular algebra.
\end{Theorem}

\begin{proof}
  Let $\alpha,\beta\in Q^+$ and $\bmu\in\Parts[\beta]$. By \autoref{T:GBasis},
  if $g\in G^\bmu$ then $ge_\alpha=\delta_{\alpha\beta}g$.  Therefore, if
  $\alpha\ne\beta$ and $\blam\in\Parts[\alpha]$ then
  $\ZHom_{\R}(G^\blam,G^\bmu)=0$. Hence, as $\mz$-modules,
  \begin{align*}
  \Sch[n]&=\ZEnd_{\R}(G^\Lambda_n)
     =\bigoplus_{\alpha,\beta\in Q^+_n}\ZHom_{\R}(G^\Lambda_\alpha,G^\Lambda_\beta)\\
   %&=\bigoplus_{\substack{\beta\in Q^+_n\\I^\beta\ne\emptyset}}
   %   \bigoplus_{\blam,\bmu\in\Parts[\beta]}\ZHom_{\R[\beta]}(G^\blam,G^\bmu)
    &=\bigoplus_{\beta\in Q^+_n}
                     \ZEnd_{\R[\beta]}(G^\Lambda_\beta)
     =\bigoplus_{\beta\in Q^+_n}\Sch.
\end{align*}
It follows that the cellular basis of \autoref{T:SCellular} restricts to give
cellular bases for the algebras $\Sch$, for $\beta\in Q^+_n$. Therefore, $\Sch$
is a quasi-hereditary graded cellular algebra for each $\beta\in Q^+_n$.

It remains to show that each of the algebras $\Sch$ is indecomposable. By the
double centralizer property (\autoref{L:DoubleCentralizer}), the algebras~$\R$
and~$\DotS$ have the same number  of blocks and~$\Sch[n]$ and~$\DotS$ have the
same number of indecomposable two-sided ideals by \autoref{L:Equiv}. By
\autoref{E:Rblocks} the blocks of $\R$ are indexed by~$Q^+_n$. As the elements
of~$Q^+_n$ also index the subalgebras~$\Sch$, the non-zero algebras~$\Sch$ must be
indecomposable giving the result.
\end{proof}

For each $\beta\in Q^+_n$ define $\SFun[\beta](M)=\SFun(M\Psi^\beta)$, for an
$\Sch[\beta]$-module $M$. Then $\SFun[\beta]$ is the subfunctor of
$\SFun$ obtained by first projecting onto the block $\Sch[\beta]$.
Hence, we have the following refinement of \autoref{P:SchurFunctor}.

\begin{Corollary}\label{C:BlockSchurFunctors}
  The functor $\SFun$ is fully faithful on projective
  $\Sch$-modules. Moreover, there is a decomposition of functors
  $\SFun\cong\bigoplus_{\beta\in Q^+_n}\SFun[\beta]$,
  where $\SFun[\beta]\map{\Sch\Mod}\R[\beta]\Mod$ is the restriction of $\SFun$ to
  $\Sch\Mod$ for $\beta\in Q^+_n$.
\end{Corollary}

\begin{proof}
  By definition,  $\SFun\cong\bigoplus_{\beta\in Q^+_n}\SFun[\beta]$ so
  we only need to check that $\SFun[\beta]$ is fully faithful on projectives. This follows because
  $\SFun=\DotFun\circ\Gomega$ and the functor $\DotFun$ is
  fully faithful on projectives by \autoref{L:DoubleCentralizer} (and
  $\Gomega$ is an equivalence of categories).
\end{proof}

\begin{Corollary}\label{E:QuasiCover}
  Suppose that $\beta\in Q^+_n$. Then $\Sch$ is a quasi-hereditary cover of
  $\R[\beta]$ in the sense of Rouquier~\cite[Definition~4.34]{Rouquier:Schur}.
\end{Corollary}

\begin{proof} Recall that $\DotG\cong\Psi^\omega\DotS$ is a projective
  $\DotS$-module. Using the graded Morita equivalence between $\DotS$
  and $\Sch[n]$, we see that $G_\beta^{\Lambda}$ is a projective
  $\Sch[\beta]$-module.  By \autoref{C:BlockSchurFunctors}, the functor
  $\dot{\Fun}_n^{\Lambda}$ is fully faithful on projectives, and so is
  $\SFun[\beta]$ because $\SFun[\beta]$ is the composition of
  $\dot{\Fun}_n^{\Lambda}$ with an equivalence of categories. This
  implies that $\Sch$ is a quasi-hereditary cover of $\R[\beta]$ in the
  sense of Rouquier~\cite[Definition~4.34]{Rouquier:Schur}.
\end{proof}

\subsection{Sign-dual quiver Schur algebras}\label{S:SignDual}
In this section we construct a twisted version of the quiver Schur
algebras by considering endomorphisms of the graded permutation modules
$G_\bmu$, for $\bmu\in\Parts$. The twisted quiver Schur algebras really
come from twisting by the sign automorphism~$\sgn$ of~$\R[\beta]$
defined in~\autoref{S:SignDual}.  The twisted Schur algebras turn out to be Ringel
dual to the algebras~$\Sch[\beta]$. We need the twisted Schur algebras
in order to understand the $\Delta$-filtration multiplicities of tilting
modules in~\autoref{S:Fock}.

Suppose that $\beta\in Q^+_n$ and recall the sign isomorphism
$\sgn\map{\R[\beta]}\Rp[\beta']$ from \autoref{E:sgn}.
Consider the $\Rp[\beta']$-module
$$G^{\beta'}_{\Lambda'}=\bigoplus_{\bmu\in\Parts[\beta]}G_{\bmu'}.$$
The \textbf{sign-dual quiver Schur algebra} of type $(\Gamma_e,\Lambda')_{\beta'}$ is
the algebra
$$\dualS=\dualS(\Gamma_e)=\ZEnd_{\Rp[\beta']}(G_{\Lambda'}^{\beta'}).$$
By \autoref{E:sgn} and \autoref{L:Gsgn} we have
\begin{align}\label{E:SchurSgnInvol}
  \dualS&=\ZEnd_{\Rp[\beta']}\(\bigoplus_{\bmu\in\Parts[\beta]}G_{\bmu'}\)
  %\cong\ZEnd_{{\Rp[\beta']}^\sgn}\(\bigoplus_{\bmu\in\Parts[\beta]}(G_{\bmu'})^\sgn\)\\
  \cong\ZEnd_{\R[\beta]}\(\bigoplus_{\bmu\in\Parts[\beta]}G^\bmu\)=\Sch.
\end{align}
That is, $\dualS\cong\Sch$ as graded algebras. For $\blam\in\Parts[\beta]$ let
$$\Tcal_{\blam'}
=\set{(\bnu',\t')|\t\in\Std_{\bnu'}(\blam')\text{ for }\bnu\in\P^{\Lambda}_{\beta}}
   =\set{(\bnu',\t')|(\bnu,\t)\in\Tcal^{\blam'}}.$$
As noted in \autoref{E:Homs},
$\ZHom_{\R[\beta]}(G_{\bmu'},G_{\bnu'})\cong G_{\bnu'}\cap G_{\bmu'}^\star$
as graded vector spaces. Therefore, by \autoref{C:intersection}, the
algebra~$\dualS$ is free as a $\mz$-module with basis
\begin{equation}\label{E:dualSBasis}
  \set{\Psi_{\bmu'\bnu'}^{\s'\t'}|(\bmu',\s'),(\bnu',\t')\in\Tcal_\blam\text{ for }
    \blam\in\P^{\Lambda'}_{\beta'}},
\end{equation}
where $\Psi_{\bmu'\bnu'}^{\s'\t'}$ is the $\Rp[\beta']$-endomorphism of
$G_{\Lambda'}^{\beta'}$ given by
$$\Psi_{\bmu'\bnu'}^{\s'\t'}(e_{\bsig'} y_{\bsig'} h)
          =\delta_{\bsig'\bnu'}\psi_{\s\t}'h,$$
for $(\bmu',\s'), (\bnu',\t')\in\Tcal_{\blam'}$ as above and
$\bsig\in\P^{\Lambda}_{\beta}$. By \autoref{P:signedPsi} and
\autoref{L:Gsgn} the isomorphism $\dualS\cong\Sch$ above sends the basis
element $\Psi^{\bmu\bnu}_{\s\t}$ of $\Sch[\beta]$ to
$\pm\Psi^{\s'\t'}_{\bmu'\bnu'}$. Therefore, by \autoref{T:SCellular} and
\autoref{T:quasi}, this basis makes $\dualS$ into a quasi-hereditary graded cellular
algebra with weight poset $(\P^{\Lambda'}_{\beta'},\ledom)$.

If $\blam\in\P^{\Lambda}_{\beta}$ let $\Delta_{\blam'}$ be the
corresponding Weyl module of $\dualS$ determined by the cellular basis
$\{\Psi^{\s'\t'}_{\bmu'\bnu'}\}$ given above and
let $L_{\blam'}=\Delta_{\blam'}/\rad\Delta_{\blam'}$ be its simple head.
%Finally, let $P_{\blam}$ be the projective cover of $L_{\blam}$.

\begin{Theorem}\label{T:DualSchur}
  Suppose that $\Lambda\in P^+$ and $\beta\in Q^+_n$. The sign isomorphism
  $\sgn\map{\R[\beta]}{\Rp[\beta']}$ induces a canonical degree persevering,
  poset reversing, isomorphism of quasi-hereditary graded cellular algebras
  $\sgn\map\Sch\dualS$. Moreover, when $\mz=K$ is a field there are isomorphisms
  $$\Delta^\bmu\cong\Delta_{\bmu'}^\sgn \quad\text{and}\quad L^\bmu\cong L_{\bmu'}^\sgn $$
  % P^\bmu\cong P_{\bmu'}^\sgn,$$
  of $\Sch$-modules, for $\bmu\in\Parts[\beta]$.
  Consequently,
  $[\Delta^\blam:L^\bmu]_q=[\Delta_{\blam'}:L_{\bmu'}]_q,$
  for all $\blam,\bmu\in\Parts[\beta]$.
\end{Theorem}

\begin{proof}By \autoref{E:SchurSgnInvol}, the
  automorphism $\sgn\map{\R[\beta]}{\Rp[\beta']}$ induces an isomorphism
  $\sgn\map{\Sch[\beta]}\dualS$ of graded algebras. This isomorphism sends the basis
  element $\Psi^{\bmu\bnu}_{\s\t}$ of~$\Sch[\beta]$ to
  $\pm\Psi^{\s'\t'}_{\bmu'\bnu'}\in\dualS$. Hence, the sign isomorphism identifies
  the cell module $\Delta^\bmu$ of $\Sch[\beta]$ with the cell module
  $\Delta_{\bmu'}$ of $\dualS$ (compare
  with \autoref{P:signedPsi} and \autoref{C:SignedSpechts}). All of the
  remaining claims now follow.
\end{proof}

To each pair of weights $(\Lambda,\beta)\in P^{+}\times Q^+$ we have now
attached four different quiver Schur algebras $\Sch[\beta]$,
$\SS^\beta_\Lambda$, $\Schp$ and $\dualS$. Each of these algebras is a
quasi-hereditary graded cellular algebra. To avoid confusion we clarify the
relationships between these four algebras. By definition,
$$\Sch[\beta]=\ZEnd_{\R[\beta]}\(\bigoplus_{\bmu\in\Parts[\beta]}G^\bmu\)
  \quad\text{and}\quad
  \SS^\beta_\Lambda=\ZEnd_{\R[\beta]}\(\bigoplus_{\bmu\in\Parts[\beta]}G_\bmu\).
$$
Both of these algebras are defined using the cyclotomic quiver Hecke
algebra $\R[\beta]$ and the multicharge~$\charge$. The cellular basis of
$\Sch[\beta]$ is obtained by lifting the $\psi$-bases of the graded
permutation modules $\{G^\bmu\}$ and the graded cellular basis
of~$\SS^\beta_\Lambda$ is defined by lifting the $\psi'$-bases
of~$\{G_\bmu\}$. Using the $\sgn$ isomorphism we can define isomorphic
versions of both of these algebras using $\Rp[\beta']$ and the
``conjugate'' multicharge~$\charge'$. That is, define
$$\dualS=\ZEnd_{\Rp[\beta']}(\bigoplus_{\bmu\in\Parts[\beta]}G_{\bmu'})
  \quad\text{and}\quad
  \Schp=\ZEnd_{\Rp[\beta']}(\bigoplus_{\bmu\in\Parts[\beta]}G^{\bmu'}).
$$
Then the isomorphism $\sgn\map{\R[\beta]}\Rp[\beta']$ induces algebra
isomorphisms $\Sch[\beta]\cong\dualS$ and
$\SS^\beta_\Lambda\cong\Schp[\beta']$. Therefore, via the algebra $\dualS$, we
can transfer all of the theory that we have
developed for $\Sch[\beta]$ to $\SS^\beta_\Lambda$. The algebras~$\Sch[\beta]$
and $\SS^\beta_\Lambda$ are not isomorphic, however, having both algebras are useful
because they give rise to different graded Schur functors
$$\SFun[\beta]\map{\Sch[\beta]\Mod}\R[\beta]\Mod\quad\text{and}\quad
  \Fun^\beta_\Lambda\map{\SS^\beta_\Lambda\Mod}\R[\beta]\Mod,$$
where $\Fun^\beta_\Lambda=\sgn\circ\Fun^{\Lambda'}_{\beta'}\circ\sgn^{-1}$.
We are abusing notation in the definition of $\Fun^\beta_\Lambda$
because the left-hand $\sgn$ is the equivalence
$\sgn\map{\Rp[\beta']\Mod}\R[\beta]\Mod$ of~\autoref{E:sgnEquiv} whereas
$\sgn^{-1}\map{\SS^\beta_\Lambda\Mod}\dualS\Mod$ is the inverse of the
equivalence given by \autoref{T:DualSchur}.

We will show in \autoref{T:Ringeldual} below that
$(\Sch[\beta],\SS^\beta_\Lambda)$ and $(\dualS, \Schp[\beta'])$ are both
pairs of Ringel dual algebras. Thus, ultimately, the $\sgn$ isomorphism
of the cyclotomic quiver Hecke algebras induces Ringel duality at the
level of the cyclotomic quiver Schur algebras.

\section{Tilting modules}
In this chapter we introduce the tilting modules for~$\Sch[n]$, and the closely
related Young modules for~$\R$. In \autoref{S:CycSchur} we use the Young modules to
prove that the cyclotomic quiver Schur algebra~$\Sch$ is Morita
equivalent to the (ungraded) cyclotomic Schur algebras
of~\cite{DJM:cyc,BK:HigherSchurWeyl} whereas the tilting modules give
one of the Kazhdan-Lusztig bases of the Fock space in~\autoref{S:Fock}.
Throughout this chapter we continue working over a field and we maintain
our standing \autoref{A:StandingAssumption} that $e=0$ or $e>n$.

\subsection{Young modules}\label{S:Young}
In this section we show that there exists a family of indecomposable
$\R$-modules indexed by $\Parts$ and that $G^\bmu$ is a direct sum of these
modules, for each $\bmu\in\Parts$.

Fix $\beta\in Q^+_n$ and recall from \autoref{E:WeightSpace} that every
$\Sch$-module has a weight space decomposition. Analogously, as a right
$\Sch$-module, the regular representation
of~$\Sch$ has a decomposition into a direct sum of left weight spaces:
\begin{equation}\label{E:Zmu}
  \Sch=\bigoplus_{\bmu\in\Parts[\beta]}Z^\bmu,
     \quad\text{where } Z^\bmu=\Psi^\bmu\Sch \text{ for } \bmu\in\Parts[\beta].
\end{equation}
Since $\Psi^\bmu$ is an idempotent in~$\Sch$, $Z^\bmu$ is a
projective $\Sch$-module. By \autoref{T:SCellular}, the module $Z^\bmu$ has basis
$\set{\Psi^{\bmu\bnu}_{\s\t}|(\bmu,\s),(\bnu,\t)\in\Tcal^\blam \text{ and
}\blam\in\Parts[\beta]}$.

Let $P^\bmu$ be the projective cover of $L^\bmu$ (in the category of graded
$\Sch$-modules). By \autoref{C:CartanSymmetric}, $P^\bmu$ has a
filtration by Weyl modules in which $\Delta^\blam$ appears with graded
multiplicity $(P^\bmu:\Delta^\blam)_q=[\Delta^\blam:L^\bmu]_q$. We now
describe an analogous $\Delta$-filtration of~$Z^\bmu$.

Fix a total ordering $\Std^\bmu(\Parts[\beta])=\{s_1,\dots,\s_z\}$ such that
$a>b$ whenever $\blam_a\gdom\blam_b$, where $\blam_c=\Shape(\s_c)$.
In particular, $\s_1=\tmu$. If $a\ge1$ let $M_a$ be the
submodule of $Z^\bmu$ spanned by
$$\set{\Psi^{\bmu\bnu}_{\s_b\t}|\t\in\Std^\bnu(\blam_b)\text{ for
}\bnu\in\Parts[\beta]\text{ and }b\ge a}.$$
By \autoref{T:SCellular} and (GC$_2$) of \autoref{D:cellular}, the cell
filtration of $\Sch[\beta]$ restricts to $Z^\bmu$ showing that
\begin{equation}\label{E:SmuFiltration}
Z^\bmu=M_1\supset M_2\supset\dots\supset M_z\supset 0
\end{equation}
is an $\Sch$-module filtration of $Z^\bmu$ with $M_a/M_{a+1}\cong
\Delta^{\blam_a}\<\deg\s_a-\deg\tmu\>$, for $1\le a\le z$. Thus, in the notation of
\autoref{S:graded}, $Z^\bmu$ has a $\Delta$-filtration in which
$\Delta^\blam$ appears with graded multiplicity
\begin{equation}\label{E:ZmuDelta}
  (Z^\bmu:\Delta^\blam)_q:=\sum_{\s\in\Std^\bmu(\blam)}q^{\deg\s-\deg\tmu}.
\end{equation}
Since $\Sch$ is quasi-hereditary $(Z^\bmu:\Delta^\blam)_q$ is
independent of the choice of~$\Delta$-filtration.

By the last paragraph $(Z^\bmu:\Delta^\bmu)_q=1$ and there is a surjection
$Z^\bmu\surjection\Delta^\bmu$. Moreover, $\Delta^\blam$ appears in~$Z^\bmu$
only if $\blam\gedom\bmu$. Therefore, since $Z^\bmu$ is projective, it follows
that
\begin{equation}\label{E:proj}
  Z^\bmu=P^\bmu\oplus\bigoplus_{\blam\gdom\bmu}z_{\blam\bmu}(q)P^\blam,
\end{equation}
for some Laurent polynomials $z_{\blam\bmu}(q)\in\N[q,q^{-1}]$. This
observation will be used later to compute the graded decomposition
numbers of~$\Sch$ in characteristic zero.

Analogously, let $Z_\bmu=\Psi_\bmu\dualSp$, where
$\Psi_\bmu=\Psi^{\tmu\tmu}_{\bmu\bmu}\in\dualSp$ is
the identity map on~$G_\bmu$.  Then, as before, $Z_\bmu$ is a projective
$\dualSp$-module so there exist Laurent polynomials
$z'_{\blam\bmu}(q)\in\N[q,q^{-1}]$ such that
$Z_{\bmu}=P_{\bmu}\oplus\bigoplus_{{\bmu}\gdom{\blam}}
         z'_{{\blam}{\bmu}}(q)P_{\blam}$,
where~$P_{\blam}$ is the projective cover of~$L_{\blam}$ in $\dualSp\Mod$. As
in \autoref{S:GradedSchurFunctors}, there is a graded Schur functor
$\SpFun\map{\dualSp\Mod}\R[\beta]\Mod$ for~$\dualSp$. Equivalently,
$\SpFun=\sgn\circ\Fun^{\Lambda'}_{\beta'}\circ\sgn^{-1}$
by \autoref{T:DualSchur},

\begin{Definition}\label{D:YoungModule}
  Suppose that $\bmu\in\Parts[\beta]$. The \textbf{graded Young modules} are the
  $\R[\beta]$-modules
  $$Y^\bmu=\SFun[\beta](P^\bmu)\qquad\text{ and }\qquad
    Y_\bmu=\SpFun(P_{\bmu}).$$
\end{Definition}

The next result gives some justification for this terminology. In
\autoref{L:UngradedYoundMods} below we will show that the graded Young modules
are graded lifts of the Young modules for~$\H$ introduced in~\cite{M:tilting}.

\begin{Proposition}\label{P:YoungProperties}
  Suppose that $\beta\in Q^+_n$ and that $\blam,\bmu\in\Parts[\beta]$. Then:
  \begin{enumerate}
    \item The Young modules $Y^\bmu$ and $Y_\bmu$ are indecomposable
    $\R[\beta]$-modules. Moreover, $Y^\bmu\cong Y_{\bmu'}^\sgn$.
    \item If $d\in\Z$ then $Y^\bmu\cong Y^\bnu\<d\>$
      if and only if~$\blam=\bmu$ and $d=0$. Similarly, $Y_\bmu\cong Y_\blam\<d\>$
      if and only if~$\blam=\bmu$ and $d=0$.
    \item
    $G^\bmu\cong Y^\bmu\oplus\bigoplus_{\blam\gdom\bmu}z_{\blam\bmu}(q)Y^\blam$
    and
    $G_\bmu\cong Y_\bmu\oplus\bigoplus_{\bmu\gdom\blam}z'_{\blam\bmu}(q)Y_\blam$.
    \item $Y^\bmu$ has a graded Specht filtration in which $S^\blam$ appears with
    graded multiplicity
      $$ (Y^\bmu:S^\blam)_q= [\Delta^\blam:L^\bmu]_q $$
      and~$Y_\bmu$ has a graded dual Specht filtration in which $S_\blam$
      appears with graded multiplicity
      $$ (Y_\bmu:S_\blam)_q=[\Delta_\blam:L_\bmu]_q. $$
    %\item if $\bmu\in\Klesh[\beta]$ then $Y^\bmu$ is the projective cover of
    %$D^\bmu$.
  \end{enumerate}
\end{Proposition}

\begin{proof}
  By \autoref{C:BlockSchurFunctors}, the functor $\SFun[\beta]$ is fully faithful
  on projective modules, so $\ZEnd_{\R[\beta]}(Y^\bmu)\cong\ZEnd_{\Sch}(P^\bmu)$
  is a local ring since $P^\bmu$ is indecomposable. Similarly,
  $\ZEnd_{\R[\beta]}(Y_\bmu)$ is a local ring. Hence,~$Y^\bmu$ and
  $Y_\bmu$ are indecomposable $\R[\beta]$-modules. Moreover, the fact
  that $\SFun[\beta]$ is fully faithful on projectives also implies~(b)
  since the $P^\bmu\<d\>$ are pairwise non-isomorphic.

  Applying the Schur functor from \autoref{P:SchurFunctor},
  $$\SFun[\beta](Z^\bmu)=\Psi^\bmu\DotS[\beta]\Psi^\omega
        \cong\ZHom_{\R[\beta]}(\R[\beta],G^\bmu)\cong G^\bmu.$$
  Hence, the first formula in part~(c) follows from~\autoref{E:proj}.
  Similarly,
  $\SpFun(Z_\bmu)\cong\ZHom_{\R[\beta]}(\R[\beta],G_\bmu)\cong G_\bmu$,
  which implies the second formula in part~(c). By \autoref{L:Gsgn},
  $G^\bmu\cong(G_{\bmu'})^\sgn$ so it follows by induction on dominance that
  $Y_\bmu\cong Y_\bmu^\sgn$, for all $\bmu\in\Parts[\beta]$, completing
  the proof of~(a).

  Now consider~(d). If $\bmu\in\Parts[\beta]$ then $L^\bmu\ne0$ by
  \autoref{T:quasi}. Therefore,
  $$(P^\bmu:\Delta^\blam)_q=[\Delta^\blam:L^\bmu]_q$$
  by \autoref{C:CartanSymmetric}, for $\blam\in\Parts[\beta]$. Therefore,
  $(Y^\bmu:S^\blam)_q=(P^\bmu:\Delta^\blam)_q=[\Delta^\blam:L^\bmu]_q$ by the
  exactness of $\SFun[\beta]$ and \autoref{P:SchurFunctor}. By the same argument,
  $(Y_\bmu:S_\blam)_q=(P_\bmu:\Delta_\blam)_q=[\Delta_\blam:L_\bmu]_q$,
  so~(d) holds. (Note that we are not claiming that the graded Specht
  filtration multiplicities for~$Y^\bmu$ and $Y_\bmu$ are independent of
  the choice of filtration.)
%
%   Finally, part~(e) follows by the exactness of $\SFun[\beta]$ and
%   \autoref{P:SchurFunctor} because~$P^\bmu$ is the projective cover
%   of~$L^\bmu$.
\end{proof}

\begin{Remark}
  The Laurent polynomials $z'_{\blam\bmu}(q)$ in part~(c) of
  \autoref{P:YoungProperties} should be computed using the analogue of
  \autoref{E:ZmuExpansion} for the algebra
  $\SS^{\Lambda}_{\beta}\cong\dualS$ whereas the graded decomposition
  number $[\Delta_\blam:L_\bmu]_q$ in part~(d) is for the algebra
  $\SS^\beta_\Lambda\cong\SS^{\Lambda'}_{\beta'}$.
\end{Remark}

\begin{Corollary}\label{C:DualYoung}
  Suppose that $\bmu\in\Parts[\beta]$, for $\beta\in Q^+_n$.
  Then, as $\R[\beta]$-modules
  $$Y^\bmu\cong(Y^\bmu)^\circledast\<2\defect\beta\>\quad\text{and}\quad
    Y_\bmu\cong(Y_\bmu)^\circledast\<2\defect\beta\>.$$
\end{Corollary}

\begin{proof}
  As both isomorphisms can be proved similarly we consider only the
  first one.  If~$\bmu$ is maximal in~$\Parts[\beta]$ then
  $Y^\bmu=G^\bmu$ by \autoref{P:YoungProperties}(c), so in this case the
  result is a special case of \autoref{T:GSelfDual}. Now,
  $G^\bmu\cong(G^\bmu)^\circledast\<2\defect\beta\>$ by
  \autoref{T:GSelfDual}. So if $\bmu$ is not maximal the result now
  follows by induction on dominance using parts~(b) and~(c) of
  \autoref{P:YoungProperties}.
\end{proof}

We want to identify the projective Young modules. Recall that
$P^\bmu$ is the projective cover of $L^\bmu$ and
$\SFun[\beta](L^\bmu)=D^\bmu$ if $\bmu\in\Klesh$ by
\autoref{P:SchurFunctor}.

\begin{Proposition}\label{P:YoungProjCover}
  Suppose that $\bmu\in\Klesh[\beta]$, for $\beta\in Q^+_n$. Then $Y^\bmu$ is the
  projective cover of~$D^\bmu$.
\end{Proposition}

\begin{proof}As $\SFun$ is exact there is a surjective map
  $Y^\bmu\twoheadrightarrow D^\bmu$. Therefore, it suffices to show that
  $Y^\bmu$ is projective since $Y^\bmu$ is indecomposable by
  \autoref{P:YoungProperties}(a).

  Recall from \autoref{S:GradedSchurFunctors}, that
  $\DotS=\End_{\R}(\DotG)$, where $\DotG=G^\Lambda_n\oplus\R$. By
  \autoref{E:SequivSS}, there is a graded Schur functor $\DotFun$ from
  $\DotS\Mod$ to $\R\Mod$ given by $\DotFun(M)=M_\omega$. In particular,
  $\DotFun(\DotG)\cong\R$ as graded $\R$-modules.

  As an $\DotS$-module, $\DotG\cong\Psi^\omega\DotS$. In particular,
  $\DotG$ is a projective $\DotS$-module. If $\blam\in\Parts$ let
  $\dot P^\blam$ be the projective cover of the irreducible
  $\DotS$-module $\dot L^\blam$. The graded multiplicity of $\dot P^\blam$ as a
  summand of~$\DotG$ is equal to
  $$\Dim\ZHom_{\DotS}(\DotG,\dot L^\blam)
         =\Dim\ZHom_{\DotS}(\Psi^\omega\DotS, \dot L^\blam)
         =\Dim\dot L^\blam\Psi^\omega=\Dim D^\blam,
   $$
   where the first equality follows because $\Psi^\omega$ is an
   idempotent and the second comes from \autoref{P:SchurFunctor}.
   Consequently,
   $\DotG\cong\bigoplus_{\blam\in\Klesh}(\Dim D^\blam)\dot P^\blam$ as an
   $\DotS$-module. By definition, $Y^\blam=\SFun(P^\blam)=\DotFun(\dot P^\blam)$, for
   all $\blam\in\Parts$. Therefore
   $$\R=\DotFun(\DotG)\cong\bigoplus_{\blam\in\Klesh}(\Dim D^\blam)Y^\blam$$
   as a right $\R$-module. The result follows.
\end{proof}

A \textbf{prinjective} module for an algebra is a module that
is both projective and injective.

\begin{Corollary} \label{P:ProjInjCover}
  Suppose that $\bmu\in\Klesh[\beta]$ and $\beta\in Q^+_n$. Then
  $P^\bmu\cong(P^{\bmu})^{\circledast}\<2\defect\beta\>$. Consequently,
  $P^\bmu$ is a prinjective $\Sch[n]$-module.
\end{Corollary}

\begin{proof} By the proof of \autoref{P:YoungProjCover}, $\DotG$ is a
projective $\DotS$-module. Moreover,
$\DotG\cong(\DotG)^\circledast\<2\defect\beta\>$ as an $\R$-module by
\autoref{T:GSelfDual} and \autoref{T:trace}. By
 a standard argument, see \cite[(1.5), (1.6)]{DJ}, this implies that
$\DotG\cong(\DotG)^\circledast\<2\defect\beta\>$ as an $\DotS$-module.
The proof of \autoref{P:YoungProjCover} shows that
$\DotG\cong\bigoplus_{\blam\in\Klesh}(\Dim D^\bmu)\dot P^\bmu$ as an
$\DotS$-module. Consequently, up to shift,
$P^\bmu$ and $(P^\bmu)^\circledast$ are both summands of
$\DotG$. Therefore, there exists $\bnu\in\Klesh$ and $d\in\Z$ such that
$(P^\bmu)^\circledast\<2\defect\beta\>\cong P^\bnu\<d\>$. Applying
the graded Schur functor and \autoref{C:DualYoung}, we deduce that
$Y^\bmu\cong(Y^\bmu)^\circledast\<2\defect\beta\>\cong Y^\bnu\<d\>$.
Hence, $d=0$ and $\bnu=\bmu$ by \autoref{P:YoungProperties}(b),
completing the proof.
\end{proof}

\subsection{Tilting modules}\label{S:Tilting}
By \autoref{T:quasi}, $\Sch$ is a quasi-hereditary algebra. An $\Sch$-module $T$ is a
(graded) \textbf{tilting module} if it has both a filtration by shifted Weyl modules
$\Delta^\blam\<k\>$, for $\blam\in\Parts[\beta]$ and $k\in\Z$, and a filtration by the
graded duals of shifted Weyl modules.

On forgetting the grading, \autoref{T:quasi} says that the ungraded
algebra $\UnS$ is quasi-hereditary. Therefore, by a famous theorem of
Ringel~\cite{Ringel}, for each $\bmu\in\Parts[\beta]$ there exists a
unique $\UnS$-module $\underline{T}^\bmu$ such that
\begin{enumerate}
  \item $\underline{T}^\bmu$ is indecomposable.
  \item $\underline{T}^\bmu$ has both a $\underline{\Delta}$-filtration and a
    $\underline{\nabla}$-filtration.
  \item $(\underline{T}^\bmu:\underline{\Delta}^\bmu)=1$ and
    $(\underline{T}^\bmu:\underline{\Delta}^\blam)\ne0$ only if $\bmu\gedom\blam$.
\end{enumerate}
Ringel's construction extends to the graded case to show that every
tilting module for $\UnS$ has a graded lift; see
\cite{Zhu:CharacteriticGradedQH,Mazorchuk:CatLinTiltingComplexes}. Since
$\underline{T}^\bmu$ is indecomposable it follows that there is a unique
graded lift $T^\bmu$ of $\underline{T}^\bmu$ such that
$(T^\bmu:\Delta^\bmu)_q=1$. The aim of this section is to show that
$T^\bmu\cong(T^\bmu)^\circledast$ is graded self-dual. To prove this we
need another description of the graded tilting modules.

Fix $\bmu\in\Parts[\beta]$ and let $\theta_\bmu\in\ZHom_{\R[\beta]}(\R[\beta],G_\bmu)$
be the map in $\DotS$ given by
$$\theta_\bmu(h) = \psi'_{\tlmu\tlmu}h,\qquad\text{for all }h\in\R[\beta].$$
We define analogues of the exterior powers for $\Sch$ using the functor $\Fomega$
from \autoref{E:SequivSS}.

\begin{Definition}\label{D:Emu}
  Suppose that $\bmu\in\Parts[\beta]$. Define
  $E^\bmu=\Fomega(\theta_\bmu\DotS)\<-\defect\beta\>$.
\end{Definition}

Observe that $E^\bmu$ is a right $\Sch$-module under composition of
maps because, by definition, $E^\bmu$ is the set of maps from
$G^\Lambda_n$ to $G_\bmu\<-\defect\beta\>$ which factor through~$\theta_\bmu$:
$$\begin{tikzpicture}
     \matrix[matrix of math nodes,row sep=1cm,column sep=16mm]{
       |(GLam)| G^\Lambda_n  & |(R)| \R[\beta] \\
                               & |(Galp)| G_\bmu\<-\defect\beta\>\\
     };
     \draw[->,dashed] (GLam) -- node[above]{$\exists\theta'$} (R);
     \draw[->] (GLam) -- node[below] {$\theta$} (Galp);
     \draw[->>] (R) -- node[right]{$\theta_\bmu$}(Galp);
  \end{tikzpicture}$$
This is similar to the description of the Weyl module $\Delta^\bmu$ given in
\autoref{R:Weyl}.

Our first aim is to give a basis for $E^\bmu$. Notice that if
$\blam\in\Parts[\beta]$, $\s\in\Std_\bmu(\blam)$ and
$\t\in\Std^\bnu(\blam)$ then $ \psi'_{\tlmu\tlmu}\psi_{\s\t}\in
G_\bmu\cap(G^\bnu)^\star$ by \autoref{C:GDualBasis}.
Therefore, we can define
$\theta^{\bmu\bnu}_{\s\t}\in\ZHom_{\R[\beta]}(G^\bnu,G_\bmu\<-\defect\beta\>)$ by
$$\theta^{\bmu\bnu}_{\s\t}(e^\bnu y^\bnu h)
        = \psi'_{\tlmu\tlmu}\psi_{\s\t}h,$$
for all $h\in\R[\beta]$. Recall that
$\Tcal^\blam=\set{(\bnu,\t)|\t\in\Std^\bnu(\blam)\text{ for }\bnu\in\Parts[\beta]}$
and that $\dot\Tcal^\blam$ is defined in the same way except that
$\bnu\in\DotParts[\beta]=\Parts[\beta]\cup\{\omega\}$.

\begin{Theorem}\label{T:EBasis}
  Suppose that $\bmu\in\Parts[\beta]$, for $\beta\in Q^+_n$. Then
  $$\set{\theta^{\bmu\bnu}_{\s\t}|\s\in\Std_\bmu(\blam)\text{ and }
    (\bnu,\t)\in\Tcal^\blam \text{ for some } \blam\in\Parts[\beta]}$$
  is a basis of $E^\bmu$. Moreover,
 $\deg\theta^{\bmu\bnu}_{\s\t}=\deg\s-\deg\tlmu+\deg\t-\deg\tnu.$
\end{Theorem}

\begin{proof}Let $\dot E^\bmu=\theta_\bmu\DotS$. Then $\dot E^\bmu$ is a
  right $\DotS$-module and $E^\bmu=\Fomega(\dot E^\bmu)\<-\defect\beta\>$. By
  \autoref{P:SSCellular}, $\dot E^\bmu$ is spanned by the maps
  $\theta_\bmu\Psi^{\bsig\bnu}_{\s\t}$, for $\bsig,\bnu\in\DotParts$,
  $(\bsig,\s), (\bnu,\t)\in\dot\Tcal^\blam$, and $\blam\in\Parts[\beta]$.
  By definition,
  $$\theta_\bmu\Psi^{\bsig\bnu}_{\s\t}(e^\bnu y^\bnu h)
              =\delta_{\bsig\omega}\psi'_{\tlmu\tlmu}\psi_{\s\t}h.$$
  Hence, applying \autoref{L:PsiDual}, $\theta_\bmu\Psi^{\bsig\bnu}_{\s\t}$ is
  non-zero only if $\bsig=\omega$, $\res(\s)=\res(\tlmu)$ and
  $\tlmu\gedom\s$, so that $\s\in\Std_\bmu(\Parts[\beta])$. Consequently, in this
  case,
  $\theta_\bmu\Psi^{\omega\bnu}_{\s\t}=\theta^{\bmu\bnu}_{\s\t}$.
  Therefore, the elements
  $$\set{\theta^{\bmu\bnu}_{\s\t}|\s\in\Std_\bmu(\blam)\text{ and }
          (\t,\bnu)\in\dot\Tcal^\blam\text{ for some }\blam\in\Parts[\beta]}$$
  span $\dot E^\bmu$. On the other hand, these elements are linearly
  independent because $\{\theta^{\bmu\bnu}_{\s\t}(e^\bnu
  y^\bnu)\}=\{\psi'_{\tlmu\tlmu}\psi_{\s\t}\}$ is a linearly independent
  subset of~$G_\bmu$ by \autoref{C:GDualBasis}(a). Hence, we have found
  a basis for $\dot E^\bmu$. Applying the functor $\Fomega$ kills the
  $\omega$-weight space of~$\dot E^\bmu$. So $\Fomega$ maps the basis
  $\{\theta^{\bmu\bnu}_{\s\t}\}$ of~$\dot E^\bmu$ to the elements in the
  statement of the theorem, or to zero if $\bnu=\omega$. Hence,
  $\{\theta^{\bmu\bnu}_{\s\t}\}$ is a basis of~$E^\bmu$.

  To complete the proof it remains to compute $\deg\theta^{\bmu\bnu}_{\s\t}$,
  for $\s\in\Std_\bmu(\blam)$, $(\bnu,\t)\in\Tcal^\blam$, and
  $\blam\in\Parts[\beta]$. Recalling the degree shifts in the definition of
  the three modules~$G^\bnu$, $G_\bmu$ and~$E^\bmu$,
  $$\deg\theta^{\bmu\bnu}_{\s\t}=\codeg\tlmu+\deg\s+\deg\t-\deg\tnu-\defect\beta.$$
  By \autoref{L:defect}, $\codeg\tlmu-\defect\beta=-\deg\tlmu$, so
  $\deg\theta^{\bmu\bnu}_{\s\t}=\deg\s-\deg\tlmu+\deg\t-\deg\tnu$
  as required.
\end{proof}

By definition, $\s\in\Std_\bmu(\Parts)$ only if $\t_\bmu\gedom\s$ which
implies that $\bmu\gedom\Shape(\s)$. Order
$\Std_\bmu(\Parts[\beta])=\{\s_1,\dots,\s_y\}$ so that
$a>b$ whenever $\blam_\a\gdom\blam_b$, where
$\blam_c=\Shape(\s_c)$, for $1\le c\le y$. (In particular, $\s_y=\tlmu$.) The
proof of \autoref{T:EBasis} shows that
$\theta^{\bmu\bnu}_{\s\t}=\theta_\bmu\Psi^{\omega\bnu}_{\s\t}$
so, as in ~\autoref{E:SmuFiltration}, the cell filtration
of~$\Sch[\beta]$ gives the following.

\begin{Corollary}\label{C:EWeylFilt}
  Suppose that $\bmu\in\Parts[\beta]$. Then
  $E^\bmu$ has a $\Delta$-filtration
  $$E^\bmu=E_1> E_2>\dots> E_y>0,$$
  such that $E_r/E_{r+1}\cong\Delta^{\blam_r}\<\deg\s_r-\deg\tlmu\>$,
  for $1\le r\le y$. In particular, $\Delta^\bmu$ is a
  submodule of $E^\bmu$, $(E^\bmu:\Delta^\bmu)_q=1$ and
  $(E^\bmu:\Delta^\blam)_q\ne0$ only if $\bmu\gedom\blam$.
\end{Corollary}

We now give a second basis of  $E^\bmu$ and use it show that $E^\bmu$ is a
tilting module.  Suppose that $\u\in\Std_\bmu(\bnu)$ and
$\v\in\Std^\blam(\bnu)$, for $\blam,\bnu\in\Parts[\beta]$.  Then
$\psi'_{\u\v}\psi_{\tlam\tlam}\in G_\bmu\cap(G^\blam)^\star$
by~\autoref{C:GSpecht}. Therefore, we can define
$\theta^{\u\v}_{\bmu\blam}\in\ZHom_{\R[\beta]}(G^\blam,G_\bmu\<-\defect\beta\>)$ by
$$\theta^{\u\v}_{\bmu\blam}(e^\blam y^\blam h)=\psi'_{\u\v}\psi_{\tlam\tlam}h,$$
for $h\in\R[\beta]$.

\begin{Lemma}\label{L:dualEBasis}
  Suppose that $\mz=K$ is a field and that $\bmu\in\Parts[\beta]$. Then
  $$\set{\theta^{\u\v}_{\bmu\bnu}|(\bmu,\u)\in\Tcal_\blam, (\bnu,\v)\in\Tcal^\blam
    \text{ for some } \blam\in\Parts[\beta]}$$
  is a basis of $E^\bmu$. Moreover,
  $\deg\theta_{\bmu\bnu}^{\u\v}=\codeg\u-\codeg\tlmu+\codeg\v-\codeg\t^{\bnu}$.
\end{Lemma}

\begin{proof}
  We first show that $\theta^{\u\v}_{\bmu\bnu}\in E^\bmu$ whenever
  $\u\in\Std_\bmu(\blam)$ and $\v\in\Std^\bnu(\blam)$, for some
  $\blam\in\Parts[\beta]$. By \autoref{T:GBasis},
  $\psi'_{\u\v}=\psi'_{\tlmu\tlmu}x$, for some $x\in\R[\beta]$. Therefore,
  $$\theta^{\u\v}_{\bmu\bnu}(e^\bnu y^\bnu h) =\psi'_{\u\v}\psi_{\t^{\bnu}\t^{\bnu}}h
           =\psi'_{\tlmu\tlmu}x\psi_{\t^{\bnu}\t^{\bnu}}h
           =\theta_\bmu(x\psi_{\t^{\bnu}\t^{\bnu}}h).$$
   That is, $\theta^{\u\v}_{\bmu\bnu}$ factors through $\theta_\bmu$ so that
   $\theta^{\u\v}_{\bmu\bnu}\in E^\bmu$ as claimed. The elements in the
   statement of the theorem are linearly independent because
   $\{\theta_{\bmu\bnu}^{\u\v}(e^\bnu y^\bnu)\}$ is a linearly
   independent subset of~$G_\bmu$ by applying $\star$ to
   \autoref{C:GDualBasis}(a).  Therefore, as we are working over a
   field, these elements give a basis of~$E^\bmu$ by counting dimensions
   using \autoref{T:EBasis}.

   Finally, as in the last paragraph of the proof of \autoref{T:EBasis}, the
   formula for the degree of $\theta_{\bmu\bnu}^{\u\v}$ follows using
   \autoref{L:defect}.
\end{proof}

Notice that, unlike \autoref{T:EBasis}, the basis of
\autoref{L:dualEBasis} does not obviously yield a
$\Delta$-filtration of~$E^\bmu$ because it is not clear how to write
the basis elements $\theta_{\bmu\bnu}^{\u\v}$ in
terms of the cellular basis of~$\Sch$. By appealing to
\autoref{T:DualSchur} it is possible to construct a
$\nabla$-filtration of~$E^\bmu$ using the basis of
\autoref{L:dualEBasis}. The existence of a $\nabla$-filtration is
also implied by the next result.

Recall the homogeneous trace form $\tau_\beta$ from \autoref{T:trace}.

\begin{Theorem}\label{T:ESelfDual}
  Suppose that $\mz=K$ is a field and that $\bmu\in\Parts[\beta]$, for $\beta\in Q^+_n$. Then
  $E^\bmu\cong(E^\bmu)^\circledast$.
\end{Theorem}

\def\form<#1|#2>{\mathbb{\Langle} #1,#2\mathbb{\Rangle}_\bmu}

\begin{proof}Using the two bases of $E^\bmu$ given by \autoref{T:EBasis} and
  \autoref{L:dualEBasis}, define $\form<\ |\ >\map{E^\bmu\times E^\bmu}K$ to be the
  unique bilinear map such that
  $$\form<\theta^{\bmu\bnu}_{\s\t}|\theta^{\u\v}_{\bmu\btau}>
                =\tau_\beta(\psi_{\s\t}\psi'_{\v\u}),$$
  for $(\bmu,\s)\in\Tcal_\blam$, $(\bnu,\t)\in\Tcal^\blam$,
  $(\bmu,\u)\in\Tcal_\bsig$ and $(\btau,\v)\in\Tcal^\bsig$ for some
  $\blam,\bsig\in\Parts[\beta]$. By \autoref{T:trace},
  $\form<\theta^{\bmu\bnu}_{\s\t}|\theta^{\s\t}_{\bmu\btau}>\ne0$ and
  $\form<\theta^{\bmu\bnu}_{\s\t}|\theta^{\u\v}_{\bmu\btau}>\ne0$ only if
  $(\u,\v)\Gedom(\s,\t)$ and $\deg(\psi_{\s\t}\psi'_{\v\u})=2\defect\beta$. Therefore,
  $\form<\ |\ >$ is a non-degenerate bilinear form.

  We claim that the bilinear form $\form<\ |\ >$ is homogeneous. To see this
  suppose that $\form<\theta^{\bmu\blam}_{\s\t}|\theta^{\u\v}_{\bmu\btau}>\ne0$,
  for basis elements as above. Then
  $\deg(\psi_{\s\t}\psi'_{\v\u})=2\defect\beta$ since $\tau_\beta$ is
  homogeneous of degree $-2\defect\beta$. Using the degree formulae in
  \autoref{T:EBasis}, together with \autoref{L:defect},
  \begin{align*}
    \deg\theta^{\bmu\blam}_{\s\t}+\deg\theta^{\u\v}_{\bmu\btau}
         &=\deg\s-\deg\t_\bmu+\deg\v-\deg\tnu
           \codeg\u-\codeg\t_\bmu+\codeg\v-\codeg\tnu\\
          &=\deg(\psi_{\s\t}\psi'_{\v\u})-2\defect\beta
           =0.
  \end{align*}
  Hence, $\form<\ |\ >$ is a homogeneous bilinear form of degree zero.

  To complete the proof it is enough to show that the form $\form<\ |\ >$ is
  associative because then the map that sends $\theta^{\bmu\blam}_{\s\t}$ to
  the function $x\mapsto \form<\theta^{\bmu\blam}_{\s\t}|x>$ is an $\Sch$-module
  homomorphism. This can be proved by repeating the argument from
  the proof of \autoref{T:GSelfDual}. We leave the details for the
  reader.
\end{proof}

\begin{Corollary}\label{C:EmuDecomp}
  Suppose that $\mz=K$ is a field and that $\bmu\in\Parts[\beta]$. Then $E^\bmu$
  is a tilting module. Moreover,
  $$E^\bmu=T^\bmu\oplus\bigoplus_{\bmu\gdom\blam}t_{\blam\bmu}(q)T^\blam,$$
  for some Laurent polynomials $t_{\blam\bmu}(q)\in\N[q,q^{-1}]$ such that
  $t_{\blam\bmu}(q)=t_{\blam\bmu}(q^{-1})$.
\end{Corollary}

\begin{proof}
  By \autoref{C:EWeylFilt} $E^\bmu$ has a $\Delta$-filtration. Therefore,
  since $E^\bmu\cong(E^\bmu)^\circledast$ it also as a $\nabla$-filtration.
  Hence, $E^\bmu$ is a tilting module so that $E^\bmu$ can be written uniquely
  as a direct sum of indecomposable tilting modules. By
  \autoref{C:EWeylFilt}, $(E^\bmu:\Delta^\bmu)_q=1$ and $(E^\bmu:\Delta^\blam)_q\ne0$
  only if $\bmu\gedom\blam$. Therefore, if~$d\in\Z$ then $T^\blam\<d\>$
  is a summand of~$E^\blam$ only if $\bmu\ge\blam$ and~$T^\bmu$ is a summand
  appearing with multiplicity~$1=(E^\bmu:\Delta^\bmu)_q$.
  Hence, $E^\bmu=T^\bmu\oplus\bigoplus_{\bmu\gdom\blam}t_{\blam\bmu}(q)T^\blam$ for
  some polynomials $t_{\blam\bmu}(q)\in\N[q,q^{-1}]$. Finally,
  $t_{\blam\bmu}(q)=t_{\blam\bmu}(q^{-1})$ because $E^\bmu$ is graded self-dual
  and because $T^\blam\cong T^\bnu\<d\>$ only if~$\blam=\bnu$ and~$d=0$.
\end{proof}

Arguing by induction on dominance we obtain the main result of this
section.

\begin{Corollary}\label{C:TiltingSelfDual}
  Suppose that $\bmu\in\Parts[\beta]$. Then $(T^\bmu)^\circledast\cong T^\bmu$.
\end{Corollary}

\begin{proof}
To start the induction note that if $\bmu$ is a minimal element
of $\Parts[\beta]$, with respect to dominance, then $E^\bmu=T^\bmu$ is
self-dual by \autoref{C:EmuDecomp}. The general case now follows by
induction using \autoref{C:EmuDecomp}.
\end{proof}

\subsection{Graded Ringel duality}
This section introduces the Ringel duality in the graded setting. The
main aim, however, is to compute the $\Delta$-filtration multiplicities
in the tilting modules.  This will allow us to identify the tilting
modules with one of the canonical bases of the Fock space in \autoref{S:Fock}.

A \textbf{full tilting module} $E^\Lambda_\beta$ for $\Sch$ is a tilting
module that contains every indecomposable tilting module, up to shift,
as a direct summand. Hence,
$$E^\Lambda_\beta=\bigoplus_{\bmu\in\Parts[\beta]}E^\bmu$$ is a full
tilting module for $\Sch$. Define the \textbf{Ringel dual} of $\Sch$ to
be the graded algebra $\ZEnd_{\Sch}(E^\Lambda_\beta)$. (Strictly
speaking, this is \textit{a} Ringel dual of~$\Sch$.)

Recall the graded Schur functor $\SFun[\beta]\map{\Sch\Mod}\R[\beta]\Mod$ from
\autoref{C:BlockSchurFunctors}.

\begin{Lemma}\label{L:SchurTilt}
  Suppose that $\bmu\in\Parts[\beta]$. Then
  $\SFun[\beta](E^\bmu)\cong G_\bmu\<-\defect\beta\>$ as an
  $\R[\beta]$-module.
\end{Lemma}

\begin{proof}
  By \autoref{P:SchurFunctor} and \autoref{L:Equiv}
   , and \autoref{D:Emu},
  \begin{align*}\SFun[\beta](E^\bmu)
    &=\SFun\big(\Fomega(\theta_\bmu\DotS)\<-\defect\beta\>\big)
     =\DotFun(\theta_\bmu\DotS\<-\defect\beta\>)
     =\theta_\bmu\DotS\Psi^\omega\<-\defect\beta\>\\
    &\cong\ZHom_{\R}(\R,G_\bmu\<-\defect\beta\>)
     \cong G_\bmu\<-\defect\beta\>,
\end{align*}
as required.
\end{proof}

\begin{Corollary}\label{C:SchurTilt}
  Suppose that $\bmu\in\Parts[\beta]$. Then
  $\SFun[\beta](T^\bmu)\cong Y_\bmu\<-\defect\beta\>$ as an
  $\R[\beta]$-module.
\end{Corollary}

\begin{proof}
  If $\bmu$ is minimal with respect to dominance in $\Parts[\beta]$ then
  $E^\bmu=T^\bmu$ by \autoref{C:EmuDecomp}. In fact,
  $E^\bmu=\Delta^\bmu=\nabla^\bmu$ by \autoref{C:EWeylFilt} and the fact
  that $T^\bmu$ is self-dual. Therefore,
  $\SFun(T^\blam)=\SFun[\beta](\nabla^\bmu)=S_\bmu\<-\defect\beta\>$ by
  \autoref{P:SchurFunctor} and \autoref{L:SchurTilt}.  On the other
  hand, $G_\bmu=S_\bmu$ by \autoref{C:GSpecht}, so
  $\SFun[\beta](T^\bmu)\cong Y_\bmu\<-\defect\beta\>$ as claimed.  If
  $\bmu$ is not minimal in $\Parts[\beta]$ the result follows by
  induction on the dominance order using \autoref{L:SchurTilt}
  since
  $E^\bmu=T^\bmu\oplus\bigoplus_{\bmu\gdom\blam}t_{\blam\bmu}(q)T^\blam$
  by \autoref{C:EmuDecomp} and
  $G_\bmu=Y_\bmu\oplus\bigoplus_{\bmu\gdom\blam}z'_{\blam\bmu}(q)Y_\blam$
  by \autoref{P:YoungProperties}(c). (Moreover,
  $t_{\blam\bmu}(q)=q^{-\defect\beta}z'_{\blam\bmu}(q)$.)
\end{proof}

% If~$A$ is an algebra let $A^{op}$ be the opposite
% algebra which is obtained by reversing the order of multiplication.

\begin{Theorem} \label{T:Ringeldual}
  Suppose that $\beta\in Q^+_n$. Then the Ringel dual of $\Sch$ is isomorphic to
  $\ddS$. In particular, $\ZEnd_{\Sch}(E^\Lambda_\beta)$ is a
  quasi-hereditary graded cellular algebra.
\end{Theorem}

\begin{proof}
  There is a natural map
  $\ZHom_{\R[\beta]}(G_\bmu,G_\bnu)\rightarrow\ZHom_{\Sch[\beta]}(E^\bmu,E^\bnu)$
  given by composition of maps, for $\bmu,\bnu\in\Parts[\beta]$. By
  \autoref{L:SchurTilt} this map is injective. On the other hand, since
  $\underline{E}^{\bmu}$ and $\underline{E}^{\bnu}$ are tilting modules for the
  quasi-hereditary algebra $\UnS[\beta]$, it is well-known (compare
  \cite[Proposition~A2.2 and Proposition A3.7]{Donkin:book}) that
  \begin{align*}
  \dim \Hom_{\Sch[\beta]}(E^\bmu,E^\bnu)
     &=\dim\Hom_{\UnS[\beta]}(\underline{E}^{\bmu},\underline{E}^{\bnu})
      =\sum_{\bsig\in\Parts[\beta]}(\underline{E}^{\bmu}:\underline{\Delta}^{\bsig})
         (\underline{E}^{\bmu}:\underline{\nabla}^{\bsig})\\
     &=\sum_{\bsig\in\Parts[\beta]}(\underline{E}^{\bmu}:\underline{\Delta}^{\bsig})
         (\underline{E}^{\bmu}:\underline{\Delta}^{\bsig})
      =\dim \Hom_{\H[\beta]}(\underline{G}_\bmu,\underline{G}_\bnu),
  \end{align*}
  where the last equality comes from \autoref{E:dualSBasis} and
  \autoref{C:EWeylFilt}.  Therefore, comparing dimensions, the Ringel dual
  of~$\Sch$ is isomorphic to $\ZEnd_{\R[\beta]}(G^\beta_\Lambda)=\ddS$, as a
  graded algebra.
\end{proof}

By \autoref{T:DualSchur}, $\ddS\cong\SS^{\Lambda'}_{\beta'}$ as graded algebras.
Note, however, that this is not an isomorphism of quasi-hereditary algebras
because the isomorphism reverses the partial ordering on the standard modules of
the two algebras.

We now identify $\ddS$ and the Ringel dual of $\Sch$. The Ringel duality functor
$\ZHom_{\Sch}(E^\Lambda_\beta,?)$, combined with \autoref{T:Ringeldual}, defines a
duality $\RFun\map{\Sch\Mod}{\ddS\Mod}$ that sends an $\Sch$-module $M$ to
$\ZHom_{\R[\beta]}(G^\beta_\Lambda, \Fun^\beta_\Lambda(M))$, where
$\Fun^\beta_\Lambda$ is the graded Schur functor for $\ddS$. It is a standard
fact that Ringel duality sends tilting modules to projective modules and
costandard modules to standard modules; for example, see
\cite[\S A4]{Donkin:book}.  That $\RFun(T^\blam)\cong P_\blam\<-\defect\beta\>$ is
immediate from the definitions whereas, using the graded Schur functor
$\Fun^\beta_\Lambda\map{\ddS\Mod}\R[\beta]\Mod$ and \autoref{R:Weyl}, shows that
$\RFun(\nabla^\bmu)\cong\Delta_\bmu\<-\defect\beta\>$.

The next result is the graded analogue of
\cite[Lemma~A4.6]{Donkin:book}. This is the result that we need in
\autoref{S:Fock} to identify the tilting modules with one of the
canonical bases of the Fock space.

\begin{Corollary}\label{L:TiltingMultiplicities}
  Suppose that $\blam,\bmu\in\Parts[\beta]$, for $\beta\in Q^+_n$. Then
  $$(T^\blam:\Delta^\bmu)_q =\overline{[\Delta_\bmu:L_\blam]_q},$$
  where $[\Delta_{\bmu}:L_{\blam}]_q$ is a graded decomposition number for the
  sign-dual quiver Schur algebra~$\ddS$.
\end{Corollary}

\begin{proof}
  Using \autoref{C:CartanSymmetric} and the remarks in the last paragraph,
  $$[\Delta_\bmu:L_\blam]_q = (P_\blam:\Delta_\bmu)_q
  =(\RFun(T^\blam):\RFun(\nabla^\bmu))_q=(T^\blam:\nabla^\bmu)_q.$$
  Therefore,
  $(T^\blam:\Delta^\bmu)_q
      =\overline{((T^\blam)^\circledast:(\Delta^\bmu)^\circledast)_q}
      =\overline{(T^\blam:\nabla^\bmu)_q}
      =\overline{[\Delta_\bmu:L_\blam]_q}$
  as required.
\end{proof}

\section{Cyclotomic Schur algebras}\label{S:CycSchur}
We are now ready to connect the quiver Schur algebras with the (ungraded)
cyclotomic Hecke algebras introduced in~\cite{DJM:cyc}
and~\cite[Theorem~C]{BK:HigherSchurWeyl}.

\subsection{Cyclotomic permutation modules}\label{S:CycPermModules}
Throughout this section we work with the ungraded Hecke algebra~$\H$.
Consequently, as in \autoref{T:BKiso}, we assume that $\mz=K$ is a suitable field. If
$w\in\Sym_n$ define $T_w=T_{i_1}\dots T_{i_k}$, where
$w=s_{i_1}\dots s_{i_k}$ is a reduced expression for~$w$. Unlike the
element $\psi_w\in\R$, $T_w$ is independent of the choice of reduced
expression for~$w$.

Suppose that $\bmu\in\Parts$. Recall that if $1\le k\le n$ and
$\t=(\t^{(1)}, \dots, \t^{(\ell)})$ is a tableau then
$\comp_\t(k)=s$ if $k$ appears in $\t^{(s)}$. Define $m^\bmu=u^\bmu
x^\bmu$ where
$$u^\bmu=\prod_{k=1}^n\prod_{s=\comp_{\tmu}(k)+1}^\ell (L_k-\xi^{(\kappa_s)})
\quad\text{and}\quad x^\bmu=\sum_{w\in\Sym_\bmu}T_w,$$
where $\xi^{(k)}$ is as defined in~\autoref{E:xi}. These definitions reduce to
\cite[Definition~3.5]{DJM:cyc} when $\xi\ne1$ and to
\cite[(6.12)--(6.13)]{BK:HigherSchurWeyl} when $\xi=1$.

\begin{Definition}[\cite{DJM:cyc,BK:HigherSchurWeyl}]\label{D:Mu}
  Suppose that $\bmu\in\Parts$ and define $\Mmu=m^\bmu\H$.
 \end{Definition}

We write $\Mmu$  rather than $M^\bmu$ to emphasize that $\Mmu$ is not
(naturally) $\Z$-graded. We will not define a graded lift of $\Mmu$. Instead,
the aim of this section is to show that~$\underline{G}^\bmu$ is a direct summand
of~$\Mmu$.

We remind the reader of our standing assumption that $e=0$ or $e>n$ from
\autoref{A:StandingAssumption}. This is crucial for the next
result -- and consequently for all of the results in this section.

\begin{Lemma}\label{L:blamCancellation}
      Suppose that $\blam\in\Parts$ and $1\ne w\in\Sym_\blam$. Then
      $e^\blam \psi_w e^\blam=0$.
\end{Lemma}

\begin{proof}
  By \autoref{D:QuiverRelations}, $\psi_w e^\blam=e(\bj)\psi_w$ where
  $\bj=w\cdot\ilam$. Now, the assumption that $e=0$ or $e>n$ implies that all
  of the nodes in row~$a$ of~$\lambda^{(l)}$ have pairwise distinct residues
  whenever $\lambda^{(l)}_a\ne0$, for $a\ge0$ and $1\le l\le\ell$. Consequently,
  $\bj\ne\ilam$ since $1\ne w\in\Sym_\blam$. Therefore, $e^\blam \psi_w e^\blam=e^\blam
  e(\bj)\psi_w=0$.
\end{proof}

\begin{Lemma}\label{L:ulamExpansion}
  Suppose that $\blam\in\Parts$. Then $e^\blam u^\blam=g^\blam(y)e^\blam y^\blam$,
  where~$g^\blam(y)$ is an invertible element of~$K[y_1,\dots,y_n]$.
\end{Lemma}

\begin{proof}We prove the Lemma only when $\xi\ne1$ and leave the case when
  $\xi=1$, which is similar, to the reader. Write $\ilam=(i_1,\dots,i_n)$ and
  let $d^\blam_1,\dots,d^\blam_n$ be as defined in \autoref{D:psis}, so
  that $d^\blam_r=\set{\comp_{\tlam}(r)<t\le\ell| i_r\equiv\kappa_t\pmod e}$,
  for $1\le r\le n$. Then, using \autoref{E:InverseBKIsoL},
  \begin{align*}
    e^\blam u^\blam
       &=\prod_{r=1}^n\prod_{t=\comp_{\tlam}(r)+1}^\ell e^\blam(L_r-\xi^{\kappa_t})\\
       &=\prod_{r=1}^n\prod_{t=\comp_{\tlam}(r)+1}^\ell e^\blam
               (\xi^{i_r}-\xi^{\kappa_t}-\xi^{i_r}y_r)\\
       &=\prod_{r=1}^n (-\xi^{i_r})^{d^\blam_r}y_r^{d^\blam_r}
            \prod_{\substack{\comp_{\tlam}(r)<t\le\ell\\i_r\not\equiv\kappa_t\pmod e}}
               e^\blam(\xi^{i_r}-\xi^{\kappa_t}-\xi^{i_r}y_r)\\
    &=e^\blam y^\blam\cdot \prod_{r=1}^n(-\xi^{i_r})^{d^\blam_r}
         \prod_{\substack{\comp_{\tlam}(r)<t\le\ell\\i_r\not\equiv\kappa_t\pmod e}}
              (\xi^{i_r}-\xi^{\kappa_t}-\xi^{i_r}y_r).
  \end{align*}
  The factor to the right of $e^\blam y^\blam$ in the last equation is a
  polynomial in $K[y_1,\dots,y_n]$ with non-zero constant term. Since each $y_r$
  is nilpotent (it has positive degree), it follows that $g(y)$ is
  invertible. All of the terms in the last equation commute, so the lemma follows.
\end{proof}

\begin{Theorem}\label{T:eme}
  Suppose that $e=0$ or $e>n$ and let $\blam\in\Parts$. Then there exists an
  invertible element $f^\blam(y)\in K[y_1,\dots,y_n]$ such that
  $$e^\blam m^\blam e^\blam = f^\blam(y) e^\blam y^\blam.$$
\end{Theorem}

\begin{proof}
  By \autoref{L:ulamExpansion}, there exists an invertible element
  $g^\blam(y)\in K[y_1,\dots,y_n]$ such that
  \begin{align*}
    e^\blam m^\blam e^\blam
      &= e^\blam u^\blam x^\blam e^\blam
       = g^\blam(y)y^\blam \sum_{w\in\Sym_\blam} e^\blam T_w e^\blam.\\
       \intertext{By \autoref{E:InverseBKIsoT}, if $w\in\Sym_n$ and
         $\bj\in I^n$ then $T_re(\bj)=(\psi_rQ_r(\bj)-P_r(\bj))e(\bj)$ so the last equation
         can be rewritten as}
    e^\blam m^\blam e^\blam
    &= g^\blam(y)y^\blam \sum_{w\in\Sym_\blam} r_w(y)e^\blam\psi_w e^\blam,
  \end{align*}
  for some $r_w(y)\in K[y_1,\dots,y_n]$. Applying \autoref{L:blamCancellation}, this
  sum collapses to give
  $$e^\blam m^\blam e^\blam=g^\blam(y) e^\blam y^\blam r_1(y)=f^\blam(y)e^\blam y^\blam,$$
  for some polynomial $f^\blam(y)\in K[y_1,\dots,y_n]$. It remains to show that
  $f^\blam(y)$ is invertible or, equivalently, that it has non-zero constant
  term. By \cite[Corollary~3.11]{HuMathas:GradedInduction}, if $1\le r\le n$ and
  $(\s,\t)\in\SStd(\Parts)$ then $y_r\psi_{\s\t}$ is a linear combination of
  terms $\psi_{\u\v}$, where $(\u,\v)\Gdom(\s,\t)$. Therefore, since $e^\blam
  y^\blam=\psi_{\tlam\tlam}$, there exist scalars $b_{\u\v}\in K$ such that
  $$f^\blam(y)e^\blam y^\blam=b_{\tlam\tlam}\psi_{\tlam\tlam}+
           \sum_{\substack{(\u,\v)\Gdom(\tlam,\tlam)\\\u,\v\in\Std^\blam(\Parts)}}
                b_{\u\v}\psi_{\u\v},$$
  where $b_{\tlam\tlam}=f^\blam(0)$ is the constant term of $f^\blam(y)$. On the
  other hand, by \cite[Theorem~3.9]{HuMathas:GradedInduction} there exist
  scalars $c_{\u\v}\in K$ such that $c_{\tlam\tlam}\ne0$ and
  $$e^\blam m^\blam e^\blam
     =e^\blam\Big( \sum_{\u,\v\Gedom\tlam}c_{\u\v}\psi_{\u\v}\Big)e^\blam
     = c_{\tlam\tlam}\psi_{\tlam\tlam}+
     \sum_{\substack{\u,\v\ne\tlam\\\u,\v\in\Std^\blam(\Parts)}}c_{\u\v}\psi_{\u\v},$$
   where the second equality follows from~\autoref{E:weights}. Hence,
   $f^\blam(0)=c_{\tlam\tlam}\ne0$ by \autoref{T:PsiBases}, and the proof is complete.
\end{proof}

\begin{Remark}\label{R:em}
  Using \autoref{T:eme} it is possible to show that $e^\blam m^\blam = f^\blam(y)
  e^\blam y^\blam+\epsilon$ where $\epsilon$ is a linear combination of
  homogeneous terms of degree strictly greater than
  $2\deg\tlam=\deg(e^\blam y^\blam)$. To see this first show that
  $e^\blam m^\blam$ is a linear combination of terms of the form
  $e^\blam m_{\blam}e(\bj)$, where $\bj\in I^\blam=\set{\bi\in
  I^n|\bi=\sigma\cdot\ilam\text{ for some }\sigma\in\Sym_\blam}$. The
  key observation is then that $\deg \psi_w e(\bj)>0$ whenever $1\ne
  w\in\Sym_\blam$ and $\bj\in I^\blam$, which can be proved by adapting
  the argument of \autoref{L:blamCancellation}. Consequently, $e^\blam
  y^\blam$ is the homogeneous component of $e^\blam m^\blam$ of minimal
  degree. Examples show that this does not always happen if we drop the
  assumption that $e=0$ or $e>n$.
\end{Remark}

Recall from \autoref{D:Mu} that $\Mmu[\blam]=m^\blam\H$.

\begin{Corollary} \label{C:emH} Suppose that  $\blam\in\Parts$. Then
  $$e^\blam\Mmu[\blam]=e^\blam m^{\blam}\H=e^\blam m^{\blam}e^\blam \H
          =e^\blam y^{\blam}\H=\underline{G}^\blam. $$
\end{Corollary}

\begin{proof} By definition, $e^\blam\Mmu[\blam]=e^\blam m^{\blam}\H$ and
  $\underline{G}^\blam=e^\blam y^{\blam}\H$ so we only need to check the two
  middle equalities. By \autoref{T:eme} there exists an invertible
  element~$f^{\blam}(y)$ such that $e^\blam
  m^{\blam}e^\blam =f^{\blam}(y)e^\blam y^{\blam}$. Consequently,
  $e^\blam m^{\blam}e^\blam \H=e^\blam y^{\blam}\H$. To complete the proof it is
  enough to show that $e^\blam m^{\blam}\in e^\blam y^{\blam}\H$. This is immediate
  because $e^\blam m^\blam=e^\blam u^\blam x^\blam\in e^\blam y^\blam\H$ by
  \autoref{L:ulamExpansion}.
\end{proof}

\begin{Definition} Suppose that $\blam\in\Parts$. Let
  $\pi^\blam\map{\Mmu[\blam]}e^\blam\Mmu[\blam]=\underline{G}^\blam$ be the surjective
  $\H$-module homomorphism given by $\pi^\blam(h)=e^\blam h$, for $h\in\Mmu[\blam]$.
\end{Definition}

\begin{Proposition} \label{P:splits} Suppose that $\blam\in\Parts$. Then the epimorphism
  $\pi^{\blam}$ splits. That is, $\pi^{\blam}$ has a right inverse $\phi^\blam$ and
  $\Mmu[\blam]\cong e^\blam\Mmu[\blam]\oplus\Ker\pi^{\blam}$.
\end{Proposition}

\begin{proof} By \autoref{T:eme}, $e^\blam m^{\blam}e^\blam =f^{\blam}(y)e^\blam y^{\blam}$
  where $f^\blam(y)$ is an invertible element of~$\H$. Define $\phi^\blam$ to be the map
  $$\phi^\blam \map{e^\blam\Mmu[\blam]}\Mmu[\blam]; e^\blam y^\blam h\mapsto m^\blam
  e^\blam f^\blam(y)^{-1} h,$$
  for $h\in\H$. To prove that $\phi^\blam$ is well-defined suppose that $e^\blam
  y^{\blam}h=0$ for some $h\in\H$. By \autoref{C:emH}, there exists $h^\blam\in\H$
  such that $e^\blam m^{\blam}=e^\blam y^{\blam}h^\blam$. Let~$\ast$ be the
  non-homogeneous anti-isomorphism of~$\H$ that fixes each of the
  non-homogeneous generators $T_r$ and $L_s$, for $1\le r<n$ and $1\le s\le n$..
  Then $(e^\blam y^{\blam}h^\blam)^\ast=(h^{\blam})^{\ast}e^\blam y^{\blam}$ because
  $e^\blam$ and $y^{\blam}$ are polynomials in $L_1,\cdots,L_n$  by
  \cite[Proposition 4.8]{HuMathas:GradedCellular} and \autoref{T:BKiso},
  respectively. Therefore,
  $$m^\blam e^\blam f^\blam(y)^{-1} h=(e^\blam y^{\blam}h^\blam)^\ast f^\blam(y)^{-1} h
     = (h^\blam)^\ast f^\blam(y)^{-1}e^\blam y^\blam h=0.$$
  That is, $\phi^\blam(e^\blam y^\blam h)=0$. Hence, $\phi^\blam$ is a
  well-defined $\H$-module homomorphism. Moreover, if $h\in\H$ then
  $$
  (\pi^{\blam}\circ\phi^\blam)(e^\blam y^{\blam}h)
        =e^\blam m^{\blam}e^\blam f^{\blam}(y)^{-1}h
        = e^\blam y^{\blam}f^{\blam}(y)f^{\blam}(y)^{-1}h=e^\blam y^{\blam}h.
  $$
  That is, $\pi^{\blam}\circ\phi^{\blam}$ is the identity map on $e^\blam\Mmu[\blam]$.
  Hence, $\pi^\blam$ splits as claimed.
\end{proof}

\begin{Corollary} \label{C:em=me} Suppose that $\blam\in\Parts$. Then $\phi^{\blam}$ induces an
  $\H$-module isomorphism $ e^\blam m^{\blam}\H\cong m^{\blam}e^\blam \H$.
\end{Corollary}

\begin{proof} This follows directly from the proof of \autoref{P:splits}. In
  fact, we have that $\phi^{\blam}\bigl(e^\blam m^{\blam}\H\bigr)=m^{\blam}e^\blam \H$.
\end{proof}

\subsection{Cyclotomic Schur algebras}
\label{S:DJMSchur}
We are now ready to show that $\Sch[n]$ is Morita
equivalent to the corresponding cyclotomic Schur algebras introduced in
\cite{DJM:cyc,BK:HigherSchurWeyl}.

\begin{Definition}[\cite{DJM:cyc,BK:HigherSchurWeyl}]\label{CycSchur}
  The \textbf{cyclotomic Schur algebra} is the algebra
  $$\ScDJM=\End_{\H}\(\bigoplus_{\bmu\in\Parts}\Mmu\).$$
\end{Definition}

Again, we write $\ScDJM$ to emphasize that $\ScDJM$ is not $\Z$-graded. Note that the algebra
$\ScDJM$  depends implicitly on the dominant weight~$\Lambda$.

By \cite[Corollary~6.18]{DJM:cyc}, $\ScDJM$ is a quasi-hereditary cellular algebra with Weyl
modules $\WDJM$ and irreducible modules $\LDJM$, for $\blam,\bmu\in\Parts$. By
\cite[Theorem~2.11]{LM:AKblocks} and \cite{Brundan:degenCentre} the blocks of $\ScDJM$ are
again labelled by~$Q^+_n$, however, the direct summands of $\Mmu$ do not necessarily belong
to the same block so it is difficult to describe the blocks of $\ScDJM$ explicitly; however,
see~\cite[Theorem~4.5]{M:seminormal}.

Recall the graded Young modules $Y^\bmu$, for $\bmu\in\Parts$, from
\autoref{D:YoungModule}.

\begin{Lemma}\label{L:UngradedYoundMods}
  Suppose that $\bmu\in\Parts$. Then $\Mmu\cong
  \underline{Y}^\bmu\oplus\bigoplus_{\blam\gdom\bmu}(\underline{Y}^\blam)^{m_{\blam\bmu}}$
  for some integers $m_{\blam\bmu}\in\N$.
\end{Lemma}

\begin{proof} By \cite[(3.5)]{M:tilting} there is a family of pairwise
  non-isomorphic (ungraded) indecomposable $\H$-modules
  $\set{\underline{y}^\bmu|\bmu\in\Parts}$ that are uniquely determined, up to
  isomorphism, by the property that
  \begin{equation}\label{E:DJMyoungmodule}\Mmu\cong \underline{y}^\bmu\oplus
     \bigoplus_{\blam\gdom\bmu}(\underline{y}^\blam)^{\oplus m_{\blam\bmu}}
  \end{equation}
  for some (in general, unknown) integers $m_{\blam\bmu}\in\N$. We show
  by induction on the dominance ordering that
  $\underline{Y}^\bnu\cong\underline{y}^\bnu$, for all $\bnu\in\Parts$.

  First suppose that $\bmu\in\Parts$ is maximal in the  dominance ordering. Then
  $\Mmu=\underline{y}^\bmu$ by~\autoref{E:DJMyoungmodule} and $G^\bmu=Y^\bmu$ by
  \autoref{P:YoungProperties}(c). Therefore,
  $\underline{Y}^\bmu\cong\underline{y}^\bmu$ since $\underline{G}^\bmu$ is a
  summand of~$\Mmu$ by \autoref{P:splits}.

  Now suppose that $\bmu$ is not maximal in the dominance ordering. Then $Y^\bmu$ is
  isomorphic to an indecomposable direct summand of~$G^\bmu$ by
  \autoref{P:YoungProperties}(c).  Therefore, there exists a
  multipartition $\blam\gedom\bmu$ such
  that $\underline{Y}^\bmu\cong\underline{y}^\blam$ by
  \autoref{P:splits} and~\autoref{E:DJMyoungmodule}. By induction, if
  $\bnu\gdom\bmu$ then $\underline{y}^\bnu\cong\underline{Y}^\bnu$, so this
  forces $\blam=\bmu$ by \autoref{P:YoungProperties}(b). That is,
  $\underline{Y}^\bmu\cong\underline{y}^\bmu$ as claimed. This completes the proof.
\end{proof}

\begin{Theorem}\label{T:DJMEquivalence}
  Suppose that $\mz$ is a field and that $e=0$ or $e>n$. Then there is an equivalence
  of highest weight categories
  $$\EDJM:\UnS\Mod \bijection \ScDJM\Mod$$
  such that $\EDJM(\underline{\Delta}^\blam)\cong\WDJM$ and
  $\EDJM(\underline{L}^\bmu)\cong\LDJM$, for all $\blam,\bmu\in\Parts$.
\end{Theorem}

\begin{proof} By \autoref{L:UngradedYoundMods}, the algebra
  $$\End_{\H}\(\bigoplus_{\bmu\in\Parts} \underline{Y}^\bmu\)$$
  is the basic algebra of $\UnS$ and it is also the basic algebra of
  $\ScDJM$. Hence the result follows because (ungraded) basic algebras are
  unique up to isomorphism, as discussed in \autoref{S:Basic}.
\end{proof}

Using the combinatorics of the cellular bases of the algebras $\ScDJM$ and
$\Sch[n]$ it is easy to see that if $\mz$ is a field then $\dim\UnS\le\dim\ScDJM$.
Moreover, this inequality is strict except for small~$n$; compare with
\autoref{R:Semistandard}. In particular, the algebras~$\UnS$ and~$\ScDJM$ are
not isomorphic in general.

\begin{Corollary}
  Suppose that $\mz$ is a field and that $e=0$ or $e>n$. Then, up to Morita
  equivalence, $\ScDJM$ depends only on~$e$, $\Lambda$ and the characteristic
  of~$\mz$.
\end{Corollary}

In particular, if $e=0$ or $e>n$ then the decomposition numbers of the
degenerate and non-degenerate cyclotomic Schur algebras depend only~$e$,
$\Lambda$ and the characteristic of the field. This generalizes
\cite[Corollary~6.3]{BK:GradedKL}, which is the analogous result for the
cyclotomic Hecke algebras (without any restriction on~$e$).

Using \autoref{L:UngradedYoundMods} it is not hard to show that the degenerate
and non-degenerate cyclotomic Schur algebras are isomorphic over any field when
$e=0$ or $e>n$. Gordon and
Losev~\cite[Proposition~6.6]{GordonLosev:CycCatOCherednik} have constructed an
explicit isomorphism between these algebras over~$\C$ when $e=0$, extending
Brundan and Kleshchev's isomorphism \autoref{T:BKiso}.

\subsection{Signed permutation modules}\label{S:SignedPerms}
The arguments in the last  two sections apply equally well to the
\textit{signed}, or \textit{twisted}, permutation modules defined in
\cite[\S4]{M:tilting}. Mirroring the definitions in
\autoref{S:CycPermModules}, for $\bmu\in\Parts$ define
$$u_\bmu=\prod_{k=1}^n\prod_{s=\comp_{\t_\bmu}(k)+1}^\ell (L_k-\xi^{(\kappa_s)})
\quad\text{and}\quad x_\bmu=\sum_{w\in\Sym_{\bmu'}}T_w.$$
Note that $\Sym_{\bmu'}$ is the \textit{column} stabiliser of $\t_\bmu$.
Let $n_\bmu=u_\bmu x_\bmu$ and define $\Nmu=n_\bmu\H$. (For similar
reason as in \autoref{R:warning}, this module is denoted $N^{\bmu'}$ in
\cite[\S4]{M:tilting}.) By \cite[Proposition~4.3]{M:tilting}, there is
an isomorphism of ungraded algebras
$\ScDJM\cong\End_{\H}(\bigoplus_\bmu\Nmu)$. In fact,
\autoref{T:DualSchur} should be considered as a graded analogue of this
result. Let $\set{\underline{\Delta}_{\blam}^{DJM}|\blam\in\Parts}$ and
$\set{\UnLSW|\bmu\in\Parts}$ be the sets of standard modules and simple
modules, respectively, of the quasi-hereditary algebra
$\End_{\H}(\bigoplus_\bmu\Nmu)$, which we consider as $\ScDJM$-modules.

By adapting the arguments leading to \autoref{L:UngradedYoundMods} we
obtain the following.

\begin{Lemma}\label{L:UngradedSignedYoundMods}
  Suppose that $\bmu\in\Parts$. Then $\Nmu\cong
  \underline{Y}_\bmu\oplus\bigoplus_{\bmu\gdom\blam}(\underline{Y}_\blam)^{n_{\blam\bmu}}$
  for some integers $n_{\blam\bmu}\in\N$.
\end{Lemma}

\subsection{Positivity}
In this section we show that, in characteristic zero, the graded decomposition numbers of
$\Sch$ are polynomials, rather than Laurent polynomials, by showing
that $\Sch$ is graded Morita equivalent to one of the cyclotomic
Quiver Schur algebras introduced by Stroppel and
Webster~\cite{StroppelWebster:QuiverSchur}. Stroppel and Webster's
results can be summarised as follows.

\begin{Theorem}[\protect{Stroppel and Webster~\cite{StroppelWebster:QuiverSchur}}]
  Suppose that $\mz=\mathbb{C}$ and $\Lambda\in P^+$. Then there exists
  a graded cellular $\mz$-algebra $\ScSW$ such that:
  \label{T:StroppelWebster}
  \begin{enumerate}
    \item As ungraded algebras $\UnScSW\cong\ScDJM$.
    \item The algebra $\ScSW$ is a quasi-hereditary graded cellular algebra with
    weight poset $(\Parts,\ledom)$,
    standard modules $\set{\WSW|\blam\in\Parts}$ and simple modules
    $\set{\LSW|\bmu\in\Parts}$ such that
    $\UnWSW\cong\underline{\Delta}_{\blam}^{DJM}$ and $\UnLSW\cong\underline{L}_{\bmu}^{DJM}$ as $\ScDJM$-modules.
    \item If $\blam,\bmu\in\Parts$ then $[\WSW:\LSW]_q\in\delta_{\blam\bmu}+q\N[q]$.
  \end{enumerate}
\end{Theorem}

By \autoref{T:DJMEquivalence} there is an ungraded equivalence between
the module categories of $\Sch$ and $\ScSW$. The next result says that
this lifts to a graded equivalence.

\begin{Theorem}\label{T:SWEquivalence}
  Suppose that $\mz$ is a field and that $e=0$ or $e>n$. Then there is a graded equivalence
  of highest weight categories
  $$\SWEquiv:\SS^n_\Lambda\Mod \bijection \ScSW\Mod$$
  such that $\SWEquiv(\Delta_\blam)\cong\WSW$ and
  $\SWEquiv(L_\bmu)\cong\LSW$, for all $\blam,\bmu\in\Parts$.
\end{Theorem}

\begin{proof}
  Recall from \autoref{S:SignedPerms} that $\Nmu$ is a signed permutation
  module and that $\ScDJM\cong\End_{\H}(\bigoplus_\bmu \Nmu)$.
  By \cite[Theorem~6.3]{StroppelWebster:QuiverSchur}, the module $\Nmu$ has a graded lift
  $N_\blam$, which is a graded $\R$-module. The Stroppel-Webster
  cyclotomic quiver Hecke algebra is (isomorphic to) the algebra
  $\ScSW=\ZEnd_{\R}(\bigoplus_\blam N_\blam)$. (Stroppel and Webster first
  define their algebra geometrically and then show that it is isomorphic
  to this algebra.) As the modules $\Nmu$ are graded, the ungraded Schur
  functor from $\UnScSW\Mod$ to $\H\Mod$ automatically lifts to a graded
  Schur functor $\SWFun\map{\ScSW\Mod}\R\Mod$. Let~$\PSW$ be the (graded) projective
  cover of~$\LSW$ and set $\YSW=\SWFun(\PSW)$, for $\bmu\in\Parts$. Then
  $\UnYSW$ is a direct summand of~$\Nmu$ by \cite[Proposition~4.4]{M:tilting}. By
  \autoref{L:UngradedSignedYoundMods}, up to shift, every indecomposable summand
  of~$N_\blam$ is isomorphic to $Y_\bmu$, for some $\bmu\in\Parts$. Arguing
  by induction on the dominance order it follows that there exist
  integers $a_\bmu\in\Z$ such that $\YSW\cong Y_\bmu\<a_\bmu\>$, for all
  $\bmu\in\Parts$.

  Let $c^{\text{SW}}_{\blam\bmu}(q)=\Dim \ZHom_{\ScSW}(\PSW,\PSW[\blam])$
  be a graded Cartan number of~$\ScSW$, for $\blam,\bmu\in\Parts[\beta]$.
  Since Schur functors are fully faithful on projective modules,
  \begin{align*} c_{\blam\bmu}(q)
    &=\Dim \ZHom_{\Sch[\beta]}(P_\bmu,P_\blam)
     =\Dim \ZHom_{\R[\beta]}(Y_\bmu,Y_\blam)\\
    &=\Dim \ZHom_{\R[\beta]}\(\YSW[\bmu]\<-a_\bmu\>,\YSW[\blam]\<-a_\blam\>)
     =q^{a_\bmu-a_\blam}\Dim \ZHom_{\R[\beta]}(\YSW[\bmu],\YSW[\blam])\\
    &=q^{a_\bmu-a_\blam}c^{\text{SW}}_{\blam\bmu}(q).
  \end{align*}
  By \autoref{T:StroppelWebster}(b), $\ScSW$ is a graded cellular algebra. Applying
  \autoref{C:CartanSymmetric}, we deduce that
  $c_{\blam\bmu}(q)=q^{a_\bmu-a_\blam}c^{\text{SW}}_{\bmu\blam}(q)
                   =q^{2(a_\bmu-a_\blam)}c_{\blam\bmu}(q)$. On the other
  hand, $\Sch[n]$ is also a graded cellular algebra, so
  \autoref{C:CartanSymmetric} now forces
  $a_\blam=a_\bmu$ whenever $c_{\blam\bmu}(q)\ne0$. For $\beta\in Q^+_n$
  the algebra $\Sch$ is an indecomposable block of $\Sch[n]$ by \autoref{T:SBlocks}. So,
  there exist integers $a_\beta\in\Z$ such that $a_\blam=a_\beta$
  whenever $\blam\in\Parts[\beta]$.  Therefore, $\ScSW$ is graded Morita
  equivalent to the algebra
  $$\bigoplus_{\beta\in Q^+_n}\ZEnd_{\R}(\bigoplus_{\bmu\in\Parts[\beta]}\YSW\)
      \cong\bigoplus_{\beta\in Q^+_n}
              \ZEnd_{\R}(\bigoplus_{\bmu\in\Parts[\beta]}Y_\bmu\<a_\beta\>\)
      \cong\bigoplus_{\beta\in Q^+_n}
              \ZEnd_{\R}(\bigoplus_{\bmu\in\Parts[\beta]}Y_\bmu\).
  $$
  The last algebra in the displayed equation is a graded basic algebra for~$\Sch$.
  Hence, there is a graded equivalence $\ScSW\Mod\longrightarrow\Sch[n]\Mod$.
  The remaining claims follow using \autoref{T:DJMEquivalence}.
\end{proof}

\begin{Corollary}\label{C:Positivity}
  Suppose that $\mz=\C$ and that $\blam,\bmu\in\Parts$. Then
  $d_{\blam\bmu}(q)\in\delta_{\blam\bmu}+q\N[q]$.
\end{Corollary}

\begin{proof}
  By \autoref{T:SWEquivalence} with \autoref{T:StroppelWebster}(c),
  $[\Delta_\blam:L_\bmu]_q=[\WSW:\LSW]_q\in\delta_{\blam\bmu}+q\N[q]$. Applying
  \autoref{T:DualSchur} gives the result.
\end{proof}

\section{Parabolic category~$\Ocal$ and the Fock space}\label{S:CatO}
%  In this chapter we prove \autoref{THEOREMB} and \autoref{THEOREMC}
%  from the introduction. That is, when $e=0$ and $\mz=\C$ we show that
%  $\Sch[n]\Mod$ is Koszul and that the gradings on $\R$ and $\Sch[n]$
%  are compatible  with previously known Koszul grading on parabolic
%  category~$\Ocal$ for~$\gl$. As this chapter is quite technical we
%  explain the strategy that we use to prove these results.

  \autoref{T:DJMEquivalence} shows that $\Sch[n]\Mod$ induces a grading on the
  category of finite dimensional modules for the cyclotomic Schur algebras
  $\ScDJM$. On the other hand, for the degenerate case ($e=0$) in characteristic zero
  Brundan and Kleshchev~\cite{BK:HigherSchurWeyl} have constructed an equivalence of
  categories $\O[n]\bijection\ScDJM\Mod$, where
  $\O[n]=\bigoplus_{\beta\in Q^+_n}\O[\beta]$ is a sum of certain integral
  blocks of the BGG parabolic category~$\Ocal$ for
  $\mathfrak{gl}_N(\C)$. Deep results of Beilinson, Ginzburg and
  Soergel~\cite[Theorem 1.1.3]{BGS:Koszul} and
  Backelin~\cite[Theorem~1.1]{Backelin:Koszul} show that $\O[n]$ admits a
  Koszul grading, so it follows that $\ScDJM$ can be endowed with a Koszul
  grading as well. The aim of this chapter is to show that the Koszul
  grading on category $\Ocal$ coincides with the grading on $\Sch[n]$.
  More precisely, we prove \autoref{THEOREMC} and show that Brundan and
  Kleshchev's equivalence can be lifted to a graded equivalence
  $\O[n]\bijection\Sch[n]\Mod$. We now give a brief outline of the main
  arguments in this chapter

  Fix $\beta\in Q^+$ and let $\SO$ be the basic algebra of
  category~$\O[\beta]$ and $\SFlat$ be the basic algebra of $\Sch$. (We
  define the category $\O[\beta]$ in \autoref{S:HigherSchurWeyl} below.)
  We define these algebras below to be the graded endomorphism algebras
  of minimal projective generators in their respective categories.  The
  algebras~$\SO$ and $\SFlat$ are both graded and, as remarked above,
  $\SO$ is a Koszul algebra and $\UnSO\cong\UnSFlat$ as ungraded
  algebras.  Unfortunately, we are not able to compare the
  algebras~$\SO$ and~$\SFlat$ directly.  Instead, the idea is to compare
  the endomorphism algebras of their minimal \textit{prinjective
  generators}.  Let $\set{\PO|\bmu\in\Parts[\beta]}$ and
  $\set{P^\bmu|\bmu\in\Parts[\beta]}$ be complete sets, up to shift, of
  the pairwise non-isomorphic projective indecomposable $\SO$-modules
  and $\Sch[\beta]$-modules, respectively. We choose the labelling
  of these modules so that
  $\bigoplus_{\bmu\in\Klesh[\beta]}\PO$ is a minimal prinjective
  generator for $\SO$ and that
  $\bigoplus_{\bmu\in\Klesh[\beta]}P^\bmu$  is a minimal prinjective
  generator for~$\Sch[\beta]$. Define
  $$\RO=\ZEnd_{\SO}\(\bigoplus_{\bmu\in\Klesh[\beta]}\PO\)^\op
    \quad\text{and}\quad
    \RFlat=\ZEnd_{\Sch[\beta]}\(\bigoplus_{\bmu\in\Klesh[\beta]}P^\bmu\)^\op.
  $$
  Then $\RO$ and $\RFlat$ are both graded basic algebras. Moreover, on
  forgetting the gradings, $\UnRO\cong\UnRFlat$ is the basic algebra of
  $\H[\beta]$.  Using~\cite{BK:GradedDecomp} we can
  determine the graded decomposition numbers of~$\R[\beta]$ and the
  graded decomposition numbers of~$\O[\beta]$ can be computed using the
  results of~\cite{Backelin:Koszul,BGS:Koszul}. In fact, it turns out
  that the decomposition matrices of $\RO$ and $\RFlat$ are equal, which implies
  that
  $$\Dim\RFlat=\sum_{\blam,\bmu\in\Klesh[\beta]}c_{\blam\bmu}(q)=\Dim\RO.$$
  The next step is to explicitly construct a homogeneous basis of $\RO$.
  The key point is that because~$\SO$ is Koszul the prinjective modules
  $\PO$, for $\bmu\in\Klesh[\beta]$, are \textit{rigidly graded} by
  \autoref{P:Rigidity}. Consequently, the socle, radical and grading
  filtrations of~$\PO$ coincide. As $\PO$ is rigid, the
  $\Sch[\beta]$-module $P^\bmu$ is also rigid, so that the radical
  and socle filtrations of $P^\bmu$ coincide. Using
  \autoref{C:Positivity} we show that the corresponding projective
  indecomposable module~$\PFlat$ for the basic algebra~$\SFlat$ is also
  rigidly graded.  This observation allows us to use the radical
  filtrations of~$\PO$ and $P^\bmu$ to construct explicit bases of~$\RO$
  and $\RFlat$. As a consequence, it follows that $\RO\cong\RFlat$ as
  graded algebras.

  The argument so far is based entirely on the prinjective modules and
  it says nothing about the projective modules $\PO$ and $P^\bmu$, for
  $\bmu\in\Parts[\beta]\setminus\Klesh[\beta]$. Nonetheless, we have
  essentially completed the proof. There is a graded Schur functor
  $\SOFun\map{\SO\Mod}\R[\beta]\Mod$ and the modules $\YO=\SOFun(\PO)$,
  for~$\bmu\in\Parts[\beta]$, are graded lifts of the Young modules
  considered in~\autoref{S:DJMSchur}. Therefore,
  $\YO\cong Y^\bmu\<k_\bmu\>$, for some integers $k_\bmu\in\Z$, because graded
  lifts of indecomposable modules are unique up to grading shift, when
  they exist. Using BGG reciprocity we show that these shifts are constant
  on~$\Klesh[\beta]$ (in fact, $k_\bmu=0$ for all~$\bmu$). As Schur
  functors are fully faithful on projectives, we therefore obtain graded
  algebra isomorphisms
  $$\SO\cong\ZEnd_{\O[\beta]}\Big(\bigoplus_{\bmu\in\Parts[\beta]}\PO[\bmu]\Big)^\op
  \cong\ZEnd_{\R[\beta]}\Big(\bigoplus_{\bmu\in\Parts[\beta]}\YO[\bmu]\Big)^\op
  \cong\ZEnd_{\R[\beta]}\Big(\bigoplus_{\bmu\in\Parts[\beta]}Y^\bmu\Big)^\op
  \cong\SFlat.$$
  Hence, $\SFlat$ is Koszul, and $\Sch\Mod$ is a Koszul category,
  as claimed.

  To make this argument work we need closely related irreducible
  and projective indecomposable modules in several categories. For the
  readers convenience, we summarise this notation now.

  $$\begin{array}{*5c}\toprule
    \text{Algebra} & \text{Irreducible}  &\text{PIM}
          &\text{Index set}\\\midrule
      \SO         &   \LO   & \PO  &\bmu\in\Parts[\beta]\\
      \Sch[\beta] &   L^\bmu & P^\bmu &\bmu\in\Parts[\beta]\\
      \SFlat      &   \LFlat  & \PFlat &\bmu\in\Parts[\beta]\\
    \bottomrule
  \end{array}$$
  We use similar notation for the functors between the module
  categories of these algebra, for the corresponding ungraded modules
  and for the cell modules of these algebras. Even though a large
  amount of notation is needed, we hope that this consistent pattern for
  labelling these modules will help the reader to understand our
  results.

\subsection{Parabolic category~$\Ocal$}\label{S:HigherSchurWeyl}
Following~\cite{BK:HigherSchurWeyl,BK:DegenAriki}, the first step in the
program outlined above is to use Brundan and Kleshchev's \textit{higher
Schur-Weyl duality} for the degenerate cyclotomic Hecke algebras to
connect the representation theory of the cyclotomic quiver Schur
algebras with the blocks of parabolic category $\Ocal$ for general
linear Lie algebras; see, for example, \cite[Chapter~9]{Humphreys:BGG}. Our
focus is somewhat different to that of
\cite{BK:HigherSchurWeyl,BK:DegenAriki} because we are interested in the
blocks of category~$\Ocal$ that correspond to $\Sch$ for a
particular~$n$. This difference of perspective makes it difficult to
extract the information that we need from
\cite{BK:HigherSchurWeyl,BK:DegenAriki}, so we are generous with the
details.

Fix $\beta\in Q^+_n$ and $\Lambda\in P^+$.
Recall from \autoref{S:KLR} that the dominant weight
$\Lambda=\Lambda_{\kappa_1}+\dots+\Lambda_{\kappa_\ell}$ is determined by our
fixed choice of multicharge $\charge=(\kappa_1,\dots,\kappa_\ell)$. For the rest
of this chapter we assume that
$$\kappa_1\ge\kappa_2\ge\dots\ge\kappa_\ell.$$
There is no loss of generality in making this assumption because we can permute
the numbers in the multicharge
%, and shift all of the generators $L_1,\dots,L_n$
%of $\H$ in \autoref{D:HeckeAlgebras} by the same scalar,
without changing the isomorphism type of~$\H$ or the graded isomorphism
type of~$\R$.

We describe a special case of the definitions and results in
\cite{BK:DegenAriki,BK:HigherSchurWeyl} that is sufficient to capture
all of the blocks of $\Sch[n]$, and hence of $\R$, when $e=0$. Set
$$J=\{\kappa_1+n-1,\kappa_1+n-2,\dots,\kappa_\ell+2-n,\kappa_\ell+1-n\}$$
and let $J_+=J\cup(J+1)$ (note that $J_+=I_+$ in the notation of
\cite{BK:DegenAriki}). The motivation for this definition is that
if $\t\in\Std(\Parts)$ then $\res_k(\t)\in J$, for $1\le k\le n$.
Consequently, if $e(\bi)\ne0$ in $\R$ then $\bi\in J^n$ by
\autoref{T:PsiBases}. This is necessary for \autoref{L:BKBijections}
below.

To help explain what Brundan and Kleshchev do, let $\gl[J_+]$ be the
general linear group of $J_+\times J_+$ matrices, where we index the
rows and columns of the matrices in $\gl[J_+]$ by $J_+$. We label the
fundamental and simple roots of $U(\gl[J_+])$ by $J_+$. In this way, we
identify $\Lambda$ and $\beta$ with weights for~$\gl[J_+]$.

Let $\pi=(\pi_1,\dots,\pi_\ell)$ be the partition defined by
$\pi_c=n+\kappa_c-\kappa_\ell$, for $1\le c\le\ell$, and let
$N=\pi_1+\dots+\pi_\ell$. Note that $\pi_c=\kappa_c+1-\inf(J)$, so this
agrees with the definitions in \cite[\S3.1]{BK:DegenAriki}. (In the notation of
\cite[\S1]{BK:HigherSchurWeyl},~$\pi$ corresponds to the partition
$(q_1\ge\dots\ge q_\ell)$.) Consider the Lie algebra $\gl$ of all
$N\times N$ matrices over $\C$.  Let $\h$ be the standard Cartan subalgebra
of diagonal matrices in~$\gl$ and let $\b\supseteq\h$ the Borel
subalgebra of upper triangular matrices. Define~$\p$ to be the standard
parabolic subalgebra of $\gl$ with Levi subalgebra
$\gl[\pi]=\gl[\pi_1]\oplus\dots\oplus\gl[\pi_\ell]$, so that
$\p=\<\b,\gl[\pi]\>$.  By definition, $\p$ depends only on~$\pi$ and hence
on~$\charge$ (or $\Lambda)$, and~$n$.

Let $\O[]=\Ocal^\p$ be the category of all finitely generated
$\gl$-modules that are locally finite dimensional over~$\p$ and
semisimple over~$\h$. This is the usual parabolic analogue of the BGG
category~$\Ocal$ except that we are only allowing modules with integral
weights or, equivalently, integral central characters. The irreducible
modules in category~$\O[]$ are naturally parameterised by highest
weights, however, following~\cite{BK:DegenAriki,BK:HigherSchurWeyl} we
will use a different labelling of the irreducible $\O$-modules that
comes from the categorification of the $\gl[J_+]$-module
$\bigwedge^\pi V=\bigwedge^{\pi_1}V\otimes\dots\otimes\bigwedge^{\pi_\ell}V$
by the blocks of $\O$. Here $V$ is the defining representation
of~$\gl[J_+]$.

Following Brundan and Kleshchev~\cite[\S2]{BK:DegenAriki}, define the
\textbf{$\Lambda$-diagram}, of column shape~$\pi$, to be the justified
array of boxes with rows indexed by
$\{\kappa_1,\kappa_1-1,\dots,\kappa_\ell+1-n\}$, in decreasing order
from top to bottom, and columns indexed by $\{1,\dots,\ell\}$, in
increasing order from left to right, with the rows left justified and
the columns bottom justified. In particular, column~$c$ of the
$\Lambda$-diagram has a node in row~$j$ if and only if $\kappa_c\ge
j\ge\kappa_\ell+1-n$.  (The $\Lambda$-diagrams should not be confused with the
Young diagrams defined in \autoref{S:tableaux}.) A
\textbf{$\Lambda$-tableau} is any filling of the $\Lambda$-diagram by
numbers in~$J_+$. The \textbf{ground state} $\Lambda$-tableau is the
$\Lambda$-tableau with a~$j$ in all of the boxes in row~$j$, whenever
$\kappa_1\ge j\in\kappa_\ell+1-n$.

\begin{Example}\label{Ex:Tab}
Let $n=3$, $\ell=4$, $\charge=(1,0,0,-2)$. Then
$\Lambda=\Lambda_1+2\Lambda_0+\Lambda_{-2}$
and  $\pi=(6,5,5,3)$. The rows of the $\Lambda$-tableaux are indexed by
$\{1,0,-1,\dots,-4\}$, so the ground state $\Lambda$-tableau is
$$\Tableau{{1},{0,0,0},{-1,-1,-1},{-2,-2,-2,-2},{-3,-3,-3,-3},{-4,-4,-4,-4}}.$$
\end{Example}

Let $\Col$ be the set of \textbf{column-strict} $\Lambda$-tableaux, which
are those $\Lambda$-tableaux with strictly decreasing entries, from
top to bottom, in each column. A $\Lambda$-tableau is \textbf{standard}
if it is column strict and its entries are weakly increasing from left
to right in each row. Brundan and Kleshchev \cite[(2.2) and
(2.3)]{BK:DegenAriki} define the \textit{weight} of a $\Lambda$-tableau
and they let $\Col[\beta]$ and $\Std^\Lambda_\beta$ be the sets of
column-strict and standard tableau of weight $\Lambda-\beta$,
respectively. Alternatively, one can use \autoref{L:BKBijections} below
to define the weight of a $\Lambda$-tableau.

Following \cite[(2.50)]{BK:DegenAriki} and \cite[Lemma
5.4]{Brundan:degenCentre}, if $\blam\in\Parts[\beta]$ define the
\textbf{$\Lambda$-tableau} of~$\blam$ to be the
$\Lambda$-tableau~$\T^\blam$ that has
$\cramped{\lambda^{(c)}_{\kappa_c-j+1}+j}$ in row~$j$ and
column~$c\in\{1,2,\dots,\ell\}$, where $\kappa_c\ge
j\ge\kappa_\ell+1-n$. That is, column~$c$ of $\T^\blam$ is obtained by
adding the parts of $\lambda^{(c)}$ to the ground
state $\Lambda$-tableau. In particular, the $\Lambda$-tableau of the
empty multipartition $(0|\dots|0)$ is the ground state
$\Lambda$-tableau. The point of these definitions is that the
$\Lambda$-tableaux naturally index a basis of~$\bigwedge^\pi V$.

\begin{Example}
  Continuing \autoref{Ex:Tab}, some $\Lambda$-tableaux are:
$$\begin{array}{*5{@{\qquad}l}}
  \blam     & (2,1||0|0) & (1|1|0|1) & (1|0|1|1) & (0|0|1^2|1)\\[2mm]
   \T^\blam
     & \Tableau{{3},{1,0,0},{-1,-1,-1},{-2,-2,-2,-2},{-3,-3,-3,-3},{-4,-4,-4,-4}}
     & \Tableau{{2},{0,1,0},{-1,-1,-1},{-2,-2,-2,-1},{-3,-3,-3,-3},{-4,-4,-4,-4}}
     & \Tableau{{2},{0,0,1},{-1,-1,-1},{-2,-2,-2,-1},{-3,-3,-3,-3},{-4,-4,-4,-4}}
     & \Tableau{{1},{0,0,1},{-1,-1,0},{-2,-2,-2,-1},{-3,-3,-3,-3},{-4,-4,-4,-4}}
\end{array}$$
Only the last two of these $\Lambda$-tableaux are standard and,
using \autoref{C:Restricted}, it is easy to see that these $\Lambda$-tableaux
correspond to the Kleshchev multipartitions in this list.
\end{Example}

If $\T$ is a $\Lambda$-tableau let $\col(\T)=(t_1,\dots,t_N)$ be the
\textit{column reading} of~$\T$, that is, the sequence obtained by
reading the entries of~$\T$ in order from top to bottom down the columns
and reading the columns in order from left to right. The symmetric group
acts from the right on such sequences by place permutations.

\begin{Lemma}[\protect{%
              Brundan and Kleshchev~\cite[(2.50) and (2.52)]{BK:DegenAriki}}]
  \label{L:BKBijections}
  Suppose that $\lambda\in P^+$ and $\beta\in Q^+_n$, for $n\ge0$. Then
  the map $\blam\mapsto\T^\blam$ defines a bijection
  $\Parts[\beta]\bijection\Col[\beta]$ that restricts to
  a bijection $\Klesh[\beta]\bijection\Std^\Lambda_\beta$.
  Moreover, if $\blam,\bmu\in\Parts[\beta]$ then
  $\col(\T^\blam)=\col(\T^\bmu)w$, for some $w\in\Sym_N$.
\end{Lemma}

\begin{proof}
  By the remarks above, $\res_k(\t)\in J$, for all $\t\in\Std(\Parts)$ and
  $1\le k\le n$. Therefore, the map $\blam\mapsto\T^\blam$ defines a bijection
  $\Parts[\beta]\bijection\Col[\beta]$ in view of
  \cite[(2.50)]{BK:DegenAriki}. Furthermore,
  by \autoref{C:Restricted} and \cite[(2.52)]{BK:DegenAriki},
  this map restricts to a bijection $\Klesh[\beta]\bijection\Std^\Lambda_\beta$.
  For the final claim, if $\blam,\bmu\in\Parts[\beta]$ then the column readings of
  $\T^\blam$ and $\T^\bmu$ belong to the same $\Sym_N$-orbit in view of
  \cite[(2.3)]{BK:DegenAriki}. Alternatively, all of the statements in
  the lemma follow easily from standard facts about abacuses
  once one realises we can identify $\T^\blam$ with an $\ell$-tuple of
  abacuses corresponding to the multipartition~$\blam$ as in
  \cite[\S3.1]{LM:AKblocks}.
\end{proof}

Brundan and Kleshchev index the irreducible representations of $\O$ by
the $\Lambda$-tableaux. To this end, let $\eps_1,\dots,\eps_N\in\h^*$ be
the standard coordinate functions on~$\h$ so that if $t=(t_{ij})\in\h$
then $\eps_i(t)=t_{ii}$ picks out the $i^{\text{th}}$ diagonal entry
of~$t$. If $\blam\in\Parts[\beta]$ let $\UnLO[\blam]$ be the
irreducible highest weight $\gl[N]$-module of highest weight
\begin{equation}\label{E:TlamWeight}
\wt(\blam)=t_1\eps_1+(t_2+1)\eps_2+\dots+(t_N+N-1)\eps_N
\end{equation}
where $\col(\T^\blam)=(t_1,\dots,t_N)$ is the column reading
of~$\T^\blam$.  (In the notation of \cite{BK:DegenAriki}, $\UnLO[\blam]$
is the module $L_{\T^\blam}$.) By construction, $\UnL[\blam]$ belongs to
$\O[\beta]$.

For $\beta\in Q^+_n$ let $\O[\beta]$ be the Serre subcategory of~$\O[]$
generated by the irreducible $\gl$-modules
$\set{\UnLO[\bmu]|\bmu\in\Parts[\beta]}$. Then $\O[\beta]$ is the full
subcategory of~$\O[]$ consisting of the modules that have all of their
composition factors in $\set{\UnLO[\bmu]|\bmu\in\Parts[\beta]}$.
Brundan~\cite[Theorem~2]{Brundan:degenCentre} shows that $\O[\beta]$ is
an indecomposable block of parabolic category~$\O[]$. All of the blocks of~$\O[]$
can be described this way, however, we are only interested in the blocks
that correspond to~$\Sch[n]$. Accordingly, set
\begin{equation}\label{E:CatOBlocks}
  \O[n]=\bigoplus_{\beta\in Q_n^+}\O[\beta].
\end{equation}
All of these categories have enough projectives. Let $\UnPO[\blam]$ be the
projective cover of $\UnLO[\blam]$ in $\O[\beta]$, for
$\blam\in\Parts[\beta]$.

Following Backelin~\cite{Backelin:Koszul} we now introduce a grading
on~$\O[\beta]$. For any module $M$ let $\iota_M$ be the identity map
on~$M$. The reader might like to recall the definition of Koszul
categories and Koszul duality from~\autoref{S:Koszul}.

\begin{Theorem}[\protect{Backelin~\cite[Theorem~1.1]{Backelin:Koszul}}]
  \label{T:CatOKoszul}
  Suppose that $\Lambda\in P^+$ and $\beta\in Q^+_n$. Then the category
  $\O[\beta]$ is Koszul. Moreover, there exists a Koszul dual
  category $\dualO$ with simple modules $\set{\dualUnLO[\bnu]|\bnu\in\Parts[\beta]}$
  such that
  $$\End_{\O[\beta]}\(\UnPO[\beta]\)^\op\cong
    \Ext_{\dualO}^\bullet\(\dualUnLO[\beta],\dualUnLO[\beta]\),\qquad
       \text{where }
            \UnPO[\beta]=\bigoplus_{\blam\in\Parts[\beta]}\UnPO[\blam]
        \text{ and }
           \dualUnLO[\beta]=\bigoplus_{\blam\in\Parts[\beta]}\dualUnLO[\blam],
  $$
  and where the algebra on the right-hand side is a positively graded
  Koszul algebra under the Yoneda product.  Moreover, this isomorphism
  can be chosen so that it sends $\iota_{\UnPO[\bnu]}$ to $\iota_{\dualLO[\bnu]}$,
  for $\bnu\in\Parts[\beta]$.
\end{Theorem}

\begin{proof}
  The existence of $\dualO$ and such an isomorphism is proved by Backelin
  in~\cite[Theorem~1.1]{Backelin:Koszul}. The isomorphism can be chosen
  so that it identifies $\iota_{\UnPO[\bnu]}$ and $\iota_{\dualLO[\bnu]}$
  by \cite[Remark~3.8]{Backelin:Koszul}.
\end{proof}

Define $\SO=\Ext_{\dualO}^\bullet\(\dualUnLO[\beta],\dualUnLO[\beta])$.
By \autoref{T:CatOKoszul}, $\SO$ is a finite dimensional Koszul algebra
and there is an equivalence of categories $\EFun:
\O[\beta]\bijection\UnSO\Mod$.

Backelin proves \autoref{T:CatOKoszul} more generally for the blocks of parabolic
category~$\Ocal$ for an arbitrary semisimple complex Lie algebra. The
proof that $\O[\beta]$ is Koszul follows easily from fundamental work of
Beilinson, Ginzburg and Soergel~\cite[Theorem~1.1.3]{BGS:Koszul}, which
considers the Borel-parabolic and singular-regular case whereas Backelin
considers the more general parabolic-parabolic and singular-singular
case. Given \cite{BGS:Koszul}, the deepest statement in
\autoref{T:CatOKoszul} is the explicit description of Koszul duality for
the categories $\O[\beta]$ and $\dualO$. In
\autoref{S:CatODecompNumbers} we give Backelin's description of $\dualO$
and use it to compute the graded decomposition numbers of $\O[\beta]$.
Mazorchuk~\cite[Theorem~7.2]{Mazorchuk:CatLinTiltingComplexes} has given
an algebraic proof of \autoref{T:CatOKoszul} starting from the
Kazhdan-Lusztig conjecture, which is known for~$\O[\beta]$. In fact,
Mazorchuk shows that the algebra $\SO$ is standard Koszul.

By \autoref{T:CatOKoszul}, $\SO=\bigoplus_{d\ge0}\SO[\beta,d]$ is Koszul
so its degree zero component $\SO[\beta,0]$ is semisimple by
\autoref{D:Koszul}. By \autoref{T:CatOKoszul}, the irreducible
$\SO$-modules are labelled by $\Parts[\beta]$ and are concentrated in
degree zero. For $\bnu\in\Parts[\beta]$ let
$\LO[\bnu]=\iota_{\PO[\bnu]}\SO[\beta,0]$ be the unique
irreducible~$\SO$-module that is concentrated in degree zero such
that $\LO[\bnu]$ is isomorphic to $\EFun(\UnLO[\bnu])$ when we forget the grading.
We are abusing notation here because $\LO[\bnu]$ is one dimensional so
that $\UnLO[\bnu]$ is \textit{not} the module obtained from~$\LO[\bnu]$
by forgetting the grading.
Let $\SO[\beta,+]=\bigoplus_{d>0}\SO[\beta,d]$. Then, by \autoref{T:CatOKoszul},
\begin{equation}\label{E:LOquotients}
  \SO/\SO[\beta,+]\cong\bigoplus_{\bnu\in\Parts[\beta]}\LO[\bnu].
\end{equation}
Let $\PO[\bnu]=\iota_{\PO[\bnu]}\SO$ be the projective cover
of~$\LO[\bnu]$ in $\SO\Mod$, so that $\PO[\bnu]$ is isomorphic to
$\EFun(\UnPO[\bnu])$ when we forget the grading. As with $\UnLO$ the
module $\UnPO$ is not obtained from $\PO$ by just forgetting the
grading.

The ungraded category $\Ocal$ has a duality $\underline\diamond$ that
is induced by the anti-isomorphism of $\gl$ that maps a matrix to its
transpose. The duality $\underline\diamond$ restricts to a duality,
also denoted $\underline\diamond$, on the (ungraded) subcategories~$\dualO$
and~$\O[\beta]$.  Since $(\dualUnLO[\bnu])^{\underline\diamond}\cong\dualUnLO[\bnu]$
taking duals induces natural isomorphisms
\begin{equation}\label{E:SOduality}
  \Ext_{\dualO}^k\(\dualUnLO[\blam], \dualUnLO[\bnu]\)
    \cong\Ext_{\dualO}^k\(\dualUnLO[\bnu], \dualUnLO[\blam]\),
\end{equation}
for $\blam,\bnu\in\Parts[\beta]$ and $k\in\Z^{\geq 0}$. Therefore,
$\underline\diamond$ induces a homogeneous anti-isomorphism $\theta$
on~$\SO$. If $M$ is a graded $\SO$-module let $M^\diamond=\Hom_\C(M,\C)$ be
the graded dual of~$M$ where the $\SO$-action is given by twisting
by~$\theta$. Note that $\Dim M^\diamond=\Dim[q^{-1}]M$.

By \autoref{T:CatOKoszul},
$\set{\iota_{\dualUnLO[\bnu]}|\bnu\in\Parts[\beta]}$ is a complete set of
pairwise orthogonal homogeneous idempotents that sum to~$1$ in~$\SO$.
For any finite dimensional $\SO$-module, define the \textbf{graded
character} of~$M$ to be
$$ \ch M:=\sum_{\bnu\in\Parts[\beta]}\dim_q(M\iota_{\dualUnLO[\bnu]})\bnu\in\Z[q,q^{-1}][\Parts[\beta]].
$$
In (ungraded) category~$\O[]$,
$(\dualUnLO[\bnu])^{\underline\diamond}\cong\dualUnLO[\bnu]$. Therefore,
$\iota_{\dualUnLO[\bnu]}$ is fixed by $\theta$, for
$\bnu\in\Parts[\beta]$, and so $\ch M^\diamond=\ch[q^{-1}]M$.

\begin{Lemma}\label{L:DualSimples}
  Suppose that $\bmu\in\Parts[\beta]$. Then
  $(\LO)^\diamond\cong\LO$ as graded $\SO$-modules.
\end{Lemma}

\begin{proof}
  As $\theta$ is homogeneous, and $\LO$ is concentrated in degree
  zero, $\(\LO)^\diamond\cong \LO[\bnu]$ for some
  $\bnu\in\Parts[\beta]$. By~\autoref{E:LOquotients}, $\ch \LO[\blam]=\blam$, for all
  $\blam\in\Parts[\beta]$. Therefore, by the last paragraph,
  $\bnu=\ch(\LO)^\diamond=\ch[q^{-1}]\LO[\bmu]=\bmu$, as required.
\end{proof}

In \autoref{P:ProjInjCover} we proved that $P^\bmu$ is a prinjective
module for all $\bmu\in\Klesh$. We now prove the corresponding result in
category~$\Ocal$.  Recall from \autoref{S:Koszul} that a module is
rigidly graded if its radical, socle and grading filtrations coincide.
The following rigidity result is crucial for the proof of
\autoref{THEOREMC}.

\begin{Proposition}\label{P:RigidPIMs}
  Suppose that $\beta\in Q^+_n$ and $\bmu\in\Klesh[\beta]$. Then $\PO$
  is a prinjective $\SO$-module. Consequently, $\PO$ is rigidly graded.
\end{Proposition}

\begin{proof}
   By \autoref{L:DualSimples}, $(\PO)^\diamond$ is an injective hull of~$\LO$.
   If $\bnu\in\Klesh[\beta]$ then
   $(\UnPO[\bnu])^{\underline\diamond}\cong\UnPO[\bnu]$ by \cite[(2.52),
   Lemma 3.2]{BK:DegenAriki}. Therefore, $\UnPO$ is a prinjective
   $\UnSO$-module. Consequently,
   $(\PO[\bmu])^\diamond\cong\PO[\bmu]\<k_\bmu\>$ for some integer
   $k_\bmu\in\Z$, because graded lifts of indecomposable modules
   are unique up to shift. Hence, $\PO$ is a prinjective $\SO$-module. It
   follows that $\soc\PO[\bmu]\cong\LO\<-k_\bmu\>$ and
   $\PO/\rad\PO\cong\LO$ are both irreducible $\SO$-modules.  By
   \autoref{T:CatOKoszul} $\SO$ is Koszul, so $\PO$ is rigidly graded
   by \autoref{P:Rigidity}(b).
\end{proof}

Let $\UnDelO$ be the parabolic Verma module of highest weight~$\wt(\blam)$
in~$\O[\beta]$.  Then there is an indecomposable graded $\SO$--module
$\DelO$ such that $\DelO\cong\EFun(\UnDelO)$ by
\cite[Proposition~3.5.7]{BGS:Koszul} and the proof of
\cite[Proposition~3.2]{Backelin:Koszul}. As remarked above.
Mazorchuk~\cite[Theorem~5.1]{Mazorchuk:CatLinTiltingComplexes} has shown
that $\SO$ is a standard Koszul algebra. In particular, $\SO$ is graded
quasi-hereditary, so using this framework we can take
$\set{\DelO|\blam\in\Parts[\beta]}$ to be the graded standard modules of
$\SO$. Since $\DelO$ is indecomposable, we fix the grading on $\DelO$ by
requiring that the surjection $\DelO\surjection\LO[\blam]$ is a
homogeneous map of degree zero in $\SO\Mod$.

Since $\SO$ is a standard Koszul algebra with a duality $\diamond$ that
fixes the simple $\SO$-modules, the category $\O[\beta]$ is a
\textit{graded highest weight category with duality} in the sense of
\cite[\S1.2]{CPS:HomDual}. As a consequence, a graded analogue of BGG
reciprocity holds in~$\O[\beta]$. The following result is essentially
\cite[Proposition~1.2.4]{CPS:HomDual}, however, we sketch a proof for
completeness.

\begin{Corollary}[Graded BGG reciprocity]\label{C:SOBGG}
  Suppose that $\blam,\bmu\in\Parts[\beta]$, for $\beta\in Q^+_n$. Then
  $$(\PO:\DelO)_q=[\DelO:\LO]_q.$$
  Consequently,
  $[\PO[\bmu]:\LO[\blam]]_q=\sum_\bnu(\PO:\DelO[\bnu])_q[\DelO[\bnu]:\LO[\blam]]_q$.
\end{Corollary}

\begin{proof}
  The \textit{costandard modules} of $\SO$ are the
  $\NabO=(\DelO)^\diamond$, for $\blam\in\Parts[\beta]$. As $\O[\beta]$
  is a highest weight category, well-known arguments show that
  $\Ext_{\SO}^1(\DelO,\NabO[\bmu]\<k\>)=0$, for all
  $\blam,\bmu\in\Parts[\beta]$ and $k\in\Z$. Therefore, the functor
  $\ZHom_{\SO}(?,\NabO)$ is exact on the subcategory of $\Delta$-filtered
  $\SO$-modules. The projective indecomposable module $\PO$ has a
  $\Delta$-filtration, so if $k\in\Z$ then
  $$(\PO:\DelO\<k\>)_q=\Dim\ZHom_{\SO}(\PO,\NabO\<k\>)
               =[\NabO\<k\>:\LO]_q=[\DelO:\LO\<k\>]_q,$$
  where the last equality follows by applying $\diamond$.  Hence,
  $(\PO:\DelO)_q=[\DelO:\LO]_q$. As $\PO$ has a $\Delta$-filtration the
  multiplicity formula for $[\PO:\LO[\blam]]_q$ now follows easily.
\end{proof}

We are now ready to make the link between parabolic category~$\Ocal$ and the
quiver Schur algebras. The following result is a reformulation of some of
Brundan and Kleshchev's main results
from~\cite{BK:HigherSchurWeyl,BK:DegenAriki}. Our \autoref{THEOREMC} from
the introduction is a graded analogue of this result.

\begin{Theorem}[Brundan and Kleshchev]\label{T:UngradedCatOEquivalence}
    Suppose that $e=0$ and $\mz=\C$. Then there is an
    equivalence of categories $\UnEFun\map{\O[n]}\UnS\Mod$.
    Moreover, $\UnEFun(\UnDelO)\cong\underline{\Delta}^\blam$ and
    $\UnEFun(\UnLO)\cong\underline{L}^\bmu$, for all $\blam,\bmu\in\Parts$.
\end{Theorem}

\begin{proof} By \cite[Theorem~C]{BK:HigherSchurWeyl} there is an equivalence of
  categories from $\O[n]$ to $\ScDJM\Mod$, which sends $\UnDelO$ to $\WDJM$
  and $\UnLO$ to~$\LDJM$, for $\blam,\bmu\in\Parts$. Hence, by
  \autoref{T:DJMEquivalence} and the remarks above, there is an equivalence
  of categories $\UnEFun\map{\O[n]}\UnS\Mod$ such that
  $\UnEFun(\UnDelO)\cong\underline{\Delta}^\blam$ and
  $\UnEFun(\UnLO)\cong\underline{L}^\bmu$.
\end{proof}

By projecting onto the blocks, there are equivalences
$\O[\beta]\bijection\UnS[\beta]\Mod$, for each $\beta\in Q^+_n$.

Ultimately we want to compare the grading on the cyclotomic quiver Schur
algebra $\Sch[\beta]$ with the Koszul grading on $\SO$ coming from
\autoref{T:CatOKoszul}.  Recall that
$d_{\blam\bmu}(q)=[\Delta^\blam:L^\bmu]_q$, for
$\blam,\bmu\in\Parts[\beta]$. Similarly, set
$d^\Ocal_{\blam\bmu}(q)=[\DelO:\LO]_q$.  Recalling
\autoref{L:KoszuldualDecomp}, define ``Koszul dual'' polynomials
$p_{\blam\bmu}(q)$ and $p^\Ocal_{\blam\bmu}(q)$ by the matrix equations
$(p_{\blam\bmu}(q))^\trans=(d_{\blam\bmu}(-q))^{-1}$ and
$(p^{\Ocal}_{\blam\bmu}(q))^\trans=(d^{\Ocal}_{\blam\bmu}(-q))^{-1}$.  Then
$d_{\blam\bmu}(q)=d^\Ocal_{\blam\bmu}(q)$ for all
$\blam\in\Parts[\beta]$ and $\bmu\in\Klesh[\beta]$ if and only if
$p_{\blam\bmu}(q)=p^\Ocal_{\blam\bmu}(q)$ for all
$\blam\in\Parts[\beta]$ and $\bmu\in\Klesh[\beta]$.
To describe these polynomial explicitly, if
$x,y\in\Sym_N$ let $P_{x,y}(t)\in\Z[t]$ be the corresponding
Kazhdan-Lusztig polynomial introduced in~\cite[Theorem~1.1]{KL}.

If $\blam\in\Parts[\beta]$ recall that $\col(\T^\blam)=(t_1,\dots,t_N)$
is the column reading of the $\Lambda$-tableau $\T^\blam$, as defined before
\autoref{L:BKBijections}. As in \cite[\S2.3]{BK:DegenAriki},  define
$w_\blam\in\Sym_n$ to be the unique permutation of minimal length such
that $\col(\T^\blam)w_\blam$ is weakly increasing. Let $Z_\blam$ be the
stabiliser of $\col(\T^\blam)w_\blam$. By \autoref{L:BKBijections},
$Z_\blam$ depends only on $\Lambda$ and $\beta$. Further, $w_\blam$ is a
minimal length left coset representative of $Z_\blam$ in $\Sym_N$ (see
\autoref{S:CatODecompNumbers} for more details.).

%For each $\blam,\bmu\in\Parts[\beta]$, let $p_{\T^\blam\T^\bmu}(q),
%d_{\T^\blam\T^\bmu}(q)$ be the polynomials given in
%\cite[(2.16),(2.17)]{BK:DegenAriki}. We have that
%$(p_{\T^\blam\T^\bmu}(q))^\trans=(d_{\T^\blam\T^\bmu}(-q))^{-1}$.
%Brundan and Kleshchev have explicitly described the polynomials
%$p_{\T^\blam\T^\bmu}(q)$ as parabolic Kazhdan-Lusztig polynomials.  In
%fact, their formula \cite[(2.18)]{BK:DegenAriki} for these polynomials
%coincides exactly with \autoref{T:CatODecompNumbers}.

We prove the following result in \autoref{S:CatODecompNumbers}.  This
proposition may be well-known to experts, however, we include a proof
because we have not found the result in the literature.

\begin{Proposition}\label{T:CatODecompNumbers}
  Suppose that $\blam,\bmu\in\Parts[\beta]$, for $\beta\in Q^+_n$. Then
  $$p^\Ocal_{\blam\bmu}(q) = q^{\ell(w_\bmu)-\ell(w_\blam)}
      \sum_{z\in Z_\blam}(-1)^{\ell(z)} P_{w_\blam z,w_\bmu}(q^{-2})\in\N[q,q^{-1}].$$
\end{Proposition}

By \autoref{L:KoszuldualDecomp}, $p^\Ocal_{\blam\bmu}(q)$ is a graded
decomposition number for the Koszul dual category $\dualO$.  As we
describe in \autoref{S:CatODecompNumbers}, $\dualO$ is also a block of
parabolic category~$\Ocal$, although it is associated to another
parabolic subalgebra~$\q$ of~$\gl$.  When $\q=\b$ this result is
equivalent to \cite[Theorem~3.11.4(i)(iv)]{BGS:Koszul}, which a
consequence of the Kazhdan-Lusztig conjecture for~$\gl$.

Given \autoref{T:CatODecompNumbers}, the next result is essentially a
reformulation of results of Brundan and
Kleshchev\cite{BK:GradedDecomp,BK:DegenAriki}, building on Ariki's
categorification theorem~\cite{Ariki:can}. This result is the key to
comparing the gradings on~$\RO$ and~$\R[\beta]$.

\begin{Theorem}\label{T:RGradedDecomp}
  Suppose that $e=0$, $\beta\in Q^+_n$ and $\mz=\C$. Then
  $d_{\blam\bmu}(q)=d^\Ocal_{\blam\bmu}(q)$
  is a parabolic Kazhdan-Lusztig polynomial, for $\blam\in\Parts[\beta]$
  and $\bmu\in\Klesh[\beta]$. In particular,
  $d_{\blam\bmu}(q)\in\delta_{\blam\bmu}+q\N[q]$.
  %, $p_{\blam,\bmu}(q)\in\Z[q]$ and $  p_{\blam,\bmu}(0)=\delta_{\blam,\bmu}=p_{\blam,\bmu}^{\Ocal}(0)$.
\end{Theorem}

\begin{proof}
%  If $\blam,\bmu\in\Parts[\beta]$ then
%  $$d_{\T^\blam\T^\bmu}(q)=(p_{\T^\bmu\T^\blam}(-q))^{-1}
%       =(p_{\bmu,\blam}^{\Ocal}(-q))^{-1}=d_{\blam\bmu}(q)^{\Ocal}$$
%  by \autoref{T:CatODecompNumbers},
%  Assume further that $\bmu\in\Klesh[\beta]$.
  By \cite[Theorem~5.14 and Corollary~5.15]{BK:GradedDecomp}, the
  Grothendieck group of the finitely generated projective $\R$-modules
  categorifes the integral highest weight module $L(\Lambda)$ for
  $U_q(\SL)$, with the projective indecomposable $\R$-modules
  corresponding to the Lusztig-Kashiwara canonical basis
  of~$L(\Lambda)$. Thus, the graded decomposition numbers of~$\R$
  give the transition matrix between the standard basis
  and the canonical basis of~$L(\Lambda)$.

  As we are assuming that $e=0$, the algebra
  $\R\cong\H$ is a degenerate cyclotomic Hecke algebra. Brundan and
  Kleshchev~\cite[Theorem~3.1 and (2.18)]{BK:DegenAriki} proved the
  analogue of Ariki's categorification theorem in the degenerate case.
  In particular, they showed that the transition matrix between the
  standard basis and the canonical basis of the Fock space is given by
  certain polynomials $\{d_{\T^\blam\T^\bmu}(q)\}$, where
  $\T^\blam\in\Col[\beta]$ and $\T^\bmu\in\Std^\Lambda_\beta$ are the
  column-strict and standard $\Lambda$-tableaux defined before
  \autoref{L:BKBijections}. The two
  papers~\cite{BK:GradedDecomp,BK:DegenAriki} use the same bar
  involution by ~\cite[\S2.5]{BK:DegenAriki} and the remarks
  after~\cite[(3.60)]{BK:GradedDecomp}. (As explained in
  \cite[\S2.2]{BK:DegenAriki}, the infinite dimensional space considered
  in \cite{BK:GradedDecomp} agrees with the finite dimensional space
  considered in \cite{BK:DegenAriki} after truncation.) Therefore, it
  follows by induction on dominance that
  $d_{\blam\bmu}(q)=d_{\T^\blam\T^\bmu}(q)$,
  for all $\blam\in\Parts[\beta]$ and $\bmu\in\Klesh[\beta]$.

  For $\blam,\bmu\in\Klesh[\beta]$ define polynomials
 $p_{\T^\blam\T^\bmu}(q)$ by the matrix equation
  $$\(p_{\T^\blam\T^\bmu}(-q)\)^\trans
        =\(d_{\T^\blam\T^\bmu}(q)\)^{-1}_{\blam,\bmu\in\Klesh[\beta]}.$$
  These polynomials coincide with the polynomials defined by Brundan and
  Kleshchev in \cite[(2.17)]{BK:DegenAriki} by
  \cite[(2.39)]{BK:DegenAriki}.  Brundan and Kleshchev explicitly
  describe the polynomials $p_{\T^\blam\T^\bmu}(q)$ as parabolic
  Kazhdan-Lusztig polynomials in~\cite[(2.18)]{BK:DegenAriki}. In fact,
  their formula for these polynomials coincides exactly with
  \autoref{T:CatODecompNumbers}. That is,
  $p_{\T^\blam\T^\bmu}(q)=p^\Ocal_{\blam\bmu}(q)$, for all
  $\blam,\bmu\in\Klesh[\beta]$. It follows that
  $d_{\blam\bmu}(q)=d_{\T^\blam\T^\bmu}(q)=d^\Ocal_{\blam\bmu}(q)$, for all
  $\blam,\bmu\in\Klesh[\beta]$.
  Therefore, $d_{\bnu\bmu}(q)=d^\Ocal_{\blam\bmu}(q)$, for all
  $\bnu\in\Parts[\beta]$ and $\bmu\in\Klesh[\beta]$, since the
  polynomials indexed by multipartitions in $\Klesh[\beta]$ uniquely
  determine the remaining polynomials under
  categorification by \cite[Corollary~5.15]{BK:GradedDecomp}.

  Finally, observe that
  $d_{\blam\bmu}(q)=d^\Ocal_{\blam\bmu}(q)\in\delta_{\blam\bmu}+q\N[q]$
  because $\SO$ is Koszul.
\end{proof}

\begin{Corollary}\label{C:SameKleshCartan}
  Suppose that $e=0$, $\beta\in Q^+_n$ and $\mz=\C$. Then
  $d_{\blam\bmu}(q)=d^\Ocal_{\blam\bmu}(q)$
  for all $\blam,\bmu\in\Klesh[\beta]$.
\end{Corollary}

\begin{proof}Applying in turn \autoref{C:SOBGG}, \autoref{T:RGradedDecomp} and
  \autoref{C:CartanSymmetric}, shows that
  $$c^\Ocal_{\blam\bmu}(q)
         =\sum_{\bnu\in\Parts[\beta]} d^\Ocal_{\bnu\blam}(q) d^\Ocal_{\bnu\bmu}(q)
         =\sum_{\bnu\in\Parts[\beta]} d_{\bnu\blam}(q) d_{\bnu\bmu}(q)
         =c_{\blam\bmu}(q),
  $$
  as required.
\end{proof}

\subsection{Comparing the KLR and category~$\Ocal$ gradings}\label{S:Comparison}
We want to lift \autoref{T:UngradedCatOEquivalence} to the graded
setting.  As a first step we show that the KLR and category~$\O[\beta]$
gradings induce the same grading on~$\H[\beta]$.  To do this we need a
graded analogue of the Hecke algebra~$\H[\beta]$ with a grading that
comes from category~$\O[\beta]$.

We start by describing~$\SO$ as an endomorphism algebra. Let
$\PO[\beta]=\bigoplus_{\bmu\in\Parts[\beta]}\PO$.

\begin{Lemma}\label{L:CanonicalIso}
  Suppose that $\beta\in Q^+_n$. Then there is an isomorphism of
  graded rings
  $$\SO\cong\ZEnd_{\SO}\(\PO[\beta]\)^\op.$$
  In particular, $\ZEnd_{\SO}(\PO[\beta])$ is Koszul.
\end{Lemma}

\begin{proof}By definition, $\PO$ is the projective cover of $\LO$ so,
  by~\autoref{E:LOquotients}, $\PO$ is a direct summand of $\SO$ by the
  universal property of projective modules.  Therefore,
  $\SO\cong\PO[\beta]$ as a right $\SO$-module and, consequently,
  $\SO\cong\ZEnd_{\SO}(\PO[\beta])^{\op}$ as graded algebras. In fact, by the
  argument of \cite[Corollary~2.5.2]{BGS:Koszul}, the
  isomorphism $\SO\cong\ZEnd_{\SO}(\PO[\beta])$ is unique.
\end{proof}

Henceforth, we identify $\SO$ and $\ZEnd_{\SO}(\PO[\beta])$ via
\autoref{L:CanonicalIso}. Motivated by \autoref{T:UngradedCatOEquivalence}, define
$\eR =\sum_{\bmu\in\Klesh[\beta]}\iota_{\PO}$.  Then $\eR\in\SO$ is a
homogeneous idempotent of degree zero. Now define
$$
\RO=\eR \SO \eR \cong \ZEnd_{\SO}\(\bigoplus_{\bmu\in\Klesh[\beta]}\PO\)^\op.
$$
By~\autoref{E:SchurFunctors} there is an exact graded functor
\begin{equation}\label{E:OSchurFunc}
  \SOFun\map{\SO\Mod}\RO\Mod; M\mapsto M\eR , \qquad\text{ for }M\in\SO\Mod.
\end{equation}
The anti-isomorphism $\theta$ of~\autoref{E:SOduality} fixes
$\SO[\beta,0]$,
the degree zero component of~$\SO$, so $\theta(\eR)=\eR$. Therefore,~$\theta$
restricts to a homogeneous anti-isomorphism of~$\RO$. Abusing
notation, let~$\diamond$ be the corresponding duality on $\R\Mod$. By
construction, there is an isomorphism of functors
$\diamond\circ\SOFun\cong\SOFun\circ\diamond$.

Let $\RO=\bigoplus_{d\in\Z}\RO[\beta,d]$ be the decomposition of~$\RO$ into its
homogeneous components.

\begin{Proposition}\label{P:ROPositive}
  Suppose that $\beta\in Q^+_n$. Then:
  \begin{enumerate}
    \item $\RO$ is a positively graded basic algebra;
    \item $\RO[\beta,0]$ is semisimple;
    \item The ungraded algebras $\UnRO$, $\UnR[\beta]$ and $\H[\beta]$ are Morita equivalent.
  \end{enumerate}
\end{Proposition}

\begin{proof}
  By \autoref{T:CatOKoszul} $\SO$ is a Koszul algebra so it is positively graded
  and its degree zero component is semisimple by \autoref{D:Koszul}. Hence,
  parts~(a) and~(b) follow because $\RO=\eR\SO\eR$. Finally,
  by \autoref{T:UngradedCatOEquivalence}, the algebra
  $\UnRO$ is Morita equivalent to the Hecke algebra $\H[\beta]$.  Hence,
  part~(c) follows because $\UnR[\beta]\cong\H[\beta]$ by
  \autoref{T:BKiso} and the remarks after~\autoref{E:Rblocks}.
\end{proof}

Using \autoref{T:RGradedDecomp}, we now to replicate the results that we
have just proved for the algebra $\RO$ for the KLR algebra~$\R[\beta]$.
To do this we need to work with a basic algebra for~$\R[\beta]$. We
start by defining
$$\SFlat=\ZEnd_{\Sch[\beta]}\(\bigoplus_{\bmu\in\Parts[\beta]}P^\blam\)^\op,$$
which a basic algebra for~$\Sch[\beta]$. Write
$\SFlat=\bigoplus_d\SFlat[\beta,d]$ for the decomposition of $\SFlat$ into its
homogeneous components.

\begin{Lemma}\label{L:SPositive}
  The algebra $\SFlat$ is a positively graded algebra. Moreover, $\SFlat[\beta,0]$
  is semisimple.
\end{Lemma}

\begin{proof}
  This is immediate from \autoref{C:Positivity}.
\end{proof}

Eventually we will show that
$\SO\cong\SFlat$ as graded algebras.
Mirroring the definition of $\RO$, define
$$\RFlat=\ZEnd_{\Sch[\beta]}\(\bigoplus_{\bmu\in\Klesh[\beta]}P^\bmu\)^\op,$$
so that $\RFlat=e^\Lambda_\beta\Sch[\beta]e^\Lambda_\beta$ where
$e^\Lambda_\beta=\sum_{\beta\in\Klesh[\beta]}\iota_{P^\bmu}$. Notice
that $\RFlat$ is positively graded by \autoref{L:SPositive} because it
is a subalgebra of $\SFlat$.

\begin{Lemma}\label{L:RFlatEquiv}
  There is a graded equivalence of categories $\R[\beta]\Mod\cong\RFlat\Mod$.
\end{Lemma}

\begin{proof}
  By \autoref{L:DoubleCentralizer}, the Schur functor $\SFun[\beta]$ is fully
  faithful on projective modules so, by \autoref{D:YoungModule},
  $$ \ZHom_{\Sch[\beta]}(P^\blam,P^\bmu)\cong\ZHom_{\R[\beta]}(Y^\blam,Y^\bmu) $$
  for all $\blam,\bmu\in\Parts[\beta]$. Therefore,
  $\RFlat\cong\ZEnd_{\R[\beta]}(\bigoplus_{\bmu\in\Klesh[\beta]}Y^\bmu)^\op$.
  By \autoref{P:YoungProjCover},
  $\bigoplus_{\bmu\in\Klesh[\beta]}Y^\bmu$ is a minimal progenerator for
  $\R[\beta]$, so $\RFlat$ is graded Morita equivalent to $\R[\beta]$ as
  claimed.
\end{proof}

By definition, $c_{\blam\bmu}(q)=\Dim\ZHom_{\R[\beta]}(Y^\bmu,Y^\blam)$.
Therefore, the proof of \autoref{L:RFlatEquiv}, together with
\autoref{C:SameKleshCartan}, implies that
$\Dim\RFlat=\sum_{\blam,\bmu\in\Klesh[\beta]}c_{\blam\bmu}(q)=\Dim\RO$.

\begin{Lemma}\label{L:DefectBound}
  Suppose that $\blam,\bmu\in\Klesh[\beta]$. Then
  $0\le\deg c_{\blam\bmu}(q)\le 2\defect\beta$ with equalities if and only if
  $\blam=\bmu$, so there homogeneous map. Moreover, $\PO$ and $P^\bmu$
  both have Loewy length~$2\defect\beta$.
\end{Lemma}

\begin{proof}
  First observe that $c_{\blam\bmu}(q)\in\delta_{\blam\bmu}+q\N[q]$ by
  \autoref{C:Positivity}. Moreover,
  $(Y^\bmu)^\circledast\<2\defect\beta\>\cong Y^\bmu$ by
  \autoref{C:DualYoung} if $\bmu\in\Klesh[\beta]$. Therefore,
  $0\le\deg c_{\blam\bmu}(q)\le 2\defect\beta$ with equalities if and only if
  $\blam=\bmu$.

  To prove the remaining claims statement first recall that $\PO$ is rigidly graded
  by \autoref{P:RigidPIMs}. Therefore, $\PO$ has Loewy length
  $\deg c_{\bmu\bmu}(q)=2\defect\beta$ because
  $c^\Ocal_{\bmu\bmu}(q)=c_{\bmu\bmu}(q)$
  by \autoref{C:SameKleshCartan}. On the other hand, $\UnEFun(\PO)\cong
  \UnP$ by \autoref{T:UngradedCatOEquivalence}. In particular, $\PO$ and
  $P^\bmu$ both have the same Loewy length, completing the proof.
\end{proof}

We would like to say that $P^\bmu$ is rigidly graded when
$\bmu\in\Klesh$. We cannot say this, however, because we have only
defined grading filtrations, and rigidly graded modules, for
positively graded algebras.  To remedy this, let
$\EFlat\map{\Sch[\beta]\Mod}\SFlat\Mod$ be the graded equivalence given by
$$M\mapsto\ZHom_{\Sch[\beta]}(\bigoplus_{\bnu\in\Parts[\beta]}P^\bnu,M).$$
Set $\PFlat=\EFlat(P^\bmu)$, for $\bmu\in\Parts[\beta]$. Then $\PFlat$
is an indecomposable projective $\SFlat$-module.

\begin{Proposition}\label{P:PFlatRigid}
  Suppose that $\bmu\in\Klesh[\beta]$. Then $\PFlat$ is rigidly graded.
\end{Proposition}

\begin{proof}As $\EFlat$ is an equivalence, $P^\bmu$ and $\PFlat$ both have Loewy length
  $2\defect\beta$ by \autoref{L:DefectBound}. The grading filtration
  $\{\Gr_d\PFlat\}$ of $\PFlat$ also has length $2\defect\beta$ because
  $\deg c_{\bmu\bmu}(q)=2\defect\beta$ by \autoref{L:DefectBound}.
  Moreover, for $0\le d\le 2\defect\beta$ the module
  $\Gr_d\PFlat/\Gr_{d+1}\PFlat$ is semisimple because it is an
  $\SFlat[\beta,0]$-module and $\SFlat[\beta,0]$ is semisimple by
  \autoref{L:SPositive}. Therefore, the grading filtration of $\PFlat$
  is a filtration with semisimple quotients and length equal to the
  Loewy length of~$\PFlat$. As $\PO$ is rigid by \autoref{P:RigidPIMs}
  and $\UnEFun(\PO)\cong \UnP$ by \autoref{T:UngradedCatOEquivalence},
  it follows that $\UnPFlat$ is rigid. Hence, the grading, radical and
  socle filtrations of~$\PFlat$ all coincide, so $\PFlat$ is rigidly
  graded as claimed.
\end{proof}

We now use the radical filtration of $\PO$ to construct a basis of
$\SO$. We could give a basis for $\RFlat$ in exactly the same way, and
we will do this implicitly in the proof of \autoref{T:MatchingGradings}
below. The key point is that if $\bmu\in\Klesh[\beta]$ then $\PO$ is
rigidly graded by \autoref{P:RigidPIMs}. In particular, the radical
filtrating of~$\PO$ is equal to its grading filtration. Therefore,
writing $c_{\blam\bmu}(q)=\sum_{d\ge0}c^{(d)}_{\blam\bmu}q^d$, for
$c^{(d)}_{\blam\bmu}\in\N$,
$$\rad^d\PO/\rad^{d+1}\PO\cong
    \bigoplus_{\blam\in\Parts[\beta]}(\LO[\blam]\<d\>)^{\oplus c^{(d)}_{\blam\bmu}},$$
for $0\le d\le z$. Fix $\blam\in\Parts[\beta]$ and $d\ge0$ with
$c^{(d)}_{\blam\bmu}\ne0$. Since $\PO[\blam]\<d\>$ is the projective cover of
$\LO[\blam]\<d\>$ there exist homogeneous maps
$\thetads\in\ZHom_{\SO}(\PO[\blam],\rad^d\PO)$ such that the diagrams commute
\begin{equation}\label{E:ThetaDefinition}
  \begin{tikzpicture}[baseline=2mm]
  \matrix[matrix of math nodes,row sep=1cm,column sep={8mm}]{
             &                          &                     &[8mm]|(Ymu)|\PO[\blam]\\
             &                          &                     &|(Ymud)|\PO[\blam]\<d\>\\
       |(0)|0&|(radd)|\rad^{d+1}(\PO)&|(rad)|\rad^d(\PO)&|(Dmu)|\LO[\blam]\<d\> &|(00)|0\\
     };
     \draw[->](Ymu)--node[right]{$\iota_{\PO[\blam]}\<d\>$}(Ymud);
     \draw[->>](Ymud) --(Dmu);
     \draw[->](Ymu)--node[above]{$\thetads$\,\,} (rad);
     \draw[->,dashed](Ymud)--(rad);
     \draw[->](0)--(radd);
     \draw[->](radd)--(rad);
     \draw[->](rad)--(Dmu);
     \draw[->](Dmu)--(00);
    \end{tikzpicture}
\end{equation}
for $1\le s\le c^{(d)}_{\blam\bmu}$ and where $\iota_{\PO[\blam]}\<d\>$
is the identity map followed by shifting all degrees by~$d$. By
embedding $\rad^d\PO$ into $\PO$ we consider $\thetads$ as a homogeneous
element of degree~$d$ in $\SO$.

\begin{Lemma}\label{L:BasicBasis}
  Suppose that $\beta\in Q^+_n$. Then
  $$\RBasis=\set{\thetads| 1\le s\le c^{(d)}_{\blam\bmu},
                                     0\le d\le2\defect\beta \text{ and }
                                     \blam,\bmu\in\Parts[\beta]}$$
  is a homogeneous basis of $\RO$.  Moreover, $\deg\thetads=d$ for all
  $\thetads\in\RBasis$.
\end{Lemma}

\begin{proof}
  Suppose that $\blam,\bmu\in\Parts[\beta]$. By construction,
  $\thetads$ is a homogeneous map of degree~$d$ and the set
  $\set{\thetads|1\le s\le c^{(d)}_{\blam\bmu}\text{ and }d\ge0}$ is linearly
  independent by construction. Counting graded dimensions,
  this set is a homogeneous basis of
  $\ZHom_{\SO}(\PO[\blam],\PO)$. If $\blam,\bmu\in\Klesh[\beta]$ then
  $c^{(d)}_{\blam\bmu}=0$ if $d>2\defect\beta$ by
  \autoref{L:DefectBound}. The lemma follows.
\end{proof}

We are ready to prove the main result of this section.

\begin{Theorem}\label{T:MatchingGradings}
  Suppose that $e=0$, $\mz=\C$ and $\beta\in Q^+_n$. Then there is a homogeneous
  algebra isomorphism $\Xi:\RO\bijection\RFlat$ of degree zero.
\end{Theorem}

\begin{proof} By \autoref{T:UngradedCatOEquivalence} there is an equivalence
  of categories $\UnEFun\map{\SO\Mod}\UnS[\beta]\Mod$ such that
  $\UnEFun(\UnPO[\blam])\cong\UnP[\blam]$, for all
  $\blam\in\Parts[\beta]$.  By \autoref{L:BasicBasis}, $\RBasis=\{\thetads\}$ is a
  basis of $\ZHom_{\SO}(\PO[\blam],\PO)$. For each $\thetads\in\RBasis$
  set $\varthetads=\UnEFun(\thetads)$. Forgetting the gradings for the
  moment, $\varthetads\in\Hom_{\UnSO}(\UnPO[\blam],\UnPO)$ and
  $\{\varthetads\}$ is a basis of
  $\Hom_{\UnS[\beta]}(\UnP[\blam],\UnP)$. Moreover, because $\UnEFun$
  preserves radical filtrations, if $\blam,\bmu\in\Klesh[\beta]$ then
  \autoref{E:ThetaDefinition} shows that there is a commutative diagram
  \begin{center}
  \begin{tikzpicture}[baseline=2mm]
  \matrix[matrix of math nodes,row sep=1cm,column sep={8mm}]{
             &                      &                 &[8mm]|(Ymu)|\UnP[\blam]\\
       |(0)|0&|(radd)|\rad^{d+1}(\UnP)&|(rad)|\rad^d(\UnP)&|(Dmu)|\UnL &|(00)|0\\
     };
     \draw[->>](Ymu) --(Dmu);
     \draw[->](Ymu)--node[above]{$\varthetads$\,\,} (rad);
     \draw[->](0)--(radd);
     \draw[->](radd)--(rad);
     \draw[->](rad)--(Dmu);
     \draw[->](Dmu)--(00);
    \end{tikzpicture}
  \end{center}
  By \autoref{P:PFlatRigid}, we can write
  $\varthetads=\sum_{k\ge d}\vartheta^{(d,s)}_{\blam\bmu,k}$, where
  $\vartheta^{(d,s)}_{\blam\bmu,k}\in\SFlat$ is homogeneous of
  degree~$k$ and $\vartheta^{(d,s)}_{\blam\bmu,d}\ne0$. By replacing
  $\varthetads$ with $\vartheta^{(d,s)}_{\blam\bmu,d}\ne0$, if
  necessary, we may assume that $\varthetads$ is homogeneous of
  degree~$d$ whenever $\blam,\bmu\in\Klesh[\beta]$. The map
  $\Xi\map{\RO}\RFlat$ given by $\Xi(\thetads)=\varthetads$, for
  $\thetads\in\RBasis$ with $\blam,\bmu\in\Klesh[\beta]$, defines an
  isomorphism of graded vector spaces $\RO\bijection\RFlat$.  As
  $\UnEFun$ respects composition of maps, and as multiplication in~$\RO$
  and in~$\RFlat[\beta]$ is homogeneous, it follows $\Xi$ is an algebra
  homomorphism. Hence, $\RO\cong\RFlat$ as graded algebras as claimed.
\end{proof}

The proof of \autoref{T:MatchingGradings} is quite subtle because we
have to work with the indecomposable prinjective modules $P^\bmu$ for
the quiver Schur algebra $\Sch[\beta]$ and use the fact that these
modules are rigid. In many ways it would be more natural to prove this
result using the graded Young modules $Y^\bmu$ but as these modules are
not known to be rigid we cannot argue this way. It would be nice to know
when the Young modules are rigid. Examples show that the algebras
$\RO$ are not always quadratic, so the ideas underpinning
\autoref{P:Rigidity} are not sufficient to answer this question.

\subsection{Graded decomposition numbers when $e=0$}\label{S:GradedDec}
Beilinson, Ginzburg and Soergel showed that a Koszul algebra has a unique
positive grading \cite[Corollary~2.5.2]{BGS:Koszul}. We now show that the
gradings on the Schur algebras $\SO$ and $\Sch[\beta]$ coincide.

Recalling the functor~$\SOFun$ from \autoref{E:OSchurFunc}, define
analogues of the Young modules for~$\RO$ by setting $\YO=\SOFun(\PO)$,
for $\bmu\in\Parts[\beta]$.  Similarly, recalling that
$\RFlat=e^\Lambda_\beta\Sch[\beta]e^\Lambda_\beta$, let
$\YFlat[\bmu]=P^\bmu e^\Lambda_\beta$ be a Young module for $\RFlat$.
Using the isomorphism $\Xi:\RO\bijection\RFlat$ from
\autoref{T:MatchingGradings} we can consider $\YO$ as an
$\RFlat$-module.

\begin{Lemma}\label{L:YoungGradings}
  Suppose that $\bmu\in\Parts[\beta]$. Then $\YFlat[\bmu]\cong\YO\<a_\bmu\>$ as
  $\RFlat$-modules, for some $a_\bmu\in\Z$.
\end{Lemma}

\begin{proof} Recall the (ungraded) Young modules from
  \cite[(3.5)]{M:tilting} that were used in the proof of
  \autoref{L:UngradedYoundMods}. By \cite[Lemma~6.11]{BK:HigherSchurWeyl} and
  \autoref{L:UngradedYoundMods}), $\EFun(\YO)$ and $Y^\bmu$ are both
  graded lifts of $\underline{y}^\bmu$. As $\underline{y}^\bmu$ is
  indecomposable, and graded lifts of indecomposable modules are unique
  up to shift, it follows that $\YFlat[\bmu]\cong\YO\<a_\bmu\>$, for
  some $a_\bmu\in\Z$.
\end{proof}

\begin{Theorem}\label{T:Koszul}
  Suppose that $e=0$, $\mz=\C$ and $\beta\in Q^+_n$. Then $\SFlat\cong\SO$ as
  graded algebras. In particular, $\SFlat$ is Koszul.
\end{Theorem}

\begin{proof}
For any $\blam,\bmu\in\Parts[\beta]$, let $c_{\blam\bmu}(q)$ and
  $c^\Ocal_{\blam\bmu}(q)$ be the graded Cartan numbers of $\SFlat$ and $\SO$,
  respectively. By construction,
  $\Dim\SFlat=\sum_{\blam,\bmu\in\Parts[\beta]}c_{\blam\bmu}(q)$. By
  definition, for any $\blam,\bmu\in\Parts[\beta]$,

  \begin{xalignat*}{3}
    c_{\blam\bmu}(q)&=\Dim \ZHom_{\Sch[\beta]}(P^\bmu,P^\blam)\\
    &=\Dim \ZHom_{\R[\beta]}(Y^\bmu,Y^\blam),&&\text{by \autoref{C:BlockSchurFunctors}},\\
    &=\Dim \ZHom_{\RFlat[\beta]}(\YFlat[\bmu],\YFlat[\blam])\\
    &=\Dim \ZHom_{\RO[\beta]}(\YO[\bmu]\<a_\bmu\>,\YO[\blam]\<a_\blam\>),
       &&\text{by \autoref{L:YoungGradings}},\\
    &=q^{a_\blam-a_\bmu}\Dim \ZHom_{\RO[\beta]}(\YO[\bmu],\YO[\blam])\\
    &=q^{a_\blam-a_\bmu}c^\Ocal_{\blam\bmu}(q).
  \end{xalignat*}
  By graded BGG reciprocity (\autoref{C:SOBGG}),
  $c^\Ocal_{\blam\bmu}(q)=c^\Ocal_{\bmu\blam}(q)$, so that
  $c_{\blam\bmu}(q)=q^{2(a_\blam-a_\bmu)}c_{\bmu\blam}(q)$. However, the
  Cartan matrix of~$\Sch[\beta]$ is symmetric by
  \autoref{C:CartanSymmetric}. Arguing as in the second last paragraph in the proof of
  \autoref{T:SWEquivalence} it follows that $a_\blam=a_\bmu$ for all
  $\blam,\bmu\in\Parts[\beta]$ since $\Sch$ is indecomposable by
  \autoref{T:SBlocks}. On the other hand, $a_\bmu=0$ for all
  $\bmu\in\Klesh[\beta]$ because if $\bmu\in\Klesh[\beta]$ then
  $\Dim\YFlat=\Dim\YO$ by \autoref{T:MatchingGradings}.  Therefore,
  $c_{\blam\bmu}(q)=c^\Ocal_{\blam\bmu}(q)\in\N[q]$ and, consequently,
  $\Dim\SFlat=\Dim\SO$.

  Using \autoref{C:BlockSchurFunctors} twice, there are homogeneous isomorphisms
  \begin{xalignat*}{3}
    \SO&\cong\ZEnd_{\SO}\(\bigoplus_{\bmu\in\Parts[\beta]}\PO\)^\op
         \cong\ZEnd_{\RO[\beta]}\(\bigoplus_{\bmu\in\Parts[\beta]}\YO\)^\op\\
    &\cong\ZEnd_{\RFlat[\beta]}\(\bigoplus_{\bmu\in\Parts[\beta]}{\YFlat[\bmu]}\)^\op
     \cong\ZEnd_{\Sch}\(\bigoplus_{\bmu\in\Parts[\beta]}P^\bmu\)^\op
     \cong\SFlat,
  \end{xalignat*}
  where the third isomorphism follows from \autoref{L:YoungGradings} using the
  fact that $a_\bmu=0$ for all $\bmu\in\Parts[\beta]$. Hence, $\SFlat\cong\SO$ is Koszul by
  \autoref{T:CatOKoszul}.
\end{proof}

Define non-negative integers $d_{\blam\bmu}^{(s)}$ by
$d_{\blam\bmu}(q)=\sum_{s\ge0}d^{(s)}_{\blam\bmu}q^s$, for
$\blam,\bmu\in\Parts[\beta]$.

\begin{Corollary}\label{C:GradedDecompPositive}
  Suppose that $e=0$, $\mz=\C$ and $\beta\in Q^+_n$ and let
  $\blam,\bmu\in \Parts[\beta]$. Then
  $d_{\blam\bmu}(q)\in\delta_{\blam\bmu}+q\N[q]$ and
  $c_{\blam\bmu}(q)\in\delta_{\blam\bmu}+q\N[q]$ and if $s\ge0$ then
   $$[\rad^s\Delta^\blam/\rad^{s+1}\Delta^\blam:L^\bmu\<s\>]_q
            =d^{(s)}_{\blam\bmu}.$$
\end{Corollary}

\begin{proof} Since $\SFlat$ is Koszul,
  $d_{\blam\bmu}(q)\in\delta_{\blam\bmu}+q\N[q]$. Consequently,
  $c_{\blam\bmu}(q)\in\delta_{\blam\bmu}+q\N[q]$ by
  \autoref{C:CartanSymmetric}. Finally, since
  $\Delta^\bmu/\rad\Delta^\bmu\cong L^\bmu$ is irreducible the last statement
  follows from \autoref{P:RigidPIMs}(a).
\end{proof}

Combining the results in this section we obtain a more precise version of
\autoref{THEOREMC} from the introduction.

\begin{Theorem}\label{T:GradedCatOEquivalence}
    Suppose that $\beta\in Q^+_n$, $e=0$ and $\mz=\C$. Then there is a
    graded equivalence of categories
    $\EFun\map{\O[\beta]}\Sch[\beta]\Mod$  such that the following
    diagram commutes:
    $$\begin{tikzpicture}
     \matrix[matrix of math nodes,row sep=1cm,column sep=16mm]{
       |(O)| \O[\beta]&|(S)|\Sch[\beta]\Mod\\
                    &|(H)|\R[\beta]\Mod\\
     };
     \draw[->](O) -- node[above]{$\EFun$} (S);
     \draw[->](S) -- node[right]{$\SFun[\beta]$} (H);
     \draw[->](O) -- node[below]{$\OFun$} (H);
    \end{tikzpicture}$$
    Moreover, $\EFun(\DelO)\cong\Delta^\blam$ and
    $\EFun(\LO)\cong L^\bmu$, for all $\blam,\bmu\in\Parts$.
\end{Theorem}

\subsection{The Fock Space}\label{S:Fock}
The aim of this subsection is to realize the projective indecomposable and
irreducible modules for~$\Sch[n]$ as the canonical and dual canonical bases of the
higher level Fock space. Throughout this subsection, we work over~$\C$
and assume that either $e=0$ or $e>n$.

Let $\Rep(\Sch[n])$ be the Grothendieck group of finitely generated
$\Sch[n]$-modules. If~$M$ is an $\Sch[n]$-module let $[M]$ be its image in
$\Rep(\Sch[n])$. Observe that $\Rep(\Sch[n])$ is naturally a $\Z[q,q^{-1}]$-module
where~$q$ acts by grading shift: $q[M]=[M\<1\>]$, for $M\in\Sch[n]\Mod$.
Similarly, let $\Proj(\Sch[n])$ be the Grothendieck group of the category of
finitely generated projective $\Sch[n]$-modules. The Cartan pairing is the
sesquilinear map (anti-linear in the first argument, linear in the second)
$$(\ ,\ )\map{\Proj(\Sch[n])\times\Rep(\Sch[n])}\Z[q,q^{-1}],\quad
         ([P],[M])=\Dim\ZHom_{\Sch[n]}(P,M),$$
for $[P]\in\Proj(\Sch[n])$ and $[M]\in\Rep(\Sch[n])$. There is a natural
embedding $\Proj(\Sch[n])\hookrightarrow\Rep(\Sch[n])$.

Define the \textbf{combinatorial Fock space} of weight~$\Lambda$ to be
$$\Fock=\bigoplus_{n\ge0}\Rep(\Sch[n]).$$
Thus,~$\Fock$ is a free $\Z[q,q^{-1}]$-module of infinite rank. Let
$\Parts[]=\bigcup_{n\ge0}\Parts$. The Fock space~$\Fock$ is equipped
with the following distinguished bases:
\begin{itemize}
  \item The irreducible modules: $\set{[L^\bmu]|\bmu\in\Parts[]}$.
  \item The standard modules $\set{[\Delta^\bmu]|\bmu\in\Parts[]}$.
  %\item The costandard modules $\set{[\nabla^\bmu]|\bmu\in\Parts[]}$.
  \item The projective indecomposable modules $\set{[P^\bmu]|\bmu\in\Parts[]}$.
  \item The tilting modules $\set{[T^\bmu]|\bmu\in\Parts[]}$.
\end{itemize}
These four sets are all bases for~$\Fock$ as a $\Z[q,q^{-1}]$-module because the
graded decomposition matrix of~$\Sch[n]$ is invertible over~$\Z[q,q^{-1}]$ by
\autoref{C:Positivity}.

The aim of this section is to clarify the relationships between these bases and
to give an algorithm for computing the graded decomposition numbers of~$\Sch[n]$.

There is a natural duality on $\Rep(\Sch[n])$ that induces an involution
on~$\Fock$. Let~$M$ be an $\Sch[n]$-module. Recall that
$M^\circledast=\ZHom_\C(M,\C)$ is the graded dual of~$M$.
Similarly, define $M^\#=\ZHom_{\Sch[n]}(M,\Sch[n])$, where $\Sch[n]$ acts
on $M^\#$ by $(f\cdot s)(x)=f(xs^{\star})$, for $f\in M^\#$ and $x\in M,
s\in\Sch[n]$. Then~$\#$ restricts to a duality on~$\Proj(\Sch[n])$.

\begin{Lemma}\label{L:involutions}
  Suppose that $M$ is an $\Sch[n]$-module. Then
         $$([P^\#],[M]) = \overline{([P],[M^\circledast])},$$
for all $[P]\in\Proj(\Sch[n])$ and $[M]\in\Rep(\Sch[n])$.
    Moreover, if $\bmu\in\Parts$ then
    $(L^\bmu)^\circledast\cong L^\bmu$,
    $(T^\bmu)^\circledast\cong T^\bmu$ and $(P^\bmu)^\#\cong P^\bmu$.
\end{Lemma}

\begin{proof}
  The first statement is well-known; see, for example,
  \cite[Lemma~2.5]{BK:GradedDecomp}. This implies that $(P^\bmu)^\#\cong P^\bmu$
  since $(L^\bmu)^\circledast\cong L^\bmu$ by \autoref{T:CellularSimples}.
  Finally, $(T^\bmu)^\circledast\cong T^\bmu$ by
  \autoref{C:TiltingSelfDual}.
\end{proof}

A map $f\map MN$ of $\Z[q,q^{-1}]$-modules is \textbf{semilinear} if it is
$\Z$-linear and $f(q^km)=q^{-k}f(m)$, for all $m\in M$ and $k\in\Z$.

\begin{Lemma}\label{L:Semilinear}
  The maps $\circledast$ and $\#$ induce semilinear involutions on~$\Fock$ such
  that
  $$(M\<d\>)^\circledast\cong M^\circledast\<-d\>\qquad\text{and}\qquad
    (N\<d\>)^\#\cong N^\#\<-d\>,$$
  for all $[M]\in\Rep(\Sch[n])$,  $[N]\in\Proj(\Sch[n])$ and $d\in\Z$.
\end{Lemma}

\begin{proof}It follows easily from the definitions that $\circledast$ is
  a duality on $\Rep(\Sch[n])$ and that $\#$ is a duality on~$\Proj(\Sch[n])$. This
  immediately implies that $\circledast$ induces an involution on~$\Fock$ with the
  required properties. Moreover, $\#$ extends to an automorphism of~$\Fock$
  because $\set{[P^\bmu]|\bmu\in\Parts[]}$ is a $\Z[q,q^{-1}]$-basis of~$\Fock$.
  The map induced by~$\#$ is an involution because
  $(P^\bmu)^\#\cong P^\bmu$ by \autoref{L:involutions}, for $\bmu\in\Parts$.
\end{proof}

We emphasize that both of these maps are semilinear -- that is, $\Z$-linear but
not $\Z[q,q^{-1}]$-linear. This is implicit in the displayed equation of
\autoref{L:Semilinear} because, for example,
$(q[M])^\circledast=[M\<1\>]^\circledast=[M^\circledast\<-1\>]=q^{-1}[M^\circledast]$.

Recall from \autoref{S:graded} that the \textbf{bar involution} on
$\Z[q,q^{-1}]$ is the $\Z$-linear automorphism of $\Z[q,q^{-1}]$ determined by
$\overline{q}=q^{-1}$. A Laurent polynomial~$f(q)$  in~$\Z[q,q^{-1}]$ is \textbf{bar
invariant} if $f(q)=\overline{f(q)}$.

\begin{Lemma}\label{L:BarInvolution}
  Suppose that $\blam\in\Parts$. Then
  $$[\Delta^\blam]^\circledast=[\Delta^\blam]
       +\sum_{\substack{\bmu\in\Parts\\\blam\gdom\bmu}}f_{\blam\bmu}(q)[\Delta^\bmu]
       \quad\text{and}\quad
  [\Delta^\blam]^\#=[\Delta^\blam]+\sum_{\substack{\bmu\in\Parts\\\bmu\gdom\blam}}
            g_{\blam\bmu}(q)[\Delta^\bmu],$$
  for some Laurent polynomials $f_{\blam\bmu}(q),g_{\blam\bmu}(q)\in\Z[q,q^{-1}]$.
\end{Lemma}

\begin{proof}
  Recall that $\(d_{\blam\bmu}(q)\)_{\blam,\bmu\in\Parts}$ is the graded
  decomposition matrix of~$\Sch[n]$. Let
  $\(e_{\blam\bmu}(q)\)_{\blam,\bmu\in\Parts}$ be the inverse graded
  decomposition matrix. Using \autoref{L:BarInvolution} we compute:
  \begin{align*}
    [\Delta^\blam]^\circledast
    &=\Big(\sum_{\substack{\bmu\in\Parts\\\blam\gedom\bmu}}
                  d_{\blam\bmu}(q)[L^\bmu]\Big)^\circledast
       =\sum_{\substack{\bmu\in\Parts\\\blam\gedom\bmu}}d_{\blam\bmu}(q^{-1})[L^\bmu]\\
  \intertext{since $(L^\bmu)^\circledast\cong L^\bmu$ by \autoref{T:quasi}. Therefore,}
    [\Delta^\blam]^\circledast
      &=\sum_{\substack{\bmu\in\Parts\\\blam\gedom\bmu}}d_{\blam\bmu}(q^{-1})
      \sum_{\substack{\bnu\in\Parts\\\bmu\gedom\bnu}}e_{\bnu\bmu}(q)[\Delta^\bnu]\\
      &=[\Delta^\blam]+\sum_{\substack{\bnu\in\Parts\\\blam\gdom\bnu}}
      \Big(\sum_{\substack{\bmu\in\Parts\\\blam\gedom\bmu\gedom\bnu}}
      e_{\bnu\bmu}(q)d_{\blam\bmu}(q^{-1})\Big)\,[\Delta^\bnu],
  \end{align*}
  where the last line follows because both the graded decomposition matrix and
  its inverse are triangular with respect to dominance by
  \autoref{C:TriangularDecomp}. The formula for~$[\Delta^\blam]^\#$ is
  proved in exactly the same way by first writing
  $[\Delta^\blam] =\sum_{\bmu\gedom\blam} e_{\bmu\blam}(q)[P^\bmu]$.
\end{proof}

By a well-known result of Lusztig~\cite[Lemma~24.2.1]{Lusztig:QuantBook},
\autoref{L:BarInvolution} implies that~$\Fock$ has several uniquely determined
`canonical bases' that are invariant under $\circledast$ and $\#$. Using
\autoref{C:Positivity} we can describe these bases explicitly.
Let~$\Fock_q(\rhd\bmu)$ (resp., ~$\Fock_q(\lhd\bmu)$) be the $\Z[q]$-sublattice of $\Fock$ with basis the images of the
standard modules $\set{[\Delta^\blam]|\bmu\lhd\blam\in\Parts[]}$ (resp., $\set{[\Delta^\blam]|\bmu\rhd\blam\in\Parts[]}$)
in~$\Fock$. Let $\Fock_{q^{-1}}(\lhd\bmu)$ be the $\Z[q^{-1}]$-sublattice of $\Fock$ spanned by the images of the
standard modules $\set{[\Delta^\blam]|\bmu\rhd\blam\in\Parts[]}$.
% Extending this notation, if $\bmu\in\Parts[]$ let
% $(\Fock_0)^{\gdom\bmu}$ be the $\Z[q]$-sublattice of $\Fock$ spanned by the
% images of the standard modules $[\Delta^{\blam}]$ with $\blam\rhd\bmu$. We
% define $(\Fock_0)^{\lhd\bmu}$ and $(\Fock_{q^{-1}})^{\rhd\bmu}$ similarly.

\begin{Theorem}\label{T:CanonicalBases}
  Suppose that $e=0$ and $\mz=\C$. Then the three bases
  $$\set{[P^\bmu]|\bmu\in\Parts[]},\quad
    \set{[L^\bmu]|\bmu\in\Parts[]}\quad\text{and}\quad
    \set{[T^\bmu]|\bmu\in\Parts[]}$$
    are ``canonical bases'' of $\Fock$ that, for $\bmu\in\Parts[]$, are
    uniquely determined by:
  \begin{enumerate}
    \item $[P^\bmu]^\#=[P^\bmu]$ and
    $[P^\bmu]\equiv[\Delta^\bmu]\pmod{q\Fock_q(\rhd\bmu)}$.%\mod q(\Fock_0)^{\gdom\bmu}$.
    \item $[L^\bmu]^\circledast=[L^\bmu]$ and
    $[L^\bmu]\equiv [\Delta^\bmu]\pmod{q\Fock_q(\lhd\bmu)}$.%\mod q(\Fock_0)^{\ldom\bmu}$.
    \item $[T^\bmu]^\circledast=[T^\bmu]$ and
    $[T^\bmu]\equiv[\Delta^\bmu]\pmod{q^{-1}\Fock_{q^{-1}}(\lhd\bmu)}$.
    \end{enumerate}
    %Moreover, $T^{\bmu}\cong T^{\bmu}$.
   \end{Theorem}

\begin{proof}
  The existence and uniqueness of bases of $\Fock$ with these properties follows
  from what is by now a standard argument
  (see~\cite[Lemma~24.2.1]{Lusztig:QuantBook}), using the triangularity of the
  involutions~$\circledast$ and~$\#$ from \autoref{L:BarInvolution}. If
  $\bmu\in\Parts$ then
  $[P^\bmu]^\#=[P^\bmu]$, $(L^\bmu)^\circledast\cong L^\bmu$
  and $(T^\bmu)^\circledast\cong T^\bmu$ by \autoref{L:involutions}.
  Furthermore, $[P^\bmu]=\sum_\blam d_{\blam\bmu}(q)[\Delta^\blam]$ and
  $[L^\bmu]=\sum_{\blam}e_{\blam\bmu}(q)[\Delta^\blam]$, where
  $d_{\blam\bmu}(q)$ and $e_{\blam\bmu}(q)$ are polynomials in~$\Z[q]$
  with constant term $d_{\blam\bmu}(0)=\delta_{\blam\bmu}=e_{\blam\bmu}(0)$ by
  \autoref{C:Positivity}. Therefore, if~$\bmu\in\Parts[]$ then
  $[P^\bmu]$ and $[L^\bmu]$ belong to~$\Fock_q(\gedom\bmu)$ and, moreover,
  $$ [P^\bmu]\equiv[\Delta^\bmu]\pmod{q\Fock_q(\gdom\bmu)}%\mod{q(\Fock_0)^{\rhd\bmu}}
  \quad\text{and}\quad
  [\Delta^\bmu]\equiv [L^{\bmu}]\pmod{q\Fock_q(\ldom\bmu)}%\mod{q(\Fock_0)^{\lhd\bmu}}.
  $$
  Hence, parts~(a) and~(b) follow. Finally,
  by \autoref{C:TiltingSelfDual} and \autoref{L:TiltingMultiplicities}, $(T^{\bmu})^{\circledast}\cong T^{\bmu}$ and $[T^\bmu]\equiv[\Delta^\bmu]\pmod{q^{-1}\Fock_{q^{-1}}(\lhd\bmu)}$. This completes the proof.
\end{proof}

We call $\set{[T^\bmu]|\bmu\in\Parts[]}$ the \textbf{canonical basis} of~$\Fock$
and $\set{[L^\blam]|\bmu\in\Parts[]}$ the \textbf{dual canonical basis}. By
\autoref{T:Ringeldual}, Ringel duality induces a automorphism of~$\Fock$
that interchanges, setwise, the canonical basis $\{[P^\bmu]\}$ and the basis
$\{[T^\blam]\}$ of tilting modules. We remark that
\autoref{T:CanonicalBases} should lift to a categorification of the
canonical bases of~$\Fock$ as a~$U_q(\mathfrak{gl}_\infty)$-module.

\begin{Remark}
  Abusing notation slightly, let $\#$ be the automorphism of $\Rep(\R)$ defined by
  $M^\#=\ZHom_{\R}(M,\R)$. Then, as noted in \cite[Remark~4.7]{BK:GradedDecomp},
  it follows from \autoref{T:trace} and \cite[Theorem~3.1]{Rouq:2KM} that
  there is an isomorphism of functors
  $\#\cong\<2\defect\beta\>\circ\circledast$. Therefore,
  $$\set{q^{\defect\beta}[D^\bmu]|\bmu\in\Klesh[\beta]\text{ for }\beta\in Q^+}$$
  is a $\#$-invariant basis of $\bigoplus_{n}\Rep(\R)$ that has similar
  uniqueness properties to the tilting module basis of~$\Fock$.
  Similarly, $\set{q^{-\defect\beta}[Y^\bmu]|\bmu\in\Klesh[\beta]}$ is a
  `canonical' $\circledast$-invariant basis of $\bigoplus_{n\ge0}\Rep(\R)$.
 \end{Remark}

% Analogously, by Lusztig's Lemma~\cite[Lemma~24.2.1]{Lusztig:QuantBook}, there
%  is a unique $\#$-invariant basis $\set{B_\bmu|\bmu\in\Parts[]}$ of~$\Fock$
%  such that $B^\bmu\equiv[\nabla^\bmu]\pmod{q\Fock_{q^{-1}}}$, however, we do not
%  know how to describe this basis.

% \begin{Corollary}
%   Yvonne's conjecture~\cite{Yvonne:Conjecture} is true when $e=0$.
% \end{Corollary}
%
% \begin{proof}
%   The parabolic KL polynomials arising in Theorem~3.26 of Uglov's paper coincide with
%   those in \autoref{C:ParabolicKLPolys}.
% \end{proof}

\subsection{An LLT algorithm for $\Sch[n]$}\label{S:LLT}
\autoref{T:CanonicalBases} is not surprising because it is a well
established mantra that the indecomposable tilting modules should correspond
to Lusztig's canonical basis, the simple modules to the dual canonical basis and
so on. The main reason for our introducing the Fock space is that we can now use
the tableau combinatorics for $\Sch$ to explicitly compute the graded decomposition
numbers for $\Sch$ and hence for parabolic category~$\O$. Thus, not only does this
combinatorics explicitly describe the grading on parabolic category~$\O$ but it
also gives an effective way of computing graded decomposition multiplicities of~$\Sch[n]$, which are certain parabolic Kazhdan-Lusztig polynomials.

If $\Lambda=\Lambda_0$ then $\H\cong\R$ is isomorphic to the
Iwahori-Hecke algebra of the symmetric group.  In this case, Lascoux,
Leclerc and Thibon~\cite{LLT} have given an efficient algorithm for
computing the canonical basis of the irreducible
$U_q(\widehat{\mathfrak{sl}})$-module $L(\Lambda_0)$. By Ariki's
Theorem~\cite{Ariki:ss,BK:GradedDecomp}, the LLT algorithm computes the
(graded) decomposition matrices of the Iwahori-Hecke algebra of the
symmetric group.

In this section we give an LLT-like algorithm for computing the canonical basis
of~$\Fock$. By \autoref{T:CanonicalBases} this gives an algorithm for computing
the graded decomposition numbers of~$\R$ and~$\Sch[n]$. To this end, if
$f(q)=\sum_{d\in\Z}f_dq^d$ is a non-zero Laurent polynomial in~$\Z[q,q^{-1}]$ let
$\mindeg f(q)=\min\set{d\in\Z|f_d\ne0}$.

Suppose that $\bmu\in\Parts[\beta]$, where $\beta\in Q^+$. Recall
from~\autoref{E:Zmu} that $Z^\bmu=\Psi^\bmu\Sch$. Similarly, recall from
\autoref{S:graded} that if~$M$ is an $\Sch$-module and
$f(q)=\sum_kf_kq^k\in\N[q,q^{-1}]$ then $f(q)M=\bigoplus_k M\<k\>^{\oplus f_k}$.

\begin{Lemma}\label{L:ZmuExpansion}
  Suppose that $\bmu\in\Parts[\beta]$. Then $(Z^\bmu)^\#\cong Z^\bmu$ and
  $$Z^\bmu=P^\bmu\oplus\bigoplus_{\blam\gdom\bmu}z_{\blam\bmu}(q)P^\blam,$$
  for some bar invariant polynomials $z_{\blam\bmu}(q)\in\N[q,q^{-1}]$.
\end{Lemma}

\begin{proof} By definition, $Z^\bmu$ is a direct summand of $\Sch$ and $(\Psi^{\bmu})^{\star}=\Psi^{\bmu}$, so
  $(Z^\bmu)^\#\cong Z^\bmu$. We already noted in~\autoref{E:proj} that
  $Z^\bmu=P^\bmu\oplus\bigoplus_\blam z_{\blam\bmu}(q)P^\blam$, for some Laurent
  polynomials $z_{\blam\bmu}(q)\in\N[q,q^{-1}]$, because $Z^\bmu$ is projective. In
  view of \autoref{L:involutions} these polynomials are bar invariant.
\end{proof}

Next observe that~\autoref{E:ZmuDelta} implies that in~$\Fock$
\begin{equation}\label{E:ZmuExpansion}
       [Z^\bmu] = [\Delta^\bmu]+\sum_{\substack{\bnu\gdom\bmu\\\s\in\Std^\bmu(\bnu)}}
                   q^{\deg\s-\deg\tmu}[\Delta^\bnu].
\end{equation}

We now show how to use \autoref{L:ZmuExpansion}  and~\autoref{E:ZmuExpansion} to
inductively compute~$[P^\bmu]$, for $\bmu\in\Parts[\beta]$, as a linear combination of
standard modules in~$\Fock$. Since $[P^\bmu]=\sum_\blam
d_{\blam\bmu}(q)[\Delta^\blam]$ this will give an algorithm for computing the
graded decomposition numbers of~$\Sch$.

If $\bmu$ is maximal in $\Parts[\beta]$, with respect to dominance, then
$Z^\bmu=P^\bmu=\Delta^\bmu$ by \autoref{L:ZmuExpansion}.
So~$[P^\bmu]=[\Delta^\bmu]$ in this case and there is nothing to do.

Now suppose that~$\bmu$ is not maximal in~$\Parts[\beta]$ and that~$[P^\blam]$ is known
whenever $\blam\in\Parts[\beta]$ and $\blam\gdom\bmu$. By~\autoref{E:ZmuExpansion} we can
write
$$[Z^\bmu]=[\Delta^\bmu]+\sum_{\substack{\bnu\in\Parts[\beta]\\\bnu\gdom\bmu}}
                                 y_{\bnu\bmu}(q)[\Delta^\bnu]$$
for some Laurent polynomials $y_{\bnu\bmu}(q)\in\N[q,q^{-1}]$ that are not all zero
since $\bmu$ is not maximal in~$\Parts[\beta]$. Let $\blam\gdom\bmu$ be any
multipartition such that $y_{\blam\bmu}(q)\ne0$ and
$$\mindeg y_{\blam\bmu}(q)\le\mindeg y_{\bnu\bmu}(q),$$
for all $\bnu\in\Parts[\beta]$. Let $d=\mindeg y_{\blam\bmu}(q)$.

If $d>0$ then $[Z^\bmu]\equiv[\Delta^\bmu]\pmod{q\Fock_q(\rhd\bmu)}$
by~\autoref{E:ZmuExpansion}. Now $[Z^\bmu]^\#=[Z^\bmu]$, by
\autoref{L:ZmuExpansion}. This forces $[Z^\bmu]=[P^\bmu]$ because $[Z^\bmu]$
satisfies the two properties that uniquely determine~$[P^\bmu]$ in
\autoref{T:CanonicalBases}(a).

Now suppose that $d\le 0$. Let
$y_{\blam\bmu}^{(d)}$ be the coefficient of~$q^d$ in~$y_{\blam\bmu}(q)$ and set
$$z_{\blam\bmu}^{(d)}=\begin{cases}
       y_{\blam\bmu}^{(d)}(q^d+q^{-d}),&\text{if $d<0$},\\
       y_{\blam\bmu}^{(d)},&\text{if }d=0.
     \end{cases}
$$
Since $[P^\bnu]\equiv[\Delta^\bnu]\pmod{q\Fock_q(\rhd\bnu)}$,
for all $\bnu\in\Parts[\beta]$, the minimally of~$d$ together with
\autoref{L:ZmuExpansion} implies that $z_{\blam\bmu}^{(d)}P^\blam$ is a direct
summand of~$Z^\bmu$. Since $[P^\blam]$ is known by induction we can now replace
$[Z^\bmu]$ with $[Z^\bmu]-p^{(d)}_{\blam\bmu}[P^\blam]$, which is still
$\#$-invariant. By repeating this process of stripping off the bar invariant
minimal degree terms we can rewrite $[Z^\bmu]$ as a linear combination of
canonical bases elements as in \autoref{L:ZmuExpansion}. This recursively
computes~$[P^\bmu]$ and so determines the graded decomposition numbers
$d_{\blam\bmu}(q)$.

Note that the Laurent polynomials $z_{\blam\bmu}(q)$ in
\autoref{L:ZmuExpansion} are given by
$z_{\blam\bmu}(q)=\sum_{d\le 0}z_{\blam\bmu}^{(d)}$. Hence, this algorithm also
decomposes~$Z^\bmu$ into a direct sum of projective modules.

\begin{Remark}
  Note that $(E^\bmu)^\circledast\cong E^\bmu$ by \autoref{T:ESelfDual}. An
  equivalent version of this algorithm computes $[T^\bmu]$ by applying the same
  ``straightening algorithm" to the element $[E^\bmu]=[E^\bmu]^\circledast$,
  where we use \autoref{C:EWeylFilt} in place of~\autoref{E:ZmuExpansion} and
  \autoref{C:EWeylFilt} in place of \autoref{L:ZmuExpansion}.
\end{Remark}

\begin{Example}
  Suppose that $e=0$, $\Lambda=3\Lambda_0$ and that
  $\beta=\alpha_{-1}+3\alpha_0+\alpha_1+\alpha_2+\alpha_3$. Then $\Sch[\beta]$ is a block
  of defect~$4$. The maximal multipartition in~$\Parts[\beta]$ is $(4,2|1|0)$ so
  $P^{(4,2|1|0)}=\Delta^{(4,2|1|0)}$. Taking $\bmu=(4,1|1|1)$ the tableaux in
  $\Std^\bmu(\Parts[\beta])$ are
  $$\begin{array}{ll}
    \Tritab({1,2,3,4},{5}|{6}|{7}) &
    \Tritab({1,2,3,4},{5,6}|{-}|{7})\\[5mm]
    \Tritab({1,2,3,4},{5,7}|{6}|{-})&
    \Tritab({1,2,3,4},{5,6}|{7}|{-})
  \end{array}$$
  Therefore,
  $[Z^\bmu]=[\Delta^{(4,1|1|1)}]+q[\Delta^{(4,2|0|1)}]+(q^2+1)[\Delta^{(4,2|1|0)}]$.
  Applying our algorithm, $[Z^\bmu]=[P^\bmu]+[P^{(4,2|1|0)}]$. Using our LLT algorithm,
  the full graded decomposition matrix of~$\Sch[\beta]$ in characteristic zero is:
  $$\begin{array}{r@{|}c@{|}l|*{16}{c}}
(0&1&4,2)&1&&&&&&&&&&&&&&&\\
(0&4,2&1)&q&1&&&&&&&&&&&&&&\\
(1&0&4,2)&q&.&1&&&&&&&&&&&&&\\
(1&1&4,1)&q^{2}&.&q&1&&&&&&&&&&&&\\
(1&1^2&4)&.&.&.&q&1&&&&&&&&&&&\\
(1&4&1^2)&.&.&.&q&.&1&&&&&&&&&&\\
(1&4,1&1)&q^{2}&q&q&q^{2}&q&q&1&&&&&&&&&\\
(1&4,2&0)&q^{3}&q^{2}&q^{2}&.&.&.&q&1&&&&&&&&\\
(1^2&1&4)&.&.&q&q^{2}&q&.&.&.&1&&&&&&&\\
(1^2&4&1)&.&.&q^{2}&q^{3}&q^{2}&q^{2}&q&.&q&1&&&&&&\\
(4&1&1^2)&.&.&q&q^{2}&.&q&.&.&.&.&1&&&&&\\
(4&1^2&1)&.&.&q^{2}&q^{3}&q^{2}&q^{2}&q&.&.&.&q&1&&&&\\
(4,1&1&1)&q^{2}&q&q^{3}+q&q^{4}&q^{3}&q^{3}&q^{2}&.&q^{2}&q&q^{2}&q&1&&&\\
(4,2&0&1)&q^{3}&q^{2}&q^{2}&.&.&.&.&.&.&.&.&.&q&1&&\\
(4,2&1&0)&q^{4}&q^{3}&q^{3}&.&.&.&q^{2}&q&.&q&.&q&q^{2}&q&1&\\
\end{array}$$
The Kleshchev multipartitions in this block are $(0|1|4,2)$ and
$(1|1|4,1)$.  If $\Lambda$ is any weight of level~$\ell=2$ then the
graded decomposition numbers of~$\Sch[\beta]$ are monomials in~$q$ by
\autoref{T:L2Positivity}. The algebra $\Sch[\beta]$ is one of the smallest
examples of a block that has a graded decomposition number that is not a
monomial.
\end{Example}

\autoref{T:L2Positivity} in the appendix shows that all of the results in this
section are valid over an arbitrary field when~$\Lambda$ is a dominant weight of
level~2 and $e=0$ or $e>n$.

%%%%%%%%%%%%%%%%%%%%%%%%%%%%%%%%%%%%%%%%%%%%%%%%%%%%%%%%%%%%%
\appendix
\def\theequation{\Alph{section}\arabic{equation}}
\def\theProposition{\Alph{section}\arabic{equation}}
\def\theLemma{\Alph{section}\arabic{equation}}
\def\theTheorem{\Alph{section}\arabic{equation}}
\def\theCorollary{\Alph{section}\arabic{equation}}

\section{Graded decomposition numbers for $\O[\beta]$}
  \label{S:CatODecompNumbers}
  This appendix proves \autoref{T:CatODecompNumbers}, which explicitly
  describes the polynomials $p^\Ocal_{\blam\bmu}(q)$ as inverse
  parabolic Kazhdan-Lusztig polynomials. Throughout this section we fix
  $\Lambda\in P^+$ and $\beta\in Q^+_n$ as in
  \autoref{S:HigherSchurWeyl}.

  To prove \autoref{T:CatODecompNumbers}, we need to
  interpolate between Brundan and Kleshchev's work
  from~\cite{BK:HigherSchurWeyl,BK:DegenAriki}, as described in
  \autoref{S:CatO}, and the more standard Lie theoretic notation used in
  in~\cite{Backelin:Koszul,BGS:Koszul}.
  %In particular, if $\nu\in\h^*$ is a dominant integral weight let~$\UnL(\nu)$ be the
  %irreducible highest weight module of highest weight~$\nu$.
  Backelin considers blocks $\Ocal^\sigma_\tau$ of parabolic
  category~$\Ocal$ (of arbitrary type), indexed by dominant integral
  weights $\tau,\sigma\in\h^*$. He defines $R^\sigma_\tau$ to be the
  endomorphism algebra of a minimal projective generator
  of~$\Ocal^\sigma_\tau$. Set $\Ocal_\tau=\Ocal^0_\tau$, a block of
  category~$\Ocal=\Ocal^0$ for $\gl$. Recall that $\Sym_N$ acts
  on~$\h^*$ via the dot-action and let
  $\Sym(\tau)=\set{w\in\Sym_N|w\cdot\tau=\tau}$ be the stabiliser of
  $\tau$ in~$\Sym_N$. Recall that the simple modules in $\Ocal_0^\tau$ are in
  bijection with the set
  $X_{\tau}^{+}:=\{\sigma\in\h^*|\text{$\<\sigma,\varepsilon_i-\varepsilon_{i+1}\>\geq 0$, whenever $(i,i+1)\in\Sym(\tau)$}\}$.
  We set
  $$\D_\tau=\set{x\in\Sym_N|\ell(wx)\ge\ell(x)\text{ for all }w\in\Sym(\tau)},\,\,\D^+_\tau=\set{x\in\Sym_N|\ell(wx)\le\ell(x)\text{ for all }w\in\Sym(\tau)}.$$
  Then $\D_\tau$ is the set of minimal length right coset representatives of
  $\Sym(\tau)$ in~$\Sym_N$, while $\D^+_\tau$ is the set of maximal length right coset representatives of $\Sym(\tau)$
  in~$\Sym_N$. It is well-known that $w\cdot 0\in X_{\tau}^{+}$ if and only if $w\in\D_\tau$.
  The reader can check that the map
  $x\mapsto xw_0$ defines a bijection $\D_\tau\bijection\D^+_\tau$ using
  the fact that $\ell(ww_0)=\ell(w_0)-\ell(w)$, for all $w\in\Sym_N$. As a result,
  the simple modules in $\Ocal_\tau$ are the modules
  $\set{\unL(x^{-1}\cdot\tau)|x\in\D_\tau^{+}}$, and
  $\set{\unL(xw_0\cdot0)|x\in\D_\tau^+}$ is a complete set of
  simple modules in~$\Ocal_0^\tau$, where $w_0$ is the (unique) element of
  longest length in $\Sym_N$.
  %\footnote{Our notations and conventions are slightly different with that in \cite{BGS:Koszul} and \cite{Backelin:Koszul} in that we are using {\it right} coset %representatives instead of {\it left} coset representatives.}

  We need to describe the block $\O[\beta]$ using the more standard Lie
  theoretic notation introduced above
  following~\cite{Backelin:Koszul,BGS:Koszul}.  Fix multipartitions
  $\blam,\bmu\in\Parts[\beta]$ and recall from \autoref{L:BKBijections}
  that $\T^\blam$ is the $\Lambda$-tableau corresponding to~$\blam$ and
  that $\col(\T^\blam)$ is the column reading of~$\T^\blam$. In
  \autoref{S:HigherSchurWeyl}, following \cite{BK:DegenAriki}, we
  defined $w_\blam$ to be the unique minimal length permutation such
  that the sequence $\col(\T^\blam) w_\blam$ is weakly
  \textit{increasing} and that $Z_\blam$ is the stabilizer
  of~$\col(\T^\blam) w_\blam$. Let $D_\blam$ be the set of minimal
  length right coset representatives of~$Z_\blam$ in~$\Sym_N$. Similarly, let
  $v_\blam\in\Sym_N$ be the unique minimal length permutation such that
  $(s_1,\dots,s_N)=\col(\T^\blam)v_\blam$ is weakly
  \textit{decreasing}. Fix the dominant integral weight
  $\phi=s_1\eps_1+(s_2+1)\eps_2+\dots+(\s_N+N-1)\eps_N\in\h^*$ for the
  rest of this appendix.  By \autoref{L:BKBijections}, $\phi$ depends
  only on~$\Lambda$ and~$\beta$. As above, let $\Sym(\phi)$ be the
  stabiliser of~$\phi$ under the dot action. Then $v_\blam^{-1}\in\D_\phi$
  and $v_\blam^{-1}w_0\in\D^+_\phi$.

  By definition the subgroups $Z_\blam$ and $\Sym(\phi)$ are conjugate
  in~$\Sym_N$.  Recall that~$P_{x,y}$ is the Kazhdan-Lusztig polynomial
  indexed by $x,y\in\Sym_N$.

  \begin{Lemma}\label{L:Conjugation}
    Suppose that $\blam\in\Parts[\beta]$ and let $\phi\in\h^*$ be as
    above. Then
    $Z_\blam=w_0\Sym(\phi)w_0$. Moreover, if $x,y\in\D_\phi$ and $z\in\Sym(\phi)$ then
    $P_{z'x',y'}=P_{zx,y}$, where
    $x'=w_0xw_0$, $y'=w_0yw_0$ and $z'=w_0zw_0$.
  \end{Lemma}

  \begin{proof}The definitions above ensure that
    $Z_\blam=w_0\Sym(\phi)w_0$.
    Finally,
   $P_{zx,y}=P_{w_0z'x'w_0,w_0y'w_0}=P_{z'x',y'}$ since
   $P_{u,v}=P_{w_0uw_0,w_0vw_0}$, for any $u,v\in\Sym_N$.
  \end{proof}

  The definitions above show that
  $\col(\T^\blam)v_\blam w_0^{\phi}w_0=\col(\T^\blam)w_\blam$,
  %Now  $w_0v_\blam\in\D_\phi^+$, so $w_0v_\blam w^\phi_0\in\D_\phi$,
  where
  $w^\phi_0$ is the element of longest length in~$\Sym(\phi)$.
  Therefore, by the minimality of $w_{\blam}$, $w_\blam=v_\blam w^\phi_0 w_0$, by
  \autoref{L:Conjugation}, and consequently
  $\wt(\blam)=v_\blam\cdot\phi$, where $\wt(\blam)$ is defined in \autoref{E:TlamWeight}. Using
  \autoref{T:UngradedCatOEquivalence} it follows that
  $\UnEFun(\unL(v_\blam\cdot\phi))\cong \unL^\blam$.
  Hence, $\O[\beta]$ is a
  subcategory of $\Ocal_\phi$.

  Recall from \autoref{S:HigherSchurWeyl} that the partition~$\pi$
  determines the parabolic subalgebra~$\p$. Let $\psi\in\h^*$ be any
  dominant integral weight such that $\Sym(\psi)=\Sym_\pi$, where
  $\Sym_\pi$ is the Young subgroup determined by~$\pi$. Using
  Backelin's notation, $\O=\Ocal^\p=\Ocal^\psi$ so that
  $\O[\beta]=\Ocal^\psi_\phi$ and $\SO=R^\psi_\phi$.  Consequently, by
  \cite[Theorem~1.1]{Backelin:Koszul}, a
  precise statement of \autoref{T:CatOKoszul} is that $R^\psi_\phi$ is
  Koszul and $E(R^\psi_\phi)=R^\phi_\xi$, where $\xi=-w_0\psi$ (no dot action!).
  Consequently, we identify the categories $\dualO$ and $R^\phi_\xi\Mod$.
%and take $\dualS=R^\phi_\xi$.

  Using \autoref{L:Conjugation} and the fact that $w_\blam=v_\blam w^\phi_0 w_0$, it is easy to see that the next result will turn out to equivalent to
  \autoref{T:CatODecompNumbers}.  Let $\q=\q(\phi)$ be the parabolic
  subalgebra of $\gl$ with Weyl group $\Sym(\phi)$. When $\q=\b$
  \autoref{P:CatODecompNumbers} is a restatement of
  \cite[Theorem~3.11.4(ii) and (iv)]{BGS:Koszul}. (When comparing
  \autoref{P:CatODecompNumbers} with \cite{BGS:Koszul} note that, in
  their formulas, $x$ and $y$ are \textit{maximal} length left coset
  representatives whereas $v_\blam$ and $v_\bmu$ are \textit{minimal}
  length {\it left} coset representatives of $\Sym(\phi)$ in $\Sym_N$.)

\begin{Proposition}\label{P:CatODecompNumbers}
  Suppose that $\blam,\bmu\in\Parts[\beta]$. Then
  $$p^\Ocal_{\blam\bmu}(q) = q^{\ell(v_\bmu)-\ell(v_\blam)}
      \sum_{z\in\Sym(\phi)}(-1)^{\ell(z)} P_{v_\blam w_0^{\phi}zw_0,v_\bmu w_0^{\phi}w_0}(q^{-2}).$$
\end{Proposition}

\begin{proof}
  As we have already noted, by \autoref{T:CatOKoszul} and
  \cite[Theorem 5.1]{Mazorchuk:CatLinTiltingComplexes}, the algebras
  $R^\psi_\phi$ and $R^\phi_\xi$ are standard Koszul algebras.
  By the remarks above,
  $\set{\unL(v_{\bmu}\cdot\phi)|\bmu\in\Parts[\beta]}$ is a complete set
  of simple modules in the category $\Ocal_\phi^\psi$. It follows that
  the Koszul dual category $\Ocal^\phi_\xi$ has (ungraded) simple
  modules $\set{\unL(w_0^{\phi}v_\bmu^{-1}w_0\cdot\xi)|\bmu\in\Parts[\beta]}$ and
  that the standard modules for $\Ocal^\phi_\xi$ correspond to the parabolic
  Verma modules
  $\set{\UnM[\q](w_0^{\phi}v_\bmu^{-1}w_0\cdot\xi)|\bmu\in\Parts[\beta]}$.  Let
  $\{L(w_0^{\phi}v_\bmu^{-1}w_0\cdot\xi)\}$ and $\{M_\q(w_0^{\phi}v_\bmu^{-1}w_0\cdot\xi)\}$
  be the corresponding graded simple and standard modules for
  $R^\phi_\xi$. The graded lifts of these modules exist because $R^\phi_\xi$ is a
  standard Koszul algebra by \cite[Theorem~5.1]{Mazorchuk:CatLinTiltingComplexes}.

  By \autoref{L:KoszuldualDecomp},
  $p^\Ocal_{\blam\bmu}(q)
             =[M_\q(w_0^{\phi}v_\blam^{-1}w_0\cdot\xi): L(w_0^{\phi}v_\bmu^{-1}w_0\cdot\xi)]_q$.
  The proposition is equivalent to the identity:
  \begin{equation}\label{E:KLIdentity}
    [M_\q(w_0^{\phi}v_\blam^{-1}w_0\cdot\xi): L(w_0^{\phi}v_\bmu^{-1}w_0\cdot\xi)]_q
         =q^{\ell(v_\blam)-\ell(v_\bmu)}
      \sum_{z\in\Sym(\phi)}(-1)^{\ell(z)} P_{v_\blam w_0^{\phi}zw_0,v_\bmu w_0^{\phi}w_0}(q^{-2}),
  \end{equation}
  for all $\blam,\bmu\in\Parts[\beta]$. Let $\levi$ be the Levi
  subalgebra of~$\q$ and set $\nu=x\cdot\xi$ for $x\in\D^+_\phi$, so
  that $\nu$ is a dominant weight for $\levi$. By definition,
  $\UnM[\q](\nu)=U(\mathfrak{g})\otimes_{U(\mathfrak{q})}\unL_\psi(\nu)$,
  where $\unL_\psi(\nu)$ is the (ungraded) irreducible integrable
  highest weight $\levi$-module of highest weight $\nu$. The irreducible
  module $\unL_\psi(\nu)$ lives in the BGG category $\Ocal(\levi)$, for
  the Lie algebra $\levi$, so it has a BGG resolution of the form
  $$0\rightarrow \underline{C}^\psi_m\rightarrow\dots\rightarrow
          \underline{C}^\psi_0 \rightarrow\unL_\psi(\nu)\rightarrow 0,
          \qquad\text{where }
      \underline{C}^\psi_k=\bigoplus_{\substack{z\in\Sym(\phi)\\\ell(z)=k}}
          \UnM[\psi](z\cdot\nu)
  $$
  and $\UnM[\psi](z\cdot\nu)$ is an ungraded Verma module in $\Ocal(\levi)$. In
  particular, $\underline{C}^\psi_0=\UnM[\psi](\nu)$ and
  $m=\ell(w^\phi_0)$. Applying the parabolic inflation-induction functor
  $U(\mathfrak{g})\otimes_{U(\q)}?$ gives an exact sequence
  $$0\rightarrow \underline{C}_{m}\rightarrow\dots\rightarrow
          \underline{C}_{0}=\unM(\nu)
          \rightarrow\UnM[\q](\nu)\rightarrow 0,
          \qquad\text{where }
      \underline{C}_{k}=\bigoplus_{\substack{z\in\Sym(\phi)\\\ell(z)=k}}
          \unM(z\cdot\nu),
  $$
  and $\unM(z\cdot\nu)$ is an ungraded Verma module.  By construction,
  all of the terms of this resolution live in the block~$\Ocal_\xi^0$ of
  (ungraded) category~$\Ocal$.  Let $M(z\cdot\nu)$ be the graded Verma
  module for~$R^0_\xi$ constructed in \cite[\S3.11]{BGS:Koszul}. We claim that this
  resolution of the ungraded parabolic Verma module $\UnM[\q](\nu)$
  lifts to a exact sequence of graded $R_\xi^0$-modules
  \begin{equation}\label{E:GradedBGG}
    0\rightarrow C_{m}\rightarrow\dots\rightarrow
          C_{0}=M(\nu) \rightarrow M_\q(\nu)\rightarrow 0,
          \qquad\text{where }
      C_{k}=\bigoplus_{\substack{z\in\Sym(\phi)\\\ell(z)=k}}
          M(z\cdot\nu),
  \end{equation}
  and all of the maps are homogeneous of degree~$1$ except for the
  map $C_0=M(\nu)\rightarrow M_\q(\nu)$, which is homogeneous of degree~$0$.
  Up to a scalar, there is a unique surjective
  map $M(\nu)\rightarrow M_\q(\nu)$, namely, the canonical quotient map which
  is homogeneous of degree zero. On the other hand, there is a non-zero
  homogeneous map $M(\sigma)\rightarrow M(\tau)$ of degree~$k$ only if
  $[M(\tau):L(\sigma)\<k\>]\ne0$. By
  \cite[Theorem~3.11.4(ii)(iv)]{BGS:Koszul} (and \autoref{P:Rigidity}),
  if $x,y\in(\D_\xi^{+})^{-1}$ then
  \begin{equation}\label{E:BGS}
    [M(x\cdot\xi): L(y\cdot\xi)]_q = q^{\ell(y)-\ell(x)}P_{x,y}(q^{-2}).
  \end{equation}
Note that $v_{\blam}\in\D_{\phi}^{-1}\cap\D_{\psi}$, because
$v_{\blam}\cdot\phi\in X_{\psi}^{+}$, and $\D_{\xi}^{+}=w_0\D_{\psi}^{+}w_0$.
If $z\in\Sym(\phi)$ then $v_{\blam}w_0^{\phi}z^{-1}\in\D_{\psi}$ because
$v_{\blam}w_0^{\phi}z^{-1}\cdot\phi=v_{\blam}\cdot\phi\in X_{\psi}^{+}$.
Hence, $v_{\blam}w_0^{\phi}z^{-1}w_0\in\D_{\psi}^{+}$. It follows that
\begin{equation}\label{maximal}
    zw_0^{\phi}v_{\blam}^{-1}w_0\in(\D_{\xi}^{+})^{-1}
    \quad\text{and}\quad
    w_0^{\phi}v_{\bmu}^{-1}w_0\in(\D_{\xi}^{+})^{-1}.
\end{equation}

%  (This is equivalent to the formula in \cite{BGS:Koszul}, which is
%  written in terms of maximal length coset representatives, because if
%  $x,y\in\D_\xi$ then $xw^\phi_0,yw^\phi_0\in\D_\phi^+$ and
%  $P_{x,y}=P_{xw^\phi_0,yw^\phi_0}$ by
%  \cite[Proposition~3.4]{deod:relII}.)

  Recall from \cite{KL} that $P_{x,y}\ne0$ only if $x\le y$ (where $\le$ is the Bruhat
  order) and $\deg_qP_{x,y}(q)\le\tfrac12(\ell(y)-\ell(x)-1)$.
  Therefore, if $\ell(y)=\ell(x)+1$ then $[M(x\cdot\xi):
  L(y\cdot\xi)]_q\ne0$ if and only if $y=sx$, for some simple reflection
  $s\in\Sym_N$, in which case $[M(x\cdot\xi): L(y\cdot\xi)]_q=q$.
  Therefore, up to a scalar, if $y=sx$ then there is at most one map
  $M(y\cdot\xi)\rightarrow M(x\cdot\xi)$, which must be homogeneous of
  degree one, and if $y\ne sx$ then
  $\Hom_{\Ocal_\xi^0}(M(y\cdot\xi),M(x\cdot\xi))=0$. Arguing by
  induction on~$k$ to determine the degree shifts on the Verma modules,
  it follows that the ungraded resolution of $\UnM[\q](\nu)$
  lifts to a graded ``BGG resolution'' of~$M_\q(\nu)$ as
  in~\autoref{E:GradedBGG}.

  We now complete the proof of the theorem. Fix $\blam,\bmu\in\Parts[\beta]$
  and set $x=w_0^{\phi}v_\blam^{-1}w_0$ and $y=w_0^{\phi}v_\bmu^{-1}w_0$
  in the formulas above. Using \autoref{E:GradedBGG} for the first
  equality, and \autoref{E:BGS} and \autoref{maximal} for the second,
  \begin{align*}
     [M_\q(w_0^{\phi}v_\blam^{-1}w_0\cdot\xi):L(w_0^{\phi}v_\bmu^{-1}w_0\cdot\xi)]_q
       &=\sum_{z\in\Sym(\phi)}(-q)^{\ell(z)} [M(zw_0^{\phi}v_\blam^{-1}w_0\cdot\xi):L(w_0^{\phi}v_\bmu^{-1}w_0\cdot\xi)]_q,\\
       &=\sum_{z\in\Sym(\phi)}(-q)^{\ell(z)}
       q^{\ell(v_\bmu^{-1})-\ell(zv_\blam^{-1})}P_{zw_0^{\phi}v_\blam^{-1}w_0,w_0^{\phi}v_\bmu^{-1}w_0}(q^{-2}),\\
       &=q^{\ell(v_\bmu)-\ell(v_\blam)}\sum_{z\in\Sym(\phi)}(-1)^{\ell(z)}
          P_{w_0 v_\blam w_0^{\phi}z,w_0v_\bmu w_0^{\phi}}(q^{-2})\\
       &=q^{\ell(v_\bmu)-\ell(v_\blam)}\sum_{z\in\Sym(\phi)}(-1)^{\ell(z)}
          P_{v_\blam w_0^{\phi}zw_0,v_\bmu w_0^{\phi}w_0}(q^{-2}),
  \end{align*}
  where the last two equalities use \autoref{L:Conjugation} and the well-known facts that
  $P_{u,v}=P_{u^{-1},v^{-1}}$ and $P_{u,v}=P_{w_0uw_0,w_0vw_0}$,
  for all $u,v\in\Sym_N$.
%and that $\ell(zd^{-1})=\ell(d)+\ell(z)$, whenever $d\in\D_\phi$ and $z\in\Sym(\phi)$.
  This completes the proof of our claim~\autoref{E:KLIdentity} and hence
  the proof of the proposition.
\end{proof}

Combining \autoref{L:Conjugation} and \autoref{P:CatODecompNumbers} we
obtain \autoref{T:CatODecompNumbers}, which is what this appendix set
out to prove.  As a corollary of \autoref{P:CatODecompNumbers}, using
\autoref{L:KoszuldualDecomp} together with \cite[Theorem~4.6]{doug} or
\cite[Proposition~3.17]{dual}, we can identify the graded decomposition
numbers of $\SO[\beta]$ with certain inverse parabolic Kazhdan--Lusztig polynomials.
By \autoref{T:Koszul} these are also the graded
decomposition numbers of~$\Sch[\beta]$.

%\begin{Corollary}\label{C:OGradedDecompNumbers}
%  Suppose that $\blam,\bmu\in\Parts$. Then
%  $d^\Ocal_{\blam\bmu}(q)=q^{\ell(w_\bmu)-\ell(w_\blam)}P_{w_\blam,w_\bmu}(q^{-2})$.
%\end{Corollary}

\section{Quiver Schur algebras of level two} In this appendix we assume that
either $e=0$ or $e>n$ and we fix a dominant weight
$\Lambda$ of level~$\ell=2$.  We will show that $\Sch[n]$ is a positively graded
basic algebra and, as a consequence, we will give a beautiful closed formula for the
graded decomposition numbers for these algebras. Our formula for these graded
decomposition numbers is new, however, by \autoref{THEOREMC} when $e=0$ these graded
decomposition numbers can be computed in parabolic category~$\Ocal$ where
different formulations of this result are already known, all going back to the
work of Lascoux and Sch{\"u}tzenberger~\cite{LasSch:KLPoly}. When $e=0$ all of
the results in this section have been obtained by Brundan and
Stroppel~\cite{BrundanStroppel:KhovanovII,BrundanStroppel:KhovanovIII} using
different methods. The extension of these results to the case $e>n$ is new.

Suppose that $\t=(\t^{(1)},\t^{(2)})$ is a standard tableau. Let
$\unt^{(c)}=\set{1\le k\le n|\comp_\t(k)=c}$ be the integers in component $c$
of~$\t$, for $c=1,2$. By assumption, $e=0$ or $e>n$ so the nodes of constant
residue in~$\t^{(c)}$ all appear on the same \textit{diagonal}
$\set{(a+d,b+d,c)\in\bmu|d\in\Z}$ in~$\t^{(c)}$. Therefore, the tableau~$\t$ is
uniquely determined by its residue sequence $\bi=\res(\t)\in I^n$ and the sets
$\unt^{(1)}$ and $\unt^{(2)}$.

Although we won't need this, the last paragraph implies that if
$i_{r+2}=i_r=i_{r+1}\pm1$, where $1\le r\le n-2$ and $\bi=\res(\t)$ for some
tableau $\t\in\Std(\Parts)$, then the permutation $d(\t)$ fixes $r$, $r+1$ and
$r+2$.  Consequently $s_rs_{r+1}s_r$ cannot appear in any reduced expression
for~$d(\t)$ so that $\psi_{d(\t)}$ depends only on~$\t$ and not on a choice of
reduced expression for~$d(\t)$. Therefore, in level two the basis
elements~$\psi_{\s\t}$ depend only on~$\s$ and~$\t$, and not on the choices of
reduced expressions. We warn the reader that this does not imply that the
$\psi_r$ satisfy the braid relations in~$\R$.

Following~\cite{M:gendeg}, define a tableau $\t$ to be \textbf{regular} if its
entries increase along the diagonals in each component. It is easy to see that
all standard tableaux are regular and that there exist regular tableaux that
are not standard. By the last paragraph, given a sequence $\bi\in I^n$ and
disjoint sets $A_1$ and $A_2$ such that $A_1\cup A_2=\{1,2,\dots,n\}$ there
exists a unique regular tableau~$\t$ such that $\res(\t)=\bi$ and
$\unt^{(c)}=A_c$, for $c=1,2$.  Note that $\t$ is not necessarily standard and,
in general, that the shape of~$\t$ need not be a bipartition.

Two nodes $(e,c,l)$ and $(r',c',l')$ are \textbf{adjacent} if $l=l'$ and either
$r=r'$ and $c=c'\pm1$, or $c=c'$ and $r=r'\pm1$. A set~$X$ of nodes is
\textbf{connected} if for any $x,y\in X$ there is a sequence $x=x_1,\dots,x_z=y$
of nodes in~$X$ such that $x_i$ and $x_{i+1}$ are adjacent, for $1\le i<z$.

\begin{Lemma}\label{L:Level2Positivity}
   Suppose that $\t\in\Std^\bmu(\blam)$. Then $\deg\t\ge\deg\tmu$ with equality
   if and only if $\t=\tmu$.
\end{Lemma}

\begin{proof}By definition, $\t\gedom\tmu$ and $\res(\t)=\res(\tmu)$. We argue
  by induction on dominance. If $\t=\tmu$ then there is nothing to prove, so
  suppose that $\t\gdom\tmu$. Let~$a$ be the smallest number in
  $\unt^{(1)}\cap\unt^{\mu^{(2)}}$. As remarked above, the tableau $\t$ is
  uniquely determined by its residue sequence $\res(\t)$ and the sets
  $\unt^{(c)}$, for $c=1,2$.  Therefore, if $1\le b\leq a$ then $b$ appears in
  exactly the same position in~$\t$ and in~$\tmu$. In particular,~$a$ is larger
  than all of the numbers in~$\unt^{\mu^{(1)}}$. Moreover, $a$ is uniquely
  determined by~$\blam$ and~$\bmu$.

  Let $X$ be the set of nodes in $\blam^{(1)}\setminus\bmu^{(1)}$ that are
  connected to $\t^{-1}(a)$ such that they are either adjacent to a node
  in~$\mu^{(1)}$ or they are in the first row or in the first column
  of~$\blam^{(1)}$. Since~$a$ is uniquely determined
  by~$\blam$ and $\bmu$, it follows that~$X$ also depends only on~$\blam$
  and~$\bmu$.

  Let $A=\set{\t(x)|x\in X}\subseteq\unt^{(1)}$. Define
  $\t_A=(\t_A^{(1)},\t_A^{(2)})$ to be the unique regular tableau with residue
  sequence $\res(\t)=\res(\tmu)$ such that $\unt^{(1)}_A=\unt^{(1)}\setminus A$
  and $\unt_A^{(2)}=\unt^{(2)}\cup A$. That is, $\t_A$ is the regular tableau
  obtained by moving the numbers in~$A$ from the first component of~$\t$ to the
  second component.

  For example, suppose that $e=0$ and $\charge=(0,1)$, so that
  $\Lambda=\Lambda_0+\Lambda_1$.  Let $\mu=(2,1^2|4^2,2)$ and
  $\blam=(5^2,2|1^2)$. Then $\t\in\Std^\bmu(\blam)$, $a=6$ and
  $$ \t=\Bigg(\ \ShadedTableau[(5,0),(4,0),(3,0),(3,-1),(2,-1),(2,-2)]%
                     {{1,2,6,7,8},{3,9,10,11,12},{4,13}}\ \Bigg|\
                     \Tableau{{5},{14}}\ \Bigg)
                     \leadsto
    \t_A=\Bigg(\ \ShadedTableau{{1,2,11,12},{3},{4}}\ \Bigg|\
    \ShadedTableau[(4,0),(3,0),(2,0),(2,-1),(1,-1),(1,-2)]{{5,6,7,8},{9,10},{13,14}}
    \ \Bigg)$$
  The shaded nodes in $\t$ mark the elements of~$A=\{6,7,8,9,10,13\}$.

  By the remarks in the first paragraph, the elements of~$A$ occupy the same
  positions in~$\t_A$ as they do in~$\tmu$. By definition, the elements of~$A$
  have distinct residues. Moreover, we can order $A=\{a_0,\dots,a_z\}$ so that
  $\res_t(a_i)=\res_t(a_0)+i$, for $0\le i\le z$. The definition of~$X$ implies
  that if $b\in A$ then $b$ is the smallest element of
  $\unt^{\mu^{(2)}}\cap\unt_1$ with residue $\res_\t(b)$. Therefore, two
  elements of~$A$ are in the same row of~$\t^{(1)}$ if and only if they are in
  the same row of~$\t_A^{(2)}$. Hence, the elements of~$A$ that are in the same
  row of~$\t$, or equivalently of~$\t_A$ or of $\tmu$, are consecutive.
  It follows that the set~$A$ is also determined by~$a$ (and $\res(\tmu)$), and
  hence that~$A$ is uniquely determined by~$\blam$ and~$\bmu$.

  By definition, $\t_A$ is obtained by moving the numbers in~$A$ from the first
  component to the second component, without changing their `shape', then
  `sliding' numbers down the diagonals in the first component to fill the gaps
  where the elements of~$A$ used to be, and then sliding numbers up the diagonal
  in the second component to make way for the elements of~$A$. Since $\t$ and
  $\tmu$ are both standard it follows that~$\t_A$ is also standard and that
  $\t\gdom\t_A\gedom\tmu$. In particular, $\t_A\in\Std^\bmu(\Parts)$.

  As remarked earlier, the elements in~$A$ occurring in a given row are
  consecutive. By definition, $a$ is the smallest element of~$A$ and
  $\t_{\downarrow(a-1)}=\t_{A\downarrow(a-1)}$. It is easy to see that
  $\deg\t_{\downarrow a}=\deg\t_{A\downarrow a}+1$.
  Adding the elements of~$A$ row by row to $\t_{\downarrow(a-1)}$ and to
  $\t_{A\downarrow(a-1)}$ it is easy to see that the only difference in the degrees
  of the tableaux~$\t$ and $\t_A$ occurs when adding~$a$, which appears in the
  `first row' of~$A$, and that the subsequent rows in~$A$ do not change the
  degrees of~$\t$ or of~$\t_A$. Hence, $\deg\t_{\downarrow
  z}=\deg\t_{A\downarrow z}+1$, where $z$ the largest element of~$A$. In view of
  the sliding construction of~$\t_A$, this implies that $\deg\t=\deg\t_A+1$.
  Therefore, by induction, $\deg\t>\deg\t_A\ge\deg\tmu$ as required.
\end{proof}

The set $A$ in the proof of \autoref{L:Level2Positivity} is uniquely determined
by the bipartitions $\blam$ and $\bmu$. Moreover,~$\t$ can be recovered from
$\t_{A}$ and $\Shape(\t_{A})$ is uniquely determined by the bipartitions $\blam$
and $\bmu$. Hence as a byproduct of the proof, we have the following.

\begin{Corollary}\label{C:L2Decomp}
  Suppose that $\blam,\bmu\in\Parts$. Then $\#\Std^\bmu(\blam)\le1$.
\end{Corollary}

If $\Std^\bmu(\blam)\ne\emptyset$ let $\tmu_\blam$ be the unique $\blam$-tableau in
$\Std^\bmu(\blam)$.

\begin{Theorem}\label{T:L2Positivity}
  Suppose that $\mz$ is a field, $e=0$ or $e>n$ and that $\Lambda\in P^+$ is
  a weight of level~$2$. Then $\Sch[n]$ is a positively graded basic
  algebra. Moreover,
  $$[\Delta^\blam:L^\bmu]_q =\begin{cases}
            q^{\deg\tmu_\blam-\deg\tmu},&\text{if }\Std^\bmu(\blam)\ne\emptyset,\\
            0,&\text{otherwise,}
          \end{cases}$$
  for $\blam,\bmu\in\Parts$.
\end{Theorem}

\begin{proof}By \autoref{T:SCellular}, $\Sch[n]$ is a quasi-hereditary
  cellular algebra with cellular basis
  $$\set{\Psi^{\bmu\bnu}_{\s\t}|(\bmu,\s),(\bnu,\t)\in\Tcal^\blam
                                      \text{ and }\blam\in\Parts}.$$
  Moreover, if $\Psi^{\bmu\bnu}_{\s\t}$ is one of these basis elements then
  $\deg\Psi^{\bmu\bnu}_{\s\t}=(\deg\s-\deg\tmu)+(\deg\t-\tnu)\ge0$ by
  \autoref{L:Level2Positivity}. Therefore, the quiver Schur algebra $\Sch[n]$ is positively
  graded.

  Now suppose that $\blam\in\Parts$. Then $\set{\Psi^\bnu_\t(\nu,\t)\in\Tcal^\blam}$
  is a basis of $\Delta^\blam$ by~\autoref{E:WeylBasis}. Moreover, by
  \autoref{L:Level2Positivity}, $\deg\Psi^\bnu_\t\ge0$ with equality if and
  only if $(\t,\bnu)=(\tlam,\blam)$. It follows that the simple module~$L^\blam$ is
  one dimensional with basis vector $\Psi^\blam_{\tlam}+\rad\Delta^\blam$
  since $\Dim L^\blam=\Dim(L^\blam)^\circledast=\overline{\Dim L^\blam}$ by
  \autoref{T:CellularSimples}. Thus, $\Dim L^\blam=1$, for all
  $\blam\in\Parts$, and $\Sch[n]$ is a basic algebra.

  Finally, since $\Sch[n]$ is positively graded, $Z^\bmu=P^\bmu$ for all
  $\bmu\in\Parts$, for example by applying our LLT algorithm from
  \autoref{S:LLT}. Therefore, by \autoref{E:ZmuExpansion},
  $$[\Delta^\blam:L^\bmu]_q=\sum_{\t\in\Std^\bmu(\blam)}q^{\deg\t-\deg\tmu}
  =\begin{cases} q^{\deg\tmu_\blam-\deg\tmu}, & \text{if $\Std^\bmu(\blam)\ne\emptyset$;}\\
  0, &\text{if $\Std^\bmu(\blam)=\emptyset$,}
  \end{cases}
  $$
  where the last equality comes from \autoref{C:L2Decomp}.
\end{proof}

\begin{Corollary}
  Suppose that $\mz$ is a field, $e=0$ or $e>n$ and that $\Lambda\in P^+$ is
  a weight of level~$2$. Then $G^\bmu=Y^\bmu$ is an
  indecomposable graded Young module, for all $\bmu\in\Parts$.
\end{Corollary}

\begin{Corollary}
  Suppose that $\mz$ is a field, $e=0$ or $e>n$ and that $\Lambda\in P^+$ is
  a weight of level~$2$. Then the graded decomposition numbers of~$\Sch[n]$ and
  $\H$, and the graded dimensions of their graded simple modules, are
  independent of the characteristic of~$\mz$.
\end{Corollary}

By \autoref{THEOREMC}, \cite[Corollary~8.20]{BrundanStroppel:KhovanovIII} and
the uniqueness of Koszul gradings~\cite[Corollary~2.5.2]{BGS:Koszul} we obtain
the link with Brundan and Stroppel's work.

\begin{Corollary}\label{C:BrundanStroppel}
  Suppose that $e=0$. Then $\Sch[n]$ is isomorphic, as a graded algebra, to the
  quasi-hereditary algebra $K^\Lambda_n$ defined by Brundan and
  Stroppel~\cite{BrundanStroppel:KhovanovII}.
\end{Corollary}

\let\u=\uold
%\bibliography{papers}

\end{document}